\documentclass[11pt,openany,leqno]{article}
\usepackage{amsmath,amsthm,amsfonts,amssymb,amscd}
\usepackage{graphics}
\usepackage[latin1]{inputenc}

\begin{document}

\theoremstyle{plain}
\newtheorem{Theo}{\bf Theorem}
\newtheorem{lema}{\bf Lemma}
\newtheorem{Propo}{\bf Proposition}
\newtheorem{Coro}{\bf Corollary}

\theoremstyle{rem}
\newtheorem{rem}{\bf Remark$\!\!$}
\renewcommand{\therem}{}

\newcommand{\Def}{\medskip\noindent
                 {\bf Definition.}~}
\newcommand{\Exem}{\medskip\noindent
                 {\bf Example.}~}

\errorcontextlines=0

\def\sectionname{Section}
\renewcommand{\thesection}{\arabic{section}}
\renewcommand{\theequation}{\arabic{equation}}
\setcounter{section}{0}

\def\Diff{\text{Diff}}
\def\Max{\text{max}}
\def\Log{\text{log}}
\def\loc{\text{loc}}
\def\inta{\text{int }}
\def\det{\text{det}}
\def\exp{\text{exp}}
\def\Re{\text{Re}}
\def\lip{\text{Lip}}
\def\leb{\text{Leb}}
\def\dom{\text{Dom}}
\def\diam{\text{diam}}

\def\gr{\operatorname{grau}}
\def\sen{\operatorname{sen}}
\def\ker{\operatorname{ker}}

\newcommand{\AL}{\hbox{\Large{$a$}}}
\newcommand{\AH}{\hbox{\Huge{$a$}}}
\newcommand{\loz}{{\lozenge}}

\newcommand{\Bbc}{{\mathbb{C}}}
\newcommand{\Bbf}{{\mathbb{F}}}
\newcommand{\Bbn}{{\mathbb{N}}}
\newcommand{\Bbp}{{\mathbb{P}}}
\newcommand{\Bbq}{{\mathbb{Q}}}
\newcommand{\Bbr}{{\mathbb{R}}}
\newcommand{\Bbz}{{\mathbb{Z}}}

\newcommand{\cA}{{\mathcal A}}
\newcommand{\cB}{{\mathcal{B}}}
\newcommand{\cC}{{\mathcal{C}}}
\newcommand{\cD}{{\mathcal D}}
\newcommand{\cE}{{\mathcal E}}
\newcommand{\cf}{{\mathcal F}}
\newcommand{\cG}{{\mathcal{G}}}
\newcommand{\cH}{{\mathcal H}}
\newcommand{\cI}{{\mathcal{I}}}
\newcommand{\cL}{{\mathcal{L}}}
\newcommand{\cM}{{\mathcal{M}}}
\newcommand{\cN}{{\mathcal N}}
\newcommand{\cO}{{\mathcal O}}
\newcommand{\cP}{{\mathcal{P}}}
\newcommand{\cR}{{\mathcal{R}}}
\newcommand{\cU}{{\mathcal{U}}}

\newcommand{\vep}{{\varepsilon}}
\newcommand{\al}{{\alpha}}
\newcommand{\be}{{\beta}}
\newcommand{\de}{{\delta}}
\newcommand{\De}{{\Delta}}
\newcommand{\ga}{{\gamma}}
\newcommand{\Ga}{{\Gamma}}
\newcommand{\ka}{{\kappa}}
\newcommand{\la}{{\lambda}}
\newcommand{\La}{{\Lambda}}
\newcommand{\Om}{{\Omega}}
\newcommand{\om}{{\omega}}
\newcommand{\pa}{{\partial}}
\newcommand{\ro}{{\rho}}
\newcommand{\vphi}{{\varphi}}
\newcommand{\te}{{\theta}}
\newcommand{\Te}{{\Theta}}
\newcommand{\sig}{{\sigma}}
\newcommand{\Sig}{{\Sigma}}

\newcommand{\wtA}{\widetilde{A}}
\newcommand{\wtB}{\widetilde{B}}
\newcommand{\wtcD}{\widetilde{\cD}}
\newcommand{\wtcN}{\widetilde{\cN}}
\newcommand{\wtcR}{\widetilde{\cR}}
\newcommand{\wtb}{\widetilde{b}}
\newcommand{\wtn}{\widetilde{n}}
\newcommand{\wtx}{\widetilde{x}}
\newcommand{\wtI}{\widetilde{I}}
\newcommand{\wtP}{\widetilde{P}}
\newcommand{\wtQ}{\widetilde{Q}}
\newcommand{\wtR}{\widetilde{R}}
\newcommand{\wtT}{\widetilde{T}}
\newcommand{\wtnu}{\widetilde{\nu}}
\newcommand{\wtbe}{\widetilde{\beta}}
\newcommand{\wtcdm}{\wtcD_{+}}
\newcommand{\wttm}{\wtT^{+}}
\newcommand{\rim}{R^{\infty}_{+}}
\newcommand{\wtrim}{\wtR^{\infty}_{+}}

\newcommand{\whI}{{\widehat{I}}}
\newcommand{\whd}{{\widehat{d}}}
\newcommand{\whh}{{\widehat{h}}}
\newcommand{\whm}{{\widehat{m}}}
\newcommand{\whn}{{\widehat{n}}}
\newcommand{\whs}{{\widehat{s}}}
\newcommand{\whu}{{\widehat{u}}}
\newcommand{\whz}{{\hat{z}}}
\newcommand{\whcC}{{\widehat{\cC}}}
\newcommand{\whA}{{\widehat{A}}}
\newcommand{\whB}{{\widehat{B}}}
\newcommand{\whP}{{\widehat{P}}}
\newcommand{\whR}{{\widehat{R}}}
\newcommand{\whQ}{{\widehat{Q}}}
\newcommand{\whbe}{{\widehat{\beta}}}
\newcommand{\whom}{{\widehat{\omega}}}
\newcommand{\whsig}{{\widehat{\sigma}}}

\newcommand{\ov}{\overline}
\newcommand{\ovbe}{{\overline{\beta}}}
\newcommand{\ovC}{{\overline{C}}}
\newcommand{\ovn}{{\overline{n}}}
\newcommand{\ovx}{{\overline{x}}}
\newcommand{\ovP}{{\overline{P}}}
\newcommand{\ovQ}{{\overline{Q}}}
\newcommand{\ovX}{{\overline{X}}}
\newcommand{\ovY}{{\overline{Y}}}
\newcommand{\ovfork}{{\overline{\pitchfork}}}
\newcommand{\ovforki}{{\overline{\pitchfork}_{I}}}
\newcommand{\whfork}{{\widehat{\pitchfork}}}
\newcommand{\whforki}{{\widehat{\pitchfork}_{I}}}

\newcommand{\una}{{\underline{a}}}
\newcommand{\unb}{{\underline{b}}}

\newcommand{\Lg}{\bigg\langle}
\newcommand{\Rg}{\bigg\rangle}
\newcommand{\lan}{\langle}
\newcommand{\ran}{\rangle}
\newcommand{\lgg}{{\left\langle\right.}}
\newcommand{\rg}{{\left.\right\rangle}}
\newcommand{\lV}{{\left\Vert \right.}}
\newcommand{\rV}{{\left.\right\Vert}}
\newcommand{\lv}{{\left\vert \right.}}
\newcommand{\rv}{{\left.\right\vert}}

\newcommand{{\fud}}{{\frac{1}{2}}}
\newcommand{{\fudt}}{{\tfrac{1}{2}}}
\newcommand{{\fut}}{{\frac{1}{3}}}
\newcommand{{\futt}}{{\tfrac{1}{3}}}
\newcommand{\fork}{{\pitchfork}}
\newcommand{\forki}{{\pitchfork_{I}}}
\newcommand{\cRim}{{\cR^{\infty}_{+}}}
\newcommand{\pii}{{(i)}}
\newcommand{\pjj}{{(j)}}
\newcommand{\pkk}{{(k)}}
\newcommand{\rl}{{r'}}
\newcommand{\pzz}{{(0)}}
\newcommand{\pnn}{{(n)}}
\newcommand{\prr}{{(r)}}
\newcommand{\pnu}{{(n-1)}}
\newcommand{\nuu}{{n-1}}
\newcommand{\pns}{{(n_s)}}
\newcommand{\An}{A^{(n)}}
\newcommand{\Bn}{B^{(n)}}
\newcommand{\dmu}{{\Delta^{-1}}}
\newcommand{\duz}{{\Delta^{-1}_{0}}}
\newcommand{\duu}{{\Delta^{-1}_{1}}}
\newcommand{\Ism}{{I_{m}}}
\newcommand{\Eio}{{E^{0}}}
\newcommand{\Io}{{I_{0}}}
\newcommand{\Isal}{{I_{\alpha}}}
\newcommand{\Qsal}{{Q_{\alpha}}}
\newcommand{\Psal}{{P_{\alpha}}}
\newcommand{\nsal}{{n_{\alpha}}}
\newcommand{\Isom}{{I_{\omega}}}
\newcommand{\nsom}{{n_{\omega}}}
\newcommand{\Psom}{{P_{\omega}}}
\newcommand{\Qsom}{{Q_{\omega}}}
\newcommand{\Is}{{I^{*}}}
\newcommand{\wtIs}{{\wtI^{*}}}
\newcommand{\As}{{A^{*}}}
\newcommand{\Bs}{{B^{*}}}
\newcommand{\Ps}{{P^{*}}}
\newcommand{\Qs}{{Q^{*}}}
\newcommand{\hs}{{h^{*}}}
\newcommand{\ns}{{n^{*}}}
\newcommand{\ys}{{y^{*}}}
\newcommand{\zs}{{z^{*}}}
\newcommand{\gas}{{\gamma^{*}}}
\newcommand{\kas}{{\kappa^{*}}}
\newcommand{\oms}{{\omega^{*}}}
\newcommand{\cRs}{{\cR^{*}}}
\newcommand{\nsr}{{n_{r}}}
\newcommand{\Psr}{{P_{r}}}
\newcommand{\Qsr}{{Q_{r}}}
\newcommand{\nsrl}{{n'_{r}}}
\newcommand{\Psrl}{{P'_{r}}}
\newcommand{\Qsrl}{{Q'_{r}}}
\newcommand{\nsu}{{n_{1}}}
\newcommand{\Psu}{{P_{1}}}
\newcommand{\Qsu}{{Q_{1}}}
\newcommand{\nsul}{{n'_{1}}}
\newcommand{\Psul}{{P'_{1}}}
\newcommand{\Qsul}{{Q'_{1}}}
\newcommand{\nsru}{{n_{r-1}}}
\newcommand{\Psru}{{P_{r-1}}}
\newcommand{\Qsru}{{Q_{r-1}}}
\newcommand{\nsd}{{n_{2}}}
\newcommand{\Psd}{{P_{2}}}
\newcommand{\Qsd}{{Q_{2}}}
\newcommand{\nsj}{{n_{j}}}
\newcommand{\Psj}{{P_{j}}}
\newcommand{\Qsj}{{Q_{j}}}
\newcommand{\nsm}{{n_{m}}}
\newcommand{\Psm}{{P_{m}}}
\newcommand{\Qsm}{{Q_{m}}}
\newcommand{\Esk}{{E_{k}}}
\newcommand{\nsk}{{n_{k}}}
\newcommand{\Psk}{{P_{k}}}
\newcommand{\Qsk}{{Q_{k}}}
\newcommand{\Vsk}{{V_{k}}}
\newcommand{\zk}{{z_{k}}}
\newcommand{\zo}{{z_{0}}}
\newcommand{\zell}{{z_{\ell}}}
\newcommand{\gask}{{\gamma_{k}}}
\newcommand{\vepk}{{\varepsilon_{k}}}
\newcommand{\nssk}{{n^{*}_{k}}}
\newcommand{\Pssk}{{P^{*}_{k}}}
\newcommand{\Qssk}{{Q^{*}_{k}}}
\newcommand{\Vssk}{{V^{*}_{k}}}
\newcommand{\Qsso}{{Q^{*}_{0}}}
\newcommand{\gassk}{{\gamma^{*}_{k}}}
\newcommand{\nskl}{{n'_{k}}}
\newcommand{\Pskl}{{P'_{k}}}
\newcommand{\Qskl}{{Q'_{k}}}
\newcommand{\gaskl}{{\gamma'_{k}}}
\newcommand{\Tim}{{T^{+}}}
\newcommand{\TiM}{{T^{-}}}
\newcommand{\Tsm}{{T_{+}}}
\newcommand{\TsM}{{T_{-}}}
\newcommand{\Ra}{{R_{a}}}
\newcommand{\Ral}{{R_{a'}}}
\newcommand{\Ld}{{L_{d}}}
\newcommand{\Lod}{{L^{0}_{d}}}
\newcommand{\hd}{{h_{d}}}
\newcommand{\hdl}{{h'_{d}}}
\newcommand{\Hd}{{H_{d}}}
\newcommand{\Hdl}{{H'_{d}}}
\newcommand{\Hdll}{{H''_{d}}}
\newcommand{\lad}{{\lambda_{d}}}
\newcommand{\ladl}{{\lambda'_{d}}}
\newcommand{\ladll}{{\lambda''_{d}}}
\newcommand{\mud}{{\mu_{d}}}
\newcommand{\mudl}{{\mu'_{d}}}
\newcommand{\nudl}{{\nu'_{d}}}

\newcommand{\ds}{{d_{s}}}
\newcommand{\dou}{{d^{0}_{u}}}
\newcommand{\dos}{{d^{0}_{s}}}
\newcommand{\dms}{{d^{+}_{s}}}
\newcommand{\dum}{{d^{+}_{u}}}
\newcommand{\dMs}{{d^{-}_{s}}}
\newcommand{\dss}{{d^{*}_{s}}}
\newcommand{\dsu}{{d^{*}_{u}}}

\def\za#1#2{#1^{+}_{#2}}
\def\zb#1#2{#1^{-}_{#2}}
\def\iim#1{#1^{i+1}}
\def\piim#1{#1^{(i+1)}}
\def\ism#1{#1_{i+1}}
\def\jim#1{#1^{j+1}}
\def\pism#1{#1_{(i+1)}}
\def\jiM#1{#1^{j-1}}
\def\jsM#1{#1_{j-1}}
\def\kiM#1{#1^{k-1}}
\def\ksM#1{#1_{k-1}}
\def\ksm#1{#1_{k+1}}
\def\dmuz#1{{\Delta^{-1}_{#1}}}
\def\qqim#1{{#1^{\infty}_{+}}}


\title{Non-Uniformly Hyperbolic Horseshoes Arising from\\
Bifurcations of Poincaré Heteroclinic Cycles\\[10pt]}
\author{Jacob PALIS\footnote{Partially support by CNPq and FAPERJ, Brazil.} ~and~ Jean-Christophe
YOCCOZ}
\date{}
\maketitle
\begin{abstract}
\noindent The purpose of this paper is to advance the knowledge of
the dynamics arising from the creation and subsequent bifurcation of
Poincaré heteroclinic cycles. The problem is central to dynamics: it
has to be addressed if, for instance, one aims at describing  the
typical orbit behaviour of a typical system, thus providing a global
scenario for the ensemble of dynamical systems - see the
Introduction and \cite{P1, P2}. Here, we shall consider smooth, i.e.
$C^\infty$, one-parameter families of dissipative, meaning
non-conservative, surface diffeomorphisms. An hetereoclinic cycle
may appear when the parameter evolves and an orbit of tangency, say
quadratic, is created between stable and unstable manifolds (lines)
of periodic orbits that belong to a basic hyperbolic set. The key
novelty is to allow this basic set, a horseshoe, to have Hausdorff
dimension bigger than one. In the present paper we do assume such a
dimension to be beyond one, but in a limited way, as explicitly
indicated in the Introduction. [A mild non-degeneracy condition on
the family of maps is assumed: at the orbit of tangency the
invariant lines, stable and unstable, cross each other with positive
relative speed]. We then prove that most diffeomorphisms,
corresponding to parameter values near the bifurcating one, are
non-uniformly hyperbolic in a neighborhood of the horseshoe and the
orbit of tangency; such diffeomorphisms display no attractors nor
repellors in such a neighborhood. A first precise formulation of our
main theorem is at the Introduction and a more encompassing version
at the end of the paper. These results were announced in \cite{PY3}.
\end{abstract}

\newpage
\setlength{\parskip}{0pt}
\tableofcontents

\renewcommand{\theequation}{\thesection.\arabic{equation}}
\setlength{\parskip}{12pt}
\section{Introduction\label{sec1}}
\subsection{The Context\label{sub1.1}}
Since Poincaré referred to a property valid ``pour la plupart des
coefficients'' of a polynomial (analytic) dynamical system, the
outstanding problem of describing at large the orbits of a
``typical'' system became the source of much creative work in
dynamics.

A stumbling block in a possible program to solve this question is
the understanding of dynamics arising from bifurcations of
homoclinic or heteroclinic cycles. Such cycles were defined by
Poincaré himself: they involve stable and unstable manifolds of
invariant sets (typically periodic orbits) that successively
intersect each other. In his classic book on Celestial Mechanics, he
prophetically stated: ``Rien n'est plus propre à nous donner une
idée de la complication de tous les problèmes de Dynamique''.

In fact, the creation and unfolding of cycles, in particular
homoclinic tangencies for surface diffeomorphisms and their
unfolding, led Newhouse to show the non-denseness of hyperbolic
dynamics, thus contradicting Smale's remarkable conjecture of the
early 60's. That is, there are systems that cannot be approximated
by one with a hyperbolic limit set. On the way, Newhouse showed that
in this context always appear surface diffeomorphisms displaying
infinitely many simultaneous sinks (periodic attractors) or sources
(periodic repellors).

Abundance of other more intricate kind of attractors, the so called
Hénon-like ones, was proved to be also present in the unfolding of
such cycles. This was another striking fact. It resulted from the
works of Benedicks-Carleson, Mora-Viana and Colli. Attractors here
mean invariant sets that attract future orbits of points of a
positive Lebesgue measure set in the phase space (space of events).

In view of all these intricacies inherent to homoclinic and
heteroclinic bifurcations, a new global conjecture has been proposed
in \cite{P1} (see also \cite{P2}) concerning a typical dynamical
system: In particular, systems with finitely many attractors should
be dense in the universe of dynamics, i.e. $C^r$ flows,
diffeomorphisms and maps, with $r\ge 1$. Also, their basins of
attraction should cover the whole phase space, except for a Lebesgue
zero measure set.

\begin{center}
\includegraphics{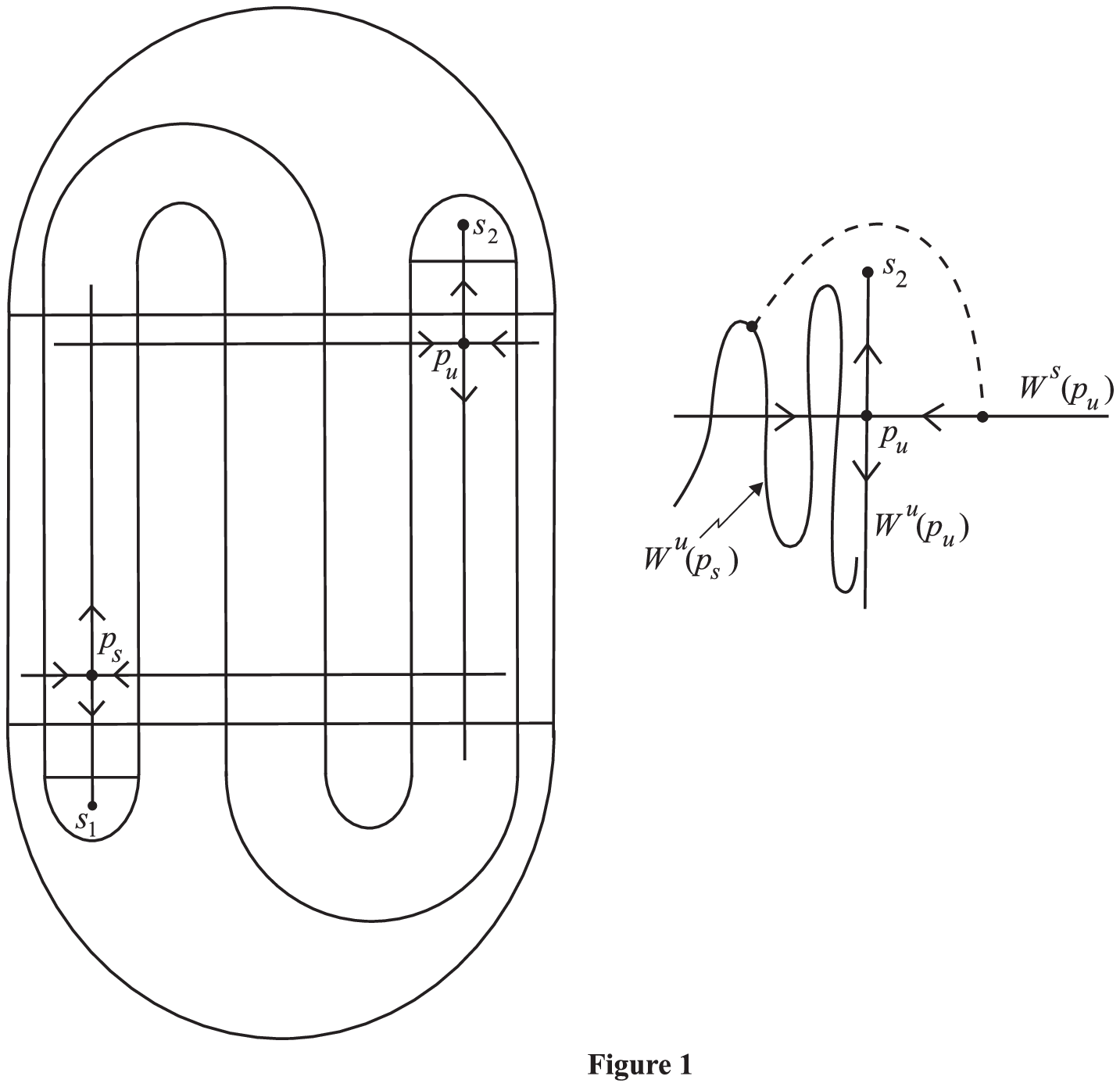}
\end{center}

An example of the creation of a heteroclinic cycle associated to a
(hyperbolic) horseshoe is indicated in Figure~1. Initially we have a
classic Smale's horseshoe map (diffeomorphism) on the two-sphere
$S^2$ with two saddle fixed points $p_s$, $p_u$ with positive
eigenvalues, a fixed point repellor outside the figure and two fixed
points attractors $s_1$ and $s_2$. The rectangle inside the figure
is sent by the map to the snake-shaped piece, while the bigger top
half-disk is sent to the small one around $s_2$ and the lower bigger
half disk is sent to the small one around $s_1$. At the right hand
side of the figure, we show how to move a small neighborhood of a
point in the stable line of $p_u$ so to create a tangency with the
unstable manifold of $p_s$. This is done through a one-parameter
family of diffeomorphisms; until we create such a tangency the
corresponding map remains hyperbolic, i.e. having a hyperbolic limit
set with no cycles among its basic sets.

Concerning Smale's conjecture on the denseness of hyperbolicity, it
is to be remarked that the question is still open for rational maps
of the Riemann sphere. On the other hand, it was known to be true in
the early 60's for flows on orientable surfaces (Peixoto), preceeded
by the case of the two-disk (Andronov-Pontryagin). Much more
recently, it was proved to be true for quadratic maps of the
interval (Graczyk-Swiatek and Lyubich) and then, more generally,
unimodal maps (Kozlovski) and, finally, multimodal maps
(Kozlovski-Shen-Van Strien).

We point out that comprehensive accounts of the results mentioned
above and many other related ones are in \cite{BDV}, \cite{P1} and
\cite{P2}.

The present paper represents a contribution to the understanding of
the dynamics arising from bifurcating a cycle of a $C^\infty$
surface diffeomorphism. The notion of ``most'' used here requires us
to consider one parameter families of diffeomorphisms $g_t$ strictly
containing the initial bifurcating one, say  $g_0$ at parameter
value $t=0$. We assume that $g_t$ is hyperbolic for $t<0$ and $|t|$
small. We suppose that the cycle is formed by a (hyperbolic)
horseshoe $K$ and an orbit of tangency $o(q)$ between stable and
unstable manifolds of different periodic orbits of $K$. We assume
the maximal invariant set in a small neighborhood of $K \cup o (q)$
to consist precisely of $K \cup o (q)$. A main novelty is that we
allow the Hausdorff dimension of $K$ to be larger than one, but not
to far from it. We show that right after the bifurcation, i.e. for
$t>0$ small, most diffeomorphisms still display no attractors nor
repellors in some neighborhood of $K \cup o (q)$. This means that
the parameter values corresponding to diffeomorphisms displaying no
attractors nor repellors should have total density at $t=0$. The
concept is again discussed in the next subsection.

Our results considerably extend those in \cite{PT}, \cite{NP}
obtained for the case when the Hausdorff dimension $HD(K)$ is
smaller than one. They were announced in \cite{PY3}.

Of course, we expect the same to be true for all cases $0 < HD (K) <
2$. For that, it seems to us that our methods need to be
considerably sharpened: we have to study deeper the dynamical
recurrence of points near tangencies of higher order (cubic, quartic
...) between stable and unstable curves. We also hope that the ideas
introduced in the present paper might be useful in broader contexts.
In the horizon lies the famous question whether for the standard
family of area preserving maps, one can find sets of positive
Lebesgue probability in parameter space such that: the corresponding
maps display non-zero Lyapunov exponents in sets of positive
Lebesgue probability in phase space. Finally, we expect our results
to be true in higher dimensions (see \cite{MPV}).

We wish to thank W. de Melo e M. Viana for fruitful conversations.

\subsection{The Setting and a First Formulation of the Main Result\label{sub1.2}}
Let $f$ be a smooth, i.e. $C^\infty$ diffeomorphism of a smooth surface $M$.

Recall that a {\it basic set} is a compact hyperbolic transitive
locally maximal invariant set. A basic set is a {\it horseshoe} if
it is infinite and is neither an attractor nor a repellor.

We assume that there exists a basic set $K$ for $f$, points $p_s$,
$p_u\in K$, $q\in M-K$ such that the following properties hold:

{\bf (H1)}~~$p_s$ and $p_u$ are periodic points and belong to
distinct periodic orbits;

{\bf (H2)}~~$W^s (p_s)$ and $W^u (p_u)$ have a quadratic tangency at
$q$;

{\bf (H3)}~~there exists a neighbourhood $U$ of $K$, a neighbourhood
$V$ of the orbit $\cO (q)$ of $q$, such that $K\cup \cO (q)$ is the
maximal invariant set in $U\cup V$.

We would like to understand, when $U$, $V$ are appropriately small
and $g$ is $C^\infty$ close to $f$, the maximal invariant set
\begin{equation}\label{eq1.1}
\La_g = \bigcap_{\Bbz} g^{-n} (U\cup V).
\end{equation}
Observe that the smaller set
\begin{equation}\label{eq1.2}
K_g = \bigcap_{\Bbz} g^{-n} (U)
\end{equation}
is a horseshoe which is the hyperbolic continuation of $K$.

Let $\cU$ be an appropriately small neighbourhood of $f$ in
$\Diff^\infty (M)$. We still denote by $p_s$, $p_u$ the continuation
of these hyperbolic periodic points in $\cU$. The condition that
$W^s (p_s)$, $W^u (p_u)$ have a quadratic tangency near $q$ defines
a codimension 1 hypersurface $\cU_0$ through $f$ in $\cU$. It
divides $\cU$ into regions $\cU_+$, $\cU_-$ such that, for $g\in
\cU_-$, $W^s (p_s)$ and $W^u (p_u)$ do not intersect near $q$ while,
for $g\in \cU_+$, $W^s (p_s)$ and $W^u (p_u)$ have two transverse
intersection points near $q$ (for obvious dynamical reasons, the
intersection is actually infinite in this case; we are really
considering here the intersection derived from the continuation of
large {\it compact} curves contained in $W^s (p_s)$ and $W^u (p_u)$.

When $g\in \cU_-$, we clearly have
\begin{equation}\label{eq1.3}
\La_g = K_g
\end{equation}
When $g\in \cU_0$, we have
\begin{equation}\label{eq1.4}
\La_g = K_g \cup \cO (q_g),
\end{equation}
where $q_g$ is the tangency point close to $q$ given by the
definition of $\cU_0$.

The interesting case is therefore $g\in \cU_+$.

It is actually not realistic to try to understand $\La_g$ for all
$g\in \cU_+$. One of the reasons is the so-called Newhouse's
phenomenon \cite{N}: there exists an open set $\cN \subset \cU_+$,
with $\cU_0 \subset {\ov{\cN}}$, such that, residually in $\cN$,
$\La_g$ has infinitely many periodic sinks or sources and so its
full dynamical description appears to be beyond reach. See also
\cite{BC}, \cite{MV}, \cite{C} for similar results involving
Hénon-like attractors.

Still, we can and shall consider most $g\in \cU_+$ in the following
sense.

We will say that a subset $\cP\subset \cU_+$ contains most $g\in
\cU_+$ if, for any smooth 1-parameter family
$(g_t)_{t\in(-t_0,t_0)}$ which is transverse to $\cU_0$ at $t=0$
(with $g_t\in \cU_+$ for $t>0$), we have
\begin{equation}\label{eq1.5}
\lim_{t\to 0} \frac{1}{t} \, \text{Leb}\, (s\in (0,t], g_s\in \cP) =
1.
\end{equation}
Denote by $W^s (K)$ (resp.~$W^u (K)$) the stable set (resp.~unstable
set) of $K$ for $f$. This is a partial foliation with a $C^{1+\al}$
Cantor transverse structure; denote by $d^0_s$ (resp.~$d^0_u$) the
transverse Hausdorff dimension of $W^s(K)$ (resp.~$W^u (K)$). The
Hausdorff dimension of $K$ is equal to $d^0_s + d^0_u$. We then
have, in some contrast to Newhouse's phenomenon:

\medskip
{\bf Theorem.~} {\it \cite{PT}, \cite{NP} Assume that $d^0_s + d^0_u
< 1$. Then, for most $g\in \cU_+$, $\La_g$ is a horseshoe.}

On the other hand, by \cite{PY1}, the same conclusion does not hold
when $d^0_s + d^0_u > 1$. The paper \cite{MY} gives substantially
more geometric information in this case, specially concerning
tangencies between stable and unstable manifolds (lines) in the
hyperbolic continuation $K_g$ of $K$. These results have been
extended to higher dimensions, as announced in \cite{MPV} and
complete proofs to appear in the near future.

In the present work, we investigate the maximal invariant set
$\La_g$, for most $g\in \cU_+$, provided that the dimensions
$d^0_s$, $d^0_u$ satisfy (see figure~2)
\begin{alignat}{1}
(d^0_s + d^0_u)^2 + (\Max (d^0_s, d^0_u))^2 < d^0_s + d^0_u + \Max
\, d^0_s, d^0_u.\tag{\bf H4}
\end{alignat}

Our results can essentially be summarized as:
\medskip

{\bf Main Theorem.~} {\it Assume that (H1), (H2), (H3), (H4) hold.
Then, for most $g\in \cU_+$, $\La_g$ is a non-uniformly hyperbolic
horseshoe.}

The meaning of a non-uniformly hyperbolic horseshoe in the present
context will be explained somewhat in the next section and more
completely in the rest of the paper. We can, however, comment that,
for most $g\in \cU_+$, $\La_g$ will be a saddle-like object in the
sense that both the stable set $W^s (\La_g)$ and the unstable set
$W^u (\La_g)$ have Lebesgue measure zero and, so, it carries no
attractors nor repellors. It will be (non-uniformly) hyperbolic in
the sense that we will construct geometric invariant measures, à la
Sinai-Ruelle-Bowen \cite{Si,Ru,BR}, on $\La_g \subset W^s (\La_g)$
and $\La_g \subset W^u (\La_g)$ with non-zero Lyapunov exponents.
Such properties of the invariant set $\La_g$ are made especially
precise in Section~\ref{sec10} and \ref{sec11}, the last ones in the
paper. They yield some rephrasing of the main result in these terms,
which is presented at the end of Section~\ref{sec11}.

\begin{center}
\includegraphics{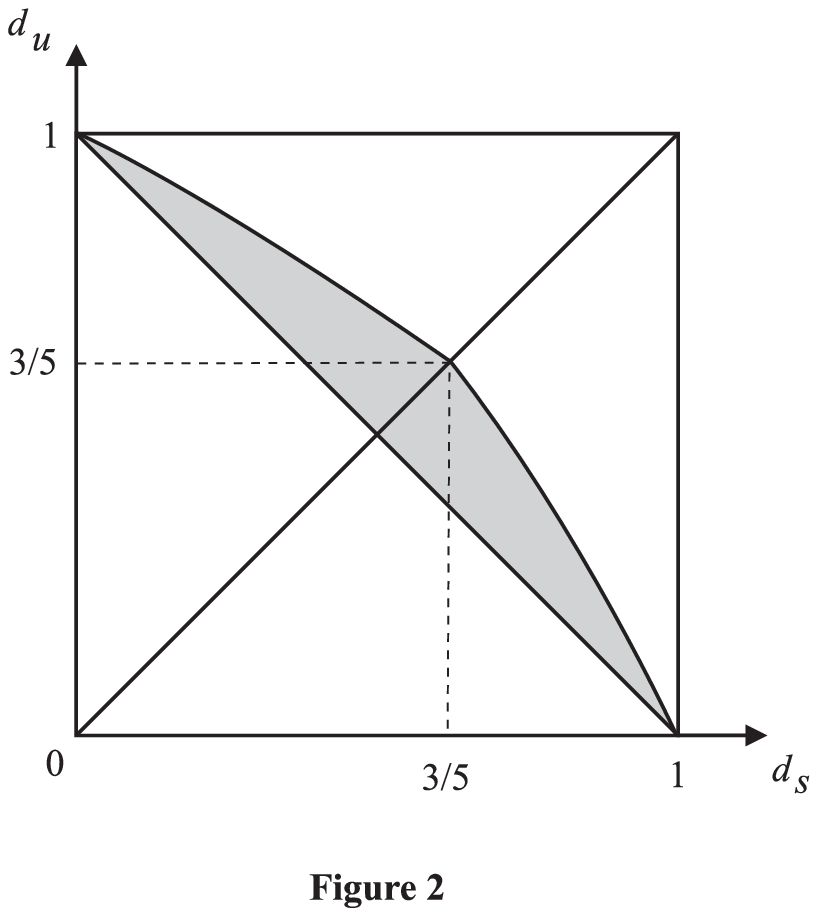}
\end{center}

\begin{rem}
In the case when $d^0_s + d^0_u < 1$, mentioned above and studied in
\cite{PT}, \cite{NP}, it is not necessary to assume that $p_s$,
$p_u$ belong to {\it distinct} periodic orbits. It is probably not
necessary in our case either, at least as far as the qualitative
statements are concerned. But, this assumption seems to make the
technicalities significantly easier in what is already a very long
construction.
\end{rem}

\subsection{A Summary of the Next Sections of the Paper\label{sub1.3}}
Sections 2--4 consist mainly of preparatory work.

In Section~2, we introduce a Markov partition by smooth disjoint
rectangles $(R_a)_{a\in\AL}$ for the horseshoe $K$. The dynamics in
the neighbourhood $U$ of $K$ is given by the transition maps from
one rectangle to another, which enjoy a nice hyperbolic behaviour.
To understand the dynamics in the larger set $U\cup V$, we need to
control the dynamics along a finite part of the orbit of $q$,
stretching from the moment this orbit goes out of $R := \sqcup R_a$
until it comes back to $R$. The region of exit of $R$ and the region
of entry into $R$ are two parabolic tongues $L_u$ and $L_s$, and the
transition map
\[
G = g^{N_0}\; : \; L_u \to L_s
\]
is a folding map which share many features with the Henon quadratic
polynomial diffeomorphisms of the plane.

Section~3 is essentially a summary of our previous work \cite{PY2},
(which was written having the present paper in mind). The important
concept of affine-like map is introduced. The basic idea, which goes
back to the early stages of the hyperbolic theory, is to describe
maps that present hyperbolic features in an implicit way exhibiting
preference for coordinates with a macroscopic range. Concretely, if
a two-dimensional diffeomorphism contracts the vertical coordinate
$y$ and expands the horizontal coordinate $x$, we use $y_0$ and
$x_1$ as independent variables associated with a point $(x_0, y_0)$
and its image $(x_1, y_1)$, writing $x_0$ and $y_1$ as functions of
$y_0$ and $x_1$.

Cone conditions are easy to formulate in this setting. A nice
feature of this implicit representation of the dynamics is that it
is time-symmetric: the map and its inverse satisfy symmetric
formulas. Another even more important feature is that this formalism
is well-adapted to the right concept of distortion (for
2-dimensional maps), yielding appropriate control on the partial
derivatives of order two.

Composition of two affine-like maps which satisfy the same cone
condition is also affine-like, and the distortion is only slightly
bigger than the distortion of the two maps.
Besides this ``simple" composition, we study ``parabolic"
compositions of the form $F_1 \circ G  \circ F_0$, where $F_0$,
$F_1$ are affine-like and $G$ is the folding map of Section~2. When
the relative positions of the parabolic strip $G (Q_0)$ (where $Q_0$
is the image of $F_0$) and $P_1$ (the domain of $F_1$) are
appropriate, the domain of $F_1 \circ G  \circ F_0$ has two
connected components and the restrictions $F^\pm$ of $F_1 \circ G
\circ F_0$ to each component is affine-like. A control of the
distortion of $F^+$ and $F^-$ is also obtained.

In Section~4, the general structure of parameter space is
introduced. The parameter coordinate is normalized by the relative
speed at the quadratic tangency of the tips of the stable and
unstable manifolds. Then, with $\vep_0$ very small, the starting
interval $I_0 := [\vep_0, 2\vep_0]$ for the parameter selection
process is introduced. A small parameter $\tau$ (with $\tau \ll 1$
but still $\vep_0^\tau \ll 1$) determines a sequence of scales
$(\vep_k)_{k\geq 0}$ in parameter space through the formula
$\vep_{k+1} = \vep_k^{1+\tau}$. At level $k$, we have disjoint
parameter intervals of length $\vep_k$ (starting from level 0 with
$I_0$). Each parameter interval of level $k$ that has been selected
is divided into $[\vep_k^{-\tau}]$ disjoint candidates of length
$\vep_{k+1}$. These candidates will pass a test to decide whether
they are selected at level $k+1$.

The test takes two forms. First, in Section~5, a property of the
parameter interval called regularity (see below) will be introduced;
candidates which do not possess this property are discarded. Such a
property is sufficient to develop in Sections~5--8 some basic
combinatorial and quantitative properties, but it is not
well-adapted to an inductive scheme. Hence, in Section~9, a stronger
property called strong regularity is introduced, and candidates have
to satisfy this property in order to be selected.

\pagebreak Sections~5--7 constitute in some sense a single logical
step: in Section~5, certain classes of restrictions of iterates of
$g_t$ are inductively defined, and the definition is only possible
because of properties that are inductively proved in Sections~6 and
7.

In Section~5, the goal is to define, for each parameter interval $I$
which is a candidate (i.e. its parent interval at the immediately
upper level has been tested as regular), a class $R (I)$ of
$I$-persistent affine-like iterates. An $I$-persistent affine-like
iterate is a triple $(P, Q, n)$ where $P$ is a vertical-like strip
in some rectangle $R_a$ depending on $t\in I$, $Q$ is a
horizontal-like strip in some rectangle $R_{a'}$, depending on $t\in
I$, and the restriction  of $g_t^n$ to $P$ is a diffeomorphism onto
$Q$ which is affine-like.

However, we do not want to have in ${\cal R}(I)$ all $I$-persistent
affine-like iterates: we will argue about them by induction (on $n$,
for instance) and in order to do this, we want to obtain them in
some explicit constructive way. Therefore, a number of Axioms,
(R1)--(R7), are introduced and together they completely determine
the class ${\cal R}(I)$. The most important feature of these Axioms is the
following: every element of ${\cal R}(I)$ consisting of more than one
iteration of $g_t$ can be obtained from simpler elements of ${\cal R}(I)$
by simple or parabolic composition; in this context, the notions of
parent and simple or non-simple child introduced here, play a
relevant role; simple composition is allowed in ${\cal R}(I)$ whenever it
makes sense; and parabolic compositions of elements of ${\cal R}(I)$ is
allowed if and only if a certain {\it transversality relation} is
satisfied.

Thus, the definition of ${\cal R}(I)$ is reduced to the definition of this
transversality relation, which is presented in Subsection~5.4. The
intuitive notion behind the formal definition is the following: an
element $(P, Q, n)$ should be $I$-transverse to an element $(P', Q',
n')$ if the distance $\de (Q, P')$ between the tip of the
parabolic-like strip $G(Q)$ and $P'$ satisfies
\[
\de (Q, P') \geq |Q|^{1-\eta} + |P'|^{1-\eta}
\]
for all $t\in I$, where $|Q|$, $|P'|$ are the widths of the strips
$Q$ and $P'$. Here $\eta$ is a small positive constant, fixed once and for all.
However, a number of properties, presented in
Section~6, are very helpful, and they require a definition of the
transversality relation that may seem quite complicated. In
Appendix~C, we explain why this seemed complication is rather
necessary.

For the starting interval $I_0$, it follows from the formal
definition that the transversality relation is never satisfied;
therefore, parabolic composition is not allowed and the class ${\cal R}
(I_0)$ is exactly the one associated with the symbolic dynamics
given by the Markov partition. We conclude Section~5 with the
definition of {\it regularity}. First, one says that a strip $P$
(from an element $(P, Q, n)\in {\cal R}(I))$ is {\it $I$-transverse} if one
can find finitely many $Q_\al$, whose union contains the unstable
set of $\La$, such that $Q_\al$ and $P$ are $I$-transverse for every
$\al$; otherwise one says that $P$ is $I$-critical. Then, given a
constant $\be > 1$, one says that the parameter interval is
$\be$-regular if any $(P, Q, n)\in {\cal R}(I)$ such that both $P$ and $Q$
are $I$-critical satisfies $|P| < |I|^\be$, $|Q| < |I|^\be$ for all
$t\in I$. Intuitively, this means that no short return to the
critical set is allowed.

In Section~6, we prove a number of properties of the transversality
relation and the classes ${\cal R}(I)$. Amongst the most important is the
following: children are born from their parents. Let us explain what
it means. Let $(P, Q, n)\in {\cal R}(I)$, and let $(\wtP, \wtQ, \wtn)$ the
element of ${\cal R}(I)$ such that $P\subset \wtP$, $P\not= \wtP$ and
$\wtP$ is the thinnest rectangle with this property; one says that
$P$ is a child of $\wtP$ and that $\wtP$ is the parent of $P$. There
are two cases; either $n = \wtn + 1$ and one says that $P$ is a
simple child; $(P, Q, n)$ is obtained by simple composition of
$(\wtP, \wtQ, \wtn)$ with an element of length~1; or $n > \wtn + 1$
and one says that $P$ is a non-simple child; one then proves that
$(P, Q, n)$ is obtained by parabolic composition of $(\wtP, \wtQ,
\wtn)$ with some element $(P_1, Q_1, n_1)$.

The most important result in Section~6 is a structure theorem for
new rectangles in Subsection~6.7. One considers an element $(P, Q,
n)$ which belongs to ${\cal R}(I)$ but not to ${\cal R}(\wtI)$, where $\wtI$ is
the parameter interval containing $I$ of the level immediately
inferior (one says that $\wtI$ is the parent of $I$). Then there is
a unique way to write $(P, Q, n)$ as the result of a sequence of $k$
parabolic compositions, possible in ${\cal R}(I)$ but not in ${\cal R}(\wtI)$, of
elements $(P_0, Q_0, n_0),\cdots,(P_k, Q_k, n_k)$. This fundamental
result has several useful corollaries.

We have grouped in Section~7 a number of calculations and estimates
related to the definitions of Section~5, and that are necessary for
the inductive construction of the classes ${\cal R}(I)$. The first
result relates the length $n$ of an element $(P, Q, n)$ to the width
$|P|$, $|Q|$ of the strips. While there is no uniform exponential
estimate as in the uniformly hyperbolic case, we are still able to
prove a stretched-exponential uniform estimate.

We prove next that a uniform cone condition is satisfied by the
affine-like iterates that we consider, and also that they have
uniformly bounded distortion. After a technical estimate related to
parabolic composition, we deal in Subsection~7.6 with the relative
speed of the strips when the parameter varies; this is clearly of
capital importance if we are to succeed. A point which is worth
mentioning is that we are not able to obtain estimates for all pairs
of strips (actually, it is easy to see that such estimates do not
exist); we have to restrict ourselves to strips satisfying certain
criticality conditions, that, fortunately, will be satisfied every
time we need some information on these speeds. In the last
Subsection~7.7, we investigate the oscillation of the widths of the
strips with the parameter. While it is just not true that the
relative oscillation is bounded (in the sense that the maximum over
a parameter interval is no greater than a constant times the
minimum), the result that we get will allow us to argue as if it
was.

At the end of Section~7, the construction of the classes ${\cal R}(I)$ is
complete, for every parameter interval $I$ whose parent $\wtI$ is
regular. But we still don't know whether a single interval $I$ is
regular.

In Section~8, we develop several quantitative estimates that will
turn out to be crucial both in the parameter selection process of
Section~9 and in the analysis of the dynamics for strongly regular
parameters in Sections~10 and 11. We first investigate, for a given
element $(\wtP, \wtQ, \wtn)$, the number of elements $(P, Q, n)$
such that $P$ is a non-simple child of $\wtP$; we show that, for
every $\vep > 0$, there are at most $\vep^{-C\eta}$ such non-simple
children with width $|P|$ larger than $\vep |\wtP|$. The constant
$\eta$ here  is small and related to the definition of the
transversality relation. The meaning of this estimate is that the
presence of non-simple children is not too significant from the
point of view of Hausdorff (or box) dimension, as it is made clear
in Subsection~8.2. In Subsection~8.3, we transfer this information
to parameter space, combining it with the result on relative speed
of strips in Subsection~7.6.

Section~9 is the longest one in the paper and deals with the
parameter selection process. The concept of regularity is very
useful to develop a number of properties of the classes $R(I)$, as
we did in Sections~5--8. Unfortunately, we are not able to prove
(and it is probably false) that, given a  $\be$-regular interval
$\wtI$, most candidates $I\subset \wtI$ at the next level are
$\be$-regular. [It is a consequence of the structure theorem of
Section~6 that all candidates are $\ovbe$-regular, where $\ovbe =
\be (1+\tau)^{-1}$ is very close to $\be$; this allow us to obtain
all qualitative consequences of regularity for all candidates; but
obviously we cannot repeat this at many successive levels of
parameter intervals, because we need to keep $\be > 1$.]  The
problem with the concept of regularity is that it is dealing with
only one scale $|\wtI|^\be$; it could happen a priori that for a
regular parameter interval $\wtI$ we have many $\wtI$-bicritical
$(P, Q, n)\in {\cal R}(I)$ with $|P|$ or $|Q|$ only slightly below the
threshold $|\wtI|^\be$ (and therefore above the next threshold
$|I|^\be$ for candidates $I\subset \wtI$); for each such $(P, Q,
n)$, we have to eliminate candidates $I$ such that $(P, Q, n)$ is
$I$-bicritical, and no candidate will survive this selection process
if there are too many $(P, Q, n)$.

The solution to this difficulty is to introduce the condition of
strong regularity, which implies regularity and gives a quantitative
control at all scales. Actually, the strong regularity condition
involves two parts. In the first, one controls the size of the
critical locus (in several slightly different ways) by a series of
eight inequalities which all amount to say that the "dimension'' of
the critical locus is not much larger than $d^0_s + d^0_u - 1$. In
this case, the parameter selection process is based on the result mentioned above in the last part of Section~8.
The second part of the strong regularity condition, by far the most
subtle one, is a quantitative estimate for the number of bicritical
elements at all scales. Because of the inductive nature of the
argument, which relies in an essential way on the structure theorem
of Section~6, we need to control the number of elements $(P, Q,
n)\in {\cal R} (I)$ such that $P$ is $I_\al$-critical, $Q$ is $I_\om$-critical and $|P|\geqslant x$ for some $t\in I$. Here,
$I_\al$ and $I_\om$ are parameter intervals containing $I$, and the
control will depend on $I_\al$, $I_\om$ and $x$. The formulas in
Subsection~9.4 present a phase transition with respect to the width
parameter $x$. Discussing this phase transition leads naturally to
the hypothesis (H4) on the transverse dimensions $d^0_s$, $d^0_u$: a
small calculation shows that (H4) is exactly what one needs to
obtain $\be$-regularity with $\be > 1$.

Having stated the strong regularity condition, the goal in the rest
of Section~9 is to prove that, given a strongly regular parameter
interval $\wtI$, most candidates $I \subset \wtI$ at the next level
are also strongly regular (the proportion of failed candidates turns
out to be not larger than $C |I|^{\tau^2}$). This requires the control
of two things. First, is to bound the number of ``new" bicritical
elements $(P, Q, n)\in {\cal R}(I)$ which did not belong to ${\cal R}(\wtI)$; this
is based on the structure theorem of Section~6 and leads to a long
but straightforward calculation. Second, is to estimate which
proportion of bicritical elements for $\wtI$ are still bicritical
for $I$; this is only necessary when $I_\al$ or $I_\om$ is equal to
$I$; when only one of the two intervals $I_\al$, $I_\om$ is equal to
$I$, the idea is simply to estimate what is the mean proportion
(over all candidates), and to discard candidates for which the
proportion is much above the mean. To compute the mean proportion,
we rely again on the result at the end of Section~8. The case where
$I = I_\al = I_\om$ is the most important and the most difficult.
When $x$ is ``large", the same argument than when $I = I_\al \not=
 I_\om$ still applies; but when $x$ is ``small", the phase
 transition of the estimate means than the argument is not
 sufficient any more. A more complicated strategy is required, which
 is explained in Subsection~9.8 and carried out in 9.9 through 9.13.

 It is worth mentioning that up to the end of Section~9, we never
 consider the dynamics for a single parameter, only for parameter
 intervals. In the last two sections, we study the dynamics for a
 strongly regular parameter value, i.e. the intersection of a
 decreasing sequence $(I_m)$ of strongly regular parameter
 intervals.

 In Section~10, we study the dynamics on the set of stable curves. A
 {\it stable curve} $\om$ is the decreasing intersection of a
 sequence of vertical-like strips $P_k$, where $(P_k, Q_k, n_k)\in {\cal R}
 = \bigcup\limits_m {\cal R} (I_m)$. The set of stable curves is denoted by
 $\cR^\infty_+$, their union by $\wtcR_+^\infty$. In order to define a map on
 $\wtcR_+^\infty$ (which is not invariant under $g$), we introduce the
 concept of prime element in $R$, i.e. one which cannot be written
 as the simple composition of shorter elements. Let then $\om$ be a
 stable curve which is not contained in infinitely many prime
 elements $P_k$, and let $(P, Q, n)$ be such that $P$ is the
 thinnest prime element containing $\om$. The image $g^n (\om)$ is
 contained in a stable curve $\om'$ and we set $T^+(\om) = \om'$,
 $\wttm/\om = g^n/\om$. This defines a map $T^+$ from a subset
 $\cD_+$ of $\cR^\infty_+$ onto $\cR^\infty_+$ which lifts to a map $\wttm$ from the
 union $\wtcD_+$ of curves in $\cD_+$ to $\wtcR_+^\infty$.

 The map $T^+$ is Bernoulli in the following sense: its domain
 $\cD_+$ splits into countably many pieces $\cR^\infty_+ (P)$  indexed by
 prime elements, and each piece is sent homeomorphically by $T^+$
 onto the intersection of $\cR^\infty_+$ with some rectangle $R_a$ of the
 Markov partition.

 The map $T^+$ is uniformly expanding (with countably many branches)
 and we introduce a one parameter family of weighted transfer
 operators in the spirit of classical uniformly hyperbolic maps. One
 has only to be careful because the presence of countably many
 branches is the source of some problems, which are dealt with in
 Subsection~10.3 using the estimates of Section~8 on the number of
 children.

 As expected, the transfer operators $L_d$, considered in the
 appropriate function space, turn out to have a positive
 eigenfunction $h_d$ associated with a dominant eingevalue
 $\la_d > 0$. There is a unique value $d_s$ such that $\la_{d_s} =
 1$. This value turns out to be, unsurprisingly, the transverse
 Hausdorff dimension of the partial foliation $\wtcR_+^\infty$ (which is
 proved in Subsection~10.5 to be transversally Lipschitzian). The
 transfer operator also allows us to identify, as usual, a measure
 $\mu_d$ with prescribed Jacobian and an invariant measure $\nu_d =
 h_d \mu_d$. For $d = d_s$, the $\mu_d$-measure (or $\nu_d$-measure)
 of the set of stable curves contained in any vertical-like strip
 $P$ is proportional to $|P|^{ds}$.

 The set $\wtcR_+^\infty - \wtcdm$ where $\wttm$ is not defined, has
 transverse dimension smaller than $d_s$, hence is negligible in a geometrical
 sense. One can lift the $T^+$-invariant measure $\nu = \nu_{d_s}$
 to a $\wttm$-invariant measure $\wtnu$ which is ergodic and then
 spread it to a $g$-invariant measure on $\La$.

 In Section~11, the last in the paper, we pursue the study of the
 dynamics of $g_t$ on $\La = \La _{g_t}$ for a strongly regular parameter $t$,
 looking now beyond the well-behaved set $\wtcR_+^\infty$ which was studied
 in Section~10. In the first part (Subsections~11.1--11.5), we study the
 intersection of the invariant set $\La$ with an unstable curve $\om^*$
 (defined as a stable curve, exchanging $P$'s and $Q$'s). The main part
 of this intersection is a countable disjoint union of dynamical copies
 of the set $\cR^\infty_+$ studied in Section~10. There are also at most
 countably many critical points, corresponding to quadratic
 tangencies between stable curves and images under $G$ of unstable
 curves. And, finally, there is an exceptional set (formed by points
 which come very close to the critical locus infinitely many times);
 but this exceptional set is small; its Hausdorff dimension is
 explicitly controlled by a value much smaller than the dimension
 $d_s$ of $\om^* \cap \La$.

 In the second part of Section~11, we prove that the invariant set
 $\La$ is a saddle-like object in the metric sense: both its stable
 set $W^s (\La)$ and its unstable set $W^u(\La)$ have Lebesgue
 measure 0. So, no attractors nor repellors are present on $\La$.
 One actually expects more: certainly the Hausdorff dimension of
 $W^s(\La)$ should be strictly less than~2, probably it is close to
 $1+ d_s$, and perhaps even equal to $1 + d_s$. However, we stick
 to the simpler, but still very meaningful result: it implies that
 $\La_g$ carries no attractor nor a repellor for most $g$.

 One has a nice combinatorial decomposition of the restricted stable set
 $W^s (\La, R)$, but to compute Lebesgue measure (or Hausdorff
 dimension), one has to transport the pieces of this decomposition
 by affine-like iterates of $g$ of high order. This is easy to do as
 far as Lebesgue measure is concerned, because bounded distortion of
 affine-like maps mean also bounded distortion of measure (bounded
 relative oscillation of Jacobians). This is much more delicate with
 respect to Hausdorff dimension: the geometry of the pieces after
 iteration can get very distorted.

 In Appendix~A, we recall all formulas related to the implicit
 representation of affine-like maps; many of them can already be
 found in \cite{PY2}, but we have also to consider the
 derivatives with respect to parameter, a setting which was not
 considered in \cite{PY2}.

 In Appendix~B, we perform some calculations related to proposition~40
 in Subsection~10.5, which generates the transversally Lipschitz
 regularity of the partial foliation $\wtcR_+^\infty$.

 In Appendix~C, we give some justification for what seems to be a
 convoluted definition of the transversality relation in
 Subsection~5.4.
\newpage

\setcounter{section}{1}
\setcounter{equation}{0}

\section{Markov Partition and Folding Map\label{sec2}}
\subsection{Markov Partition and Related Charts\label{sub2.1}}
We will choose once and for all a finite system of smooth charts
\[
I^s_a \times I^u_a \; \overset{\thickapprox}{\longrightarrow}\; R_a
\subset M,\qquad a\in \AL
\]
indexed by a finite alphabet $\AL$. Each chart depends smoothly on
$g\in \cU$; the intervals $I^s_a$, $I^u_a$ are compact; the
rectangles $R_a$ are disjoint.

Let $R = \bigcup\limits_{\AL} R_a$. We choose the charts in order to
have:

{\bf (MP1)}~~for each $g\in\cU$, $K_g$ is the maximal invariant set
in int~$R$; for each $g\in\cU$, $a\in\AL$, one has
\begin{eqnarray}
&g (\pa  I^s_a \times I^u_a) \cap R = \emptyset,\label{eq2.1}\\
&g^{-1} (I^s_a \times \pa I^u_a) \cap R = \emptyset;\label{eq2.2}
\end{eqnarray}

{\bf (MP2)}~~for each $g\in\cU$, the family $(R_a \cap
K_g)_{a\in\AL}$ induces a Markov partition for the horseshoe $K_g$.
Let
\begin{equation}
\cB = \{ (a, a')\in \AL^2, \;\; f (R_a) \cap R_{a'} \not=
\emptyset\}.\label{eq2.3}
\end{equation}
The Markov partition provides a coding which is a topological
conjugacy between the horseshoe $K_g$ and the subshift of finite
type of $\AL^\Bbz$ defined by $\cB$.

\subsection{The Parabolic Tongues $L_u$, $L_s$\label{sub2.2}}
Denote by $a_s$, $a_u\in\AL$ the letters such that $p_s\in R_{a_s}$,
$p_u\in R_{a_u}$. We choose the corresponding charts in order to
have:

{\bf (MP3)}~~for each $g\in \cU$, the equation of the local stable
manifold $W^s_{\loc} (p_s)$ is $\{ x_{a_s} = 0\}$, the equation of
the local unstable manifold $W^u_{\loc} (p_u)$ is $\{ y_{a_u} =
0\}$.

We have written $x_a$ (resp.~$y_a$) for the coordinate in $I^s_a$
(resp.~$I^u_a$). We also choose the rectangles $R_a$ in order to
have, for some integer $N_0\geqslant 2$:

{\bf (MP4)}~~for each $g\in \cU_0$, there are points $q_s$, $q_u$ in
the orbit of $q$ such that

--~~for $n\ge 0$, $g^n (q_s)$ and $g^n (p_s)$ belong to the interior
of the same rectangle;

--~~for $n\le 0$, $g^n (q_u)$ and $g^n (p_u)$ belong to the interior
of the same rectangle;

--~~$q_s = g^{N_0} (q_u)$ and $g^i (q_u)$ does not belong to $R$ for
$0 < i < N_0$.

Consider small pieces of $W^s (p_s)$, $W^u (p_u)$ which are tangent
at $q_u$ for $g\in \cU_0$. When $g\in \cU_+$, these pieces will meet
in two points and bound a compact lenticular region $L_u \subset
\inta R_{a_u}$. Taking the image under $g^{N_0}$, we get another
lenticular region $L_s \subset \inta R_{a_s}$. These regions are
called {\it parabolic tongues}. See figure~3.

\begin{center}
\includegraphics{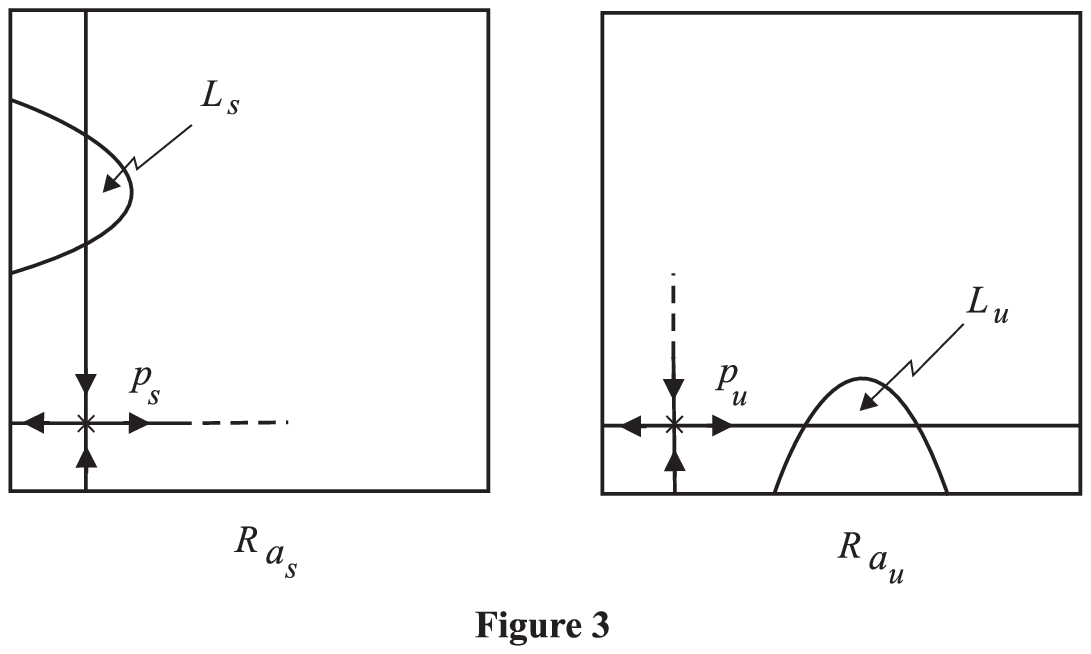}
\end{center}

Define then, for $g\in \cU_+$
\begin{equation}\label{eq2.4}
\whR = R \bigcup_{0<i<N_0} g^i (L_u)
\end{equation}
The maximal invariant set we are interested in is
\begin{equation}\label{eq2.5}
\La_g = \bigcap_\Bbz g^{-n} (\whR).
\end{equation}
We also define
\begin{eqnarray}
W^s (\La_g, \whR) &= &\bigcap_{n\geqslant 0} g^{-n}
(\whR),\label{eq2.6}\\
W^u (\La_g, \whR) &= &\bigcap_{n\leqslant 0} g^{-n}
(\whR)\label{eq2.7}
\end{eqnarray}
The dynamics in $\whR$ are generated by

--~~the transition maps related to the Markov partition:
\[
g: R_a \cap g^{-1} (R_{a'}) \to g (R_a) \cap R_{a'}, for \; (a, a')
\in \cB;
\]
--~~the folding map $G := g^{N_0}$ from $L_u$ onto $L_s$.

\subsection{The Folding Map $G$\label{sub2.3}}
For simplicity, we write $(x_s, y_s)$, $(x_u, y_u)$ for the
coordinates in $R_{a_s} \supset L_s$, $R_{a_u} \supset L_u$.

The folding map $G$ is {\it parabolic} in the sense of \cite{PY2};
let us recall this definition.

Consider the graph $\Ga_G$ of the restriction $G$ of $g^{N_0}$ to
the component of $R_{a_u} \cap g^{-N_0} (R_{a_s})$ which contains
$L_u$ (for $g\in \cU_+$; we then follow this component in the rest
of $\cU$). Using the corresponding charts, we can view $\Ga_G$ as a
surface in $I^s_{a_u} \times I^u_{a_u} \times  I^s_{a_s} \times
I^u_{a_s}$. Denote by $\pi$ the projection from $I^s_{a_u} \times
I^u_{a_u} \times  I^s_{a_s} \times  I^u_{a_s}$ onto $I^u_{a_u}
\times  I^s_{a_s}$.  For $\cU$ small enough, from the quadratic
tangency at $q$ and (MP3) we deduce that:

{\bf (P1)}~~the restriction of $\pi$ to $\Ga_G$ is a fold map (in
the sense of singulary theory).

Denote by $\Ga_0 \subset I^u_{a_u} \times I^s_{a_s}$ the smooth
curve which is the image of the critical locus of this fold map. It
divides $I^u_{a_u}\times I^s_{a_s}$ into two regions $\Ga_+$,
$\Ga_-$ such that $\Ga_+ \cup \Ga_0$ is the image of the fold map.
We can reformulate (P1) as:

{\bf (P$'$1)} (i)~~for $(y^0, x^0)\in \Ga_0$, the image $G (\{ y_u =
y^0\})$ meets $\{ x_s = x^0\}$ in a single point, interior to both
curves, at which the curves have a quadratic tangency;

(ii)~~for $(y^0, x^0) \in \Ga_-$, the curves $G (\{ y_u = y^0\})$
and $\{ x_s = x^0\}$ do not intersect;

(iii)~~for $(y^0, x^0) \in \Ga_+$, the curves $G ( \{ y_u = y_0\})$
and $\{ x_s = x^0\}$ intersect transversally in two points.

As $G$ is a diffeomorphism, the tangents to $\Ga_0$ are never
vertical or horizontal. Therefore, we can and will choose a smooth
function $\te$ on $I^u_{a_u} \times I^s_{a_s}$ such that

{\bf (P2)}~~$\te \equiv 0$  on  $\Ga_0$, $\te > 0$ on $\Ga_+$, $\te
< 0$ on $\Ga_-$;

{\bf (P3)}~~the partial derivatives $\te_y$, $\te_x$ of $\te$ do not
vanish on $I^u_{a_u}\times I^s_{a_s}$.

\begin{rem}
The choice of $\te$ is far from unique. One could for instance
choose $\te$ of the form
\begin{equation}\label{eq2.8}
\te (y_u, x_s) = \vep_u y_u + \vep_s \chi (x_s),
\end{equation}
with $\vep_s$, $\vep_u \in \{ -1, +1\}$ and $\chi$ monotone
increasing. We prefer not to specify a particular choice in order to
keep a time-symmetric setting between positive and negative
iterations.
\end{rem}
From $\te$, we define a smooth function $w$ on $\Ga_G$ by

{\bf (P4)}~~$w^2 = \te \circ \pi$

(there are two choices for $w$; the other is $-w$).

Then, from (P3) we obtain smooth maps $Y_u$, $X_s$ implicitly
defined by
\begin{alignat}{1}
w^2 &= \;\te (Y_u (w, x_s), x_s)\tag{\bf P5}\\
    &= \;\te (y_u, X_s (w, y_u)\notag
\end{alignat}
On the graph $\Ga_G$, we can use either $(x_u, y_u)$ or $(x_s, y_s)$
or $(w, y_u)$ or $(x_s, w)$ as coordinates; therefore we can
factorize $G$ as $G_+ \circ G_0 \circ G_-$:
\begin{equation}
(x_u, y_u) \overset{G_-}{\longrightarrow} (w, y_u)
\overset{G_0}{\longrightarrow} (x_s, w)
\overset{G_+}{\longrightarrow} (x_s, y_s)\label{eq2.9}
\end{equation}
with
\begin{alignat}{3}
&G_0 (w, y_u)      \;&=  &\;(X_s (w,y_u), w),\tag{\bf P6}\\
&G_0^{-1} (x_s, w) \;&=  &\;(w, Y_u (w, x_s)),\notag\\
&G_+ (x_s, w)      \;&=  &\;(x_s, Y_s (w, x_s)),\notag\\
&G_-^{-1} (w, y_u) \;&=  &\;(X_u (w,y_u), y_u).\notag
\end{alignat}
The last two formulas define smooth maps $Y_s$, $X_u$ and the
partial derivatives $Y_{s,w}$, $X_{u,w}$ do not vanish as $G_+$,
$G_-$ are diffeomorphisms. Observe that the map $G_0$ is very
similar to a quadratic Hénon-like map.

\newpage

\setcounter{section}{2}
\setcounter{equation}{0}
\section{Affine-like Maps\label{sec3}}
This section is essentially a summary of \cite{PY2}.

\subsection{Definition and Implicit Representation\label{sub3.1}}
Let $I^s_0$, $I^u_0$, $I^s_1$, $I^u_1$ be non trivial compact
intervals, $x_0$, $y_0$, $x_1$, $y_1$ the corresponding coordinates.
Consider a smooth diffeomorphisms $F$ whose domain is a {\it
vertical strip}
\[
P = \{ \vphi^- (y_0) \leqslant x_0 \leqslant \vphi^+ (y_0) \}
\subset I^s_0 \times I^u_0
\]
and whose image is a {\it horizontal strip}
\[
Q = \{ \psi^- (x_1) \leqslant y_1 \leqslant \psi^+ (x_1) \} \subset
I^s_1 \times I^u_1.
\]
We say that $F$ is {\it affine-like} if

{\bf (AL1)}~~the restriction to the graph of $F$ of the projection
onto $I^u_0 \times I^s_1$ is a diffeomorphism onto $I^u_0 \times
I^s_1$.

This allow us to define smooth maps $A$, $B$ on $I^u_0\times I^s_1$
such that
\begin{equation}
F (x_0, y_0) = (x_1, y_1) \iff
\begin{cases}
x_0 = A (y_0, x_1)\\
y_1 = B (y_0, x_1).
\end{cases}\label{eq3.l}
\end{equation}
The pair $(A,B)$ is the {\it implicit representation} (or
definition) of the affine-like map $F$. See figure~4. In the
formulas below, we shall most of the time omit the arguments of the
functions considered, which should be obvious from the context. We
will write $A_x$, $A_y$, $A_{xx}$, $B_x$, $B_y\cdots$ for partial
derivatives.

\begin{center}
\includegraphics{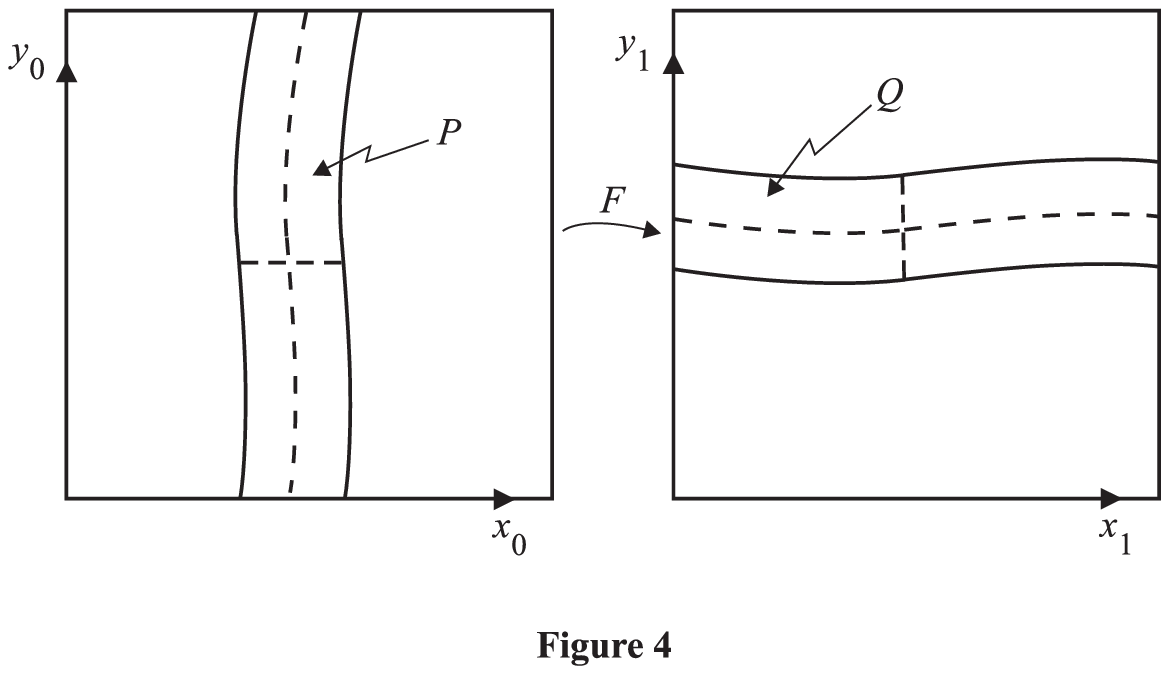}
\end{center}

On the graph of $F$, we have
\begin{eqnarray}
&dx_0 = A_y\, dy_0 + A_x\, dx_1,\label{eq3.2}\\
&dy_1 = B_y\, dy_0 + B_x\, dx_1,\nonumber
\end{eqnarray}
which leads to
\begin{eqnarray}
DF &= &A^{-1}_x
\begin{pmatrix}
1   &-A_y\\
B_x &A_x B_y - A_y B_x,\label{eq3.3}
\end{pmatrix}\\
DF^{-1} &= &B^{-1}_y
\begin{pmatrix}
A_x B_y - A_y B_x &A_y\\
-B_x &1,\label{eq3.4}
\end{pmatrix}\\
\det\; DF &= &A^{-1}_x B_y.\label{eq3.5}
\end{eqnarray}
The main advantage of the implicit representation is the symmetry
between positive and negative iteration.

\subsection{Cone Condition and Distortion\label{sub3.2}}
Let $\la$, $u$, $v > 0$ satisfy
\begin{equation}
1 < uv \leqslant \la^2.\label{eq3.6}
\end{equation}
Let $(X_0, Y_0)$ be a tangent vector at some point in the domain of
$F$, and let $(X_1, Y_1)$ be its image under $TF$.  The usual cone
condition with parameters $(\la, u, v)$ is:

{\bf (AL2)} (i)~~if $|Y_0| \leqslant u |X_0|$, then $|Y_1| \leqslant
v^{-1} |X_1|$ and $|X_1| \geqslant \la |X_0|$;

(ii)~~if $|X_1| \leqslant v |Y_1|$, then $|X_0| \leqslant u^{-1}
|Y_0|$ and $|Y_0| \geqslant \la |Y_1|$.

This is readily seen to be equivalent to
\begin{alignat}{1}
&\la |A_x| + u |A_y| \leqslant 1,\tag{\bf AL$'$2}\\
&\la |B_y| + v |B_x| \leqslant 1,\notag
\end{alignat}
everywhere on $I^u_0\times I^s_1$.

We will also need to control partial derivatives of second order of
$A$, $B$. By (\ref{eq3.5}), the partial derivatives $A_x$, $B_y$ do
not vanish on $I^u_0 \times I^s_1$. It turns out that the right way
to look at partial derivatives of second order is to consider the
six functions
\begin{eqnarray*}
&\pa_x\, \Log |A_x|,\, \pa_y\, \Log |A_x|,\, A_{yy},\notag\\
&\pa_y\, \Log |B_y|,\, \pa_x\, \Log |B_y|,\, B_{xx}.\notag
\end{eqnarray*}
We define the {\it distortion} of an affine-like map $F$, and denote
by $D(F)$, the maximal absolute value attained by any one of these
six functions on $I^u_0 \times I^s_1$.

We also define the {\it width} of the domain $P$ of $F$ by
\begin{equation}
|P|\; := \;\Max |A_x|,\label{eq3.7}
\end{equation}
and the width of the image $Q$ by
\begin{equation}
|Q|\; := \;\Max |B_y|.\label{eq3.8}
\end{equation}

\subsection{Simple Composition\label{sub3.3}}
The composition of two affine-like maps is not always affine-like.
However, the composition of two affine-like maps which also satisfy
the same cone condition (AL2) will again be affine-like and satisfy
the same cone condition (actually a better one).

More precisely, let $I^s_0$, $I^u_0$, $I^s_1$, $I^u_1$,  $I^s_2$,
$I^u_2$ be compact intervals. Let $F: P\to Q$ and $F' : P' \to Q'$
be affine-like maps with domains $P\subset I^s_0 \times I^u_0$,
$P'\subset I^s_1 \times I^u_1$ and images $Q \subset I^s_1 \times
I^u_1$, $Q' \subset I^s_2 \times I^u_2$.  We assume that both $F$
and $F'$ satisfy (AL2) (or (AL'2)) with parameters $\la$, $u$, $v$.
The composition $F'' = F' \circ F$ has domain $P'' = P \cap F^{-1}
(P')$ and image $Q'' = Q' \cap F' (Q)$. It satisfies (AL1) and (AL2)
with parameters $\la^2$, $u$, $v$ (cf.~\cite{PY2}). See figure~5.

\begin{center}
\includegraphics{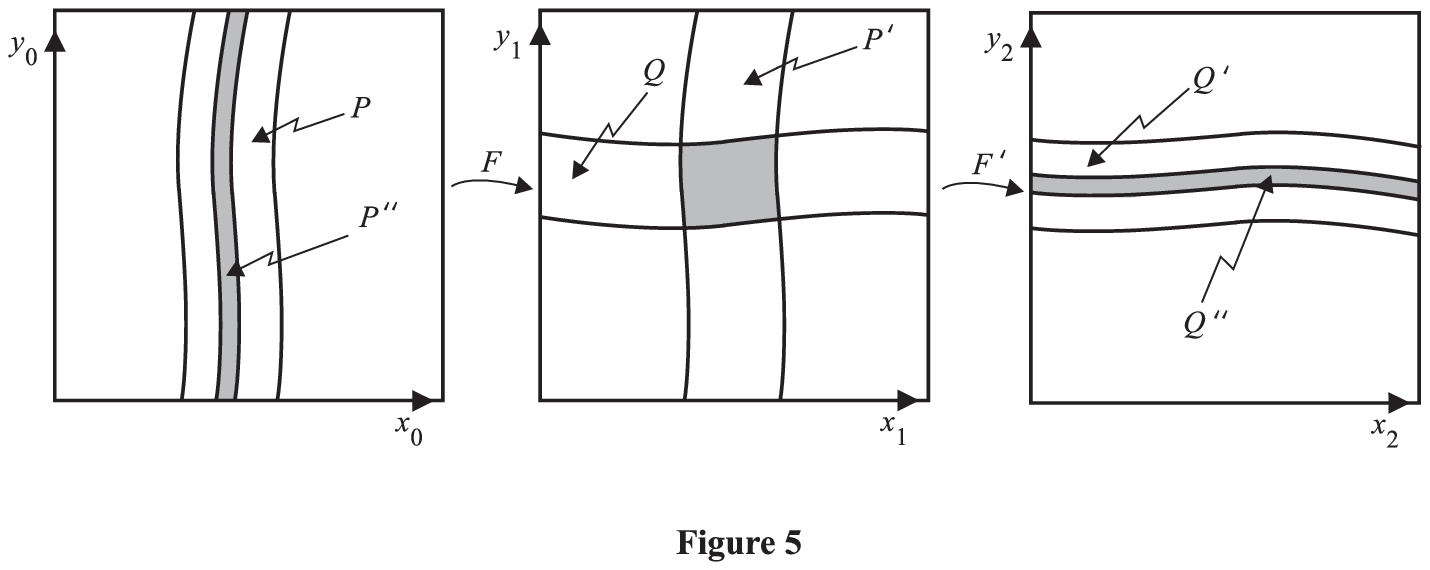}
\end{center}

Let $(A,B)$, $(A',B')$, $(A'',B'')$ be the implicit representations
of $F$, $F'$, $F''$ respectively.  Define
\begin{equation}
\De\; := \;1 - A'_y B_x > 1- u^{-1} v^{-1} > 0.\label{eq3.9}
\end{equation}
The partial derivatives of first order of $A''$, $B''$ are given by
\begin{alignat}{1}
A''_x  &= A_x A'_x \De^{-1},\label{eq3.10}\\
B''_y  &= B_y B'_y \De^{-1},\notag
\end{alignat}
\begin{alignat}{1}
A''_y  &= A_y + A'_y A_x B_y \De^{-1},\label{eq3.11}\\
B''_x  &= B'_x + B_x A'_x B'_y \De^{-1}.\notag
\end{alignat}
From (\ref{eq3.10}), we get
\begin{alignat}{2}
&C^{-1} \;\leqslant \;\frac{|P''|}{|P|\;|P'|} \; &\leqslant\; C,\label{eq3.12}\\
&C^{-1} \;\leqslant \;\frac{|Q''|}{|Q|\;|Q'|} \; &\leqslant\;
C,\notag
\end{alignat}
where the constants are uniform once $u$, $v$ are fixed and the
distortions are uniformly bounded.

The formulas for the partial derivatives of second order are derived
in \cite{PY2} and recalled in Appendix~A. They lead to the following
estimate for the distortion:
\begin{equation}
D (F'') \leqslant \Max \Bigl\{ D (F) + C|Q| (D(F) + D(F')), \; D(F')
+ C|P'| (D(F) + D (F' ))\Bigr\},\label{eq3.13}
\end{equation}
where $C$ depends only on $u$, $v$.

\subsection{Properties of the Markov Partition\label{sub3.4}}
We choose charts for the Markov partition discussed in
Subsection~\ref{sub2.1} in order to have the following property, for
some $\la$, $u$, $v$ satisfying (\ref{eq3.6}):

{\bf (MP5)}~~for any $(a,a')\in \cB$, any $g\in\cU$, the transition
map $g_{a,a'}$ from $P_{aa'} = R_a \cap g^{-1} (R_{a'})$ onto
$Q_{aa'} = R_{a'} \cap g (R_a)$ is affine-like and also satisfies
the cone condition (AL2).

These values of $(\la, u, v)$ will be fixed in  what follows.

To any finite word $\una = (a_0,\ldots,a_n)$ with transitions in
$\cB$, we have a composition
\[
g_{\una} = g_{a_{n-1} a_n} \circ \ldots \circ g_{a_{0} a_{1}}
\]
which satisfies also (AL1) and (AL2).

Moreover, as the widths decrease exponentially with the number of
iterations, it follows from (\ref{eq3.13}) that there exists $D_0 >
0$ such that all $g_\una$ satisfy

{\bf (MP6)}~~$D (g_\una) \le D_0$.

\subsection{Parabolic Composition\label{sub3.5}}
Let $G$ be the folding map of Subsection~\ref{sub2.3}, satisfying
properties (P1)--(P6).

Let also $I^s_0$, $I^u_0$, $I^s_1$, $I^u_1$ be compact intervals ;
let $F_0$ be an affine-like map from a vertical strip $P_0 \subset
I^s_0\times I^u_0$ to a horizontal strip $Q_0 \subset
I^s_{a_u}\times I^u_{a_u}$; let $F_1$
be an affine-like map from a vertical strip $P_1 \subset
I^s_{a_s}\times I^u_{a_s}$ to a horizontal strip $Q_1 \subset
I^s_1\times I^u_1$.

We recall from \cite{PY2} how, under appropriate hypotheses, the
composition $F_1\circ G\circ F_0$ defines two affine-like maps
$F^\pm$ with domains $P^\pm \subset P_0$ and image $Q^\pm \subset
Q_1$. See figure~6.

\begin{center}
\includegraphics{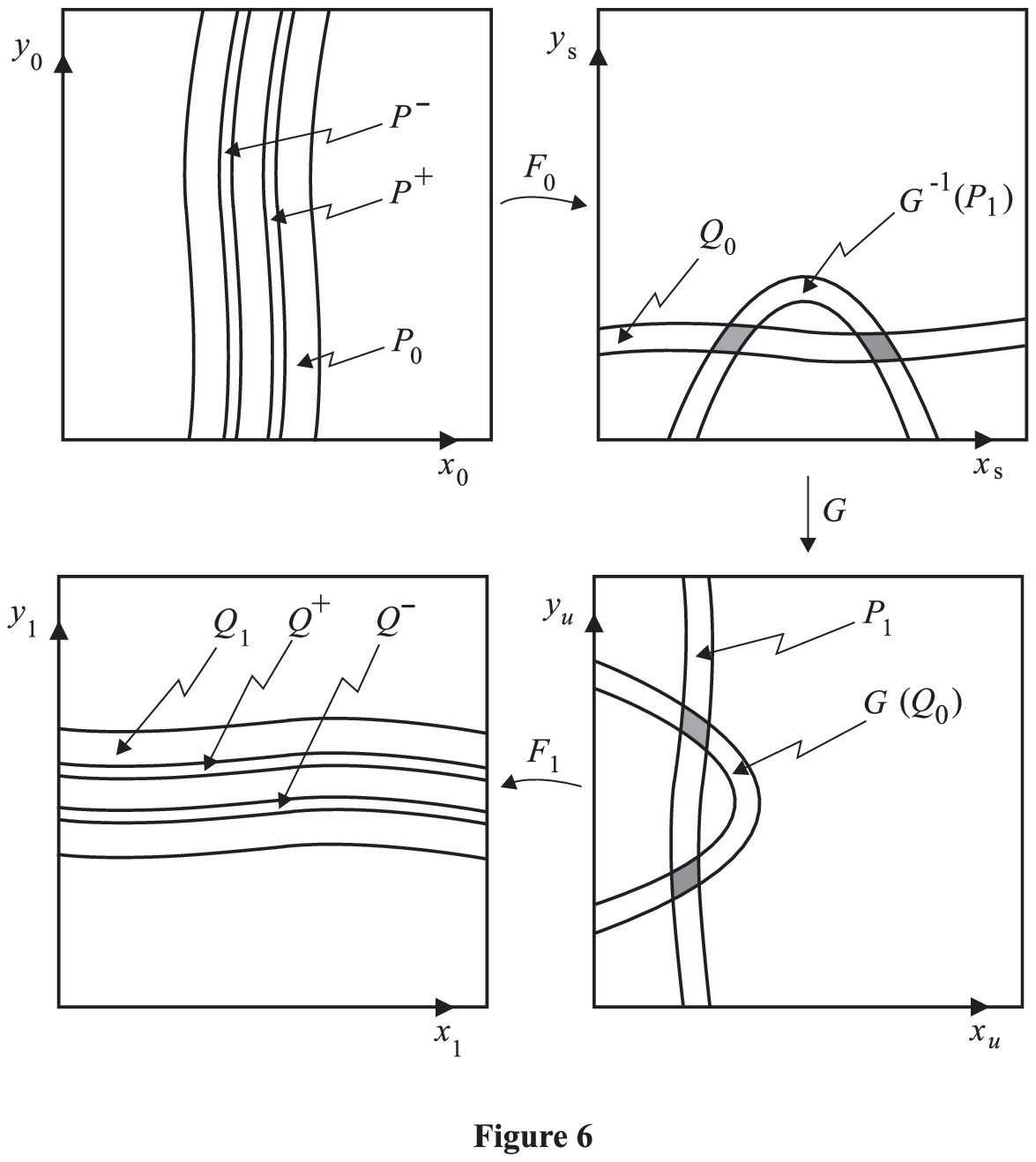}
\end{center}

Let $(A_0, B_0)$, $(A_1, B_1)$ be implicit representations of $F_0$,
$F_1$, respectively. We assume that
\begin{alignat}{2}
&|A_{1,y}|< b,\qquad &|A_{1,yy}|< b,\tag{\bf PC1}\\
&|B_{0,x}|< b,\qquad &|B_{0,xx}|< b,\notag
\end{alignat}
with $b \ll 1$. In the system
\begin{eqnarray}
&x_u = X_u (w, y_u),\label{eq3.14}\\
&y_u = B_0 (y_0, x_u),\notag
\end{eqnarray}
we can, as $|B_{0,x}| \ll 1$, eliminate $y_u$ and solve for $x_u$ to
define
\begin{equation}
x_u = X (w, y_0)\label{eq3.15}
\end{equation}
Similarly, in the system
\begin{eqnarray}
&y_s = Y_s (w, x_s),\label{eq3.16}\\
&x_s = A_1 (y_s, x_1),\notag
\end{eqnarray}
we eliminate $x_s$, solve for $y_s$ to get
\begin{equation}
y_s = Y (w, x_1).\label{eq3.17}
\end{equation}
The next step is to define
\begin{equation}
C (w, y_0, x_1) := w^2 - \te \Bigl( B_0 (y_0, X (w, y_0)), A_1 (Y
(w, x_1), x_1)\Bigr).\label{eq3.18}
\end{equation}
This quantity has the following geometrical interpretation. Fix
values $y^*_0$, $x^*_1$ for $y_0$, $x_1$. The image $G_- \circ F_0
(\{ y_0 = y^*_0\})$ is the graph
\begin{equation}
\ga_0 = \Bigl\{ y_u = B_0 (y^*_0, X (w,y^*_0))\Bigr\};\label{eq3.19}
\end{equation}
symmetrically, $G_+^{-1} \circ F_1^{-1} (\{ x_1 = x^*_1\})$ is the
graph
\begin{equation}
\ga_1 = \Bigl\{ x_s = A_1 (Y (w,x^*_1), x^*_1)\Bigr\}.\label{eq3.20}
\end{equation}
Then, $C (w,y^*_0,x^*_1)$ gives the relative position of the two
curves $\ga_0$ and $G_0^{-1} (\ga_1)$ (or equivalently $G_0 (\ga_0)$
and $\ga_1$).  More precisely, it is positive for all $w$ if the two
curves do not intersect; it vanishes at the intersection points and
is negative between the intersection points.

It follows from (PC1) just above that
\begin{alignat}{2}
|C_w - 2w| &\ll &1,\label{eq3.21}\\
|C_{ww}-2| &\ll &1.\label{eq3.22}
\end{alignat}
Therefore, for fixed values of $y_0$ and $x_1$, $C$ has a unique
minimum as a function of $w$; we denote by $\ovC (y_0, x_1)$ the
corresponding minimum value. We have $\ovC (y^*_0, x^*_1) > 0$
(resp.~$=0$, resp.~$<0$) if and only if the curves $\ga_0$ and
$G_0^{-1} (\ga_1)$ do not intersect (resp. are tangent, resp. have
two transverse intersection points).

In order to consider parabolic compositions, we shall require that
$\ovC (y_0, x_1) < 0$ everywhere on $I^u_0\times I^s_1$. Setting
\begin{equation}
\de = \min_{y_0, x_1} - \ovC (y_0, x_1)\label{eq3.23}
\end{equation}
we actually want to have
\begin{alignat}{1}
\de > b^{-1} (|P_1| + |Q_0|).\tag{\bf PC2}
\end{alignat}
The geometric interpretation of this requirement is clear: the
displacement of one of the rectangles and the image of the other
should be much bigger than the sum of their widths. In other words,
the distance between the tip of the parabolic strip $G^{-1}_0 (P_1)$
and the horizontal strip $Q_0$ should be much bigger than the widths
of these strips.

Assume now that (PC1) and (PC2) are satisfied; the equation $C (w,
y_0, x_1) = 0$ defines two smooth functions
\begin{equation}
w = W^\pm (y_0, x_1)\label{eq3.24}
\end{equation}
with $w^+ > w^-$.  One then defines
\begin{alignat}{2}
A^\pm (y_0, x_1)\; &:= &\;A_0 \Bigl(y_0, X (W^\pm (y_0, x_1),
y_0)\Bigr),\label{eq3.25}\\
B^\pm (y_0, x_1)\; &:= &\;B_1 \Bigl( Y (W^\pm (y_0, x_1),
x_1),x_1\Bigr).\label{eq3.26}
\end{alignat}
As shown in \cite{PY2}, the pair $(A^+, B^+)$ (resp.~$(A^-,B^-)$)
implicitly defines an affine-like map $F^+$ (resp.~$F^-$).

Denote by $P^+$ (resp.~$P^-$) the domain of $F^+$ and by $Q^+$
(resp.~$Q^-$) the domain of $F^-$. Then $P^+$ and $P^-$ are the two
components of $P_0 \cap (G\circ F_0)^{-1} (P_1)$, $Q^+$ and $Q^-$
are the two components of $Q_1 \cap (F_1\circ G)(Q_0)$; $F^+$
(resp.~$F^-$) is the restriction of $F_1 \circ G \circ F_0$ to $P^+$
(resp.~$P^-$).

The formulas for the partial derivatives of $A^\pm$, $B^\pm$ are
derived in \cite{PY2} and recalled in Appendix~A. They provide the
following estimate for the widths:
\begin{alignat}{2}
C^{-1} &\leqslant\; &\frac{|P^\pm|}{|P_0|\;|P_1| \de^{-\fudt}} \;
\leqslant\; C,\label{eq3.27}\\
C^{-1} &\leqslant\; &\frac{|Q^\pm|}{|Q_0|\;|Q_1| \de^{-\fudt}} \;
\leqslant\; C,\label{eq3.28}
\end{alignat}
where the constants are uniform once $b$ is fixed and the
distortions are uniformly bounded.

From [PY2, Theorem~3.7], we also have the following estimate for the
distortion of $F^\pm$: assuming that $b$ is small enough (in terms
of the partial derivatives of first order of $X_u$, $Y_s$, $\te$),
we have
\begin{equation}
D (F^\pm) \le \Max \Bigl\{ D (F_0) + C |Q_0| \de^{-1}, D (F_1) + C
|P_1| \de^{-1}\Bigr \},\label{eq3.29}
\end{equation}
provided that $D (F_0) + D (F_1) \le \de^{-\fudt}$. The constant $C$
in (\ref{eq3.29}) depends only on the partial derivatives of first
order of $X_s$, $Y_u$, $\te$.

We also recall from [PY2, formula~3.50] the estimate:
\begin{alignat}{2}
|A^\pm_y - A_{0,y}|\; &\le\; &C|P_0|\; |Q_0| \de^{-\fudt},\label{eq3.30}\\
|B^\pm_x - B_{1,x}|\; &\le\; &C|P_1|\; |Q_1|
\de^{-\fudt},\label{eq3.31}
\end{alignat}
where the right-hand terms must be small by (PC2).

As a concluding remark for this section, let us observe that, while
conditions (PC1), (PC2) are {\it necessary} in order to consider
parabolic composition, they are not {\it sufficient}: in
Section~\ref{sec5}, the requirement for parabolic composition will
be much more restrictive than (PC2).
\newpage

\setcounter{section}{3}
\setcounter{equation}{0}

\section{Structure of Parameter Space\label{sec4}}
\subsection{Some Important Constants\label{sub4.1}}
Throughout the rest of the paper, we will use four main constants
$\vep_0$, $\eta$, $\tau$, $\be$ which satisfy
\begin{equation}
0 < \vep_0 \ll \eta \ll \tau \ll \be-1 < 1.\label{eq4.1}
\end{equation}
We roughly explain the meaning of each constant:

--~~$\vep_0$ is the maximal width of the parabolic tongues $L_u$,
$L_s$. It is also the size of the parameter interval we start with.

--~~$\eta$ is involved in the transversality relation (defined in
Section~5) which allows parabolic composition: instead of the
condition (PC2) of Subsection~\ref{sub3.5}, roughly speaking we will
ask that
\begin{equation}
\de \geqslant (|P_1| + |Q_0|)^{1-\eta}.\label{eq4.2}
\end{equation}
--~~$\tau$ relates the successive scales of the parameter intervals
we will consider through the formula $\vep_{k+1} = \vep_k^{1+\tau}$.

--~~$\be$ will actually be given in Section~\ref{sec9} by an
explicit formula in terms of $d^0_s$, $d^0_u$; the condition (H4)
involving $d^0_s$, $d^0_u$ in Section~\ref{sec1} is required because
we need $\be>1$. It appears in the definition of regularity in
Section~\ref{sec5}, which controls the recurrence of the ''critical
locus''.

\subsection{One-Parameter Families\label{sub4.2}}
From now on, we fix a one-parameter family $(g_t)_{t\in (-t_0,t_0)}$
in $\cU$. We assume that the family is transverse to $\cU_0$ at
$t=0$, with $g_t\in \cU_+$ for $t>0$ and $g_t\in \cU_-$ for $t<0$.

Observe that $g_0$ satisfies exactly the same assumptions as $f$,
provided $\cU$ is small enough. Therefore, we may and shall, assume
that $g_0 = f$.

We will first reparametrize the family in order to make some
computations simpler. Consider the folding map $G_t = g_t^{N_0}$ of
Subsection~\ref{sub2.3}. If $t_0$ is small enough, $G_t$ is a fold
map for all values of $t\in (-t_0,t_0)$. Moreover, we can in
properties (P2), (P3) of Subsection~\ref{sub2.3} choose a function
$\te$ which depends smoothly on $t$.

From (MP3), Subsection~\ref{sub2.2}, the values $y_u=0$, $x_s=0$ of
the arguments of $\te$ correspond to $W^u_\loc (p_u)$ and $W^s_\loc
(p_s)$ respectively. Therefore, the transversality of the family to
$\cU_0$ is equivalent to
\begin{equation}
\frac{\pa}{\pa t}\; \te_t (0,0) \mid_{t=0}\; > \; 0.\label{eq4.3}
\end{equation}
Taking $t_0$ small enough, we can therefore reparametrize our family
in order to have
\begin{equation}
\te_t (0,0) \equiv t,\qquad t\in (-t_0, t_0).\label{eq4.4}
\end{equation}

\subsection{Parameter Intervals\label{sub4.3}}
The starting parameter interval will be
\begin{equation}
I_0\; := \;[\vep_0, 2\vep_0],\label{eq4.5}
\end{equation}
where, as explained above, $\vep_0$ will be taken very small. This
is the only parameter interval at level~0.

At level $k$, we will deal with parameter intervals of length
$\vep_k$, where the sequence of scales $\vep_k$ is defined
inductively by
\begin{equation}
\vep_{k+1} = \vep_k^{1+\tau}.\label{eq4.6}
\end{equation}
The constant $\tau$ is small, but $\vep_0$ is much smaller and in
particular we will have $\vep_0^{\tau^2} \ll 1$. Every parameter
interval of level $k$ is divided into $[\vep_k^{-\tau}]$ parameter
intervals of level $k+1$.

The remaining part, if any, is discarded; it is of length $<
\vep_{k+1}$; the total length discarded in this way is smaller than
$\vep_1 \ll \vep_0$.

Let $\wtI$ be a parameter interval of level $k$ and $I$ be a
parameter interval of level $k+1$ contained in $\wtI$. We say that
$\wtI$ is the {\it parent} of $I$ and that $I$ is a {\it child} of
$\wtI$.

\subsection{The Selection Process\label{sub4.4}}
In Section~\ref{sec5}, we will define what it means for a parameter
interval to be {\it regular}. The starting interval $I_0$ will be
regular.

Given a regular parameter interval $\wtI$ of level $k$, we divide it
into its children: these parameter intervals of level $k+1$ are the
{\it candidates}. We then test each candidate for regularity and
discard those which are not regular. We then proceed to level $k+1$
with each surviving candidate.

The {\it regular} parameters are those which are the intersection of
a decreasing sequence of regular parameter intervals. For such
parameters, we are able to carry out some analysis of the maximal
invariant set $\La_{g_t}$.

\subsection{Strongly Regular Parameters\label{sub4.5}}
The regularity property is, in some sense, the minimal requirement
that is needed to keep control on the geometry and dynamics of the
maximal invariant set. However, this requirement is of an
essentially qualitative character and this leads in particular to
the following difficulty: we are not able to estimate which
proportion of the children of a regular parameter interval are also
regular.

To circumvent this problem, we define in Section~\ref{sec9} a
stronger property for parameter intervals, called {\it strong
regularity}. It implies regularity, and is better adapted to the
inductive selection process. It also gives additional geometric
information on the maximal invariant set.

When $\wtI$ is a strongly regular parameter interval of level $k$,
we will show in Section~\ref{sec9} that most candidates of level
$k+1$ contained in $\wtI$ are also strongly regular. The proportion
of discarded candidates is less than $\al_k$, with
\begin{equation}
\sum_{k\geqslant 0} \al_k \ll 1,\label{eq4.7}
\end{equation}
the $\ll$ sign means that the sum gets arbitrarily small as $\vep_0$
goes to zero. Then we can conclude that most parameters are strongly
regular in the sense that they are equal to the intersection of
decreasing sequences of strongly regular parameter intervals.

The non-uniformly hyperbolic horseshoes that are the subject of our
study are exactly the maximal invariant set $\La_g$ for strongly
regular $g\in \cU_+$.
\newpage

\setcounter{section}{4}
\setcounter{equation}{0}

\section{Classes of Affine-like Iterates and the Transversality Relation\label{sec5}}
\subsection{Affine-like Iterates\label{sub5.1}}
Let $I$ be a parameter interval of some level.

\Def An {\it $I$-persistent affine-like iterate} is a triple
$(P,Q,n)$ such that

--~~$P$ is a vertical strip in some $R_a$, depending smoothly on
$t\in I$;

--~~$Q$ is a horizontal strip in some $R_{a'}$, depending smoothly
on $t\in I$;

--~~$n$ is a nonnegative integer;

--~~for each $t\in I$, the restriction of $g^n_t$ to $P_t$ is an
affine-like map onto $Q_t$, i.e. property (AL1) of
Subsection~\ref{sub3.1} holds;

--~~for each $t\in I$, each $m\in [0,n]$, we have $g_t^m
(P_t)\subset \whR$.

\medskip\noindent
{\bf Examples.}

1.~~For $n=0$, the $I$-persistent affine-like iterates are the
$(R_a, R_a, 0)$, $a\in \AL$.

2.~~For $n=1$, the $I$-persistent affine-like iterates are the
$(P_{aa'}, Q_{aa'}, 1)$, $(a,a')\in \cB$.

3.~~More generally, for any finite word $\una = (a_0, \ldots, a_n)$
with transitions in $\cB$, the map $g_\una$ of
Subsection~\ref{sub3.4} defined an $I$-persistent affine-like
iterate $(P_\una, Q_\una, n)$.

\medskip\noindent
{\bf Notation.~} If $P$ is a vertical strip $\{\vphi_- (y)\leqslant
x \leqslant \vphi_+ (y)\}$ we denote by $\pa P$ the vertical part of
the boundary, i.e. the two graphs $\{ x = \vphi^\pm (y)\}$.
Similarly for horizontal strips.

If $(P, Q, n)$ is an $I$-persistent affine-like iterate and $I'$ is
a parameter interval contained in $I$, $(P, Q, n)$ also defines by
restriction an $I'$-persistent affine-like iterate. A slightly less
trivial property is given by

\begin{Propo}\label{propo1}
Let $(P,Q,n)$, $(P',Q',n')$ be $I$-persistent affine-like iterates.
We have

a)~~if $n=n'$, then either $P=P'$ and $Q=Q'$ for all $t\in I$ or
$P\cap P'=\emptyset$, $Q\cap Q'=\emptyset$ for all $t\in I$.

b)~~if $n<n'$, then either $P\supset P'$, $\pa P\cap P'=\emptyset$
for all $t\in I$ or $P\cap P'=\emptyset$, for all $t\in I$.
\end{Propo}

\begin{rem}
Throughout the paper, except in Section~\ref{sec9} (where we break
the symmetry assuming $d^0_s\geqslant d^0_u$), we will keep a
time-symmetric setting. Thus every property stated for the domains
$P'$s is also valid for the images $Q'$s. This apply for instance to
part b) of the proposition.
\end{rem}

\medskip\noindent {\it Proof.~} By the definition of an
$I$-persistent affine-like iterate, for all $t\in I$, $P$ is a
connected component of $R\cap g_t^{-n} (R)$ and also of
$\bigcap\limits_{0\leqslant m\leqslant n} g_t^{-m} (\whR)$.

a)~~If $n=n'$ and $P\cap P' \ne \emptyset$ for some $t_0\in I$, we
must have $P=P'$ at $t_0$ and hence, $P\cap P'\ne \emptyset$ for $t$
close to $t_0$. It follows that $P=P'$ for all $t\in I$, and also
$Q=Q'$ for all $t\in I$.

b)~~Assume that $n<n'$ and $P\cap P'\ne\emptyset$ for some $t_0\in
I$, then $P'\subset P$ at $t_0$ (since $\bigcap\limits_{0\leqslant
m\leqslant n'} g_{t_0}^{-m} (\whR)$ is contained in
$\bigcap\limits_{0\leqslant m\leqslant n} g_{t_0}^{-m} (\whR)$),
hence $P'\cap P\ne \emptyset$ for $t$ close to $t_0$ and $P'\subset
P$ for all $t\in I$.

Let $t\in I$, $z\in \pa P$; then, $g_t^n (z)$ belongs to the
vertical boundary of $R$ and $g_t^{n+1}(z)\not\in \whR$; therefore,
$z\not\in P'$. This proves that $\pa P\cap P'$ is empty for all $t\in
I$.\hfill $\square$

\subsection{The Classes $\cR(I)$: General Overview\label{sub5.2}}
It would be nice to work with the class of all $I$-persistent
affine-like iterates, but with this approach one faces two problems:

--~~$I$-persistent affine-like iterates do not satisfy a uniform
cone condition, and they do not have uniformly bounded distortion;

--~~even if we force such uniformity in the definition, a major
problem is that we lack some control on the way in which long
$I$-persistent affine-like iterates are constructed from shorter
ones by simple or parabolic composition.

To overcome these problems, we will define, for every parameter
interval (whose parent is regular, see Subsection~\ref{sub4.4} and
the end of the present section on the definition of regularity) a
subset $\cR(I)$ of the set of all $I$-persistent affine-like
iterates; all elements of $\cR(I)$ with $n>1$ can be, from the very
definition, obtained from shorter ones by simple or parabolic
composition; the elements of $\cR(I)$ will turn out to satisfy a
uniform cone condition and have uniformly bounded distortion.

The main ingredient in the definition of $\cR(I)$ is a {\it
transversality relation} which is an appropriate strengthening of
condition (PC2) in Subsection~\ref{sub3.5}. Simple composition is
allowed whenever it makes sense, but parabolic composition is only
allowed when this transversality relation holds.

The definition of the transversality relation, given later in this
section, is quite involved; this is because we want some
combinatorial properties proved in Section~\ref{sec6} to be
satisfied. Such properties make our later work much easier.

All the process of this definition is based on a double induction:

--~~an induction on the level of the parameter interval, starting
with $I_0 = [\vep_0, 2\vep_0]$ at level 0;

--~~for a given parameter interval, an induction on the length $n$
of the $I$-persistent affine-like iterates which are considered.

In order for the definition of $\cR(I)$ to be consistent, a number
of properties (uniform cone condition, uniform bound on distortion
and many others) must hold; unfortunately, the proof of any of these
properties is inductive, using several {\it other} properties for
shorter iterates. Our presentation, is therefore, organized through
the following scheme: we will first define $\cR(I)$ and list all the
required properties conditioning the definition. Then we will
proceed to the proofs of these properties, and of other related
ones. The whole process will occupy the rest of the present section
and Sections~\ref{sec6} and \ref{sec7}.

\subsection{Definition of the Special Class of Affine-Like Iterates $\cR(I)$\label{sub5.3}}
Let $I$ be a parameter interval of some level.

We claim that there exists a (unique) class of $I$-persistent
affine-like iterates which satisfy the properties (R1)--(R7) below.

{\bf (R1)}~~For any word $\una = (a_0,\ldots,a_n)$ with transitions
in $\cB$, the element $(P_\una, Q_\una, n)$ (see example~3 above)
belongs to $\cR(I)$.

For the starting interval $I_0 = [\vep_0, 2\vep_0]$, it will turn
out that one obtains in this way all elements of $\cR(I_0)$.

Recall from (MP5), (MP6) in Subsection~\ref{sub3.4} that all
$(P_\una, Q_\una, n)$ with $n > 0$ satisfy for all $t\in I_0$ a uniform cone
condition (AL2) with parameters $\la$, $u$, $v$ (satisfying
$1<uv\leqslant \la^2)$, and have distortion bounded by $D_0$. Let
$u_0 = u^\fud$, $v_0 = v^\fud$.

{\bf (R2)}~~All $(P,Q,n)\in \cR(I)$ satisfy for all $t\in I$ the
cone condition (AL2) with parameters $\la,u_0,v_0$ and have
distortion bounded by $2D_0$ for all $t\in I$.

Let $(P,Q,n)$, $(P',Q',n')$ be elements of $\cR(I)$ such that
$Q\subset R_a$, $P'\subset R_a$ for some $a\in\AL$. As both iterates
satisfy the cone condition (AL2) with parameters $\la$, $u_0$,
$v_0$, we know from Subsection~\ref{sub3.3} that the simple
composition defined by
\begin{equation}
P'' = P \cap g^{-n} (P'), \qquad Q'' = Q' \cap g^{n'}(Q),\qquad n''
= n+n',\label{eq5.1}
\end{equation}
is an ($I$-persistent) affine-like iterate.

The next condition states that it should also belong to $\cR(I)$.

{\bf (R3)}~~The class $\cR(I)$ is stable under simple composition.

We now turn to parabolic composition.

We first define two special elements which belong to $\cR(I)$
according to (R1): define $(P_s, Q_s, n_s)$ (resp.~$(P_u, Q_u,
n_u)$) to be the element $(P_\una, Q_\una, n)$ with maximal length
$n$ such that $L_s\subset P_s$ for all $t\in I_0$ (resp.~$L_u\subset
Q_u$ for all $t\in I_0$). We have that $p_s\in P_s$ and $p_u\in
Q_u$. See figure~7.

\begin{center}
\includegraphics{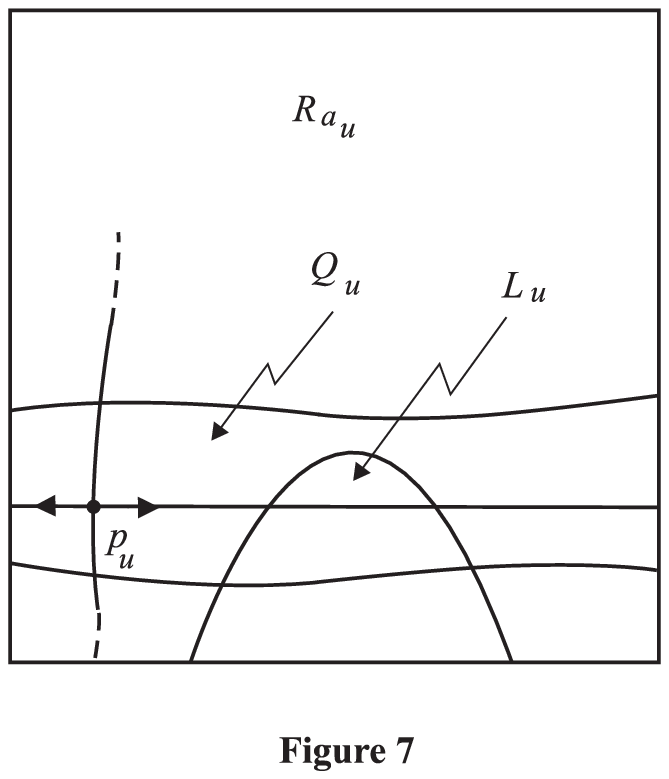}
\end{center}

We obviously have, for all $t\in I_0$
\begin{alignat}{2}
C^{-1} \vep_0 &\leqslant \;|P_s| &\leqslant \;C \vep_0,\label{eq5.2}\\
C^{-1} \vep_0 &\leqslant \;|Q_u| &\leqslant \;C \vep_0.\notag
\end{alignat}

The next condition guarantees that property (PC1) in
Subsection~\ref{sub3.5} is satisfied.

{\bf (R4)}~~Let $(A,B)$ be the implicit representation of an
affine-like iterate $(P,Q,n)\in \cR(I)$.

a)~~If $P\subset P_s$, then for all $t\in I$ we have
\[
|A_y| \leqslant C\vep_0,\qquad |A_{yy}|\leqslant C\vep_0.
\]
b)~~If $Q\subset Q_u$, then for all $t\in I$ we have
\[
|B_x| \leqslant C\vep_0,\qquad |B_{xx}|\leqslant C\vep_0.
\]
Here and in the sequel, the letter $C$ denotes various constants
which depend only on our initial diffeomorphism $f$, but {\it not}
on $\eta$, $\tau$, $\vep_0$.

Let $(P_0,Q_0,n_0)$, $(P_1,Q_1,n_1)$ be elements in $\cR(I)$ with
$Q_0\subset Q_u$, $P_1\subset P_s$. In these circunstances, we will
define in Subsection~\ref{sub5.4} a {\it transversality relation}
denoted by $Q_0 \pitchfork_I P_1$ which may or may not hold. When it
holds, it implies condition (PC2) of Subsection~\ref{sub3.5} for all
$t\in I$ (see (R7) below).

{\bf (R5)}~~If $(P_0,Q_0,n_0)$, $(P_1,Q_1,n_1)$ as above satisfy
$Q_0\pitchfork_I P_1$, then both $I$-persistent affine-like iterates
obtained from the parabolic composition $g_t^{n_1} \circ G_t \circ
g_t^{n_0}$ belong to $\cR(I)$.

Writing $(P^+,Q^+,n)$ and $(P^-,Q^-,n)$ for these two iterates, we
have $n = n_0 + n_1 + N_0$. The domains $P^+$ and $P^-$ are the two
components of $g^{-n_0} (Q_0\cap G_t^{-1} (P_1))$; the images $Q^+$
and $Q^-$ are the two components of $g^{n_1} (P_1\cap G_t (Q_0))$.
See figure~6, Subsection~\ref{sub3.5}.

When $(P_0,Q_0,n_0)$, $(P_1,Q_1,n_1)$ satisfy $Q_0\subset Q_u$,
$P_1\subset P_s$, $Q_0\forki P_1$, we say that their parabolic
composition is allowed in $\cR(I)$.

{\bf (R6)}~~Any $(P,Q,n)\in \cR(I)$ with $n>1$ can be obtained from
shorter elements by simple composition or (allowed) parabolic
composition.

Typically, an element of $\cR(I)$ can be obtained in many ways by
composition of shorter ones. We say that an element of $\cR(I)$ is
{\it prime} if it cannot be obtained by {\it simple} composition of
shorter ones. Prime elements play a key role in the description of
the dynamics for regular parameters in Section~\ref{sec10}.

It is pretty clear from conditions (R1), (R3), (R5), (R6) alone that
there is at most one class $\cR(I)$ satisfying these conditions. The
existence of $\cR(I)$, i.e. the proof of the consistency of
conditions (R1)--(R6), is much more delicate. There is actually a
seventh property (R7) formulated in the next subsection and related
to the condition (PC2) for parabolic composition.

\medskip\noindent
{\bf Parent-Child Terminology and Notations for Compositions.}

Let $(P,Q,n)$, $(\wtP,\wtQ,\wtn)\in \cR(I)$ with $P\subset\wtP$,
$n>\wtn$. If there is no $(\whP,\whQ,\whn)\in\cR(I)$ with $P\subset
\whP \subset \wtP$ and $n>\whn>\wtn$, we say that $P$ is a {\it
child} of $\wtP$ and that $\wtP$ is the {\it parent} of $P$; if
moreover $n = \wtn + 1$, we say that $P$ is a {\it simple} child;
otherwise it is a non-simple child.

Let $(P_0,Q_0,n_0)$, $(P_1,Q_1,n_1)\in \cR(I)$. If $Q_0$, $P_1$ are
contained in a same rectangle $R_a$, the simple composition
$(P,Q,n)\in \cR(I)$ of these elements will be written as
\begin{equation}
(P,Q,n) = (P_0,Q_0,n_0)\;\; * \;\;(P_1,Q_1,n_1),\label{eq5.3}
\end{equation}
If $Q_0\subset Q_u$, $P_1\subset P_s$ and $Q_0\pitchfork_I P_1$, any
of the two elements $(\whP, \whQ, \whn)$ obtained by the
corresponding allowed parabolic composition will be written as
\begin{equation}
(\whP,\whQ,\whn) \in (P_0,Q_0,n_0)\;\; \square \;\;
(P_1,Q_1,n_1).\label{eq5.4}
\end{equation}

\subsection{Definition of the Transversality Relation\label{sub5.4}}
Let $I$ be a parameter interval of some level, and let
$(P_0,Q_0,n_0)$, $(P_1,Q_1,n_1)$ be elements of $\cR(I)$ which
satisfy $Q_0\subset Q_u$, $P_1\subset P_s$.

From (R4) the condition (PC1) of Subsection~\ref{sub3.5} is
satisfied provided $\vep_0$ small enough. Denote by $(x_0, y_0)$
(resp.~$(x_1,y_1)$) the coordinates in the rectangle containing
$P_0$ (resp.~$Q_1$). A function $\ovC (y_0,x_1)$ was defined in
Subsection~\ref{eq3.5}, together with
\begin{equation}
\de (Q_0, P_1)\; =\; \min_{y_0} \; \min_{x_1}\; -\; \ovC (y_0,
x_1).\label{eq5.5}
\end{equation}
In Subsection~\ref{sub3.5}, we were asking for $\de$ to be much
larger than $|P_0|$ and $|Q_1|$. From the formulas of Appendix~A, we
have
\begin{alignat}{2}
C^{-1} |P_1| &\leqslant \;|\ovC_x|\; &\leqslant \;C|P_1|,\label{eq5.6}\\
C^{-1} |Q_0| &\leqslant \;|\ovC_y|\; &\leqslant
\;C|Q_0|.\label{eq5.7}
\end{alignat}
We will introduce
\begin{alignat}{2}
\de_L (Q_0, P_1) &:= &\;\;\max_{y_0} \; \min_{x_1} - \;\ovC (y_0,x_1),\label{eq5.8}\\
\de_R (Q_0, P_1) &:= &\;\;\min_{y_0} \; \max_{x_1} - \;\ovC (y_0,x_1),\label{eq5.9}\\
\de_{LR} (Q_0, P_1) &:= &\;\;\max_{y_0} \; \max_{x_1} - \;\ovC
(y_0,x_1).\label{eq5.10}
\end{alignat}
All together, $\de$, $\de_L$, $\de_R$, $\de_{LR}$ are the values of
$-\ovC$ at the four corners of the rectangle of definition of
$\ovC$. We have from (\ref{eq5.6}), (\ref{eq5.7}) that
\begin{alignat}{2}
C^{-1} |Q_0| &\leqslant \;\de_L (Q_0,P_1) - \de (Q_0, P_1)\; &\leqslant C |Q_0|,\label{eq5.11}\\
C^{-1} |P_1| &\leqslant \;\de_R (Q_0,P_1) - \de (Q_0, P_1)\; &\leqslant C |P_1|,\label{eq5.12}\\
C^{-1} |Q_0| &\leqslant \;\de_{LR} (Q_0,P_1) - \de_R (Q_0, P_1)\; &\leqslant C |Q_0|,\label{eq5.13}\\
C^{-1} |P_1| &\leqslant \;\de_{LR} (Q_0,P_1) - \de_L (Q_0, P_1)\;
&\leqslant C |P_1|.\label{eq5.14}
\end{alignat}

{\bf Preliminary Definition.~} {\it We write $Q_0\ovfork_I P_1$ if
the following holds}

{\bf (T1)}~~for {\it all} $t\in I$,
\[
\de_{LR} (Q_0, P_1) \geqslant 2 |I|,
\]
{\bf (T2)}~~for {\it some} $t_0\in I$,
\[
\de_{R} (Q_0, P_1) \geqslant 2 |Q_0|^{1-\eta},
\]
{\bf (T3)}~~for {\it some} $t_1\in I$,
\[
\de_{L} (Q_0, P_1) \geqslant 2 |P_1|^{1-\eta}.
\]

\medskip\noindent
{\bf Definition.}

{\it We say that $Q_0$, $P_1$ are $I$-{\it transverse} and write
$Q_0\forki P_1$ if there exist a parameter interval $\wtI\supset I$,
elements $(\wtP_0, \wtQ_0, \wtn_0)$, $(\wtP_1, \wtQ_1, \wtn_1)\in
\cR(\wtI)$ with $\wtP_1 \supset P_1$, $\wtQ_0 \supset Q_0$ such that
$\wtQ_0 \ovfork_{\wtI}\, \wtP_1$.}

\medskip\noindent
{\bf Remark.}

1.~~{\it Taking $\wtI = I$, $\wtP_0 = P_0$, $\wtQ_1 = Q_1$, it is
obvious that if $Q_0 \ovforki P_1$, then $Q_0\forki P_1$.

2.~~In view of our inductive procedure, all $(\wtP_0, \wtQ_0,
\wtn_0)$, $(\wtP_1, \wtQ_1, \wtn_1)$ which have to be considered
have been constructed before $(P_0, Q_0, n_0)$, $(P_1, Q_1, n_1)$.

3.~~As mentioned before, the definition of the transversality
relation is quite involved. Some justification for the choice of
quantifiers in (T1), (T2), (T3) can be found in Appendix~C.}

At first sight, it appears that properties (T2), (T3) above are not
quite sufficient to guarantee condition (PC2) of parabolic
composition (Subsection~\ref{sub3.5}), because they involve only one
value of the parameter. The next property takes care of this
problem.

{\bf (R7)}~~If $(P_0,Q_0,n_0)$, $(P_1,Q_1,n_1)\in \cR(I)$ satisfy
$Q_0\subset Q_u$, $P_1\subset P_s$ and $Q_0\forki P_1$ holds, then,
for all $t\in I$, we have
\[
\de (Q_0, P_1)\; \geqslant\; C^{-1} \Bigl( |P_1|^{1-\eta} +
|Q_0|^{1-\eta}\Bigr).
\]

\subsection{The Class $\cR(I_0)$\label{sub5.5}}
Recall that in Subsection~\ref{sub4.2} we had the normalization
\begin{alignat}{1}
\te (0, 0, t) \equiv t,\qquad |t|\leqslant t_0.\tag{4.4}
\end{alignat}
As $\te$ is monotonous in both variables, it follows from the
definition of $C$ (formula (\ref{eq3.18}) of
Subsection~\ref{sub3.5}) that for all $t\in (-t_0, t_0)$:
\begin{equation}
-\ovC (y_0, x_1, t)\; \leqslant\; t\label{eq5.15}
\end{equation}
Therefore, for the starting interval $I_0 = [\vep_0, 2\vep_0]$,
condition (T1) above is never satisfied.

Parabolic composition is never allowed in $\cR(I_0)$. Thus (see
Subsection~\ref{sub5.3}) the class $\cR(I_0)$ consists only of the
elements $(P_\una, Q_\una, n)$ given by (R1) which are associated
with the horseshoe $K$.

Condition (R2) is satisfied because of the choice of $u_0$, $v_0$,
$D_0$. Condition (R3) is obviously satisfied, as are conditions (R5)
and (R6). Condition (R4) will be checked in Section~\ref{sec7}.

\subsection{Definition of the Regularity Property\label{sub5.6}}
We introduce first some terminology and some concepts related to the
transversality relation.

{\bf 5.6.1~} Let $I$ be a parameter interval of some level, and let
$(P_0,Q_0,n_0)$, $(P_1,Q_1,n_1)$ be elements of $\cR(I)$ such that
$Q_0\subset Q_u$ and $P_1\subset P_s$. When $Q_0$ and $P_1$ are not
$I$-transverse, two cases may happen:

--~~if $G_t(Q_0)\cap P_1 = \emptyset$ for all $t\in I$, we say that
$Q_0$ and $P_1$ are $I$-{\it separated};

--~~otherwise, we say that $Q_0$ and $P_1$ are $I$-{\it critically
related}.

{\bf 5.6.2~} Let $(P,Q,n)\in \cR(I)$. An $I$-{\it decomposition} of
$P$ is a finite family $(P_\al,Q_\al,n_\al)$ of elements of $\cR(I)$
such that the $P'_\al$s are disjoint, contained in $P$ and satisfy
\begin{equation}
W^s (\La,\whR) \cap P\; = \; \bigsqcup_\al\;(W^s (\La,\whR) \cap
P_\al),\label{eq5.16}
\end{equation}
where $W^s (\La,\whR)$ was defined in Subsection~\ref{sub2.2}. We
say that $P$ is $I$-{\it decomposable} if it admits a non-trivial
$I$-{\it decomposition}. Then, there is a coarsest non-trivial
$I$-{\it decomposition}, namely by the children of $P$, which is
called the {\it canonical $I$-decomposition.}

\begin{rem}
We will see in Section~\ref{sec8} that any $P$ has only finitely
many children.
\end{rem}

{\bf 5.6.3~} Let $(P,Q,n)\in \cR(I)$. We say that $Q$ is {\it
$I$-transverse} if either $Q\cap Q_u = \emptyset$ or $Q\subset Q_u$
and there exists an {\it $I$-decomposition} $(P_\al, Q_\al, n_\al)$
of $P_s$ such that, for any $\al$, $Q$ and $P_\al$ are either
$I$-transverse or $I$-separated.

We say that $Q$ is {\it $I$-critical} when it is not $I$-transverse.
This is always the case if $Q\supset Q_u$.

We also define in a symmetric way an $I$-decomposition for $Q$, and
$I$-transversality or $I$-criticality for $P$.

{\bf 5.6.4~} We say that $(P,Q,n)\in \cR(I)$ is {\it $I$-bicritical}
if both $P$ and $Q$ are $I$-critical. The corresponding iterate
should be seen as describing some recurrence to the ``critical
locus''.

\Def Let $\be > 1$. We say that the parameter interval $I$ is {\it
$\be$-regular} (or just regular when the value of $\be$ is fixed) if
any $I$-bicritical element $(P,Q,n)\in \cR (I)$ satisfies, for all
$t\in I$:
\begin{equation}
|P| < |I|^\be,\qquad |Q| < |I|^\be.\label{eq5.18}
\end{equation}

\bigskip\noindent
{\bf Remark.}

1.~~{\it When $d^0_s\geqslant d^0_u$, in Section~\ref{sec9} we will
take for $\be$ a number satisfying
\begin{equation}
1 < \be  < \be_{max} := \frac{(1-d^0_u)(d^0_s + d^0_u)}{d^0_s
(2d^0_s + d^0_u -1)}.\label{eq5.19}
\end{equation}
Condition (H4) in Section~\ref{sub1.2} is actually equivalent to
$\be_{max} > 1$. When $d^0_u \geqslant d^0_s$, we exchange $d^0_s$
and $d^0_u$ in the definition of $\be_{max}$.

2.~~We will see in Section~\ref{sec6} that if $\wtI$ is a
$\be$-regular parameter interval, and $I\subset \wtI$ is a candidate
at the next level, then $I$ is at least $\ovbe$-regular with}
\begin{equation}
\ovbe = \be (1+\tau)^{-1}.\label{eq5.20}
\end{equation}
This will allow us to apply to candidates, with only slightly worse
constants, all results that have been proven for regular parameter
intervals. But, obviously, we cannot let the regularity exponent to
deteriorate too much (it must stay $> 1$), which explains why
candidates have to pass the $\be$-regularity test.
\newpage

\setcounter{section}{5}
\setcounter{equation}{0}

\section{Some Properties of the Classes $\cR(I)$\label{sec6}}
We recall that all parameter intervals that we consider in the
sequel are assumed to be regular or candidates (meaning in this case
that their parent is regular).

\subsection{Transversality is Hereditary\label{sub6.1}}
The following obvious but fundamental property was forced into the
definition of the transversality relation.

\begin{Propo}\label{propo2}
Let $\wtI \supset I$ be parameter intervals. Let $(P_0,Q_0,n_0)$,
$(P_1,Q_1,n_1)\in \cR(I)$ and $(\wtP_0, \wtQ_0, \wtn_0)$, $(\wtP_1,
\wtQ_1, \wtn_1)\in \cR(\wtI)$. Assume that $Q_0\subset \wtQ_0
\subset Q_u$ and $P_1 \subset \wtP_1 \subset P_s$. If $\wtQ_0$ and
$\wtP_1$ are $\wtI$-transverse, then $Q_0$ and $P_1$ are
$I$-transverse.
\end{Propo}

\begin{Coro}\label{coro1}
Let $\wtI \supset I$ be parameter intervals, and let $(P_0, Q_0,
n_0)$, $(P_1, Q_1, n_1)\in \cR (\wtI) \cap \cR(I)$ be such that $Q_0
\subset Q_u$, $P_1 \subset P_s$. If their parabolic composition is
allowed in $\cR (\wtI)$, it is also allowed in $\cR(I)$.
\end{Coro}

\begin{proof}
This is the case $\wtQ_0 = Q_0$, $\wtP_1 = P_1$ of the
proposition.
\end{proof}

\begin{Coro}\label{coro2}
Let $\wtI \supset I$ be parameter intervals. Then $\cR (\wtI)
\subset \cR(I)$.
\end{Coro}

\begin{rem}
This is a slight abus de langage of no consequence: properly
speaking, we mean that the restriction to $I$ of any $(P,Q,n)\in
\cR(\wtI)$ belongs to $\cR(I)$.
\end{rem}

\begin{proof}
As simple composition does not depend on the parameter interval,
Corollary~\ref{coro2} is an immediate consequence by induction on
length of property (R6).
\end{proof}

\begin{Coro}\label{coro3}
Let $\wtI \supset I$ be parameter intervals, and let $(P,Q,n)\in
\cR(\wtI)$. If $Q$ is $\wtI$-transverse, then it is also
$I$-transverse.
\end{Coro}

\begin{proof}
If $Q\cap Q_u = \emptyset$, this is obvious. Assume therefore that
$Q\subset Q_u$. Then there exists an $\wtI$-decomposition $(P_\al,
Q_\al, n_\al)_\al$ of $P_s$ by elements of $\cR (\wtI)$ such that
for all $\al$, $Q$ and $P_\al$ are either $\wtI$-transverse or
$\wtI$-separated.

First observe that $(P_\al, Q_\al, n_\al)\in \cR (I)$ and therefore
this is also an $I$-decomposition of $P_s$. By
Corollary~\ref{coro1}, if $Q$ and $P_\al$ are $\wtI$-transverse,
they are also $I$-transverse. On the other hand, it is obvious from
the definition that if $Q$ and $P_\al$ are $\wtI$-separated they are
also $I$-separated. The result follows.
\end{proof}

\subsection{Criticality and Decomposability\label{sub6.2}}
\begin{Propo}\label{propo3}
Let $I$ be a parameter interval, and let $(P,Q,n)\in \cR(I)$. If $Q$
is $I$-transverse, then $P$ is $I$-decomposable.
\end{Propo}

\begin{proof}
Let us first assume that $Q\cap Q_u = \emptyset$. Let $a\in \AL$ be
such that $Q\subset R_a$. We have
\begin{equation}
R_a \cap W^s (\La, \whR) \subset \bigcup_{(a,a')\in\cB} \Bigl
(P_{a,a'} \cap W^s (\La,\whR)\Bigr) \cup L_u;\label{eq6.1}
\end{equation}
for each $a'\in \AL$ such that $(a,a')\in \cB$, we have the simple
child of $P$:
\begin{equation}
\Bigl( P (a'), Q (a'), n+1\Bigr) = (P,Q,n) \; * \; (P_{aa'},
Q_{aa'}, 1),\label{eq6.2}
\end{equation}
and together they form by (\ref{eq6.1}) an $I$-decomposition of $P$
(the canonical one).

Let us now assume that $Q\subset Q_u$. As $Q$ is $I$-transverse,
there is an $I$-decomposition $(P_\al, Q_\al, n_\al)_\al$ of $P_s$
such that, for each $\al$, $Q$ and $P_\al$ are not $I$-critically
related. For each $\al$ such that $Q$ and $P_\al$ are
$I$-transverse, let $(P^\pm_\al, Q^\pm_\al, n_\al+n + N_0)$ be the
two elements produced by the allowed parabolic composition. Together
with the simple children defined by (\ref{eq6.2}), they form an
$I$-decomposition of $P$.
\end{proof}

\begin{Coro}\label{coro4}
Let $I$ be a $\be$-regular parameter interval and let $(P,Q,n)\in
\cR(I)$. If $P$ is $I$-critical and $|P| > |I|^\be$ or $|Q|
> |I|^\be$ for some $t\in I$, then $P$ is $I$-decomposable.
\end{Coro}

\begin{proof}
Indeed, by the very definition of regularity, $Q$ cannot be
$I$-critical.
\end{proof}

The decomposability of "fat" critical rectangles is crucial to our
analysis. As mentioned before, Corollary~\ref{coro4} will apply to
candidate intervals with $\ovbe = \be (1+\tau)^{-1}$ instead of
$\be$.

\subsection{Concavity\label{sub6.3}}
The following result is a partial converse to
Proposition~\ref{propo2}. We call {\it concavity} the corresponding
property of $\cR(I)$.

\begin{Propo}\label{propo4}
Let $I$ be a parameter interval. Let $(P_0,Q_0,n_0)$,
$(P'_0,Q'_0,n'_0)$, $(P_1,Q_1,n_1)$ $(P'_1,Q'_1,n'_1)$ be elements
of $\cR(I)$ such that
\[
Q_0 \subset Q'_0 \subset Q_u,\qquad P_1 \subset P'_1 \subset P_s.
\]
If both $Q_0\forki P'_1$ and $Q'_0\forki P_1$ hold, then $Q'_0\forki
P'_1$ also holds.
\end{Propo}

\begin{proof}
From the definition (\ref{eq3.18}) of the function $C$ and the
monotonicity of $\te$ with respect to each variable, it follows that
we have, for each $t\in I$:
\begin{alignat}{1}
&\de (Q_0, P_1)\; \geqslant \;\max \Bigl( \de (Q'_0, P_1), \de (Q_0,
P'_1)\Bigr),\label{eq6.3}\\
&\de_L (Q'_0, P_1)\; \geqslant\; \de_L (Q_0, P_1)\; \geqslant
\;\de_L (Q_0,
P'_1),\label{eq6.4}\\
&\de_R (Q_0, P'_1) \;\geqslant \de_R (Q_0, P_1) \;\geqslant \de_R
(Q'_0, P_1),\label{eq6.5}\\
&\min \Bigl( \de_{LR} (Q'_0, P_1), \de_{LR} (Q_0, P'_1)\Bigr)
\;\geqslant \;\de_{LR} (Q_0, P_1).\label{eq6.6}
\end{alignat}

By definition of the transversality relation, there exist parameter
intervals $\wtI_1$, $\wtI_2$ containing $I$, elements $(\wtP_0,
\wtQ_0, \wtn_0)$, $(\wtP'_1, \wtQ'_1, \wtn'_1)$ in $\cR (\wtI_1)$,
elements $(\wtP'_0, \wtQ'_0, \wtn'_0)$, $(\wtP_1, \wtQ_1, \wtn_1)$
in $\cR (\wtI_2)$ such that $Q_0 \subset \wtQ_0 , Q'_1 \subset \wtQ'_1, Q'_0 \subset \wtQ'_0 , P_1\subset \wtP_1$ and
\begin{eqnarray}
&\wtQ_0 \;\ovfork_{\wtI_1}\; \wtP'_1,\label{eq6.7}\\
&\wtQ'_0 \;\ovfork_{\wtI_2}\; \wtP_1.\label{eq6.8}
\end{eqnarray}

If we have either $Q'_0 \subset \wtQ_0$ or $P'_1\subset \wtP_1$, we
can already conclude that $Q'_0 \forki P'_1$. Assume therefore that
$\wtQ_0 \subset Q'_0$, $\wtP_1\subset P'_1$. Assume also for
instance that $\wtI_1 \subset \wtI_2$. We have $(\wtP_1, \wtQ_1,
\wtn_1)\in \cR (\wtI_2)$ and $\wtP_1 \subset P'_1$. By the coherence
property (Proposition~\ref{propo6} in Subsection~\ref{sub6.5}), the
element $(P'_1, Q'_1, n'_1)$ belongs to $\cR (\wtI_2)$. We will show
that
\begin{equation}
\wtQ'_0\; \ovfork_{\wtI_2} \;P'_1\label{eq6.9}
\end{equation}
which implies that $Q'_0$ and $P'_1$ are $I$-transverse. We check
properties (T1)--(T3) of Subsection~\ref{sub5.5}. First, by
(\ref{eq6.6}) and (\ref{eq6.8}), we have, for all $t\in\wtI_2$
\begin{equation}
\de_{LR} (\wtQ'_0, P'_1) \;\geqslant \;\de_{LR} (\wtQ'_0, \wtP_1)
\geqslant 2 |\wtI_2|.\label{eq6.10}
\end{equation}

Next, by (\ref{eq6.8}) and (\ref{eq6.5}), there exists $t_0\in
\wtI_2$ such that
\begin{equation}
\de_{R} (\wtQ'_0, P'_1)\; \geqslant\; \de_{R} (\wtQ'_0, \wtP_1)
\geqslant 2 |\wtQ'_0|^{1-\eta}.\label{eq6.11}
\end{equation}

Finally, by (\ref{eq6.7}), there exists $t_1\in \wtI_1 \subset
\wtI_2$ such that
\begin{equation}
\de_{L} (\wtQ_0, \wtP'_1)\; \geqslant\; 2
|\wtP'_1|^{1-\eta}.\label{eq6.12}
\end{equation}

For this value $t_1$, we then have, by (\ref{eq6.4}):
\begin{equation}
\de_L (\wtQ'_0, P'_1) \geqslant \de_L (\wtQ_0, \wtP'_1) \geqslant
2|\wtP'_1|^{1-\eta} \geqslant 2 |P'_1|^{1-\eta}.\label{eq6.13}
\end{equation}
We have proved (\ref{eq6.9}) and this concludes the proof of the
proposition.
\end{proof}

\begin{rem}
The concavity property is very helpful in the sequel. The proof of
the proposition shows why the definition of the transversality
relation had to be complicated.
\end{rem}

\begin{Coro}\label{coro5}
Let $I$ be a parameter interval, $(P_0, Q_0, n_0)$, $(P_1, Q_1,
n_1)$, $(P'_1, Q'_1, n'_1)$ be elements of $\cR(I)$ such that $Q_0
\subset Q_u$, $P_1 \subset P'_1 \subset P_s$. If $P'_1$ is
$I$-transverse and $Q_0\forki P_1$ holds, then $Q_0 \forki P'_1$
also holds.
\end{Coro}

\begin{proof}
There exists an $I$-decomposition $(P_\al, Q_\al, n_\al)$ of $Q_u$
such that, for any $\al$, $Q_\al$ and $P'_1$ are either
$I$-transverse or $I$-separated. There exists $\al$ such that
$Q_\al$ and $Q_0$ do intersect. As $Q_0 \forki P_1$ holds, $Q_\al$
and $P'_1$ must be $I$-transverse. If $Q_\al \supset Q_0$, it
follows from Proposition~\ref{propo2} that $Q_0$ and $P'_1$ are
$I$-transverse. If $Q_\al \subset Q_0$, the same conclusion follows
from Proposition~\ref{propo4}.
\end{proof}

\subsection{Children are Born From Their Parent\label{sub6.4}}
Let $I$ be a parameter interval, $(\wtP, \wtQ, \wtn)$ be an element
of $\cR(I)$ and $a\in\AL$ the letter such that $\wtQ\subset \cR_a$.
The simple children of $\wtP$ are given by formula (\ref{eq6.2}) and
parametrized by the elements $a'\in \AL$ such that $(a,a')\in \cB$.
All children are simple unless $\wtQ\subset Q_u$.

\begin{Propo}\label{propo5}
Assume $\wtQ \subset Q_u$ and let $(P,Q,n)$ be a non-simple child of
$\wtP$. Then there exists $(P_1, Q_1, n_1)\in \cR (I)$ such that
$\wtQ\;\forki P_1$ and
\[
(P, Q, n) \in (\wtP, \wtQ, \wtn) \;\; \square \;\; (P_1, Q_1, n_1).
\]
Moreover, the parent $\wtP_1$ of $P_1$ is $I$-critical.
\end{Propo}

\bigskip\noindent
{\it Proof.}

1.~~By (R6) in Subsection~\ref{eq5.3}, $(P,Q,n)$ is obtained from
shorter elements by simple or parabolic composition. Let us first
assume that we can write
\begin{equation}
(P, Q, n) = (P_0, Q_0, n_0)\;\; * \;\; (P_1, Q_1,
n_1),\label{eq6.14}
\end{equation}
with $n_0$, $n_1 > 0$. As $P$ is a non-simple child, we must have
$n_1 > 1$. Let $(\wtP_1, \wtQ_1, \wtn_1)$ be the element of $\cR(I)$
such that $\wtP_1$ is the parent of $P_1$. If $P_1$ was a simple
child, one would be able to write
\begin{equation}
(P_1, Q_1, n_1) = (\wtP_1, \wtQ_1, \wtn_1)\;\; * \;\; (P_{aa'},
Q_{aa'}, 1)\label{eq6.15}
\end{equation}
for some $(a,a')\in \cB$, and then
\begin{equation}
(P, Q, n) = (\wtP_0, \wtQ_0, \wtn_0)\;\; * \;\; (P_{aa'}, Q_{aa'},
1)\label{eq6.16}
\end{equation}
where
\begin{equation}
(\wtP_0, \wtQ_0, \wtn_0) = (P_0, Q_0, n_0)\;\; * \;\; (\wtP_1,
\wtQ_1, \wtn_1).\label{eq6.17}
\end{equation}
Thus $P$ would be a simple child, in contradiction with the
assumption of the proposition. Therefore, $P_1$ is a non-simple
child; by induction on the length, we can write
\begin{equation}
(P_1, Q_1, n_1) \in (\wtP_1, \wtQ_1, \wtn_1)\;\; \square \;\; (P_2,
Q_2, n_2)\label{eq6.18}
\end{equation}
for some $(P_2,Q_2, n_2)\in \cR(I)$. We must have
\begin{equation}
\wtQ_1 \; \forki\; P_2.\label{eq6.19}
\end{equation}
We define
\begin{equation}
(\wtP_0, \wtQ_0, \wtn_0) := (P_0, Q_0, n_0)\;\; * \;\; (\wtP_1,
\wtQ_1, \wtn_1).\label{eq6.20}
\end{equation}
We then have $\wtQ_0 \subset \wtQ_1$ and hence, from (\ref{eq6.19})
and Proposition~\ref{propo1}
\begin{equation}
\wtQ_0 \; \forki\; P_2.\label{eq6.21}
\end{equation}
Thus, the parabolic composition of $(\wtP_0, \wtQ_0, \wtn_0)$ and
$(P_2, Q_2, n_2)$ is allowed; we obviously have:
\begin{equation}
(P, Q, n)\in (\wtP_0, \wtQ_0, \wtn_0) \;\; \square \;\; (P_2, Q_2,
n_2).\label{eq6.22}
\end{equation}

\medskip\noindent
2.~~We have shown so far that it is always possible to write
\begin{equation}
(P, Q, n)\in (P_0,Q_0,n_0) \;\; \square \;\; (P_1, Q_1,
n_1)\label{eq6.23}
\end{equation}
for some allowed parabolic composition in $\cR(I)$. We take $n_0$
maximal in (\ref{eq6.23}), and assume by contradiction that $n_0 <
\wtn$. Let $(\wtP_1, \wtQ_1, \wtn_1)$ be the element of $\cR(I)$
such that $\wtP_1$ is the parent of $P_1$; let $(P'_0, Q'_0, n'_0)$
be the element of $\cR(I)$ such that $P'_0$ is the child of $P_0$
containing $P$. As $n_0 < \wtn$, we have $n'_0 < n$. As $g_t^{n_0}
(P) \subset L_u$ by (\ref{eq6.23}), $P'_0$ must be a non-simple
child. Then, from the induction hypothesis, we can write
\begin{equation}
(P'_0, Q'_0, n'_0) \in (P_0,Q_0,n_0)\;\; \square \;\; (P'_1, Q'_1,
n'_1)\label{eq6.24}
\end{equation}
for some $(P'_1,Q'_1,n'_1)\in \cR(I)$ with $P'_1 \supset P_1$. We
have thus $P'_1 \supset \wtP_1$. As $Q_0$ is $I$-transverse to
$P'_1$, it is also $I$-transverse to $\wtP_1$
(Proposition~\ref{propo2}), and we can define $(\wtP_0, \wtQ_0,
\wtn_0)\in \cR(I)$ by
\begin{equation}
(\wtP_0, \wtQ_0, \wtn_0) \in (P_0, Q_0, n_0)\;\; \square \;\;
(\wtP_1, \wtQ_1, \wtn_1),\label{eq6.25}
\end{equation}
and $\wtP_0 \supset P$. If $P_1$ was a simple child of $\wtP_1$, $P$
would be a simple child of $\wtP_0$. Therefore, by the induction
hypothesis, we can write
\begin{equation}
(P_1, Q_1, n_1) \in (\wtP_1, \wtQ_1, \wtn_1)\;\; \square \;\; (P_2,
Q_2, n_2).\label{eq6.26}
\end{equation}
for some $(P_2, Q_2, n_2)$ with $\wtQ_1 \forki P_2$. But then
$\wtQ_0 \forki P_2$ (Proposition~\ref{propo2}) and we have
\begin{equation}
(P,Q,n) \in (\wtP_0, \wtQ_0, \wtn_0)\;\; \square \;\; (P_2, Q_2,
n_2)\label{eq6.27}
\end{equation}
which contradicts the maximality of $n_0$. We have proven the first
part of the proposition.

\medskip\noindent
3.~~Assume by contradiction that the parent $\wtP_1$ of $P_1$ is
$I$-transverse. As $Q_0\forki P_1$, we have $\wtP_1\subset P_s$.
There exists, therefore, an $I$-decomposition $(P_\al,Q_\al,n_\al)$
of $Q_u$ such that $Q_\al$ and $\wtP_1$ are never critically
related. By definition of an $I$-decomposition, there exists $\al$
such that $g^{\wtn} (P) \cap Q_\al \ne \emptyset$; then $Q_\al$ and
$\wtP_1$ are not $I$-separated and they must be $I$-transverse.

As $\wtQ\cap Q_\al \ne \emptyset$, we have either $\wtQ\subset
Q_\al$, and $\wtQ\forki \wtP_1$ by Proposition~\ref{propo2}, or
$\wtQ\supset Q_\al$, and $\wtQ \forki \wtP_1$ again by concavity
(Proposition~\ref{propo4}), because both $\wtQ_\al \forki \wtP_1$
and $\wtQ \forki P_1$ hold. Thus, we can define $(\whP, \whQ, \whn)$
by
\begin{equation}
(\whP, \whQ, \whn) \in (\wtP, \wtQ, \wtn)\;\; \square \;\; (\wtP_1,
\wtQ_1, \wtn_1),\label{eq6.28}
\end{equation}
and $P\subset \whP$. We then have $P\subset \whP \subset \wtP$,
$\whP\ne P$, $\whP\not= \wtP$; thus, $P$ would not be a child of
$\wtP$. \hfill $\square$

\subsection{Coherence and Parametric Concavity\label{sub6.5}}
The coherence property, asserted in the next proposition, means that
larger rectangles are constructed before thinner ones.

\begin{Propo}\label{propo6}
Let $I\subset \wtI$ be parameter intervals, and let
$(P,Q,n)\in\cR(I)$, $(\wtP,\wtQ,\wtn)\in \cR(\wtI)$. If $\wtP\subset
P$, then $(P,Q,n)\in \cR(\wtI)$.
\end{Propo}
The next property, called parametric concavity, is another partial
converse to Proposition~\ref{propo2}; it is formally very similar to
Proposition~\ref{propo4}.

\begin{Propo}\label{propo7}
Let $I\subset I'$ be parameter intervals and let $(P_0,Q_0,n_0)$,
$(P'_0, Q'_0, n'_0)$, $(P_1,Q_1,n_1)$ be elements of $\cR (I')$ such
that $Q_0 \subset Q'_0 \subset Q_u$ and $P_1 \subset P_s$. If both $Q_0 \fork_{I'} P_1$ and $Q'_0 \forki P_1$
hold, then $Q'_0 \fork_{I'} P_1$ also holds.
\end{Propo}
Obviously there is a similar statement exchanging $P'$s and $Q'$s.
In the proof of Proposition~\ref{propo6}, we will use the following
result, which is of independent interest.

\begin{Propo}\label{propo8}
Let $I$ be a parameter interval, and let $(P,Q,n)$, $(P', Q', n')$,
$(P'',Q'',n'')$ be elements $\cR(I)$ with $P\subset P'\subset P''$,
$P'\ne P''$. We have:

a)~~If $(P,Q,n) = (P'',Q'',n'')*(P_1,Q_1,n_1)$ for some
$(P_1,Q_1,n_1)\in \cR(I)$, then
\[
(P',Q',n') = (P'',Q'',n'')\;\;*\;\; (P'_1, Q'_1, n'_1)
\]
for some $(P'_1, Q'_1, n'_1)\in \cR(I)$.

b)~~If $(P,Q,n) \in (P'',Q'',n'')\;\square\;(P_1,Q_1,n_1)$ for some
$(P_1,Q_1,n_1)\in \cR(I)$, then
\[
(P',Q',n') \in (P'',Q'',n'')\;\;\square\;\; (P'_1, Q'_1, n'_1)
\]
for some $(P'_1, Q'_1, n'_1)\in \cR(I)$.
\end{Propo}

\medskip\noindent
{\it Proof of Proposition~\ref{propo8}.} We only consider case b),
case a) being similar but easier. It is sufficient to consider the
case where $P'$ is the parent of $P$. Let $(\wtP_1, \wtQ_1, \wtn_1)$
be the element of $\cR(I)$ such that $\wtP_1$ is the parent of
$P_1$. We claim that $Q''$ and $\wtP_1$ are $I$-transverse.

Indeed, let $(\ovP, \ovQ, \ovn)\in \cR(I)$ such that $\ovP$ is the
child of $P''$ containing $P'$. As $g^{n''}(P)\subset L_u$, $\ovP$
is a non-simple child; by Proposition~\ref{propo5}, there exists
$(\ovP_1, \ovQ_1, \ovn_1)$ such that
\begin{equation}
(\ovP, \ovQ, \ovn) \in (P'', Q'', n'')\;\; \square \;\; (\ovP_1,
\ovQ_1, \ovn_1).\label{eq6.29}
\end{equation}
We then have $Q''\forki \ovP_1$ and $\ovP_1\supset \wtP_1$; the
claim then follows from Proposition~\ref{propo2}.

We can therefore define $(\wtP, \wtQ, \wtn)\in \cR(I)$ by
\begin{equation}
(\wtP, \wtQ, \wtn) \in (P'', Q'', n'')\;\; \square \;\; (\wtP_1,
\wtQ_1, \wtn_1)\label{eq6.30}
\end{equation}
and $P\subset \wtP$. Let us show that $\wtP = P'$. Otherwise, we
have $P'\subset \wtP$, and $P' \ne \wtP$. Let $(\whP,\whQ,\whn)$ be
the element of $\cR(I)$ such that $\whP$ is the child of $\wtP$
containing $P'$.

If $P_1$ was a simple child of $\wtP_1$, we would have $\wtn_1 = n_1
- 1$ and $\wtn = n-1$, forcing $\wtP = P'$. Hence, $P_1$ is a
non-simple child; this implies that $g^{\wtn}(P)\subset L_u$ and
that $\whP$ is also a non-simple child. By Proposition~\ref{propo5},
there exist $(P_2,Q_2,n_2)$, $(\whP_2, \whQ_2, \whn_2)$ in $\cR(I)$
such that $P_2 \subset \whP_2$ and
\begin{eqnarray}
(P_1, Q_1, n_1) &\in &(\wtP_1, \wtQ_1, \wtn_1) \;\;\square \;\;
(P_2,Q_2,n_2),\label{eq6.31}\\
(\whP, \whQ, \whn) &\in &(\wtP, \wtQ, \wtn) \;\;\square \;\;
(\whP_2,\whQ_2,\whn_2).\label{eq6.32}
\end{eqnarray}
Then both $\wtQ_1 \forki P_2$ and $\wtQ \forki \whP_2$ hold. By
concavity (Proposition~\ref{propo4}), we also have $\wtQ_1 \forki
\whP_2$.  We then define $(\whP_1, \whQ_1, \whn_1)$ by:
\begin{equation}
(\whP_1, \whQ_1, \whn_1) \in (\wtP_1, \wtQ_1, \wtn_1) \;\; \square
\;\; (\whP_2,\whQ_2,\whn_2)\label{eq6.33}
\end{equation}
and $P_1 \subset \whP_1$. Then $P_1 \subset \whP_1 \subset
\wtP_1$ and $\whP_1$ is distinct from $P_1$ and $\wtP_1$, which
contradicts that $\wtP_1$ is the parent of $P_1$. \hfill $\square$

\medskip\noindent
{\it Proof of Proposition~\ref{propo6}.} It is sufficient to
consider the case where $P$ is the parent of $\wtP$ in $\cR(I)$. Let
$(P_0,Q_0,n_0)$ be the element of $\cR(\wtI)$ such that $P_0$ is the
parent of $\wtP$ in $\cR (\wtI)$. We want to show that $P_0 = P$. This is clear if $\tilde P$ is a simple child of $P_0$.  We assume therefore that $\tilde P$ is a non-simple child of $P_0$. We
know that $P\subset P_0$. By Proposition~\ref{propo5}, there exists
$(\wtP_1, \wtQ_1, \wtn_1)\in \cR(\wtI)$ such that
\begin{equation}
(\wtP, \wtQ, \wtn) \in (P_0, Q_0, n_0) \;\; \square \;\;
(\wtP_1,\wtQ_1,\wtn_1)\label{eq6.34}
\end{equation}
with $Q_0\fork_{\wtI} \wtP_1$. By Proposition~\ref{propo8}, there
exists $(P_1,Q_1,n_1)\in \cR(I)$ such that
\begin{equation}
(P, Q, n) \in (P_0, Q_0, n_0) \;\; \square \;\;
(P_1,Q_1,n_1)\label{eq6.35}
\end{equation}
with $Q_0\forki P_1$. By induction on the length, as $\wtP_1\subset
P_1$ we must have $(P_1,Q_1,n_1)\in\cR(\wtI)$. By parameter
concavity (Proposition~\ref{propo7}), as both $Q_0 \fork_{\wtI}
\wtP_1$ and $Q_0 \forki P_1$ hold, $Q_0 \fork_{\wtI} P_1$ must also
hold, which implies $(P,Q,n)\in \cR(\wtI)$. \hfill $\square$

\medskip\noindent
{\it Proof of Proposition~\ref{propo7}.} By definition of the
transversality relation, there exist parameter intervals
$\wtI\supset I$, $\wtI'\supset I'$ and elements $(\wtP_0, \wtQ_0,
\wtn_0)$, $(\wtP_1, \wtQ_1, \wtn_1)\in \cR (\wtI')$, $(\wtP'_0,
\wtQ'_0, \wtn'_0)$, $(\wtP'_1, \wtQ'_1, \wtn'_1)\in \cR(\wtI)$ such
that $Q_0 \subset \wtQ_0$, $Q'_0 \subset \wtQ'_0$, $P_1\subset
\wtP_1$, $P_1\subset \wtP'_1$ and
\begin{eqnarray}
&\wtQ_0 \; \fork_{\wtI'}\; \wtP_1,\label{eq6.36}\\
&\wtQ'_0 \; \fork_{\wtI}\; \wtP'_1\label{eq6.37}
\end{eqnarray}
both hold. If either $Q'_0\subset \wtQ_0$ or $I'\subset \wtI$, we
conclude immediately that $Q'_0 \fork_{I'} P_1$ holds. Assume
therefore that $\wtQ_0 \subset Q'_0$ and $\wtI\subset I'$. Let
$P^*_1$ be the largest of $\wtP_1$, $\wtP'_1$.

If $\wtP_1 \subset \wtP'_1$, $(\wtP'_1, \wtQ'_1, \wtn'_1)\in \cR
(\wtI')$ by coherence (Proposition~\ref{propo6}), hence
$(P^*_1,Q^*_1,n^*_1)$ always belong to $\cR(\wtI')$. We will show
that
\begin{equation}
\wtQ'_0 \; \ovfork_{\wtI'}\; P^*_1\label{eq6.38}
\end{equation}
holds, which implies $Q'_0 \fork_{I'} P_1$.

We check properties (T1)--(T3) of Subsection~\ref{sub5.4}. For all
$t\in \wtI'$, we have by (\ref{eq6.36}) and (\ref{eq6.6})
\begin{equation}
\de_{LR} (\wtQ'_0, P^*_1)\;\geqslant\; \de_{LR} (\wtQ_0,\wtP_1)\;
\geqslant\; 2|\wtI'|.\label{eq6.39}
\end{equation}
By (\ref{eq6.37}), there exists $t_0\in \wtI\subset I'\subset \wtI'$
such that
\begin{equation}
\de_{R} (\wtQ'_0, \wtP'_1)\;\geqslant\; 2 |\wtQ'_0|^{1-\eta}.
\label{eq6.40}
\end{equation}
Then, by (\ref{eq6.5}), for the same $t_0$, we have
\begin{equation}
\de_{R} (\wtQ'_0, P^*_1)\;\geqslant\; \de_R  (\wtQ'_0, \wtP'_1)\;
\geqslant\; 2 |\wtQ'_0|^{1-\eta}. \label{eq6.41}
\end{equation}
When $P^*_1 = \wtP'_1$, it follows directly from (\ref{eq6.37}) that
we have
\begin{equation}
\de_{L} (\wtQ'_0, P^*_1)\;\geqslant\; 2
|P^*_1|^{1-\eta}\label{eq6.42}
\end{equation}
for some $t_1\in \wtI\subset \wtI'$.

When $P^*_1 = \wtP_1$, we use (\ref{eq6.36}) and (\ref{eq6.4}) to
find $t_1 \in \wtI'$ such that
\begin{equation}
\de_{L} (\wtQ'_0, P^*_1)\;\geqslant\; \de_L  (\wtQ_0, \wtP_1)\;
\geqslant\; 2 |P^*_1|^{1-\eta}.\label{eq6.43}
\end{equation}
We have thus proved (\ref{eq6.38}).\hfill $\square$

\bigskip\noindent
{\bf Remark.}

1.~~The coherence property means in particular that the parent-child
relation does not depend on the parameter interval (once both parent
and child are defined).

2.~~We have presented together Propositions~\ref{propo6},
\ref{propo7} and \ref{propo8} because the proofs are interconnected.
But, actually, coherence (Proposition~\ref{propo6}) was already used
in the proof of Proposition~\ref{propo4}. Logically speaking, all
properties in Sections~\ref{sec5}--\ref{sec7} should be proved
together (as was already mentioned in Section~\ref{sec5}).

\subsection{Further Criteria for Transversality\label{sub6.6}}
In this subsection, we give other sufficient conditions for
transversality that can be seen as partial converses to
Proposition~\ref{propo2}.

\begin{Propo}\label{propo9}
Let $I$ be a parameter interval and let $(P,Q,n)$, $(P_0,Q_0,n_0)$,
$(P_1,Q_1,n_1)$ be elements of $\cR(I)$ such that $Q\subset Q_u$ and
$P_0 \subset P_1 \subset P_s$. Assume that $Q\forki P_0$ holds and
that $2|P_1|^{1-\eta} \leqslant |I|$ for some $t_1\in I$. Then $Q$
and $P_1$ are also $I$-transverse.
\end{Propo}

{\it Proof.~} By definition of the transversality relation, there
exist $\wtI\supset I$, $(\wtP,\wtQ,\wtn)$, $(\wtP_0,\wtQ_0,\wtn_0)\in
\cR(\wtI)$ such that $Q\subset \wtQ$, $P_0\subset \wtP_0$ and
$\wtQ\ovfork_{\wtI} \wtP_0$.

If $P_1\subset \wtP_0$ this already implies that $Q\forki P_1$. Let
us assume that $\wtP_0\subset P_1$. We will show that
$\wtQ\ovfork_{\wtI} P_1$ holds. By coherence
(Proposition~\ref{propo4}), we have $(P_1,Q_1,n_1)\in \cR(\wtI)$.
Let us check (T1)--(T3).

By (T1) for $\wtQ$, $\wtP_0$ and (\ref{eq6.6}), we have, for all
$t\in \wtI$:
\begin{equation}
\de_{LR} (\wtQ, P_1)\;\geqslant\; \de_{LR} (\wtQ,\wtP_0)\;
\geqslant\; 2|\wtI|.\label{eq6.44}
\end{equation}
By (T2) for $\wtQ$, $\wtP_0$ and (\ref{eq6.5}), there exists $t_0\in
\wtI$ such that
\begin{equation}
\de_{R} (\wtQ, P_1)\;\geqslant\; \de_{R} (\wtQ,\wtP_0)\; \geqslant\;
2|\wtQ|^{1-\eta}.\label{eq6.45}
\end{equation}
Finally, we have, for all $t\in I$, by (\ref{eq5.14})
\begin{eqnarray}
\de_L (\wtQ,P_1) &\geqslant &\de_{LR} (\wtQ, P_1) - C|P_1|\label{eq6.46}\\
&\geqslant &2|I| - C|P_1|.\notag
\end{eqnarray}
But, for $t=t_1$, we have, if $\vep_0$ is small enough
\begin{equation}
2 |I| - C |P_1| \geqslant 4 |P_1|^{1-\eta} - C |P_1| \geqslant 2
|P_1|^{1-\eta},\label{eq6.47}
\end{equation}
\hfill $\square$

\begin{Propo}\label{propo10}
Let $I$ be a parameter interval and let $(P,Q,n)$, $(P_0,Q_0,n_0)$,
$(P_1,Q_1,n_1)$ be elements of $\cR(I)$ such that $Q\subset Q_u$,
$P_0 \subset P_1 \subset P_s$. Assume that $Q\forki P_0$ holds and
that $|P_1| \leqslant \fud \;|Q|$ for all $t\in I$. Then $Q$ and
$P_1$ are also $I$-transverse.
\end{Propo}

\begin{proof}
The argument is the same as in Proposition~\ref{propo9}, with a
slight difference to check (T3) for $\wtQ$, $P_1$. We use
(\ref{eq5.13}), (\ref{eq5.14}) to obtain, for the value $t_0$ of the
parameter given by (T2) for $\wtQ$, $\wtP_0$:
\begin{eqnarray}
\de_L (\wtQ,P_1) &\geqslant &\de_{R} (\wtQ, P_1) - C|P_1|,\label{eq6.48}\\
&\geqslant &2 |Q|^{1-\eta} - C|P_1|,\notag\\
&\geqslant &2 (2 |P_1|)^{1-\eta} - C|P_1|,\notag\\
&\geqslant &2 |P_1|^{1-\eta},\notag
\end{eqnarray}
if $\vep_0$ is small enough.
\end{proof}

\begin{Propo}\label{propo11}
Let $I\subset I'$ be parameter intervals and let $(P_0,Q_0,n_0)$,
$(P_1,Q_1,n_1)$ be elements of $\cR(I')$ such that $Q_0\subset Q_u$,
$P_1 \subset P_s$. Assume that $Q_0\forki P_1$ holds and that we
have $2 |I'| < |P_1|^{1-\eta}$, for all $t\in I$. Then $Q_0$ and
$P_1$ are also $I'$-transverse.
\end{Propo}

\begin{proof}
By definition of the transversality relation, there exists
$\wtI\supset I$, $(\wtP_0, \wtQ_0, \wtn_0)$, $(\wtP_1, \wtQ_1,
\wtn_1)\in \cR(\wtI)$ such that $\wtQ_0 \supset Q_0$, $\wtP_1\supset
P_1$ and $\wtQ_0 \ovfork_{\wtI} \wtP_1$ holds. If $\wtI\supset I'$,
this already implies that $Q_0 \fork_{I'} P_1$ holds. Assume thus
that $\wtI\subset I'$. We show that $\wtQ_0 \ovfork_{I'} \wtP_1$.

Conditions (T2), (T3) in $I'$ follow from the same conditions in
$\wtI$. For the value $t_1$ given by (T3), we have
\begin{equation}
\de_{LR} (\wtQ_0, \wtP_1)\; \geqslant\; \de_L (\wtQ_0, \wtP_1)
\;\geqslant \; 2 |P_1|^{1-\eta} \;\geqslant\; 4 |I'|.\label{eq6.49}
\end{equation}

By Corollary 11 in Subsection~\ref{sub7.6}, this implies
\begin{equation}
\de_{LR} (\wtQ_0, \wtP_1) \geqslant 2|I'|,\qquad \forall t\in
I'.\label{eq6.50}
\end{equation}
which is (T1).
\end{proof}

\subsection{A Structure Theorem for New Rectangles\label{sub6.7}}
{\bf 6.7.1~Associativity of Parabolic Composition.~} Let $I$ be a
parameter interval, and let $(P_0,Q_0,n_0)$, $(P_1,Q_1,n_1)$,
$(P_2,Q_2,n_2)$ be elements in $\cR(I)$ such that $Q_0\subset Q_u$,
$Q_1\subset Q_u$, $P_1\subset P_s$, $P_2\subset P_s$. We assume that
both $Q_0 \forki P_1$ and $Q_1 \forki P_2$ hold.

Parabolic composition of $(P_0, Q_0, n_0)$, $(P_1, Q_1, n_1)$
produces two elements $(\za{P}{01}$, $\za{Q}{01}$, $\za{n}{01})$,
$(\zb{P}{01}, \zb{Q}{01}, \zb{n}{01})$.

As $\za{Q}{01}$ and $\zb{Q}{01}$ are contained in $Q_1$, it follows
from Prop.~2 that both $\za{Q}{01} \forki P_2$ and $\zb{Q}{01}
\forki P_2$ hold.

In the same way, parabolic composition of $(P_1,Q_1,n_1)$,
$(P_2,Q_2,n_2)$ produces two elements $(\za{P}{12}$, $\za{Q}{12},
\za{n}{12})$, $(\zb{P}{12}, \zb{Q}{12}, \zb{n}{12})$ such that both
$Q_0 \forki \za{P}{12}$ and $Q_0 \forki \zb{P}{12}$ hold.

It is clear that the four elements of $\cR(I)$ obtained by parabolic
composition of $(\za{P}{01}$, $\za{Q}{01}, \za{n}{01})$ or
$(\zb{P}{01}, \zb{Q}{01}, \zb{n}{01})$ with $(P_2, Q_2, n_2)$ are
the same as the four elements obtained by the parabolic composition
of $(P_0, Q_0, n_0)$ with $(\za{P}{12}$, $\za{Q}{12}$, $n^+_{12})$
or $(\zb{P}{12}$, $\zb{Q}{12}$, $\zb{n}{12})$. Their domains are the
components of $P_0\cap (G_t\circ g_t^{n_0})^{-1} P_1 \cap (G_t \circ
g_t^{n_1} \circ G_t \circ g_t^{n_0})^{-1} P_2$. If $(P,Q,n)$ is any
of these four elements, we will write
\begin{equation}
(P,Q,n) \in (P_0, Q_0, n_0) \; \square \; (P_1, Q_1, n_1) \; \square
\; (P_2, Q_2, n_2).\label{eq6.51}
\end{equation}
The same considerations extend immediately, by induction on $k$, to
the case of elements $(P_0, Q_0, n_0)$, $(P_k, Q_k, n_k)$ such that
$P_i\subset P_s$ for $0 < i < k$, $Q_i \subset Q_u$ for $0\le i <
k$, and $Q_i \forki P_{i+1}$ holds for $0\le i < k$. Then the
successive parabolic compositions of $(P_0, Q_0, n_0)$, $\cdots$,
$(P_k,Q_k,n_k)$ produce $2^k$ elements and we will write for any
such element $(P,Q,n)$:
\begin{equation}
(P,Q,n) \in (P_0, Q_0, n_0) \; \square \; \cdots \; \square \; (P_k,
Q_k, n_k).\label{eq6.52}
\end{equation}

{\bf 6.7.2~Statement of the Structure Theorem.~} We have seen in
Subsection~\ref{sub5.5} that parabolic composition is never allowed
in the class $\cR (I_0)$ associated to the starting interval $I_0 =
[\vep_0, 2\vep_0]$. This class consists exactly of the affine-like
iterates associated to the Markov partition of the initial horseshoe
$K_{g_t}$.

On the other hand, for elements $(P,Q,n)$ belonging to some class
$\cR(I)$ but which are not (restrictions of) an element of
$\cR(I_0)$, parabolic composition must occur. The following theorem
gives some rather precise information on this process.

\begin{Theo}\label{theo1}
Let $I$, be a parameter interval of level $>0$, $\wtI$ be the parent
interval, and let $(P,Q,n)$ be an element of $\cR(I)$ which is not
(the restriction of) an element of $\cR(\wtI)$. Then there exists
$k>0$, elements $(P_0,Q_0,n_0)$, $\cdots$, $(P_k,Q_k,n_k)$ of
$\cR(\wtI)$ such that $Q_i \subset Q_u$ for $0\le i < k$,
$P_i\subset P_s$ for $0 < i \le k$, $Q_i \forki P_{i+1}$ holds for
$0\le i < k$, $Q_i \fork_{\wtI} P_{i+1}$ does not hold for $0\le i <
k$ and
\[
(P,Q,n) \in (P_0, Q_0, n_0) \; \square \; \cdots \; \square \; (P_k,
Q_k, n_k).
\]
Moreover, these elements are uniquely determined by these
conditions, $P_i$ is $\wtI$-critical for $0< i \le k$ and $Q_i$ is
$\wtI$-critical for $0\le i < k$.
\end{Theo}

\bigskip\noindent
{\bf 6.7.3~} We will first introduce a concept, relative to an
element $(P,Q,n)$ as in the theorem above, that leads to the
determination of the $(P_i, Q_i, n_i)$.

Let $m$, $p$ be integers such that $0\le m\le p\le n$. We say that
$[m,p]$ is an $\wtI$-interval if there exists $(\wtP, \wtQ, \wtn)
\in \cR (\wtI)$ such that
\[
g^m_t (P) \subset \wtP \; \text{ for all }\; t\in I\; \text{ and }\;
\wtn = p - m.
\]

\begin{lema}\label{lema1}
The union of two $\wtI$-intervals with non empty intersection is an
$\wtI$-interval.
\end{lema}

\bigskip\noindent
{\it Proof.~} Let $[m,p]$, $[m',p']$ be these $\wtI$-intervals, and
let $(\wtP, \wtQ, \wtn)$, $(\wtP',\wtQ', \wtn')$ be the
corresponding elements of $\cR(\wtI)$. Without loss of generality,
we may assume that $m< m' \le p < p'$. Replacing if necessary $\wtP$
by a larger rectangle, we also assume that the element $(\whP, \whQ,
\whn)$ of $\cR(\wtI)$ such that $\whP$ is the parent of $\wtP$
satisfies $m + \whn < m'$. There are now two cases:

a)~~$p = m'$.

Let $R_a$ be the rectangle containing $\wtQ$. Then $R_a \supset \wtQ
= g^{\wtn} (\wtP) \supset g^p (P) = g^{m'} (P)$; thus $\wtP'$ is
also contained in $R_a$ and the simple composition
\begin{equation}
(\wtP'', \wtQ'', \wtn'') := (\wtP, \wtQ, \wtn) \; * \; (\wtP',
\wtQ', \wtn')\label{eq6.53}
\end{equation}
is defined. We have $m + \wtn'' = p'$ and $g_t^m (P) \subset P''$.

b)~~$p > m'$.

Then, $\wtP$ is not a simple child of $\whP$, because otherwise we
would have $\whn = \wtn - 1 \ge m' - m$. By
Proposition~\ref{propo5}, there exists $(\wtP_0, \wtQ_0, \wtn_0)$ in
$\cR (\wtI)$ such that
\begin{equation}
(\wtP, \wtQ, \wtn) \in (\whP, \whQ, \whn) \; \square \; (\wtP_0,
\wtQ_0, \wtn_0).\label{eq6.54}
\end{equation}
The element $(\wtP_0, \wtQ_0, \wtn_0)$ of $\cR(\wtI)$ is associated
to the $\wtI$-interval $[\whm, p]$, where $\whm =  m + \whn + N_0$.
We have $m + \whn < m'$ and $g_t^{m+\whn} (P) \subset L_u$, hence
also $m + \whn + N_0 = \whm \le m'$.

\medskip
To conclude the proof, we argue by induction on the total length $p'
- m$ of the interval considered. The case $p' - m = 0$ is trivial.
In the other case, we have the $\wtI$-intervals $[\whm,p]$ and
$[m',p']$ with $m < \whm \le m'$ and hence by induction $[\whm, p']$
is an $\wtI$-interval. Let $(\wtP_1, \wtQ_1, \wtn_1)$ be the
corresponding element of $\cR (\wtI)$; we have $g_t^\whm (P) \subset
\wtP_1 \subset \wtP_0$. From (\ref{eq6.54}), $\whQ \fork_{\wtI}
\wtP_0$ holds, hence also does $\whQ \fork_{\wtI} \wtP_1$ by
Proposition~\ref{propo2}. Then, the parabolic composition of $(\whP,
\whQ, \whn)$ and $(\wtP_1, \wtQ_1, \wtn_1)$ is allowed and defines
an element of $\cR (\wtI)$ which guarantees that $[m,p']$ is an
$\wtI$-interval.\hfill $\square$

\bigskip\noindent
{\bf 6.7.4~} We will now show that the $(P_i, Q_i, n_i)$ in the
theorem are uniquely determined by their properties. Indeed, define
$m_0 = 0$, $p_0 = n_0$ and for $i>0$:
\begin{equation}
m_i = p_{i-1} + N_0,\qquad  p_i = m_i + n_i.\label{eq6.55}
\end{equation}

\begin{lema}\label{lema2}
The maximal $\wtI$-intervals are exactly the $[m_i,p_i]$, $0\le i\le
k$, with associated elements $(P_i,Q_i,n_i)$.
\end{lema}

\begin{proof}
First, the $[m_i,p_i]$ are indeed $\wtI$-intervals with associated
elements $(P_i,Q_i,n_i)$. To complete the proof, it is sufficient to
show that no $\wtI$-interval $[m,p]$ can intersect a gap $(p_i, \ism
m)$. Assume by contradiction that there exists such a $[m,p]$ with
associated element $(\wtP,\wtQ,\wtn)$ and minimal $\wtn = p-m$. As
$g^\ell(P)\subset g^{\ell-p_i}(L_u)$ does not intersect $R$ for
$p_i< \ell < \ism m$, we must have $m\leqslant p_i$ and $\ism
m\leqslant p$. By property (R6) of $\cR(\wtI)$
(Subsection~\ref{eq5.3}) and the minimality of $\wtn$, there exists
$(\wtP_0,\wtQ_0,\wtn_0)$, $(\wtP_1, \wtQ_1, \wtn_1)$ in $\cR (\wtI)$
such that
\begin{equation}
(\wtP,\wtQ,\wtn) \; \in \; (\wtP_0, \wtQ_0, \wtn_0) \; \square\;
(\wtP_1, \wtQ_1, \wtn_1)\label{eq6.56}
\end{equation}
with $\wtn_0 = p_i - m \leqslant n_i$, $\wtn_1 = p - \ism m
\leqslant \ism n$. But then, from $\wtQ_0 \supset Q_i$, $\wtP_1
\supset \ism P$ and $\wtQ_0 \fork_{\wtI} \wtP_1$, we deduce from
Proposition~\ref{propo2} that $Q_i \fork_{\wtI} P_{i+1}$ holds, a
contradiction.
\end{proof}

\medskip\noindent
{\bf 6.7.5~} Lemma~\ref{lema2} allow us to {\it define} $k$ as being
the number of maximal $\wtI$-intervals minus one, and to define the
$(P_i,Q_i,n_i)\in \cR(\wtI)$ as the elements of $\cR (\wtI)$
associated to the successive maximal $\wtI$-intervals. Observe that
the maximal $\wtI$-intervals $[m_i,p_i]$, $(0\leqslant i\leqslant
k)$ must indeed satisfy $m_0 = 0$, $\ism m = p_i + N_0$ for
$0\leqslant i < k$: ~~every $\ell\in [0,n]$ not contained in an
$\wtI$-interval is such that $g_t^{\ell-N} (P)\subset L_u$ for some
$0 < N < N_0$ and then no $\wtI$-interval intersects with $(\ell-N,
\ell-N+N_0)$, while $\{[\ell-N, \ell-N + N_0]\}$ are
$\wtI$-intervals. We observe also that $Q_i \fork_{\wtI} P_{i+1}$ does not hold because otherwise $[m_i, \ism p]$ would be an
$\wtI$-interval.

\medskip\noindent
{\bf 6.7.6~} Let $0\leqslant i < k$. Let us assume by induction over
$i$ that $P_j$ is $\wtI$-critical for $0 < j \leqslant i$, $Q_j$ is
$\wtI$-critical for $0\leqslant j < i$, $Q_j \forki P_{j+1}$ holds
for $0\leqslant j < i$ and that we have an element of $\cR (I)$:
\begin{equation}
(P^\pii, Q^\pii, p_i) \in (P_0, Q_0, n_0) \; \square \cdots \;
\square (P_i, Q_i, n_i)\label{eq6.57}
\end{equation}
such that $P\subset P^\pii$. The assumption is vacuously true for
$i=0$. We will prove it at step $i+1$. For $i=k$, it gives the
properties stated in the theorem for the $(P_i,Q_i,n_i)$.

\medskip\noindent
{\bf 6.7.7~} We first prove that $Q_i$ is $\wtI$-critical. Assume by
contradiction that $Q_i$ is $\wtI$-transverse. Then $P_i$ is
$\wtI$-decomposable. Let $(\whP_i, \whQ_i, \whn_i)$ be an element of
$\cR (\wtI)$ such that $\whP_i$ is a child of $P_i$ intersecting
$g^{m_i} (P\cap \La)$. On the other hand, let $(\whP^\pii,
\whQ^\pii, \whn^\pii)$ be the element of $\cR (I)$ such that
$\whP^\pii$ is the child of $P^\pii$ containing $P$. We apply
Proposition~\ref{propo5} twice. We find $(\ism\wtP, \ism\wtQ,
\ism\wtn)$ in $\cR (I)$, $(\ism\wtP', \ism\wtQ', \ism\wtn')$ in $\cR
(\wtI)$ such that
\begin{eqnarray}
&(\whP_i, \whQ_i, \whn_i) \in (P_i, Q_i, n_i) \; \square \;
(\ism\wtP', \ism\wtQ', \ism\wtn'),\label{eq6.58}\\
&(\whP^\pii, \whQ^\pii, \whn^\pii) \in (P^\pii, Q^\pii, p_i) \;
\square \; (\ism\wtP, \ism\wtQ, \ism\wtn).\label{eq6.59}
\end{eqnarray}

If we had $\ism\wtn' > \ism\wtn$, from $Q_i\fork_{\wtI} \ism\wtP'$
and $Q^\pii \forki \ism\wtP$, we would deduce first by
Proposition~\ref{propo2} that $Q_i \forki \ism\wtP'$, then by
concavity (Proposition~\ref{propo4}) that $Q_i \forki \ism\wtP$;
parabolic composition would yield an element $(\ovP_i, \ovQ_i,
\ovn_i)\in \cR (I)$ with $P_i {\supsetneqq} \ovP_i {\supsetneqq}
\hat P_i$, in contradiction with coherence (Proposition~\ref{propo6})
and the definition of $\whP_i$.

Therefore, we must have $\ism\wtn' \leqslant \ism\wtn$; as
$\whn^\pii \leqslant n$, we have also $m_i + \whn_i \leqslant n$,
and $[m_i, m_i + \whn_i]$ is an $\whI$-interval strictly larger than
$[m_i, p_i]$; this contradiction shows that $Q_i$ is
$\wtI$-critical.

\medskip\noindent
{\bf 6.7.8~} The proof that $\ism P$ is $\wtI$-critical is rather
similar. We assume by contradiction that it is $\wtI$-transverse.
Then $\ism Q$ is $\wtI$-decomposable and we can find $(\ism P^*,
\ism Q^*, \ism n^*) \in \cR (\wtI)$ such that $\ism Q^*$ is a child
of $\ism Q$ intersecting $g^{\ism p} (P\cap \La)$. By
Proposition~\ref{propo5}, there exists $(P^\pii_*, Q^\pii_*,
n^\pii_*) \in \cR (\wtI)$ such that
\begin{equation}
(\ism P^*, \ism Q^*, \ism n^*) \in (P^\pii_*, Q^\pii_*, n^\pii_*) \;
\square \; (\ism P, \ism Q, \ism n).\label{eq6.60}
\end{equation}
If we had $n^\pii_* > p_i$, we would derive a contradiction as
follows: we should have $Q^\pii \forki \ism P$ either (if $\ism \wtn
\leqslant \ism n$) from $Q^\pii \forki \ism\wtP$ by
Proposition~\ref{propo2} or (if $\ism \wtn \geqslant \ism n$) from
$Q^\pii \forki \ism\wtP$ and $Q^\pii_* \forki \ism P$ by
Proposition~\ref{propo4}; then, parabolic composition of $(P^\pii,
Q^\pii, p_i)$ and $(\ism P,\newline \ism Q, \ism n)$ produces an
element $(P^{(i+1)}, Q^{(i+1)}, \ism p)$ in $\cR (I)$ with $\ism Q
{\supsetneqq} Q^{(i+1)} {\supsetneqq} \ism Q^*$, which is not
compatible with coherence (Proposition~\ref{propo6}).

Thus, we must have $n^\pii_* \leqslant p_i$; then we also have $\ism
n^* \leqslant \ism p$ and $[\ism p - \ism n^*, \ism p]$ is an
$\wtI$-interval strictly larger than $[\ism m, \ism p]$. This
contradiction shows that $\ism P$ is $\wtI$-critical.

\medskip\noindent
{\bf 6.7.9~} We now prove that $Q^\pii$ and $P_{i+1}$ are
$I$-transverse. If $\ism \wtn \le \ism n$, we have $Q^\pii \forki
\ism\wtP$ by (\ref{eq6.59}) and thus also $Q^\pii\forki \ism P$ by
Proposition~\ref{propo2}. Let  us assume that $\ism\wtn > \ism n$.
In this case, we claim that $\ism Q$ is $\wtI$-critical. Indeed, if
it was $\wtI$-transverse, $\ism P$ would be $\wtI$-decomposable
and we would find an element $(\ism\whP, \ism\whQ, \ism\whn)$ of
$\cR(\wtI)$ such that $\ism\whP$ is a child of $\ism P$ intersecting
$g^{\ism m} (P\cap \La)$. By coherence (Proposition~\ref{propo6}),
we should have $\ism \whP \supset \ism\wtP$ and $[\ism m, \ism m +
\ism \whn]$ would be an $\wtI$-interval larger than $[\ism m, \ism
p]$, a contradiction.

As $(\ism P, \ism Q, \ism n)$ is $\wtI$-bicritical, and the parent
interval $\wtI$ is always assumed to be $\be$-regular, we have, for
all $t\in \wtI$
\begin{equation}
|\ism P| \; < \; |\wtI|^\be\label{eq6.61}
\end{equation}
and thus also (with $\vep_0$ small enough)
\begin{equation}
2 |\ism P|^{1-\eta} \; < \; |I|.\label{eq6.62}
\end{equation}
It now follows from Proposition~\ref{propo9} and $Q^\pii \forki
\ism\wtP$ that $Q^\pii$ and $\ism P$ are $I$-transverse.

\medskip\noindent
{\bf 6.7.10~} When $i=0$, $Q^\pzz = Q_0$ and we have already shown
that $Q_i$ and $\ism P$ are $I$-transverse.

When $i>0$, $(P_i, Q_i, n_i)$ is $\wtI$-bicritical and, therefore,
we have for all $t\in \wtI$:
\begin{equation}
|Q_i|\; < \; |\wtI|^\be,\label{eq6.63}
\end{equation}
and thus also
\begin{equation}
2 |Q_i|^{1-\eta} \; < \; |I|.\label{eq6.64}
\end{equation}
It follows from Proposition~\ref{propo9} and $Q^\pii \forki \ism P$
that $Q_i$ and $\ism P$ are $I$-transverse.

To conclude the induction step of 6.7.6, we simply observe that the
parabolic composition of $(P^\pii, Q^\pii, n^\pii)$ and $(\ism P,
\ism Q, \ism n)$ is allowed in $\cR (I)$; it produces an element
$(P^{(i+1)}, Q^{(i+1)}, \ism p)\in \cR (I)$ such that $P^{(i+1)}$ intersects
$P$ and therefore contains $P$.

The proof of the theorem is now complete.

\medskip\noindent
{\bf 6.7.11~} In the next two Corollaries, the setting and notations
are those of Theorem~\ref{theo1}.

\begin{Coro}\label{coro6}
For all $t\in I$, we have
\[
|P| \; \leqslant \; C^k |P_0| \; |P_1| \cdots |P_k|\;
|I|^{-\tfrac{k}{2}}
\]
\end{Coro}

\begin{proof}
We have $P = P^{(k)}$, with $P^\pii$ defined in (\ref{eq6.57}); we
prove that, for all $t\in I$
\begin{equation}
|P^\pii| \; \leqslant \; C^i |P_0| \cdots |P_i|\;
|I|^{-\tfrac{i}{2}}.\label{eq6.65}
\end{equation}
As $P^\pzz = P_0$, this is true for $i=0$. The induction step is a
consequence of the key estimate (\ref{eq3.27}) for parabolic
composition if we know that, for all $t\in I$:
\begin{equation}
\de (Q^\pii, \ism P) \; \geqslant \; |I|.\label{eq6.66}
\end{equation}
As $Q_i$ and $\ism P$ are $I$-transverse, there exists $I^*\supset
I$, $(P^*_i, Q^*_i, n^*_i)$, $(\ism P^*, \ism Q^*, \ism n^*) \in \cR
(I^*)$ such that $Q_i \subset Q^*_i$, $\ism P\subset \ism P^*$ and
$Q^*_i \ovfork_{I^*} \ism P^*$ holds. From (T1) in
Subsection~\ref{sub5.4}, we have, for all $t\in I^*$
\begin{equation}
\de_{LR} (Q^*_i, \ism P^*) \; \geqslant \; 2 |I^*|.\label{eq6.67}
\end{equation}
From (\ref{eq6.3}), (\ref{eq5.11})--(\ref{eq5.14}) and (R7), we get,
for all $t\in I$:
\begin{alignat}{1}
\de (Q^\pii, \ism P) \; &\geqslant \; \de (Q^*_i, \ism
P^*)\label{eq6.68}\\
&\geqslant \; \de_{LR} (Q^*_i, \ism P^*) - c (|\ism P^*| + |Q^*_i|)\notag\\
&\geqslant \; \fud\;\de_{LR} (Q^*_i, \ism P^*)\notag\\
&\geqslant \; |I^*| \; \geqslant \; |I|,\notag
\end{alignat}
which concludes the proof.
\end{proof}

\begin{Coro}\label{coro7}
For all $0 < i < k$, and $t\in \wtI$, we have
\[
|P_i| \; <  \; |\wtI|^{\be}.
\]
For all $t\in I$, we also have
\[
|P_k| \; <  \; c |\wtI|^{(1-\eta)^{-1}}.
\]
\end{Coro}

\begin{proof}
The first assertion is an immediate consequence of the regularity of
the parameter interval $\wtI$, because $(P_i, Q_i, n_i)$ is
$\wtI$-bicritical for $0 < i < k$.

For the second assertion, we first observe that in the proof of
Corollary~\ref{coro6}, we must have $I^* = I$ because $\ksM Q
\fork_{\wtI} P_k$ does not hold; for the same reason, there exists
$t^*\in \wtI$ such that
\begin{equation}
\de_{LR} (\ksM Q^*, P^*_k) \; < \; 2 |\wtI|.\label{eq6.69}
\end{equation}
But then, from (\ref{eq6.6}) and Corollary 11 in
Subsection~\ref{sub7.6}, we have, for all $t\in \wtI$:
\begin{alignat}{1}
\de (\ksM Q, P_k) \; &<  \; \de_{LR} (\ksM Q, P_k)\label{eq6.70}\\
&< \; \de_{LR} (\ksM Q^*, P^*_k)\notag\\
&< \; c |\wtI|.\notag
\end{alignat}
From (R7), we have, for all $t\in I$
\begin{equation}
\de (\ksM Q, P_k) \; > \; c^{-1} |P_k|^{1-\eta}\label{eq6.71}
\end{equation}
and the second assertion of the Corollary follows.
\end{proof}

\medskip\noindent
{\bf 6.7.12~}
\begin{Coro}\label{coro8}
Any $(P,Q,n)$ in $\cR (I)$ but not in $\cR (\wtI)$ satisfies, for
all $t\in I$:
\[
|P| \; \leqslant \; |\wtI|^{\fud},\qquad\qquad |Q|\; \leqslant \;
|\wtI|^{\fudt}.
\]
\end{Coro}

\begin{proof}
From (\ref{eq6.70}), (\ref{eq6.71}) and (\ref{eq3.27}), we actually
have, for all $t\in I$
\begin{equation}
|P|\; \leqslant \; |P_0|\; |\wtI|^{\fudt},\label{eq6.72}
\end{equation}
and we must have $|P_0| \ll 1$ as $Q_0 \subset Q_u$.
\end{proof}

\medskip\noindent
{\bf 6.7.13~}
\begin{Coro}\label{coro9}
Any candidate interval $I$ is $\ovbe$-regular, with $\ovbe = \be
(1+\tau)^{-1}$.
\end{Coro}

\begin{proof}
Let $(P,Q,n)$ be an $I$-bicritical element of $\cR (I)$. If
$(P,Q,n)$ belongs to $\cR (\wtI)$, then it is $\wtI$-bicritical and
we have, for all $t\in \wtI$:
\begin{equation}
\max (|P|, |Q|) \; < \; |\wtI|^\be = |I|^\ovbe.\label{eq6.73}
\end{equation}
Assume now that $(P,Q,n)$ does not belong to $\cR (\wtI)$. We apply
Theorem~\ref{theo1}. The element $(P_0, Q_0, n_0)\in \cR (\wtI)$ is
$\wtI$-bicritical: $Q_0$ is $\wtI$-critical by Theorem~\ref{theo1},
and $P_0$ is $\wtI$-critical because $P$ is $I$-critical. Applying
Corollaries~\ref{coro6} and \ref{coro7} gives, for all $t\in I$
\begin{alignat}{1}
|P| &\leqslant \; C |P_0|\; |P_k|\; |I|^{-\fudt}\label{eq6.74}\\
&\leqslant \; C |\wtI|^{\be + 1 - \fudt\, (1+\tau)},\notag
\end{alignat}
with $\be + 1 - \fud\, (1+\tau) > \be + \fut$, so $I$-bicritical
elements are much thinner in this case.
\end{proof}

\begin{rem}
The same phenomenon will be important again in Section~\ref{sec9}:
''fat" bicritical elements were created much earlier.
\end{rem}

\medskip\noindent
{\bf 6.7.14~} The last result in this section is a complement to
Proposition~\ref{propo5} in Subsection~\ref{sub6.4}.

Let $I$ be a parameter interval, and let $(\wtP, \wtQ, \wtn)$,
$(P,Q,n)$ be elements of $\cR (I)$ such that $P$ is a non-simple
child of $\wtP$. From Proposition~\ref{propo5}, we know that there
exists  $(P_1, Q_1, n_1)\in \cR (I)$ such that $\wtQ \forki P_1$
holds and
\begin{equation}
(P,Q,n) \; \in \; (\wtP, \wtQ, \wtn) \; \square \; (P_1, Q_1,
n_1).\label{eq6.75}
\end{equation}

\begin{Propo}\label{propo12}
Assume moreover that $\wtP$ is $I$-critical. Then, for all $t\in I$
we have
\[
\de (\wtQ, P_1)\; \geqslant \; |\wtP|^{\tfrac{1}{\ovbe}}.
\]
\end{Propo}

\begin{proof}
Consider first the case where $\wtQ$ is $I$-critical. Then $(\wtP,
\wtQ, \wtn)$ is $I$-bicritical. From Corollary~\ref{coro9}, we know
that $I$ is $\ovbe$-regular and therefore we have, for all $t\in I$:
\begin{equation}
|\wtP| \; < \; |I|^{\ovbe}.\label{eq6.76}
\end{equation}
On the other hand, as $\wtQ \forki P_1$ holds, we can find $(\wtP^*,
\wtQ^*, \wtn^*)$, $(P^*_1, Q^*_1, n^*_1)\in \cR (I)$ with $\wtQ
\subset \wtQ^*$, $P_1 \subset P^*_1$ such that, for all $t\in I$
\begin{equation}
\de_{LR} (\wtQ^*, P^*_1)\; \geqslant\; 2 |I|.\label{eq6.77}
\end{equation}
From (R7) (cf.~Subsection~\ref{eq5.4}) and
(\ref{eq5.11})--(\ref{eq5.14}), it now follows that
\begin{alignat}{1}
\de (\wtQ, P_1) \;
&\geqslant \;\de (\wtQ^*, P^*_1)\label{eq6.78}\\
&\geqslant \;\fudt\;\de_{LR} (\wtQ^*, P^*_1)\notag\\
&\geqslant \;|I|\notag\\
&\geqslant \;|\wtP|^{\tfrac{1}{\ovbe}}\notag.
\end{alignat}

From now on we assume that $\wtQ$ is $I$-transverse. Let $I^*
\supset I$ be the largest parameter interval such that $(\wtP, \wtQ,
\wtn) \in \cR (I^*)$ and $\wtQ$ is $I^*$-transverse. As the
transversality relation never holds for the starting interval
$[\vep_0, 2\vep_0] = I_0$, $I^*$ is not equal to $I_0$. Let $\wtI^*$
be the parent of $I^*$. We claim that for all $t\in I^*$, we have:
\begin{equation}
|\wtP| \; \leqslant\; |\wtI^*|^\be = |I^*|^{\ovbe}.\label{eq6.79}
\end{equation}

Indeed, from the definition of $I^*$, we have that either $(\wtP,
\wtQ, \wtn)\not\in \cR (\wtI^*)$, or $(\wtP,\wtQ, \wtn)\in \cR
(\wtI^*)$ and $\wtQ$ is $\wtI^*$-critical. In the second case, $(\wtP,
\wtQ, \wtn)$ is $\wtI^*$-bicritical and $\wtI^*$ is $\be$-regular,
which gives (\ref{eq6.79}). In the first case, we obtain from
(\ref{eq6.74}) in the proof of Corollary~\ref{coro9} (where only the
$I$-criticality of $P$ was used) that, for all $t\in I^*$
\begin{equation}
|\wtP| \; \leqslant\; C |\wtI^*|^{\be+\tfrac{1}{3}} <
|I^*|^{\ovbe}.\label{eq6.80}
\end{equation}
The claim is proved.

As $\wtQ$ is $I^*$-transverse, there exists an $\wtI$-decomposition
$(P_\al, Q_\al, n_\al)$ of $P_s$ in $\cR (I^*)$ such that, for any
$\al$, $\wtQ$ and $P_\al$ are either $I^*$-separated or
$I^*$-transverse. There exists $\al_0$ such that $P_1$ and
$P_{\al_0}$ intersect; $\wtQ$ and $P_{\al_0}$ must be
$I^*$-transverse. This implies, as for (\ref{eq6.77}) above, that we
have, for all $t\in I^*$:
\begin{equation}
\de (\wtQ, P_{\al_0}) \; \geqslant \; \frac{3}{2}\;
|I^*|.\label{eq6.81}
\end{equation}
If $P_1 \subset P_{\al_0}$, we have
\begin{equation}
\de (\wtQ, P_1) \; \geqslant \; \de (\wtQ, P_{\al_0}).\label{eq6.82}
\end{equation}
If $P_{\al_0} \subset P_1$, we have from
(\ref{eq5.11})--(\ref{eq5.14}) and (R7), for all $t\in I$:
\begin{alignat}{1}
\de (\wtQ, P_1)\;
&\geqslant \;\frac{2}{3}\; \de_R (\wtQ, P_1)\label{eq6.83}\\
&\geqslant \;\frac{2}{3}\; \de_R (\wtQ, P_{\al_0})\notag\\
&\geqslant \;\frac{2}{3}\; \de (\wtQ, P_{\al_0})\notag.
\end{alignat}
In all cases, combining this with (\ref{eq6.81}) and (\ref{eq6.79})
gives the required estimate.
\end{proof}

\newpage

\setcounter{section}{6}
\setcounter{equation}{0}

\section{Estimates for the Classes $\cR(I)$\label{sec7}}
\subsection{A Stretched Exponential Estimate for Widths\label{sub7.1}}
The next proposition is a substitute for the uniform exponential
estimates for widths that are characteristic of the uniformly
hyperbolic dynamics. We denote by $\ga$ the constant
\begin{equation}
\ga \; := \; \frac{\Log\, \frac{3}{2}}{\Log \,2}\; \in \;
(0,1).\label{7.1}
\end{equation}

\begin{Propo}\label{propo13}
Let $I$ be a parameter interval and let $(P,Q,n)$ be an element of
$\cR (I)$. For all $t\in I$, we have
\[
|P| \; \leqslant \; C\; \exp (-n^\ga)
\]
with the stronger estimate
\[
|P|\; \leqslant \; C\; \exp (-2 n^\ga)
\]
when the parent of $P$ is $I$-critical.
\end{Propo}

\begin{proof}
If $I$ is the starting interval $I_0 = [\vep_0, 2\vep_0]$, a
stronger exponential bound actually holds. The proof is by induction
on $n$, and the estimates are therefore valid when $n = O (\log\;
\vep_0^{-1})$ (in which case $(P,Q,n) \in \cR (\Io)$). Let $(\wtP,
\wtQ, \wtn)$ be the element of $\cR (I)$ such that $\wtP$ is the
parent of $P$. In the case where $P$ is a simple child of $\wtP$,
the bound for $P$ easily follows from the bound for $\wtP$. Let us
therefore assume that $P$ is a non-simple child of $\wtP$, in which
case there exists by Proposition~\ref{propo5} an element $(P_1, Q_1,
n_1)\in \cR (I)$ such that $\tilde Q \forki P_1$ holds, the parent of $P_1$
is $I$-critical and
\begin{equation}
(P,Q,n) \; \in \; (\wtP, \wtQ, \wtn) \; \square \; (P_1, Q_1,
n_1).\label{7.2}
\end{equation}
By the induction hypothesis, we have
\begin{alignat}{1}
|\wtP|\; &\leqslant \; C \, \exp (-\wtn^\ga),\label{eq7.3}\\
|P_1|\; &\leqslant \; C \, \exp (-2 n_1^\ga).\label{eq7.4}
\end{alignat}

From (R7) and (\ref{eq3.27}), we have, for all $t\in I$:
\begin{equation}
|P|\; \leqslant \; C \, |\wtP| |P_1| ^{\tfrac{1+\eta}{2}} \; \ll \;
|\wtP|\; |P_1|^{\fudt},\label{eq7.5}
\end{equation}
and this gives
\begin{equation}
|P| \; \leqslant \; \exp (-\wtn^\ga - n_1^\ga).\label{eq7.6}
\end{equation}
This proves the required estimate in the general case, because, when
$\wtn$ and $n_1$ are large, we have
\begin{equation}
\wtn^\ga + n_1^\ga \; \geqslant \; (\wtn + n_1 + N_0)^\ga\; = \;
n^\ga.\label{eq7.7}
\end{equation}
Assume now that $\wtP$ is $I$-critical. Instead of (\ref{eq7.3}), we
have
\begin{equation}
|\wtP| \; \leqslant \; C \, \exp (-2 \wtn^\ga),\label{eq7.8}
\end{equation}
and then from (\ref{eq7.5}) we obtain
\begin{equation}
|P| \; \leqslant \; \exp (-2 \wtn^\ga - n_1^\ga).\label{eq7.9}
\end{equation}
This will be useful when $\wtn \geqslant n_1$. When $\wtn \leqslant
n_1$, we prefer to rely on Proposition~\ref{propo12} which gives
\begin{equation}
\de (\wtQ, P_1) \; \geqslant \; |\wtP|^{1/\ovbe} \; \gg \;
|\wtP|.\label{7.10}
\end{equation}

When we use this in (\ref{eq3.27}) and combine with (\ref{eq7.4})
and (\ref{eq7.8}), we get
\begin{equation}
|P| \; \leqslant \; \exp (-\wtn^\ga - 2 n_1^\ga).\label{eq7.11}
\end{equation}
To get the required estimate from (\ref{eq7.9}), (\ref{eq7.11}), we
have only to observe that the function $u\mapsto u^\ga + 2
(1-u)^\ga$ is concave on $[0, \fud]$ and equal to 2 for $u=0$ and
$u=\fud$.
\end{proof}

\subsection{Uniform Cone Condition\label{sub7.2}}
In this subsection, we will check that all elements $(P,Q,n)\in \cR
(I)$ satisfy the cone condition (AL2) of Subsection~\ref{sub3.2} for
the parameters $\la$, $u_0$, $v_0$ of Subsection~\ref{sub5.3}: we
have $u_0 = u^{1/2}$, $v_0 = v^{1/2}$, and all $(P,Q,n)\in \cR
([\vep_0, 2\vep_0])$ satisfy (AL2) with parameters $\la$, $u$, $v$.

Let $(A,B)$ be the implicit representation of the affine-like
iterate $(P,Q,n)$; we have to prove that
\begin{alignat}{1}
\la |A_x| + u_0 |A_y|\; \leqslant \; 1,\tag{\bf AL2}\\
\la |B_y| + v_0 |B_x|\; \leqslant \; 1.\notag
\end{alignat}

Let $u_1 = u^{\frac{3}{4}}$, $v_1 = v^{\frac{3}{4}}$. We will prove
that, for all $t\in I$, we have
\begin{equation}
|A_y| \; < \; u_1^{-1},\qquad\qquad |B_x| \; < \;
v_1^{-1}.\label{eq7.12}
\end{equation}
This is sufficient to obtain (AL2): we already know that if
$(P,Q,n)\in \cR (\Io)$ them (AL2) is satisfied; on the other hand,
if $(P,Q,n)\not\in \cR (I_0)$, then, for all $t\in I$, we have from
Corollary~\ref{coro8}
\begin{equation}
|A_x| \; < \; C\,\vep_0^{\fudt},\qquad\qquad |B_y| \; < \;
C\,\vep_0^{\fudt};\label{eq7.13}
\end{equation}
with $\vep_0$ small enough, (\ref{eq7.13}) and (\ref{eq7.12}) give
(AL2).

Let us now proceed with the proof of (\ref{eq7.12}). When
$(P,Q,n)\in \cR (I_0)$, we have the stronger estimate:
\begin{equation}
|A_y| \; < \; u^{-1},\qquad\qquad |B_x| \; < \;
v^{-1}.\label{eq7.14}
\end{equation}
Let $(\wtP, \wtQ, \wtn)$ be the element of $\cR (I)$ such that
$\wtP$ is the parent of $P$. Denote by $(\wtA, \wtB)$ the implicit
representation of the corresponding affine-like iterate.

If $P$ is a simple child, we use formula (\ref{eq3.11}) of
Subsection~\ref{sub3.3} to obtain
\begin{equation}
|A_y - \wtA_y| \; \leqslant \; C\, |\wtP| \; |\wtQ|.\label{eq7.15}
\end{equation}
If $P$ is a non-simple child, it is obtained
(Proposition~\ref{propo5}) by the parabolic composition of $(\wtP,
\wtQ, \wtn)$ with some $(P_1, Q_1, n_1)$; we use formula
(\ref{eq3.30}) of Subsection~\ref{sub3.5} to obtain
\begin{equation}
|A_y - \wtA_y| \; \leqslant \; C\,|\wtP| \; |\wtQ| \, (\de (\wtQ,
P_1))^{-\fudt}.\label{eq7.16}
\end{equation}
From (R7), $\de (\wtQ, P_1)$ is much larger than $|\wtQ|$. In all
cases, we have
\begin{equation}
|A_y - \wtA_y| \; \leqslant \; C\,|\wtP| \; \leqslant \; C \exp
(-\wtn^\ga),\label{eq7.17}
\end{equation}
where we have used Proposition~\ref{propo13} in the last inequality.
We only need (\ref{eq7.17}) when $\wtn$ is large (because we already
have (\ref{eq7.14}) otherwise), and the series $\sum\,exp (-n^\ga)$
is convergent. Therefore (\ref{eq7.12}) is a consequence of
(\ref{eq7.14}) and (\ref{eq7.17}).

The proof of (AL2), i.e., the first part of condition (R2) in
Subsection~\ref{sub5.3}, is now complete.

\subsection{Bounded Distortion\label{sub7.3}}
We now check the second half of property (R2) in
Subsection~\ref{sub5.3}. We have to prove that, for all $(P,Q,n)\in
\cR(I)$, we have the following estimate on distortion:
\begin{equation}
D ( g^n_t / P)  \; \leqslant \; 2 D_0.\label{eq7.18}
\end{equation}
Here, the constant $D_0$ corresponds to the stronger estimate we
obtain from (MP6) when $(P,Q,n)\in \cR (I_0)$:
\begin{equation}
D ( g^n_t / P) \; \leqslant \; D_0.\label{eq7.19}
\end{equation}
For $m>0$, define
\begin{equation}
D (m) = \sup_{\substack{(P,Q,n)\in \cR(I)\\ n\leqslant m}} \;
\sup_{t\in I} \; D (g^n_t/P).\label{eq7.20}
\end{equation}

To obtain, by induction on $m$, a bound for the non-decreasing
sequence $D(m)$, we combine (\ref{eq7.19}), (which gives
$D(m)\leqslant D_0$ for $m= O \;(\log \,\vep_0^{-1})$),
Proposition~\ref{propo13} and the bounds (\ref{eq3.13}) in
Subsection~\ref{sub3.3} (for simple composition) and (\ref{eq3.29})
in Subsection~\ref{sub3.5} (for parabolic composition). We set
\begin{equation}
D_s (m) = \max_{\substack{n>0,\,n'>0\\n+n'\leqslant m}} \; D_s
(n,n')\label{eq7.21}
\end{equation}
with
\begin{equation}
D_s (n,n') = D (n) + c\; \exp (-n^\ga)(D(n) + D(n')).\label{eq7.22}
\end{equation}
We also set
\begin{equation}
D_p (m) = \max_{\substack{n\gg 0,\,n'\gg 0\\n+n'+N_0\leqslant m}} \;
D_p (n,n')\label{eq7.23}
\end{equation}
with
\begin{equation}
D_p (n,n') = D (n) + c\; \exp (-\eta\, n^\ga).\label{eq7.24}
\end{equation}
Then, we have, as long as $D(m)$ is not too large (cf.~the condition
for (\ref{eq3.29}) to hold)
\begin{equation}
D (m) \; \leqslant \; \max (D_s (m), D_p (m)).\label{eq7.25}
\end{equation}
The reason for formula (\ref{eq7.24}) to hold is that the term
$c\,|P_1| \de^{-1}$ of (\ref{eq3.29}) is smaller than $|P_1|^\eta$
by (R7); then one uses Proposition~\ref{propo13}.

It is now clear that (\ref{eq7.18}) follows from (\ref{eq7.19}),
(\ref{eq7.21})--(\ref{eq7.25}).

\subsection{Estimates for the Special Rectangles $P_s$ and $Q_u$\label{sub7.4}}
In the next subsection, we will check the estimates contained in
condition (R4) of Subsection~\ref{sub5.3} concerning the class $\cR
(I)$.

These estimates, which are related to parabolic composition, are
valid for an element $(P,Q,n)$ of $\cR (I)$ which satisfies
$Q\subset Q_u$ (or $P\subset P_s$).

In the present section, we will be concerned with the affine-like
iterates which are directly associated with the elements
$(P_s,Q_s,n_s)$ and $(P_u,Q_u,n_u)$.

We will make the computations for $(P_s,Q_s,n_s)$ the other case is
obviously symmetric. We will assume that the periodic point $p_s$ is
{\it fixed}: the general case is completely similar, but the
notations are more awkward.

In this subsection, we just write $(x,y)$ for the coordinates in the
rectangle $R_{a_s}$ containing $p_s$; we denote by $(A,B)$ the
implicit representation of the affine-like iterate
\begin{equation}
G_t\; : \; (R_{a_s}) \cap g_t^{-1} (R_{a_s}) \; \longrightarrow\;
g_t\; (R_{a_s}) \cap R_{a_s}.\label{eq7.26}
\end{equation}

For $n\geqslant 0$, we denote by $(\An, \Bn)$ the implicit
representation of the $n^{th}$ iterate of this restriction.

As the equation of $W^s_{\loc} (p_s)$ is $\{x=0\}$ (cf.~(MP3) in
Subsection~\ref{sub2.2}), we have
\begin{equation}
A (y, 0, t) \; \equiv \; 0,\label{eq7.27}
\end{equation}
from which we deduce
\begin{alignat}{1}
|A_y (y, x, t)| \; &\leqslant \; c |x|,\label{eq7.28}\\
|A_t (y, x, t)| \; &\leqslant \; c |x|,\notag\\
|A_{yy} (y, x, t)| \; &\leqslant \; c |x|,\notag\\
|A_{yt} (y, x, t)| \; &\leqslant \; c |x|.\notag
\end{alignat}
Denote by $\mu = \mu (t)$ the unstable eigenvalue of $D_{g_t}$ at
$p_s$. For all $t$, $x$, $y$, $n$, we have
\begin{equation}
c^{-1} \mu^{-n} \; \leqslant\;  |A_x^\pnn (x,y,t)| \; \leqslant \; c
\,\mu^{-n}.\label{eq7.29}
\end{equation}
Let $(x_i, y_i)_{0\leqslant i\leqslant n}$ be an orbit of $g_t$ in
$R_{a_s}$. For all $0\leqslant \ell \leqslant m \leqslant n$, we
have:
\begin{equation}
c^{-1} \mu^{m-\ell} |x_\ell| \; \leqslant\;  |x_m| \; \leqslant \; c
\,\mu^{m-\ell} |x_\ell|.\label{eq7.30}
\end{equation}

\begin{Propo}\label{propo14}
The following estimates hold:
\begin{eqnarray}
|A_y^\pnn (y_0, x_n, t)| \; &\leqslant \; c |x_0| \; &\leqslant \; c
\mu^{-n} |x_n|,\label{eq7.31}\\
|A_t^\pnn (y_0, x_n, t)| \; &\leqslant \; cn |x_0| \; &\leqslant \;
cn \mu^{-n} |x_n|,\label{eq7.32}\\
|A_{yy}^\pnn (y_0, x_n, t)| \; &\leqslant \; c |x_0| \; &\leqslant
\; c \mu^{-n} |x_n|.\label{eq7.33}
\end{eqnarray}
\end{Propo}

\medskip\noindent
{\it Proof of \ref{eq7.31}:} From formula (\ref{eq3.11}) in
Subsection~\ref{sub3.3}, we have:
\begin{equation}
A_y^\pnn (y_0,x_n,t) = A_y^{\pnu} (y_0, x_{\nuu},t) + A_y A_x^{\pnu}
B_y^{\pnu} \De^{-1},\label{eq7.34}
\end{equation}
with $B_y^({\nuu})\De$ exponentially small with $n$ and, using
(\ref{eq7.28})--(\ref{eq7.30}):
\begin{equation}
\lv A_y (y_{\nuu},x_n,t) A_x^{\pnu} (y_0, x_{\nuu},t) \rv \;
\leqslant \; c |x_0|.\label{eq7.35}
\end{equation}

The inequality (\ref{eq7.31}) is now clear.

\medskip\noindent
{\it Proof of \ref{eq7.32}:} We use here formulas (A6), (A10) of
Appendix~A which give
\begin{eqnarray}
&\lv A_t^\pnn (y_0,x_n,t) - A_t^{\pnu} (y_0, x_{\nuu},
t)\rv\label{eq7.36}\\
&\leqslant \; C\mu^{-n} \Bigl( \lv A_t(y_{\nuu}, x_n,t)\rv + \lv
B_t^{\pnu} (y_0, x_{\nuu},t)\rv\; \lv A_y (y_{\nuu}, x_,
t)\rv\Bigr),\notag
\end{eqnarray}
\begin{equation}
\lv B_t^\pnn (y_0,x_n,t) - B_t^{\pnu} (y_0,x_{n-1},t)\rv \;\leqslant \;
C \lv B_y^{\pnu}\rv_{C^0} \Bigl( 1 + \lv
A_t^{\pnu}\rv_{C^0}\Bigr).\label{eq7.37}
\end{equation}
As $|B_y^\pnn|_{C_0}$ is exponentially small, we deduce from
(\ref{eq7.37}) that
\begin{equation}
\lv B_t^\pnn (y_0, x_n,t) \rv \; \leqslant \; C,\label{eq7.38}
\end{equation}
and then, from (\ref{eq7.36}), (\ref{eq7.28}) that (\ref{eq7.32})
holds.

\medskip\noindent
{\it Proof of \ref{eq7.33}:} We use formulas (A6), (A18), (A20) of
Appendix~A to obtain
\begin{equation}
A_{yy}^\pnn = A_{yy}^\pnu + 2A_{xy}^\pnu X_y + A_{xx}^\pnu X^2_y +
A_{x}^\pnu X_{yy},\label{eq7.39}
\end{equation}
with
\begin{eqnarray}
X_y &= &A_y B_y^\pnu \dmu,\label{eq7.40}\\
\De &= &1 - A_y B_x^\pnu,\label{eq7.41}\\
X_{yy} &= &B_y^\pnu \dmu \left( A_{yy} B_y^\pnu \dmu + A_y
\pa_y \log\, |B_y^\pnu|+ \right.\label{eq7.42}\\
&+ &A_y X_y \pa_x \log \, |B_y^\pnu| - \left. A_y^\pnu \De_y
\dmu\right),\notag\\
-\De_y &= &A_{yy} B_y^\pnu B_x^\pnu \dmu + A_y B_{xy}^\pnu + A_y
B_{xx}^\pnu X_y.\label{eq7.43}
\end{eqnarray}
In these formulas, $A^\pnu$, $B^\pnu$ and their derivatives are
taken at $(y_0, x_\nuu, t)$, $A$, $B$ and their derivatives are
taken at $(y_\nuu, x_n, t)$. The terms $B_x^\pnu$, $B_{xx}^\pnu$,
$\pa_x \log |B_y^\pnu|$, $\pa_y \log |B_y^\pnu|$, $\dmu$ are bounded
by the uniform cone condition and the uniform distortion; the terms
$B_y^\pnu$, $B_{xy}^\pnu$ are exponentially small. Also, from
(\ref{eq7.28}) we have:
\begin{eqnarray}
\lv A_y (y_\nuu,x_n,t)\rv \; &\leqslant \; &C |x_n|,\label{eq7.44}\\
\lv A_{yy} (y_\nuu,x_n,t)\rv \; &\leqslant \; &C |x_n|.\notag
\end{eqnarray}

We conclude that we can write
\begin{equation}
A_{yy}^\pnn (y_0,x_n,t) =  A_{yy}^\pnn (y_0, x_\nuu,t) + \mu^{-n}
x_n r_n\label{eq7.45}
\end{equation}
with $r_n$ exponentially small; this leads to (\ref{eq7.34}).\hfill
$\square$

\begin{Coro}\label{coro10}
For the special rectangle $(P_s, Q_s, n_s)$, we have:
\begin{eqnarray*}
&|A_y^{\pns}|_{C^0} \; &\leqslant \; C\vep_0,\\
&|A_{yy}^{\pns}|_{C^0} \; &\leqslant \; C\vep_0,\\
&|A_t^{\pns}|_{C^0} \; &\leqslant \; C\vep_0 \log \, \vep_0^{-1}.
\end{eqnarray*}
\end{Coro}

\begin{proof}
We have only to observe that $\mu^{-n_s}$ is of order $\vep_0$.
\end{proof}

Obviously, the same estimates hold for the other special element
$(P_u, Q_u, n_u)$.

\subsection{Proof of the Property (R4) of Section~5.3\label{sub7.5}}
We have to show that, for an element $(P,Q,n)$ in $\cR(I)$, with
associated implicit representation $(A,B)$, we have
\begin{equation}
|A_y|\; \leqslant \; C\vep_0,\qquad |A_{yy}| \; \leqslant \; C\vep_0
\label{eq7.46}
\end{equation}
whenever $P\subset P_s$ and
\begin{equation}
|B_x|\; \leqslant \; C\vep_0,\qquad |B_{xx}| \; \leqslant \; C\vep_0
\label{eq7.47}
\end{equation}
whenever $Q\subset Q_u$. We will deal only with (\ref{eq7.46}), the
other case being symmetric.

We have already proved (\ref{eq7.46}) when $P=P_s$, cf.~the first
two inequalities of Corollary~\ref{coro10}. We will prove
(\ref{eq7.46}) by induction on the length $n$. Let us denote by
$(\wtP,\wtQ,\wtn)$ the element of $\cR(I)$ such that $\wtP$ is the
parent of $P$, by $(\wtA,\wtB)$ the implicit representation of this
affine-like iterate. We assume now that $P\subsetneqq P_s$, i.e.,
$\wtP \subset P_s$. The are two possibilities.

First, we consider the easier case when $P$ is a simple child of
$\wtP$. Let $(P^*, Q^*, 1)$, with associated implicit representation
$(A^*, B^*)$, be the element such that
\begin{equation}
(P,Q,n) = (\wtP,\wtQ,\wtn)\;*\; (P^*,Q^*,1).\label{eq7.48}
\end{equation}
From formula (\ref{eq3.11}) in Subsection~\ref{sub3.3}, we have:
\begin{equation}
A_y = \wtA_y + \wtA_x \wtB_y A^*_y \dmu,\label{eq7.49}
\end{equation}
with $|\wtA_x|\le C\vep_0$, $A^*_y$, $\dmu$ bounded and $\wtB_y$
satisfying the stretched exponential estimate of
Proposition~\ref{propo13}: this is fully in line with the first
inequality in (\ref{eq7.46}).

Next, we have, as in (\ref{eq7.39}) above:
\begin{equation}
A_{yy} = \wtA_{yy} + 2\wtA_{xy} X_y + \wtA_{xx} X^2_y + \wtA_x
X_{yy},\label{eq7.50}
\end{equation}
with
\begin{eqnarray}
X_y &= &A^*_y \wtB_y \dmu,\label{eq7.51}\\
\De &= &1 - A^*_y \wtB_x,\label{eq7.52}\\
X_{yy} &= &\dmu \wtB_y \Bigl( A^*_{yy} \wtB_y \dmu + A^*_y \pa_y
\log |\wtB_y| + A^*_y X_y \pa_x \log |\wtB_y| - \wtA_y \De_y \dmu\Bigr),\label{eq7.53}\\
-\De_y &= &\As_{yy} \wtB_y \wtB_x \dmu + A^*_y \wtB_{xy} + A^*_y
\wtB_{xx} X_y.\label{eq7.54}
\end{eqnarray}

Here, we have:
\begin{eqnarray}
|\wtA_{xy}| &= &|\wtA_x|\; \lv \pa_y \log |\wtA_x|\rv\; \leqslant \;
C\vep_0,\label{eq7.55}\\
|\wtA_{xx}| &= &|\wtA_x|\; \lv \pa_x \log |\wtA_x|\rv\; \leqslant \;
C\vep_0,\label{eq7.56}\\
|\wtA_{x}| &\leqslant &C\vep_0.\label{eq7.57}
\end{eqnarray}
Moreover, $\wtB_y$ satisfies the stretched exponential estimate of
Proposition~\ref{propo13}; the uniform cone condition and uniform
distortion imply that the same stretched exponential estimate holds
for $X_y$ and $X_{yy}$. Therefore, the estimate for $A_{yy} -
\wtA_{yy}$ in (\ref{eq7.50}) is fully in line with the second
inequality in (\ref{eq7.46}). This concludes the case where $P$ is a
simple child.

We now consider the more difficult case where $P$ is a non-simple
child of $\wtP$. By Proposition~\ref{propo5}, we can find an element
$(P^*, Q^*, n^*)\in \cR(I)$, with associated implicit representation
$(A^*, B^*)$, such that:
\begin{equation}
(P,Q,n)\; \in \; (\wtP, \wtQ, \wtn) \; \square \; (\Ps, \Qs, \ns).
\label{eq7.58}
\end{equation}

From the formulas in Appendix~A, we have:
\begin{eqnarray}
A_y &= &\wtA_y + \wtA_x (X_y + X_w W_y),\label{eq7.59}\\
X_y &= &X_{u,y} \wtB_y \dmuz 0,\label{eq7.60}\\
X_w &= &X_{u,w} \dmuz 0,\label{eq7.61}\\
\De_0 &= &1 - X_{u,y} \wtB_x,\label{eq7.62}\\
W_y &= &-C_y C_w^{-1},\label{eq7.63}\\
-C_y &= &\te_y \wtB_y \dmuz 0.\label{eq7.64}
\end{eqnarray}
Here $\te$, $X_u$, $W$ and $C$ are associated with the fold map $G$
as in Subsections~\ref{sub2.3}, \ref{sub3.5}.

As $|C_{ww}-2|$ is small (cf.~\ref{eq3.22})), the quantity $C_w$ is
related to $\de (\wtQ, P^*)$ by
\begin{equation}
C^{-1} \de (\wtQ, \Ps)^{\fudt}\; \leqslant \; |C_w| \; \leqslant \;
C \de (\wtQ, \Ps)^{\fudt}. \label{eq7.65}
\end{equation}
From (R7), $\de (\wtQ, \Ps)$ is always much bigger than $|\wtB_y|$;
therefore, from (\ref{eq7.63}), (\ref{eq7.64}) we have:
\begin{equation}
|W_y| \; \leqslant \; C |\wtB_y|^{\fudt}. \label{eq7.66}
\end{equation}
We then obtain from Proposition~\ref{propo13}, that the term $X_y +
X_w W_y$ in (\ref{eq7.59}) satisfies a stretched exponential
estimate; in view of (\ref{eq7.57}), this concludes the proof of the
first inequality in (\ref{eq7.46}).

For the second inequality, we write, following the formulas in
Appendix~A:
\begin{eqnarray}
A_{yy} &= &\wtA_{yy} + 2\wtA_{xy} (X_y + X_w W_y) + \wtA_{xx} (X_y +
X_w W_y)^2 +\label{eq7.67}\\
&+ &\wtA_x (X_{yy} + 2 X_{wy} W_y + X_{ww} W^2_y + X_w
W_{yy}),\notag\\
W_{yy} &= &-C^{-1}_w (C_{ww} W^2_y + 2C_{wy} W_y + C_{yy}),\label{eq7.68}\\
-C_{wy} &= &\te_{yy} \wtB_x \wtB_y X_w \dmuz 0 + \te_{xy} \wtB_y
A^*_y Y_{s,w} \dmuz 0 \dmuz 1 + \te_y \ovY_{wy},\label{eq7.69}\\
-C_{yy} &= &\te_{yy} \wtB^2_y \De^{-2}_0 + \te_{y} \ovY_{yy},\label{eq7.70}\\
\ovY_{wy} &= &\wtB_{xy} X_w + \wtB_{xx} X_w X_y + \wtB_x X_{wy},\label{eq7.71}\\
\ovY_{yy} &= &\wtB_{yy} + 2\wtB_{xy} X_y + \wtB_{xx} X^2_y + \wtB_x X_{wy},\label{eq7.72}\\
X_{ww} &= &\dmuz 0 (X_{u,ww} + 2\wtB_x X_w X_{u,wy} + \wtB^2_x
X_{u,yy} X^2_w + \wtB_{xx} X_{u,y} X^2_w),\label{eq7.73}\\
X_{wy} &= &\dmuz 0 \Bigl[ (X_{u,wy} + \wtB_x X_{u,yy}
X_w)(\wtB_y + \wtB_x X_y)\label{eq7.74}\\
&+ &X_{u,y} X_w (\wtB_{xy} + \wtB_{xx} X_y)\Bigr],\notag\\
X_{yy} &= &\dmuz{0} \Bigl[ X_{u,yy} (\wtB_y + \wtB_x X_y)^2 +
X_{u,y} (\wtB_{yy} + 2\wtB_{xy} X_y + \wtB_{xx}
X^2_y)\Bigr],\label{eq7.75}
\end{eqnarray}
where $\dmuz 0$, $\dmuz 1$ are bounded. In formula (\ref{eq7.67}),
we have that

--~~$|\wtA_{xy}|$, $|\wtA_{xx}|$, $|\wtA_x|$ are bounded by
$C\vep_0$, cf.~(\ref{eq7.55})--(\ref{eq7.57}) above;

--~~$X_w$, $X_{ww}$ are bounded, cf.~(\ref{eq7.61}) and
(\ref{eq7.73});

--~~$X_y$, $W_y$ satisfy a stretched exponential estimate, as we
have already seen earlier;

--~~$X_{wy}$, $X_{yy}$ also satisfy a stretched exponential
estimate, cf.~(\ref{eq7.74}) and (\ref{eq7.75}).

The remaining term in (\ref{eq7.67}) is $\wtA_x X_w W_{yy}$, and we
have to estimate $W_{yy}$ from (\ref{eq7.68}). From (\ref{eq7.75}),
we have, using bounded distortion:
\begin{equation}
|X_{yy}| \; \leqslant \; C |\wtB_y|,\label{eq7.76}
\end{equation}
from which we deduce, by (\ref{eq7.72}), that
\begin{equation}
|\ovY_{yy}| \; \leqslant \; C |\wtB_y|.\label{eq7.77}
\end{equation}
In the same way, we obtain
\begin{eqnarray}
|X_{wy}| \; &\leqslant \; &C |\wtB_y|,\label{eq7.78}\\
|\ovY_{wy}| \; &\leqslant \; &C |\wtB_y|.\label{eq7.79}
\end{eqnarray}
Plugging this into (\ref{eq7.69}), (\ref{eq7.70}) yields
\begin{eqnarray}
|C_{wy}| \; &\leqslant \; &C |\wtB_y|,\label{eq7.80}\\
|C_{yy}| \; &\leqslant \; &C |\wtB_y|.\label{eq7.81}
\end{eqnarray}
We can now conclude, using (\ref{eq7.65}), (\ref{eq7.66}), that
\begin{eqnarray}
|W_{yy}| \; &\leqslant \; &C |\wtB_y|\, \de (\wtQ, \Ps)^{-\fudt},\label{eq7.82}\\
&\leqslant \; &C |\wtB_y|^{\fudt}.\notag
\end{eqnarray}
All terms in (\ref{eq7.67}) are now under control, $A_{yy} -
\wtA_{yy}$ being bounded by $C\vep_0$ times a stretched
exponentially small term: this concludes the proof of the second
inequality in (\ref{eq7.46}), and so of condition (R4).

\subsection{Relative Speeds of Critical Rectangles\label{sub7.6}}
Let $I$ be a parameter interval, and let $(P_0, Q_0, n_0)$, $(P_1,
Q_1, n_1)$ be elements of $\cR (I)$ such that $Q_0 \subset Q_u$,
$P_1 \subset P_s$.

The displacements $\de (Q_0, P_1)$, $\de_L (Q_0, P_1)$, $\de_R (Q_0,
P_1)$, $\de_{LR} (Q_0, P_1)$ were introduced in formulas
(\ref{eq5.5}) and (\ref{eq5.8})--(\ref{eq5.10}) of
Subsection~\ref{sub5.4} (see also (\ref{eq3.23}) in
Subsection~\ref{sub3.5}) and are the values at the four corners of
the rectangle of definition of the function $\ovC (y_0, x_1)$
introduces in Subsection~\ref{sub3.5} as
\begin{equation}
\ovC (y_0, x_1) \; = \; \min_w \; C (w, y_0, x_1).\label{eq7.83}
\end{equation}
All these quantities also depend on the parameter $t$, and we want
in this section to estimate the variation with the parameter of the
displacements, which amounts to estimate the partial derivative
$C_t$.

Let $(A_0, B_0)$, $(A_1, B_1)$ be the implicit representations for
$(P_0, Q_0, n_0)$, $(P_1, Q_1, n_1)$ respectively. As will be seen
below, an estimate for $C_t$ depends very much on estimates for the
partial derivatives $A_{1,t}$, $B_{0,t}$. Good estimates for these
two quantities are not available for all $(P_0, Q_0, n_0)$, $(P_1,
Q_1, n_1)$. We will only consider elements satisfying conditions
$(*s)$, $(*u)$ below; fortunately, these conditions will always be
satisfied whenever we are interested in the variation of the
displacements.

Let $(P,Q,n)$ be an element of $\cR (I)$, such that $P\subset P_s$.
Let $\Is$ be the largest parameter interval such that $(P,Q,n)$ is
(the restriction of) an element of $\cR (\Is)$. If $\Is$ is not the
starting interval $I_0 = [\vep_0, 2\vep_0]$, let $\wtIs$ be the
parent interval and let $(\wtP, \wtQ, \wtn)$ be the longest element
of $\cR (\wtIs)$ such that $P\subset \wtP$; $(\wtP, \wtQ, \wtn)$ is
the element which is denoted by $(P_0, Q_0, n_0)$ in the structure
theorem of Subsection~\ref{sub6.7}, as $[0,\wtn]$ is the maximal
initial $\wtIs$-interval.

We say that $(P,Q,n)$ satisfies condition $(*s)$ if either $\Is =
I_0$ or $\Is \ne I_0$ and $\wtP$ is $\wtIs$-critical.

We define in a symmetric way a condition $(*u)$ (when $Q\subset
Q_u$).

\begin{Propo}\label{propo15}
Let $(P,Q,n)$ be an element of $\cR (I)$ with $P\subset P_s$. Let
$(A,B)$ be the implicit representation of $(P,Q,n)$. If $(P,Q,n)$
satisfies condition $(*s)$, we have
\[
|A_t| \; \leqslant \; \vep_0^{\fudt}.
\]
If $(P,Q,n)\in \cR (I_0)$, we have the stronger estimate:
\[
|A_t| \; \leqslant \; C \vep_0\; \log \; \vep_0^{-1}.
\]
\end{Propo}

\begin{proof}
We first show for $(P,Q,n)\in \cR (I_0)$ and $P\subset P_s$, that:
\begin{equation}
|A_t| \; \leqslant \; C \vep_0\; \log \; \vep_0^{-1}, |B_t| \;
\leqslant \; C.\label{eq7.84}
\end{equation}
When $P = P_s$, the first inequality is part of
Corollary~\ref{coro10} and the second is (\ref{eq7.38}). If
$P\subsetneqq P_s$, we set
\begin{equation}
(P,Q,n) \; = \; (\wtP, \wtQ, \wtn)\; * \; (\Ps, \Qs,
1).\label{eq7.85}
\end{equation}
Writing $(\wtA, \wtB)$ and $(\As, \Bs)$ for the implicit
representations of $(\wtP, \wtQ, \wtn)$ and $(\Ps, \Qs, 1)$,
respectively, we have, from the formulas in Appendix~A:
\begin{eqnarray}
A_t &= &\wtA_t + \wtA_x (A^*_t + A^*_y \wtB_t)\, \dmu,\label{eq7.86}\\
B_t &= &B^*_t + B^*_y (\wtB_t + \wtB_x A^*_t)\, \dmu,\label{eq7.87}
\end{eqnarray}
with, as usual, $\De = 1 - \wtB_x A^*_y$.

This gives, assuming (\ref{eq7.84}) for $(\wtA, \wtB)$,
\begin{eqnarray}
&|A_t - \wtA_t| &\leqslant \; C |P|,\label{eq7.88}\\
&|B_t - \wtB_t\, B_y \, \wtB_y^{-1}| &\leqslant \; C,\label{eq7.89}
\end{eqnarray}
where we have used $B_y = \wtB_y  B^*_y \dmu$. Clearly,
(\ref{eq7.88}) and (\ref{eq7.89}) imply, by iteration,
(\ref{eq7.84}).

We now turn to the case where $\Is \ne I_0$. In this case, we
introduce the integer $k\geqslant 1$ and the elements $(P_i, Q_i,
n_i)$, $0\leqslant i \leqslant k$, of $\cR (\wtIs)$ given by the
structure theorem of Subsection~\ref{sub6.7}. We also denote by
$(P^\pii, Q^\pii, n^\pii)$ the element of $\cR (\wtIs)$ such that $P
\subset P^\pii$ and
\begin{equation}
(P^\pii, Q^\pii, n^\pii) \; \in \; (P_0, Q_0, n_0)\; \square \;
\cdots \; \square\; (P_i, Q_i, n_i).\label{eq7.90}
\end{equation}
Our proof will be by induction, on the level of the parameter
interval and on the integer $k$ (for a fixed parameter interval).

We first observe that if $P$ is $I$-critical, then $P_0$ is
$\wtIs$-critical and $P$ satisfies $(*s)$. Therefore, if $P$
satisfies $(*s)$, all $P^\pii$ also satisfy $(*s)$; moreover,
$Q^\pii$ satisfy $(*u)$ and $\ism P$ satisfy $(*s)$ by
Theorem~\ref{theo1}.

Fix $0\leqslant i < k$ and let $(\piim A, \piim B)$, $(A^\pii,
B^\pii)$ and $(\As, \Bs)$ be the implicit representations of $(\piim
P, \piim Q, \piim n)$, $(P^\pii, Q^\pii, n^\pii)$ and $(\ism P, \ism
Q, \ism n)$ respectively. From the formulas in Appendix~A, we have
\begin{eqnarray}
\piim A_t &= &A^\pii_t + A^\pii_x (X_t + X_w W_t),\label{eq7.91}\\
X_t &= &(X_{u,t} + X_{u,y} B^\pii_t) \dmuz 0,\label{eq7.92}\\
X_w &= &X_{u,w} \dmuz 0,\label{eq7.93}\\
W_t &= &-C_t\; C^{-1}_w.\label{eq7.94}
\end{eqnarray}
In these formulas, $X_{u,t}$, $X_{u,y}$, $X_{u,w}$, $\dmuz 0$ are
uniformly bounded and $B^\pii_t$ is bounded by $\vep_0^{1/2}$ from
the induction hypothesis. The term $C_w^{-1}$ is estimated by
(\ref{eq7.65}) above and $C_t$ is uniformly bounded by the induction
hypothesis and Corollary~\ref{coro11} below. We conclude that
\begin{eqnarray}
|\piim A_t - A^\pii_t| \; &\leqslant \; &C |P^\pii| \de (Q^\pii, \ism P)^{-\fudt}\label{eq7.95}\\
&\leqslant \; &C |P^\pii|\; |\Is|^{-\fudt},\notag
\end{eqnarray}
as $Q^\pii$ and $\ism P$ are $\Is$-transverse. We have for
$0\leqslant j \leqslant i$
\begin{equation}
|P_j| \; \leqslant \; |\wtIs|^\be\label{eq7.96}
\end{equation}
because $(P_j, Q_j, n_j)$ is $\wtIs$-bicritical (we use here, for
$j=0$, the hypothesis that $P_0$ is $\wtIs$-critical).

Then, by Corollary~\ref{coro6}, we obtain
\begin{equation}
|P^\pii| \; \leqslant \; C^i \, |\wtIs|^{(i+1)\be - \tfrac{i}{2}\,
(1+\tau)}\label{eq7.97}
\end{equation}
which implies
\begin{equation}
|\piim A_t - A^\pii_t| \; \leqslant \; |\wtIs|^{\tfrac{\be}{2}\,
(1+i)}.\label{eq7.98}
\end{equation}
Summing (\ref{eq7.98}) over $i$ and the levels of parameter
intervals leads to the estimate of the proposition.
\end{proof}

\begin{Coro}\label{coro11}
Let $(P_0, Q_0, n_0)$, $(P_1, Q_1, n_1)$ be elements of $\cR (I)$
with $Q_0 \subset Q_u$, $P_1 \subset P_s$. Assume that $Q_0$
satisfies $(*u)$ and that $P_1$ satisfies $(*s)$. Then, the function
$C$ introduced in Subsection~\ref{sub3.5} satisfy
\[
|C_t + 1| \; \leqslant \; C\, \vep_0^{\fud}.
\]
\end{Coro}

\begin{proof}
From formula (A35) in Appendix~A, using the notations there, we have
\begin{eqnarray}
-C_t &= &\te_x \ovX_t + \te_y \ovY_t + \te_t,\label{eq7.99}\\
\ovX_t &= &(A_{1,t} + A_{1,y}\, Y_{s,t})\, \dmuz 1,\label{eq7.100}\\
\ovY_t &= &(B_{0,t} + B_{0,x}\, X_{u,t})\, \dmuz 0,\label{eq7.101}
\end{eqnarray}
with $\dmuz 0$, $\dmuz 1$ uniformly bounded. The value of $\te_t$ is
taken at $(\ovX, \ovY, t)$, with
\begin{equation}
|\ovX| \; \leqslant \; C \, \vep_0,\qquad |\ovY|\; \leqslant \; C\,
\vep_0.\label{eq7.102}
\end{equation}
On the other hand, we have, in Subsection~\ref{sub4.2}, normalized
the parameter in order to have
\begin{equation}
\te_t (0, 0, t) \; \equiv \; 1.\label{eq7.103}
\end{equation}
We, therefore, have
\begin{equation}
|\te_t (\ovX, \ovY, t) - 1| \; \leqslant \; C\vep_0.\label{eq7.104}
\end{equation}
In (\ref{eq7.100}) and (\ref{eq7.101}), we have $|A_{1,y}| <
C\vep_0$, $|B_{0,x}| < C\vep_0$, by (R4) and $|A_{1,t}| <
\vep_0^{\fud}$, $|B_{0,t}| < \vep_0^{\fud}$ by
Proposition~\ref{propo15}. The Corollary follows, as $\te_x$,
$\te_y$, $Y_{s,t}$, $X_{u,t}$ are uniformly bounded.
\end{proof}

\subsection{Variation of Width of Critical Rectangles\label{sub7.7}}
Our main purpose now is to prove property (R7) of
Subsection~\ref{sub5.4}:

{\bf (R7)}~~If $(P_0, Q_0, n_0)$, $(P_1, Q_1, n_1)\in \cR (I)$
satisfy $Q_0\subset Q_u$, $P_1\subset P_s$ and $Q_0 \forki P_1$
holds, then, for all $t\in I$, we have
\[
\de (Q_0, P_1) \; \geqslant \; C^{-1} \Bigl( |P_1|^{1-\eta} +
|Q_0|^{1-\eta}\Bigr).
\]
A priori, the transversality condition gives some control through
(T2), (T3) in Subsection~\ref{sub5.4} only for {\it some} values of
the parameter. However, from (T1) in Subsection~\ref{sub5.4} and
Corollary~\ref{coro11} above, we know that the order of magnitude of
$\de (Q_0, P_1)$ is the same through out $I$. Therefore, to obtain
(R7), we do need to control how the widths $|P_1|$ and $|Q_0|$ vary
through $I$. Good estimates will be obtained under the same
conditions $(*s)$ or $(*u)$ used to obtain
Proposition~\ref{propo15}. Again, the estimates are even better for
an element $(P,Q,n)$ in $\cR (I_0)$, involving only simple
composition.

\begin{Propo}\label{propo16}
Let $(P,Q,n)$ be an element of $\cR (I_0)$, and let $(A,B)$ be the
implicit representation of $(P,Q,n)$. We have
\begin{eqnarray*}
\lv\pa_t \,\log |A_x|\rv &\leqslant \; Cn,\qquad &|A_{yt}| \;
\leqslant\; C\\
\lv\pa_t \,\log |B_y|\rv &\leqslant \; Cn,\qquad &|B_{xt}| \;
\leqslant\; C.
\end{eqnarray*}
\end{Propo}

\begin{proof}
We first observe that we have, by the same proof as for the second
inequality in (\ref{eq7.84}):
\begin{equation}
|A_t| \; \leqslant \; C,\qquad |B_t|\; \leqslant \;
C.\label{eq7.105}
\end{equation}
We write
\begin{equation}
(P,Q,n) \; = \; (\wtP, \wtQ, \wtn)\; * \; (\Ps, \Qs,
1),\label{eq7.106}
\end{equation}
and denote by $(\wtA, \wtB)$, $(\As, \Bs)$ the implicit
representations of $(\wtP, \wtQ, \wtn)$, $(\Ps,\Qs,1)$,
respectively. By the formulas of Appendix~A, we have
\begin{eqnarray}
\pa_t \log |A_x| &= &\pa_t \log |\wtA_x| + \pa_t \log |A_x^*| + X_t
\pa_x \log |\wtA_x|\label{eq7.107}\\
&+ &Y_t \pa_y \log |A^*_x| - \De_t \dmu,\notag\\
X_t &= &(A^*_t + A^*_y \wtB_t) \dmu,\label{eq7.108}\\
Y_t &= &(\wtB_t + A^*_t \wtB_x) \dmu,\label{eq7.109}\\
\De &= &1 - A^*_y \wtB_x,\label{eq7.110}\\
-\De_t &= &\wtB_{xt} A^*_y + \wtB_{xx} A^*_y X_t + \wtB_x A^*_{yt} +
\wtB_x A^*_{yy} Y_t.\label{eq7.111}
\end{eqnarray}
As $\wtB_t$, $\wtB_x$, $A^*_t$, $A^*_y$, $\dmu$ are uniformly
bounded, the same is true for $X_t$, $Y_t$. Using bounded distortion
and the cone condition then leads to
\begin{equation}
|\pa_t \log |A_x| - \pa_t \log |\wtA_x|| \; \leqslant \; C (1 +
|\wtB_{xt}|).\label{eq7.112}
\end{equation}
Still from Appendix~A, we have
\begin{eqnarray}
A_{yt} - \wtA_{yt} &= &X_t \wtA_{xy} + X_y \wtA_{xt} + X_t X_y
\wtA_{xx} + \wtA_x X_{yt},\label{eq7.113}\\
X_y &= &A^*_y \wtB_y \dmu,\label{eq7.114}\\
X_{yt} &= &\wtB_y \dmu (A^*_{yy} Y_t + A^*_y \pa_t \log |\wtB_y| +
A^*_{yt}\label{eq7.115}\\
&+ &A^*_y X_t \pa_x \log |\wtB_y| - A^*_y \De_t \dmu).\notag
\end{eqnarray}
As $X_t$, $X_y$, $Y_t$ are bounded and also using bounded
distortion, we have
\begin{eqnarray}
|X_t \wtA_{xy}|\; &\leqslant\; &C |\wtA_x|,\label{eq7.116}\\
|X_y \wtA_{xt}|\; &\leqslant\; &C |\wtA_x|\; |\pa_t\log |\wtA_x||,\label{eq7.117}\\
|X_t X_y \wtA_{xx}|\; &\leqslant\; &C |\wtA_x|,\label{eq7.118}\\
|X_{yt}|\; &\leqslant\; &C |\wtB_y|(1+|\pa_t\log \wtB_y| +
|\wtB_{xt}|).\label{eq7.119}
\end{eqnarray}
We, therefore, obtain
\begin{equation}
|A_{yt} - \wtA_{yt}|\; \leqslant\; C |\wtA_x| (1+ |\pa_t\log |A_x||)
+ C |\wtA_x|\; |\wtB_y| (|\pa_t \log |\wtB_y|| +
|\wtB_{xt}|).\label{eq7.120}
\end{equation}
We have symmetric estimates for $B_{xt}$ and $\pa_t \log |B_y|$,
writing now
\begin{equation}
(P,Q,n) \; = \; (\whP^*, \whQ^*, 1)\; * \; (\whP, \whQ,
n-1).\label{eq7.121}
\end{equation}
As $|\wtA_x|$, $|\wtB_y|$ are exponentially small, the estimates
(\ref{eq7.112}), (\ref{eq7.120}) and the other two estimates for
$B_{xt}$, $\pa_t \log |B_y|$ lead by summation to the bounds of the
proposition.
\end{proof}

For the last estimate of this section, the setting is the same as
that of Proposition~\ref{propo15}.

\begin{Propo}\label{propo17}
Let $(P,Q,n)$ be an element of $\cR (I)$ with $P\subset P_s$. Let
$(A,B)$ be the implicit representation of $(P,Q,n)$. If $(P,Q,n)$
satisfies condition $(*s)$, we have
\begin{equation}
|\pa_t \log |A_x|| \; \leqslant \; C_0 \; \frac{\log |P|}{|I|\log
|I|},\label{eq7.122}
\end{equation}
\begin{equation}
|A_{yt}| \; \leqslant \; C_0.\label{eq7.123}
\end{equation}
\end{Propo}

\begin{proof}
The method is the same as in Proposition~\ref{propo15}, but, as now
we have to deal with second instead of first order partial
derivatives, calculations are more
complicated.

When $(P,Q,n)\in \cR(I_0)$, $n$ and $|\log |P||$ are of the same
order; as $|I|\,|\log |I||$ is always smaller than $\vep_0\log
\vep_0^{-1}$, the estimates in Proposition~\ref{propo16} are
stronger in the present case than the ones in (\ref{eq7.122}),
(\ref{eq7.123}).

We will now assume that $(P,Q,n)\not\in \cR (I_0)$; the proof of
(\ref{eq7.122}), (\ref{eq7.123}) is by induction on the level of the
parameter interval. Let $\Is\ne I_0$ be the largest parameter
interval such that $(P,Q,n)\in \cR (\Is)$; as
$|I|\;|\log\,|I||\leqslant |\Is|\,|\log\,|\Is||$, (\ref{eq7.122})
for $\cR (I)$ follows from (\ref{eq7.122}) for $\cR (\Is)$. We can
therefore assume that $\Is = I$ and denote by $\wtI$ the parent
interval.

We apply the structure theorem of Subsection~\ref{sub6.7} and use
the same notations as in the proof of Proposition~\ref{propo15}: we
have an integer $k$, elements $(P_i,Q_i,n_i)$ in $\cR(\wtI)$ for
$0\leqslant i\leqslant k$ and partial compositions of $(P^\pii,
Q^\pii, n^\pii)$ in $\cR (I)$ for $0\leqslant i\leqslant k$. We
denote by $(A,B)$, $(\wtA, \wtB)$, $(\As,\Bs)$ the implicit
representations of $(\piim P, \piim Q, \piim n)$, $(P^\pii, Q^\pii,
n^\pii)$, $(\ism P, \ism Q, \ism n)$, respectively, for some fixed
integer $0\leqslant i < k$. We know that both $(P^\pii, Q^\pii,
n^\pii)$, $(\ism P, \ism Q, \ism n)$ satisfy $(*s)$, and $(P^\pii,
Q^\pii, n^\pii)$ also satisfies $(*u)$. We have from formulas (A47),
(A49) in Appendix~A:
\begin{eqnarray}
\pa_t \log |A_x| - \pa_t \log |\wtA_x| &= &\pa_x \log |\wtA_x|\,
(X_t + X_w W_t) + \pa_t \log |X_w| +\label{eq7.124}\\
&+ &W_t \pa_w \log |X_w| + \pa_t \log |W_x|,\notag
\end{eqnarray}
\begin{eqnarray}
A_{yt} - \wtA_{yt} &= &\wtA_{xy}\, (X_t + X_w W_t) + \wtA_{xt}\,
(X_y + X_w W_y) +\label{eq7.125}\\
&+ &\wtA_{xx} \,(X_t + X_w W_t)(X_y + X_w W_y) +\notag\\
&+ &\wtA_x (X_{yt} + X_{wy} W_t + X_{wt} W_y + X_{ww} W_y W_t + X_w
W_{yt}).\notag
\end{eqnarray}
We need to estimate all terms in the right-hand sides of these
formulas.

\medskip\noindent
{\bf Terms involving $\wtA$.} By bounded distortion, $|\pa_x\log
|\wtA_x||$ is bounded, $\wtA_x$, $\wtA_{xy}$, $\wtA_{xx}$ are
bounded by $C |\wtA_x|$ and $\wtA_{xt}$ is part of the induction
\begin{equation}
|\wtA_{xt}| \; = \; |\wtA_x|\; |\pa_t \log |\wtA_x||.\label{eq7.126}
\end{equation}

\medskip\noindent
{\bf Terms involving the first order partial derivatives of $X$.} We
have dealt earlier with $X_t$, $X_w$, $X_y$ (cf.~(\ref{eq7.92}),
(\ref{eq7.93}), (\ref{eq7.60}) and these formulas easily give
\begin{eqnarray}
|X_t| \; &\leqslant \; &C,\label{eq7.127}\\
C^{-1}\; &\leqslant \; &|X_w|\; \leqslant \; C,\label{eq7.128}\\
|X_y| \; &\leqslant \; &C |\wtB_y|.\label{eq7.129}
\end{eqnarray}

\medskip\noindent
{\bf Terms involving the first order partial derivatives of $W$.} We
have dealt earlier with  $W_y = -C_y C^{-1}_w$ and $W_t = -C_t
C_w^{-1}$; the term $C_w^{-1}$ is estimated by (\ref{eq7.65}), the
term $C_y$ from formula (\ref{eq7.64}) and the term $C_t$ from
Corollary~\ref{propo11}. One obtains
\begin{equation}
C^{-1} \de (Q^\pii, \ism P)^{-\fudt} |\wtB_y| \leqslant |W_y|
\leqslant C\de (Q^\pii, \ism P)^{-\fudt} |\wtB_y|,\label{eq7.130}
\end{equation}
\begin{equation}
C^{-1} \de (Q^\pii, \ism P)^{-\fudt} \leqslant |W_t| \leqslant C\de
(Q^\pii, \ism P)^{-\fudt}.\label{eq7.131}
\end{equation}

\medskip\noindent
{\bf Terms involving the second order partial derivatives of $X$.}
The formulas (A63) of Appendix~A express the second order partial
derivatives of $X$ in terms of partial derivatives of first and
second order of $X_u$ (which are bounded), partial derivatives of
second order of $\wtB$ (which are controlled by the distortion or by
the induction hypothesis), partial derivatives of first order of $X$
itself (see above), the bounded quantity $\dmuz 0$, and partial
derivatives of first order of $\ovY$ (defined in (A31)). These
partial derivatives given by (A33) are easy to estimate: we have
$\ovY_w = \wtB_x X_w$, $\ovY_y = \wtB_y \dmuz 0$, $\ovY_t = (\wtB_t
+ \wtB_x X_{u,t}) \dmuz 0$, hence
\begin{eqnarray}
|\ovY_w| \; &\leqslant \; &C,\label{eq7.132}\\
|\ovY_y| \; &\leqslant \; &C |\wtB_y|,\label{eq7.133}\\
|\ovY_t| \; &\leqslant \; &C,\label{eq7.134}
\end{eqnarray}
where we have used Proposition~\ref{propo15} to bound $\wtB_t$.
Plugging these estimates in the formulas (A63) gives:
\begin{eqnarray}
|X_{ww}| \; &\leqslant \; &C,\label{eq7.135}\\
|X_{wy}| \; &\leqslant \; &C |\wtB_y|,\label{eq7.136}\\
|X_{yy}| \; &\leqslant \; &C |\wtB_y|,\label{eq7.137}\\
|X_{wt}| \; &\leqslant \; &C (1 + |\wtB_{xt}|),\label{eq7.138}\\
|X_{yt}| \; &\leqslant \; &C |\wtB_y|\, (1 + |\wtB_{xt}| + \pa_t
\log |\wtB_y|).\label{eq7.139}
\end{eqnarray}

\medskip\noindent
{\bf Terms involving the second order partial derivatives of $W$.}
The formulas (A55) of Appendix~A express the second order partial
derivatives of $W$ in terms of $C_w^{-1}$ (controlled by
(\ref{eq7.65})), the partial derivatives of first order of $W$ (see
above) and the partial derivatives of second order of $C$. These
partial derivatives of second order of $C$ are in turn expressed in
formulas (A56)--(A60) in terms of partial derivatives of $\te$
(which are bounded) and partial derivatives of first and second
order of $\ovX$ and $\ovY$. The partial derivatives of first order
of $\ovY$ have been estimated above, those of $\ovX$ satisfy in the
same way the inequalities
\begin{eqnarray}
|\ovX_{w}| \; &\leqslant \; &C,\label{eq7.140}\\
|\ovX_{x}| \; &\leqslant \; &C |\wtA_x^*|,\label{eq7.141}\\
|\ovX_{t}| \; &\leqslant \; &C.\label{eq7.142}
\end{eqnarray}
The partial derivatives of second order of $\ovX$ and $\ovY$ are
expressed in formulas (A61), (A62). The formulas (A61) contain
partial derivatives of first and second order of $Y$, which are
firstly estimated in the same way, through formulas (A29), (A64), as
those of $X$:
\begin{eqnarray}
|Y_t| \; &\leqslant \; &C,\label{eq7.143}\\
C^{-1}\; &\leqslant \; &|Y_w|\; \leqslant \; C,\label{eq7.144}\\
|Y_x| \; &\leqslant \; &C |A^*_x|,\label{eq7.145}\\
|Y_{ww}| \; &\leqslant \; &C,\label{eq7.146}\\
|Y_{wx}| \; &\leqslant \; &C |A^*_x|,\label{eq7.147}\\
|Y_{xx}| \; &\leqslant \; &C |A^*_x|,\label{eq7.148}\\
|Y_{wt}| \; &\leqslant \; &C (1 + |A^*_{yt}|),\label{eq7.149}\\
|Y_{xt}| \; &\leqslant \; &C |A^*_x|\, (1 + |A^*_{yt}| + |\pa_t \log
|A^*_x||).\label{eq7.150}
\end{eqnarray}
We can then estimate the second order partial derivatives of $\ovX$
and $\ovY$:
\begin{eqnarray}
|\ovX_{ww}| \; &\leqslant \; &C,\label{eq7.151}\\
|\ovX_{wx}| \; &\leqslant \; &C |A^*_x|,\label{eq7.152}\\
|\ovX_{wt}| \; &\leqslant \; &C (1 + |A^*_{yt}|,\label{eq7.153}\\
|\ovX_{xx}| \; &\leqslant \; &C |A^*_x|,\label{eq7.154}\\
|\ovX_{xt}| \; &\leqslant \; &C |A^*_x|\, (1 + |A^*_{yt}| + |\pa_t
\log |A^*_x||),\label{eq7.155}
\end{eqnarray}
\begin{eqnarray}
|\ovY_{ww}| \; &\leqslant \; &C,\label{eq7.156}\\
|\ovY_{wy}| \; &\leqslant \; &C |\wtB_y|,\label{eq7.157}\\
|\ovY_{wt}| \; &\leqslant \; &C (1 + |\wtB^*_{xt}|),\label{eq7.158}\\
|\ovY_{yy}| \; &\leqslant \; &C |\wtB_y|,\label{eq7.159}\\
|\ovY_{yt}| \; &\leqslant \; &C |\wtB_y|\, (1 + |\wtB_{xt}| + |\pa_t
\log |\wtB_y||).\label{eq7.160}
\end{eqnarray}

The next step is to estimate the partial derivatives of second order
of $C$ (besides $C_{ww}$ which is already known to the close to 2)
\begin{eqnarray}
|C_{wx}| \; &\leqslant \; &C |A^*_x|,\label{eq7.161}\\
|C_{wy}| \; &\leqslant \; &C |\wtB_y|,\label{eq7.162}\\
|C_{wt}| \; &\leqslant \; &C (1 + |\wtB^*_{xt}| + |A^*_{yt}|),\label{eq7.163}\\
|C_{xx}| \; &\leqslant \; &C |A^*_x|,\label{eq7.164}\\
|C_{xy}| \; &\leqslant \; &C |A^*_x|\,|\wtB_y|,\label{eq7.165}\\
|C_{yy}| \; &\leqslant \; &C |\wtB_y|,\label{eq7.166}\\
|C_{xt}| \; &\leqslant \; &C |A^*_x|\,(1 + |A^*_{yt}| + |\pa_t \log
|A^*_x||),\label{eq7.167}\\
|C_{yt}| \; &\leqslant \; &C |\wtB_y|\,(1 + |\wtB_{xt}| + |\pa_t
\log |\wtB_y||).\label{eq7.168}
\end{eqnarray}
Finally, we are able to estimate the partial derivatives of second
order of $W$:
\begin{equation}
|W_{xt}| \; \leqslant \; C \de^{-\fudt} |A^*_x|\,(\de^{-1} +
\de^{-\fudt} (|\wtB_{xt}| + |A^*_{yt}|) + |\pa_t \log
|A^*_x||),\label{eq7.169}
\end{equation}
\begin{equation}
|W_{yt}| \; \leqslant \; C \de^{-\fudt} |\wtB_y|\,(\de^{-1} +
\de^{-\fudt} (|\wtB_{xt}| + |A^*_{yt}|) + |\pa_t \log
|\wtB_y||).\label{eq7.170}
\end{equation}
We have written $\de$ for $\de (Q^\pii, \ism P)$.

We are now ready to come back to formulas (\ref{eq7.124}),
(\ref{eq7.125}) above.  We get:
\begin{equation}
|\pa_t \log |A_x|- \pa_t \log |\wtA_x|| \leqslant C (\de^{-1} +
\de^{-\fudt} (|\wtB_{xt}| + |A^*_{yt}|) + |\pa_t \log
|A^*_x||),\label{eq7.171}
\end{equation}
\begin{equation}
|A_{yt}-\wtA_{yt}| \leqslant C \de^{-\fudt} |\wtA_x| (1 + |\wtB_y|
K),\label{eq7.172}
\end{equation}
with
\begin{equation}
K = \de^{-1} + \de^{-\fudt} (|\wtB_{xt}| + |A^*_{yt}|) + |\pa_t \log
|\wtA_x| + |\pa_t \log |\wtB_y||.\label{eq7.173}
\end{equation}
By the induction hypothesis, we have
\begin{eqnarray}
|\wtB_{xt}| \; &\leqslant \; &C_0,\label{eq7.174}\\
|A^*_{yt}| \; &\leqslant \; &C_0,\label{eq7.175}\\
|\pa_t \log\,|A^*_x|| \; &\leqslant \; &C_0\; \frac{\log\,
|A^*_x|}{|\wtI|\,\log\,|\wtI|},\label{eq7.176}\\
|\pa_t \log\,|\wtA_x|| \; &\leqslant \; &C_0\; \frac{\log\,
|\wtA_x|}{|I|\,\log\,|I|},\label{eq7.177}\\
|\pa_t \log\,|\wtB_y|| \; &\leqslant \; &C_0\; \frac{\log\,
|\wtB_y|}{|I|\,\log\,|I|}.\label{eq7.178}
\end{eqnarray}
Here $C_0$ is large but independent of $\vep_0$. This means that the
term $\de^{-\fudt} (|\wtB_{xt}| + |A^*_{yt}|)$ in (\ref{eq7.171})
and (\ref{eq7.173}) is dominated by $\de^{-1}$. As $|I| =
|\wtI|^{1+\tau}$, in order to prove (\ref{eq7.122}) by induction, we
need to have, in view of (\ref{eq7.171}):
\begin{equation}
C |I|\, |\log |I||\, \de^{-1} + C\,C_0 |\wtI|^\tau |\log |A^*_x|| +
C_0 |\log |\wtA_x|| \; \leqslant \; C_0 |\log |A_x||.\label{eq7.179}
\end{equation}

We have here $\de\geqslant 2 |I|$ from the definition of the
transversality relation and, by (\ref{eq3.27}):
\begin{equation}
|\log |A_x|| \; \geqslant \; |\log |\wtA_x|| + |\log |A^*_x|| -
\fud\, |\log |I||.\label{eq7.180}
\end{equation}
Therefore, (\ref{eq7.179}) will hold as far as
\begin{equation}
\Bigl( \frac{C_0}{2} + C\Bigr)\,|\log |I|| \; \leqslant \; C_0 (1 -
C |\wtI|^\tau) \, |\log |A^*_x||.\label{eq7.181}
\end{equation}
From (R7), we know that $|A^*_x|$ is much smaller than $\de$. On the
other hand, as $Q_i$ and $\ism P$ are not $\wtI$-transverse, $\de$
cannot be much larger than $\wtI$. Therefore, we must have
\begin{equation}
|\log |A^*_x|| \; \geqslant \; \log\, |\wtI| = (1 + \tau) \, \log
|I|,\label{eq7.182}
\end{equation}
from which (\ref{eq7.181}) follows if we take $C_0\geqslant 3C$.
This completes the proof of the induction step for (\ref{eq7.122}).

To do the same for (\ref{eq7.123}), we estimate the right-hand side
of (\ref{eq7.172}). From the proof of Proposition~\ref{propo15},
formula (\ref{eq7.97}), we have
\begin{equation}
|\wtA_x| \; \leqslant \; C^i \, |\wtI|^{(i+1)\be - \tfrac{i}{2}\,
(1+\tau)}\label{eq7.183}
\end{equation}
From (R7), we have
\begin{equation}
|\wtB_y| \; \leqslant \; C \,\de^{(1-\eta)^{-1}}.\label{eq7.184}
\end{equation}
The displacement $\de = \de (Q^\pii, \ism P)$ satisfies
\begin{equation}
2 |I| \; \leqslant \; \de \; \leqslant \; C |\wtI|.\label{eq7.185}
\end{equation}
This gives (as $\be > 1$)
\begin{eqnarray}
\de^{-\fudt} |\wtA_x| \; &\leqslant \; &|\wtI|^{\be/2\, (i+1)},\label{eq7.186}\\
\de^{-\tfrac{3}{2}} |\wtA_x|\,|\wtB_y| \; &\leqslant \; &|\wtI|^{\be/2\, (i+1)},\label{eq7.187}\\
\de^{-\fudt} |\wtA_x|\,|\wtB_y|\,|\pa_t\log |\wtA_x|| \; &\leqslant
\; &|\wtI|^{\be/2\, (i+1)},\label{eq7.188}\\
\de^{-\fudt} |\wtA_x|\,|\wtB_y|\,|\pa_t\log |\wtB_y|| \; &\leqslant
\; &|\wtI|^{\be/2\, (i+1)}.\label{eq7.189}
\end{eqnarray}
This leads to:
\begin{equation}
|A_{yt} - \wtA_{yt}| \; \leqslant \; C\, |\wtI|^{\be/2\,
(i+1)}.\label{eq7.190}
\end{equation}
We can now sum over $i$ and then over the different levels of
parameter intervals to obtain (\ref{eq7.123}). The proof of
Proposition~\ref{propo17} is complete.
\end{proof}

Let us see that property (R7) of Subsection~\ref{sub5.4} is a
consequence of Proposition~\ref{propo17}. Consider $(P_0, Q_0,
n_0)$, $(P_1, Q_1, n_1)$ to be elements of $\cR (I)$ such that $Q_0
\subset Q_u$, $P_1 \subset P_s$ and $Q_0 \forki P_1$ holds.

We can assume that there are no $\whI \supset I$, $(\whP_0, \whQ_0,
\whn_0)$, $(\whP_1, \whQ_1, \whn_1)\in \cR (\whI)$ with $Q_0\subset
\whQ_0$, $P_1 \subset \whP_1$, $(Q_0, P_1, I) \ne (\whQ_0, \whP_1,
\whI)$ and $\whQ_0 \fork_{\whI} \whP_1$: otherwise, as
$|Q_0|\leqslant |\whQ_0|$, $|P_1|\leqslant |\whP_1|$ and $\de (Q_0,
P_1) \geqslant \de (\whQ_0, \whP_1)$, property (R7) for $(Q_0, P_1,
I)$ would be inherited from $(\whQ_0, \whP_1, \whI)$.

In view of this maximality property, we claim that $Q_0$ must
satisfy condition $(*u)$, and $P_1$ must satisfy condition $(*s)$.

Let us first finish the proof of (R7) assuming the claim to be true.
Indeed, we have $\de (Q_0, P_1)\geqslant 2 |I|$ for all $t\in I$ by
(T1) of Subsection~\ref{sub5.4}. If $|\log |Q_0||$ is much larger
than $|\log |I||$ for all $t\in I$, we obviously have $\de (Q_0,
P_1)\geqslant |Q_0|^2$; but if $|\log |Q_0||\leqslant C |\log |I||$
for some $t\in I$, we obtain from Proposition~\ref{propo17} that
\begin{equation}
\max_I\; |Q_0| \; \leqslant \; C\, \min_I\; |Q_0|.\label{eq7.191}
\end{equation}
We also know from Proposition~\ref{propo15} that
\begin{equation}
\max_I\; \de (Q_0, P_1) \; \leqslant \; C\, \min_I\; \de (Q_0,
P_1).\label{eq7.192}
\end{equation}
It then follows from (T2) in Subsection~\ref{sub5.4} that
\begin{equation}
\de (Q_0, P_1) \; \geqslant \; C^{-1} |Q_0|^{1-\eta}.\label{eq7.193}
\end{equation}
for all $t\in I$. We argue with $P_1$ in a symmetric way.

Finally, we prove the claim. Let us show, for instance, that $P_1$
satisfies condition $(*s)$. Let $\Is$ be the largest parameter
interval such that $(P_1, Q_1, n_1)\in \cR(\Is)$. If $\Is$ is the
starting interval $I_0$, $P_1$ satisfies condition $(*s)$. Assume
therefore that $\Is\ne I_0$; let $\wtI^*$ be the parent interval,
$(\wtP_1, \wtQ_1, \wtn_1)$ the element of $\cR (\wtI^*)$ such that
$\wtP_1$ is the thinnest rectangle containing $P_1$. We have to show
that $\wtP_1$ is $\wtI^*$-critical. Assume by contradiction that
$\wtP_1$ is $\wtI^*$-transverse. Then, there exists an
$\wtI^*$-decomposition $(P_\al, Q_\al, n_\al)$ of $Q_u$ such that,
for every $\al$, $Q_\al$ and $\wtP_1$ are $\wtI^*$-separated or
$\wtI^*$-transverse. Let $\al_0$ be such that $Q_{\al_0}$ and $Q_0$
intersect. Then, $Q_{\al_0}$ and $\wtP_1$ must be
$\wtI^*$-transverse. But then $\wtP_1$ and $Q_0$ must be
$I$-transverse, either from Proposition~\ref{propo2} if $Q_0 \subset
Q_{\al_0}$, or from concavity (Proposition~\ref{propo4}) if
$Q_{\al_0}\subset Q_0$; this contradicts the maximality of $(Q_0,
P_1, I)$ and proves the claim.

The proof of property (R7) is complete.

The existence and properties of the classes $\cR (I)$ are now fully
justified. What we do not know at this moment is whether there
exists any regular parameter interval at all! This will be the
subject of Section~\ref{sec9}. Before, we develop in the next
section some results that will turn out to be essential in
Sections~\ref{sec9} and \ref{sec10}.

\newpage

\setcounter{section}{7}
\setcounter{equation}{0}

\section{Number of Children and Dimension Estimates\label{sec8}}
\subsection{Estimates on the Number of Children\label{sub8.1}}
We start with some preliminary results.

\begin{Propo}\label{propo18}
Let $I\subset \wtI$ be parameter intervals, and let $(\wtP, \wtQ,
\wtn)$ be an element of $\cR (\wtI)$. We assume that $\wtQ$ is
$\wtI$-transverse. Then, any element $(P,Q,n)$ in $\cR (I)$ such
that $P$ is a child of $\wtP$ is already an element of $\cR (\wtI)$.
\end{Propo}

\begin{proof}
We can assume that $P$ is a non-simple child. Then $(P,Q,n)$ is
obtained by parabolic composition in $\cR (I)$ of $(\wtP, \wtQ,
\wtn)$ with some $(P_1, Q_1, n_1)\in \cR (I)$. As $\wtQ$ is
$\wtI$-transverse, there exists an $\wtI$-decomposition $(P_\al,
Q_\al, n_\al)$ of $P_s$ such that each $P_\al$ is $\wtI$-separated
or $\wtI$-transverse with $\wtQ$. Let $\al_0$ be such that
$P_{\al_0}$ and $P_1$ intersect. Then, $\wtQ \fork_{\wtI} P_{\al_0}$
holds, and also $\wtQ\forki P_1$; if we had $P_1 \subsetneqq
P_{\al_0}$, this would imply that $\wtQ$ would be $I$-transverse to
the parent $\wtP_1$ of $P_1$ and $P$ would not be a child of $\wtP$.
Therefore, we must have $P_{\al_0} \subset P_1$. By coherence
(Proposition~\ref{propo6}), we have that $(P_1, Q_1, n_1)\in \cR
(\wtI)$. By parametric concavity (Proposition~\ref{propo7}), from
$\wtQ \forki P_1$ and $\wtQ \fork_{\wtI} P_{\al_0}$, we deduce that
$\wtQ \fork_{\wtI} P_1$ also holds and $(P,Q,n)\in \cR (\wtI)$.
\end{proof}

\begin{Propo}\label{propo19}
Let $I$ be a parameter interval, and let $I_1 \supset I$ be the
largest parameter interval such
that
\begin{equation}
|I_1|^\be \; < \; \Bigl( \fud \;
|I|\Bigr)^{\tfrac{1}{1-\eta}}.\label{eq8.1}
\end{equation}
Let $(\wtP, \wtQ, \wtn)$, $(P,Q,n)$ be elements of $\cR (I)$ such
that $P$ is a non-simple child of $\wtP$. Let $(P_1, Q_1, n_1)$,
$(\wtP_1, \wtQ_1, \wtn_1)$ be the elements of $\cR (I)$ such that
\[
(P, Q, n) \; \in \; (\wtP, \wtQ, \wtn)\; \square\; (P_1, Q_1, n_1)
\]
and $\wtP_1$ is the parent of $P_1$.

Then, $(P_1, Q_1, n_1)$ belongs to $\cR (I_1)$, $\wtP_1$ is
$I$-critical, $\wtQ_1$ is $I_1$-transverse and we have
\begin{equation}
2 |\wtP_1|^{1-\eta} \; > \; |I|\label{eq8.2}
\end{equation}
for all $t\in I$.
\end{Propo}

\begin{rem}
As parabolic composition is possible, we have $I\ne I_0$; then, as
$\be > (1-\eta)^{-1}$, we must have $I_1 \supsetneqq I$ and $I_1$ is $\beta$-regular.
\end{rem}

\begin{proof}
That $\wtP_1$ is $I$-critical has already been proved in
Proposition~\ref{propo5}. Also, as $\wtQ \forki P_1$ holds but
$\wtQ\forki \wtP_1$ does not hold (because $P$ is a non-simple child
of $\wtP$), we deduce (\ref{eq8.2}) from Proposition~\ref{propo9}.
Then, by definition of $I_1$, we have:
\begin{equation}
|\wtP_1| \; > \; |I_1|^\be\label{eq8.3}
\end{equation}
for all $t\in I$. Let us show that $(\wtP_1, \wtQ_1, \wtn_1)$
belongs to $\cR (I_1)$. Otherwise, there would exist $I_2\supset I$,
with parent interval $\wtI_2 \subset I_1$, such that $(\wtP_1,
\wtQ_1, \wtn_1)$ belongs to $\cR (I_2)$ but not to $\cR (\wtI_2)$.
We apply the inequality (\ref{eq6.74}) in the proof of
Corollary~\ref{coro9} (Subsection~\ref{sub6.7}) to get
\begin{equation}
|\wtP_1| \; \leqslant \; C |\wtI_2|^{\be +
\tfrac{1}{3}},\label{eq8.4}
\end{equation}
in contradiction with (\ref{eq8.3}). Therefore, $(\wtP_1, \wtQ_1,
\wtn_1)$ belongs to $\cR (I_1)$. As $I_1\supsetneqq I$ , $I_1$ is
$\be$-regular ; $(\wtP_1, \wtQ_1, \wtn_1)$ cannot be $I_1$-bicritical in view of (\ref{eq8.3}) ; $\wtP_1$ is $I_1$-critical and hence
$\wtQ_1$ is $I_1$-transverse. Proposition~\ref{propo18} then shows
that $(P_1, Q_1, n_1)\in \cR (I_1)$.
\end{proof}

\begin{Coro}\label{coro12}
Let $I$ be a parameter interval and let $(\wtP, \wtQ, \wtn)$ be an
element of $\cR (I)$. The number of $(P,Q,n)\in \cR (I)$ such that
$P$ is a child of $\wtP$ is finite.
\end{Coro}

\begin{proof}
We argue by induction on the level of the parameter interval.

If $I$ is the starting interval $I_0$, $\wtP$ has only simple
children and the assertion is obvious. Assume that $I\subsetneqq
I_0$. The number of simple children is finite, and we have to show
that the same is true for the number of non-simple children. For
every non-simple child $P$ of $\wtP$, let $I_1$, $(P_1, Q_1, n_1)$,
$(\wtP_1, \wtQ_1, \wtn_1)\in \cR (I_1)$ be as in
Proposition~\ref{propo19}. By the induction hypothesis, there is for
each fixed $\wtP_1$ only a finite number of  possibilities for
$P_1$. On the other hand, in view of relation (\ref{eq8.2}), there
are obviously only a finite number of possibilities for $\wtP_1$.
The induction step is complete, and this completes the proof.
\end{proof}

We want to make the finiteness assertion quantitative, and will do
that in two distinct ways. In each case, we have to estimate in the
proof of Corollary~\ref{coro12} the number of possibilities for
$\wtP_1$, and the number of possibilities for $P_1$ once $\wtP_1$ is
fixed.

\begin{Propo}\label{propo20}
Let $I$ be a parameter interval, and let $(\wtP, \wtQ, \wtn)$ be an
element of $\cR (I)$. The number of $(P,Q,n)\in \cR (I)$ such that
$P$ is a child of $\wtP$ is at most $|I|^{-c\eta}$, where $c$ is a
constant depending only on $\be$.
\end{Propo}

\begin{proof}
We argue again by induction on the level of $I$, following the proof
of Corollary~\ref{coro12}. When $I = I_0$, the number of (simple)
children is at most the number of rectangles in the Markov
partition, which is much smaller than $\vep_0^{-c\eta}$ when
$\vep_0$ is small enough.

When $I\ne I_0$, the number of possibilities for $P_1$ when $\wtP_1$
is fixed in at most $|I_1|^{-c\eta}$ by the induction hypothesis. We
have to estimate the number of possibilities for $\wtP_1$. We know
that $\wtQ \forki \wtP_1$ does not hold, but $\wtQ \forki P_1$
holds.

As $\wtP_1$ is $I$-critical, it satisfies condition $(*s)$, defined
just before Proposition~\ref{propo15} in Subsection~\ref{sub7.5}. We
have from (\ref{eq8.2}) and Proposition~\ref{propo17}
\begin{equation}
\max_I \; |\wtP_1| \; \leqslant \; C \,\min_I\,
|\wtP_1|.\label{eq8.5}
\end{equation}

\begin{lema}\label{lema3}
We have for all $t\in I$
\begin{equation}
\de (\wtQ, \wtP_1) \; \leqslant \; C \,
|\wtP_1|^{1-\eta}.\label{eq8.6}
\end{equation}
\end{lema}

\begin{proof}
Let $(\whP, \whQ, \whn)$ be the element of $\cR (I)$ with smallest
$\whn$ such that $\wtQ\subset \whQ$ and $\whQ \forki P_1$ holds. The
parent $Q^*$ of $\whQ$ is $I$-critical: otherwise, from an
$I$-decomposition of $P_s$, one would find $(P_\al, Q_\al, n_\al)$
with $Q^* \forki P_\al$ and $P_\al \cap P_1\ne \emptyset$; one would
conclude that $Q^* \forki P_1$ from Proposition~\ref{propo2} (if
$P_1 \subset P_\al)$ or Proposition~\ref{propo4} (if $P_\al \subset
P_1$). Thus, $\whQ$ satisfies condition $(*u)$. From
Corollary~\ref{coro11} and Proposition~\ref{propo17}, we get
\begin{equation}
\max_I \; |\whQ| \; \leqslant \; C \,\min_I\, |\whQ|,\label{eq8.7}
\end{equation}
\begin{equation}
\max_I \; \de_{LR} (\whQ, \wtP_1) - \min_I\,\de_{LR} (\whQ, \wtP_1)
\; \leqslant \; 2 |I|.\label{eq8.8}
\end{equation}
As $\wtQ \forki \wtP_1$ does not hold, $\whQ \ovforki \wtP_1$ does
not hold either and at least one of the following three inequalities
must hold:
\begin{eqnarray}
\de_{LR} (\whQ, \wtP_1) \, &< \, &2 |I|\; \text{ for some $t_0\in
I$};\label{eq8.9}\\
\de_{R} (\whQ, \wtP_1) \, &< \, &2 |\whQ|^{1-\eta}\; \text{ for all
$t\in I$};\label{eq8.10}\\
\de_{L} (\whQ, \wtP_1) \, &< \, &2 |\wtP_1|^{1-\eta}\; \text{ for
all $t\in I$}.\label{eq8.11}
\end{eqnarray}
By Proposition~\ref{propo10}, as $\whQ \forki P_1$ holds but $\whQ
\forki \wtP_1$ does not hold, we have, for some $t_1\in I$:
\begin{equation}
|\wtP_1| \; > \; \fud\, |\whQ|.\label{eq8.12}
\end{equation}
We can now prove (\ref{eq8.6}).

If (\ref{eq8.11}) holds, we have, for all $t\in I$:
\begin{equation}
\de (\wtQ, \wtP_1) \; \leqslant \; \de_L (\wtQ, \wtP_1) \; \leqslant
\; \de_L (\whQ, \wtP_1) \; < \; 2 |\wtP_1|^{1-\eta}.\label{eq8.13}
\end{equation}
If (\ref{eq8.12}) holds, we have from (\ref{eq5.12}), (\ref{eq8.5}),
(\ref{eq8.7}), (\ref{eq8.12}):
\begin{eqnarray}
\de (\wtQ, \wtP_1) \; &\leqslant \; &\de_{LR} (\whQ,
\wtP_1)\label{eq8.14}\\
&\leqslant \; &\de_{R} (\whQ, \wtP_1) + c |\whQ|\notag\\
&\leqslant \; &3\, |\whQ|^{1-\eta}\; \leqslant \; C
|\wtP_1|^{1-\eta},\notag
\end{eqnarray}
for all $t\in I$.

Finally, if (\ref{eq8.9}) holds, we have from (\ref{eq8.5}),
(\ref{eq8.8}), (\ref{eq8.2}), for all $t\in I$:
\begin{eqnarray}
\de (\wtQ, \wtP_1) \; &\leqslant \; &\de_{LR} (\whQ,
\wtP_1)\label{eq8.15}\\
&\leqslant \; &4\, |I|\; \leqslant \; 8 |\wtP_1|^{1-\eta}.\notag
\end{eqnarray}
\end{proof}

We are now able to estimate the number of possibilities for $\wtP_1$
and show that this number is at most
\begin{equation}
C\, |I|^{-\tfrac{\eta}{1-\eta}}.\label{eq8.16}
\end{equation}
This indeed follows from (\ref{eq8.6}), (\ref{eq8.2}) and the fact
that if two distinct $\wtP_1$ are not disjoint, the ratio of their
widths is bounded away from 1 (so the $\wtP_1$ at a given scale are
disjoint; one then sums over scales). The total number of children
is thus bounded by
\begin{equation}
C + 2C\, |I|^{-\tfrac{\eta}{1-\eta}}\; |I_1|^{-c\eta},\label{eq8.17}
\end{equation}
where $I_1$ was the largest parameter interval
satisfying
\begin{equation}
|I_1|^\be \; < \; \Bigl( \fud \;
|I|\Bigr)^{\tfrac{1}{1-\eta}}.\label{eq8.18}
\end{equation}
If $|I| > 2\vep_0^{\be (1-\eta)}$, we have $I_1 = I_0$; in this
case, the term $|I_1|^{-c\eta}$ in (\ref{eq8.17}) is unnecessary
because $\wtP_1$ has only simple children. If $|I| \leqslant 2
\vep_0^{\be(1-\eta)}$, we have
\begin{equation}
|I_1|^{\be (1+\tau)^{-1}} \; \geqslant \; \Bigl( \fud \;
|I|\Bigr)^{\tfrac{1}{1-\eta}}\label{eq8.19}
\end{equation}
and the term in (\ref{eq8.17}) is bounded by $|I|^{-c\eta}$ provided
that
\begin{equation}
c\eta \; > \; \frac{\eta}{1-\eta} + c\eta \; \frac{1+\tau}{1-\eta}\;
\be^{-1}. \label{eq8.20}
\end{equation}
As $\eta$, $\tau$ are very small, any choice of $c >
\frac{\be}{\be-1}$ yields (\ref{eq8.20}). Then, as $c>1$, such a
choice is also convenient when $|I| > 2\, \vep_0^{\be(1-\eta)}$, and
this concludes the proof of Lemma~\ref{lema3}.
\end{proof}

In Proposition~\ref{propo20}, we have estimated the total number of
children in terms of the level of the parameter interval.

When $\wtQ$ is $I$-transverse, Proposition~\ref{propo18} guarantees
that there will not be any new child of $\wtQ$ when we consider
parameter intervals $\whI \subset I$ of higher level.

The same is true when for all $t\in I$
\begin{equation}
|\wtQ|^{1-\eta} \; \geqslant \; C |I|.\label{eq8.21}
\end{equation}
with some large enough constant $C$. Indeed, let $\Is$ be a
parameter interval strictly smaller than $I$, with parent interval
$\wtI^*$, and let $(P,Q,n)$ be an element of $\cR (\Is)$ such that
$P$ is a non-simple child of $\wtP$. Let $(P_1, Q_1, n_1)$,
$(\wtP_1, \wtQ_1, \wtn_1)$ be as in the previous propositions. By
Proposition~\ref{propo19}, the element $(P_1, Q_1, n_1)$ belongs to
$\cR (\wtI^*)$. It is then easy to deduce from $\wtQ \fork_{\Is}
P_1$ and (\ref{eq8.21}), using (R7) and $\wtI^* \subset I$, that
$\wtQ \fork_{\wtI^*} P_1$ also holds. This proves by induction that
$(P,Q,n)$ belongs to $\cR (I)$.

In the next proposition, we are interested, not in the total number
of children, but in the number of children of a given width. The estimate is independent on the
level of the parameter interval.

\begin{Propo}\label{propo21}
Let $I$ be a parameter interval, and let $(\wtP, \wtQ, \wtn)$ be an
element of $\cR (I)$. For any $\vep>0$, the number of elements
$(P,Q,n)\in \cR (I)$ such that $P$ is a non-simple child of $\wtP$
satisfying $|P|\geqslant \vep |\wtP|$ for some $t\in I$, is at most
$\vep^{-c'\eta}$, where $c'$ is a constant depending only on the
regularity parameter $\be$.
\end{Propo}

\begin{proof}
Let $\vep > 0$, and let $(P,Q,n)$ be an element of $\cR (I)$ such
that $P$ is a non-simple child of $\wtP$. We assume, for some
$t_0\in I$, that:
\begin{equation}
|P| \; \geqslant \; \vep \, |\wtP|.\label{eq8.22}
\end{equation}
Let $(P_1, Q_1, n_1)$, $(\wtP_1, \wtQ_1, \wtn_1)\in \cR (I)$ be as
in Proposition~\ref{propo19}. From (\ref{eq3.27}), we have, for all
$t\in I$:
\begin{equation}
|P| \; \leqslant \; C
\,|\wtP|\;|P_1|\,\de\,(\wtQ,P_1)^{-\fudt}.\label{eq8.23}
\end{equation}
Property (R7) guarantees that, for all $t\in I$
\begin{equation}
\de (\wtQ, P_1) \; \geqslant \; C^{-1}\,
|P_1|^{1-\eta}.\label{eq8.24}
\end{equation}
Combining (\ref{eq8.22}), (\ref{eq8.23}), (\ref{eq8.24}), we have,
for some $t_0\in I$
\begin{equation}
\de (\wtQ, P_1) \; \geqslant \; C^{-1}\,
\vep^{2\;\tfrac{1-\eta}{1+\eta}}.\label{eq8.25}
\end{equation}
As we always have
\begin{equation}
\de (\wtQ, P_1) \; \leqslant \; C \vep_0,\label{eq8.26}
\end{equation}
there is no non-simple child satisfying (\ref{eq8.22}) unless $\vep
< \vep_0^{\fud}$; we will assume that this holds in the sequel.

From Lemma~\ref{lema3} above, we have, for all $t\in I$:
\begin{equation}
\de (\wtQ, \wtP_1) \; \leqslant \; C\,
|\wtP_1|^{1-\eta},\label{eq8.27}
\end{equation}
and thus, from (\ref{eq5.11}), also
\begin{eqnarray}
\de (\wtQ, P_1) \; &\leqslant \; &\de_L (\wtQ,\wtP_1)\label{eq8.28}\\
&\leqslant \; &\de (\wtQ, \wtP_1) + C\,|\wtP_1|\notag\\
&\leqslant \; &C\, |\wtP_1|^{1-\eta}.\notag
\end{eqnarray}
Combining (\ref{eq8.25}) and (\ref{eq8.28}), we get, for some
$t_0\in I$
\begin{equation}
|\wtP_1| \; \geqslant \; C^{-1}\,
\vep^{\tfrac{2}{1+\eta}},\label{eq8.29}
\end{equation}
an inequality which actually holds for all $t\in I$ in view of
(\ref{eq8.5}).

As in the proof of Proposition~\ref{propo20}, the number of
$(\wtP_1, \wtQ_1, \wtn_1)$ for which both (\ref{eq8.28}),
(\ref{eq8.29}) hold is easily seen to be at most
\begin{equation}
C\; \vep^{-\tfrac{2\eta}{1+\eta}}.\label{eq8.30}
\end{equation}
On the other hand, let $\whI = I$ if $|I|\geqslant \vep^{2
(\be+1/3)^{-1}}$; otherwise, let $\whI$ be the largest parameter
interval which contains $I$ and satisfies $|\whI| < \vep^{2(\be +
1/3)^{-1}}$. As $\wtP_1$ is $I$-critical, the same argument as in
the proof of Corollary~\ref{coro9} in Subsection~\ref{sub6.7} shows,
from (\ref{eq8.29}), that we must have
\begin{equation}
(\wtP_1, \wtQ_1, \wtn_1) \; \in \; \cR (\whI).\label{eq8.31}
\end{equation}
Now, $\wtP_1$ is $\whI$-critical. If $\whI\ne I$, we have, from
(\ref{eq8.29}) and the definition of $\whI$, that
\begin{equation}
|\wtP_1|\; > \; |\whI|^\be,\;\;\; \text{ for all $\;t\in
I$}.\label{eq8.32}
\end{equation}
Therefore, $(\wtP_1, \wtQ_1, \wtn_1)$ cannot be $\whI$-bicritical,
$Q_1$ is $\whI$-transverse and we conclude from
Proposition~\ref{propo18} that $(P_1, Q_1, n_1)$ also belongs to
$\cR (\whI)$. The same is also obviously true when $\whI = I$. We
apply Proposition~\ref{propo20}: for each fixed $(\wtP_1, \wtQ_1,
\wtn_1)$, the number of children $P_1$ is at most $|\whI|^{-c\eta}$.
But we have
\begin{alignat}{3}
|\whI| &= &|I|\, \geqslant \, \vep^{2 (\be + 1/3)^{-1}}\;\;\;
&\text{ if $\;|I|\; \geqslant \; \vep^{2(\be+1/3)^{-1}}$}\label{eq8.33}\\
|\whI| &= &\vep_0\qquad\qquad\qquad\; &\text{ if
$\;\vep^{2 (\be + 1/3)^{-1}} > \vep_0$}\label{eq8.34}\\
|\whI| &\geqslant &\vep^{2 (1+\tau)(\be + 1/3)^{-1}}\;\;\; &\text{
if $\;\vep_0 \geqslant \vep^{2 (\be + 1/3)^{-1}} >
|I|$}\label{eq8.35}
\end{alignat}

In all cases, this gives
\begin{equation}
|\whI|^{-c\eta} \; \geqslant\; \vep^{-c_0\eta}.\label{eq8.36}
\end{equation}
Combining Proposition~\ref{propo20} with the previous estimate in (\ref{eq8.30}) for the
number of possibilities for $\wtP_1$ gives therefore the required estimate.
\end{proof}

\subsection{A Dimension Estimate\label{sub8.2}}
The goal of this subsection is to obtain a bound on the number of
elements $(P,Q,n)$ in $\cR(I)$ with width $|P|$ bounded from below.
This is a first step towards estimating the transverse dimension of
the stable set $W^s (\La)$, which is necessary in order to achieve
our parameter selection in Section~\ref{sec9}.

Let $I$ be a parameter interval, and let $(\Ps, \Qs, \ns)$ be an
element of $\cR (I)$. We introduce, in the spirit of Laplace,
Dirichlet and Poincaré, the series
\begin{equation}
\Te (\Ps, I, s) \; = \; \sum \; |P|^s,\label{eq8.37}
\end{equation}
where the sum runs over elements $(P,Q,n)\in \cR (I)$ such that
$P\subset \Ps$. Here $s$ is a complex variable and the series is at
first a formal object, but we will soon see that it is uniformly
convergent in a half-plane $\{ \Re\,s > \sig_0\}$. The goal of this
subjection is to obtain a nice estimate for $\sig_0$ and for $\Te$
in this half-plane.

The width $|P|$ and therefore also the series $\Te$, depend on the
parameter $t\in I$; but the estimate that we will get is uniform
with respect to the parameter. The dependence of the estimate on
$\Ps$ is also quite straightforward, through a simple scaling
factor.

Let us recall that we denote by $\dos$ the transverse Hausdorff
dimension of the stable foliation $W^s (K)$ of the initial
horseshoes $K$ for the value $0$ of the parameter. It is well-known
that this transverse Hausdorff dimension depends smoothly on the
parameter, and it controls in a precise way the number of cylinders
(for the Markov partition) of a given size; more precisely, as these
cylinders correspond exactly to the elements of $\cR (\Io)$, we know
that, for all $t\in \Io$ and all $\vep>0$ the number of $(P,Q,n)\in
\cR (\Io)$ such that $|P|\geqslant \vep$ is at most
\begin{equation}
C\, \vep^{-(\dos + C\vep_0)}.\label{eq8.38}
\end{equation}
This shows that for $\Te (\Ps, \Io; s)$ we could take $\sig_0 = \dos
+ C\vep_0$. For smaller parameter intervals, we have to allow a
slightly larger margin with relation to the initial value $\dos$.

\begin{Propo}\label{propo22}
The series $\Te (\Ps, I, s)$ is uniformly convergent in the
half-plane\linebreak $\{ \Re\, s\geqslant \dos +
\vep_0^{\frac{1}{3}\;\dos}\}$. When $(\Ps, \Qs, \ns)\in \cR (\Io)$,
we have for $\Re\, s\geqslant \dos + \vep_0^{\frac{1}{5}\;\dos}$
\[
\lv \Te (\Ps, I; s) - \Te (\Ps, \Io; s)\rv \; \geqslant \; |\Ps|^s
\vep_0^{\frac{1}{20}\,\dos}.
\]
\end{Propo}

\begin{proof}
Let $(P,Q,n)$ be an element of $\cR (I)$ with $P\subset \Ps$.
Consider the intermediary rectangles
\begin{equation}
\Ps = P (0) \subset P (1) \subset \cdots \subset P (\ell) =
P\label{eq8.39}
\end{equation}
with $P (i)$ the parent of $P (i+1)$. Let
\begin{equation}
\ell_0 \; < \; \ell_1 \; < \; \cdots \; < \;
\ell_{k-1}\label{eq8.40}
\end{equation}
be the indices such that $P (\ell_j + 1)$ is a non-simple child of
$P (\ell_j)$.

We also define for $0\leqslant j\leqslant k$ elements $(P^\pjj,
Q^\pjj, n^\pjj)\in \cR (\Io)$ by the following properties
\begin{equation}
(P (\ell_0), Q (\ell_0), n (\ell_0)) = (P(0), Q(0), n(0) \; *\;
(P^\pzz, Q^\pzz, n^\pzz)),\label{eq8.41}
\end{equation}
\begin{equation}
(P(\ell_j), Q(\ell_j), n(\ell_j)) = (P(\jsM\ell + 1), Q(\jsM\ell+1),
n(\jsM\ell + 1) * (P^\pjj, Q^\pjj, n^\pjj),\label{eq8.42}
\end{equation}
\begin{equation}
(P,Q,n) = (P(\ksM\ell + 1), Q(\ksM\ell+1), n(\ksM\ell + 1)
* (P^\pkk, Q^\pkk, n^\pkk).\label{eq8.43}
\end{equation}

We now estimate the widths from (\ref{eq3.10}):
\begin{eqnarray}
|P(\ell_0)| \; &\leqslant \; &C |\Ps|\; |P^\pzz|,\label{eq8.44}\\
|P(\ell_j)| \; &\leqslant \; &C |P (\jsM\ell + 1)|\;|P^\pjj|,\label{eq8.45}\\
|P| \; &\leqslant \; &C |P (\ksM\ell + 1)|\;|P^\pkk|.\label{eq8.46}
\end{eqnarray}
From (\ref{eq3.27}) and property (R7), we also have:
\begin{equation}
|P(\ell_j+1)| \; < \; \vep_0^{\fud} |P (\ell_j)|.\label{eq8.47}
\end{equation}

Define $m_j$ for $0\leqslant j < k$ to be the largest integer such
that, for all $t\in I$
\begin{equation}
|P(\ell_j+1)| \; \leqslant \; 2^{-m_j} \vep_0^{\fud} |P
(\ell_j)|.\label{eq8.48}
\end{equation}
From Proposition~\ref{propo21}, for each fixed $P(\ell_j)$, the
number of non-simple children $P(\ell_j+1)$ satisfying
(\ref{eq8.48}) is at most
\begin{equation}
\Bigl ( 2^{m_j+1}\; \vep_0^{-\fud} \Bigr)^{c'\eta}.\label{eq8.49}
\end{equation}
Combining (\ref{eq8.44}), (\ref{eq8.45}), (\ref{eq8.46}) and
(\ref{eq8.48}), we also have
\begin{equation}
|P|\; \leqslant \; C^{k+1} |\Ps| \Bigl( \prod^k_0 |P^\pjj|\Bigr)\;
\vep_0^{\frac{k}{2}} 2^{-\Sigma m_j}.\label{eq8.50}
\end{equation}

We will take this to the power $s$ and sum over $P$. We introduce
(corresponding to the term $|P^\pjj|^s$)
\begin{eqnarray}
\Te_0 (s) :&= &\sum_{(\whP,\whQ,\whn)\in \cR(\Io)}
|\whP|^s\label{eq8.51}\\
&= &\sum_{\AL} \; \Te (R_a, \Io, s),\notag
\end{eqnarray}
and also
\begin{equation}
\te (s) := \sum_{m\geqslant 0} (C\vep_0^{\fud} \; 2^{-m})^s (2^{m+1}
\vep_0^{-\fud})^{c'\eta}.\label{eq8.52}
\end{equation}
The function $\Te_0$ is controlled by (\ref{eq8.38}), while $\te$
satisfies
\begin{equation}
\te (s) = 2^{c'\eta}\, C^s \, \vep_0^{\fud (s-c'\eta)} \; \Bigl(
1-2^{-(s-c'\eta)}\Bigr)^{-1},\label{eq8.53}
\end{equation}
and therefore, for $C^{-1} < s < C$:
\begin{equation}
C^{-1} \, \vep_0^{\fud (s-c'\eta)} \; \leqslant \; \te (s)\;
\leqslant \; C\, \vep_0^{\fud (s-c'\eta)}.\label{eq8.54}
\end{equation}
From (\ref{eq8.38}), we have, for $s > \dos + C\vep_0$:
\begin{equation}
\Te_0 (s) \; \leqslant \; C\, (s-\dos - C\vep_0)^{-1}.\label{eq8.55}
\end{equation}
In particular, for $s\geqslant \dos + \vep_0^{1/3\,\dos}$, we have
\begin{eqnarray}
\Te_0 (s)\;\leqslant \; C \vep_0^{-\frac{1}{3}\,\dos},\label{eq8.56}\\
\Te_0 (s)\,\te (s) \; \leqslant \; C
\vep_0^{\frac{1}{10}\,\dos}.\label{eq8.57}
\end{eqnarray}
But, from (\ref{eq8.50}), we have for real $s$
\begin{equation}
\Te (\Ps,I,s)\; \leqslant \; C^s |\Ps|^s \sum_{k\geqslant 0}
\Te_0^{k+1} (s)\, \te^k (s),\label{eq8.58}
\end{equation}
and therefore we deduce from (\ref{eq8.57}) that the series defining
$\Te$ is uniformly convergent in the half plane $\{\Re\, s\geqslant
\dos + \vep_0^{\frac{1}{3}\,\dos}\}$.

Assume now that $(\Ps,\Qs,\ns)\in \cR(\Io)$; the difference $\Te
(\Ps,I,s) - \Te (\Ps,\Io,s)$ consists of the sum of $|P|^s$ over
those $P$ for which $k>0$. We get, for real $s$
\begin{equation}
\Te (\Ps,I,s) - \Te (\Ps,\Io,s) \; \leqslant \; C^s |\Ps|^s
\sum_{k>0} \Te_0^{k+1} (s)\, \te^k (s).\label{eq8.59}
\end{equation}
For $s>\dos + \vep_0^{\frac{1}{5}\,\dos}$, we have, from
(\ref{eq8.55}), (\ref{eq8.54}):
\begin{eqnarray}
\Te_0 (s)\;&\leqslant \; &C \vep_0^{-\frac{1}{5}\,\dos},\label{eq8.60}\\
\Te_0^2 \,\te (s) \; &\leqslant \; &C
\vep_0^{\frac{1}{15}\,\dos},\label{eq8.61}
\end{eqnarray}
which gives the second part of the proposition.
\end{proof}

\subsection{Transfer to Parameter Space\label{sub8.3}}
\noindent {\bf 8.3.1~} Our goal in this subsection will be to prove
the following result, which expresses a transfer of the dimension
estimate of Subsection~\ref{sub8.2} to parameter space.

\begin{Propo}\label{propo23}
Let $\wtI$ be a regular parameter interval. Let $(\Ps,\Qs,\ns)$ be
an element of $\cR (\wtI)$ such that $\Qs$ is $\wtI$-critical and
\begin{equation}
|\Qs| \; \leqslant \; \fud\;
|\wtI|^{(1+\tau)(1-\eta)^{-1}}\label{eq8.62}
\end{equation}
for all $t\in\wtI$. The, the number of candidates $I\subset \wtI$ of
the next level, such that $\Qs$ is $I$-critical, is at most
$|\wtI|^{-\tau\,\dms}$, where $\dms = \dos + C\eta\tau^{-1}$ can be
made arbitrarily close to $\dos$.
\end{Propo}

Recall that the total number of candidates is $|\wtI|^{-\tau}$.
Proposition~\ref{propo23} is the key estimate that will allow us in
Section~\ref{sec9} to proceed with the selection process for
parameters. The rest of the section is devoted to the proof of
Proposition~\ref{propo23}.

\medskip\noindent
{\bf 8.3.2~} We start with some general observations, that could
have been made much earlier, but are only useful now.

Let $I$ be a parameter interval. Let $(P,Q,n)$, $(P_0,Q_0,n_0)$,
$(P_1,Q_1,n_1)$, be elements of $\cR (I)$ such that $Q\subset Q_u$,
$P_0\subset P_s$, $P_1 \subset P_s$ and $P_0\cap P_1 = \emptyset$.
Associated to the pair $(Q, P_0)$ (resp.~$(Q,P_1)$), we have defined
in Subsection~\ref{sub3.5} a function $C_0 (w,y,x)$ (resp.~$C_1
(w,y,x)$) which is the basis for the definition of $\de (Q, P_0)$,
$\de_{LR}(Q,P_0),\cdots$ (resp.~$\de (Q,P_1)$,
$\de_{LR}(Q,P_1),\cdots$). It follows immediately from the
monotonicity of the function $\te$ that, since $P_0$ and $P_1$ are
disjoint, we have either
\begin{eqnarray}
C_0 (w,y,x_0) \; &> \; &C_1 (w,y,x_1)\;\; \text{ for all $\;w,y,x_0,
x_1,t$},\;\text{ or}\label{eq8.63}\\
C_0 (w,y,x_0) \; &< \; &C_1 (w,y,x_1)\;\; \text{ for all $\;w,y,x_0,
x_1,t$}.\label{eq8.64}
\end{eqnarray}
In the first case, we have, from the definitions in
Subsection~\ref{sub5.4}:
\begin{eqnarray}
\de_R (Q, P_0) \; &< \; &\de (Q, P_1),\label{eq8.65}\\
\de_{LR} (Q, P_0) \; &< \; &\de_L (Q, P_1),\notag
\end{eqnarray}
while in the second the same inequalities hold exchanging $P_0$ and
$P_1$.

We see in particular that if (\ref{eq8.63}) holds and $Q,P_1$ are
$I$-separated, then $Q$, $P_0$ are also $I$-separated.

On the other hand, we will prove (see figure~8).

\begin{center}
\includegraphics{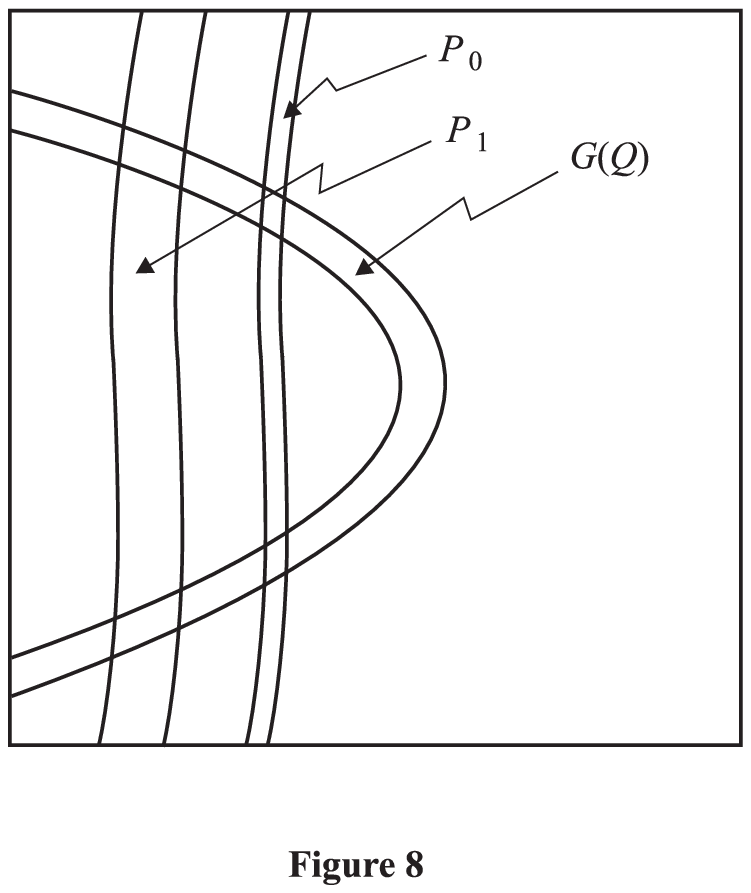}
\end{center}

\begin{Propo}\label{propo24}
Assume that (\ref{eq8.63}) holds and that $|P_1|^{1-\eta}\leqslant
|I|$ for some $t\in I$. Then, if $Q$ and $P_0$ are $I$-transverse,
$Q$ and $P_1$ are also $I$-transverse.
\end{Propo}

\begin{proof}
Let $\wtI\supset I$ and $(\wtP,\wtQ,\wtn)$, $(\wtP_0, \wtQ_0,
\wtn_0)\in \cR (\wtI)$ be such that $\wtQ\supset Q$, $\wtP_0 \supset
P_0$ and $\wtQ \ovfork_{\wtI} \wtP_0$ holds. If $P_1 \subset
\wtP_0$, we immediately conclude that $Q$ and $P_1$ are
$I$-transverse. We assume, therefore, that $P_1 \cap \wtP_0 =
\emptyset$;  replacing $(P_0, Q_0, n_0)$ by $(\wtP_0, \wtQ_0,
\wtn_0)$, and $(P,Q,n)$ by $(\wtP,\wtQ,\wtn)$, we can also assume
that $(P,Q,n)$, $(P_0, Q_0, n_0)\in \cR (\wtI)$ and $Q\ovfork_{\wtI}
P_0$ holds.

Let $(\wtP_1, \wtQ_1, \wtn_1)$, be the element of $\cR (\wtI)$ with
$P_1\subset \wtP_1$ and smallest $\wtP_1$. We will prove that $Q$
and $\wtP_1$ are $\wtI$-transverse. We have, for all $t\in \wtI$,
that
\begin{equation}
\de_{LR} (Q, \wtP_1) \; > \; \de_{LR} (Q, P_0) \; \geqslant \; 2\,
|\wtI|,\label{eq8.66}
\end{equation}
and also for some $t_0\in \wtI$,
\begin{equation}
\de_{R} (Q, \wtP_1) \; > \; \de_{R} (Q, P_0) \; \geqslant \; 2\,
|Q|^{1-\eta}.\label{eq8.67}
\end{equation}
Therefore, assuming that $Q \fork_{\wtI} \wtP_1$ does not hold, we
must have, for all $t\in \wtI$, that
\begin{equation}
\de_{L} (Q, \wtP_1) \; < \; 2\, |\wtP_1|^{1-\eta}.\label{eq8.68}
\end{equation}
We cannot have in this case $\wtP_1 = P_1$, because, for all $t\in
I$,
\begin{equation}
\de_{L} (Q, P_1) \; > \; \de_{LR} (Q, P_0) \; \geqslant \; 2\,
|\wtI|,\label{eq8.69}
\end{equation}
and (\ref{eq8.68}), (\ref{eq8.69}) together would contradict the
hypothesis of the proposition. Therefore, $\wtP_1$ strictly contains
$P_1$ and $\wtI$ strictly contains $I$. But, then, applying the
structure theorem of Subsection~\ref{sub6.7} to the child of
$\wtP_1$ which contains $P_1$, we  obtain that $\wtQ_1$ is
$\wtI$-critical. As $\wtI$ is $\be$-regular, it then follows from
(\ref{eq8.68}), (\ref{eq8.69}) and (\ref{eq5.14}) that $\wtP_1$ is
$\wtI$-transverse. This implies that there exists $(P',Q',n')\in \cR
(\wtI)$ with $Q\cap Q' \ne \emptyset$ such that $Q'\fork_{\wtI}
\wtP_1$ holds. If $Q\subset Q'$, it follows that $Q\fork_{\wtI}
\wtP_1$ holds. When $Q'\subset Q$ we use both $Q \fork_{\wtI} P_0$
and $Q' \fork_{\wtI} \wtP_1$ to conclude, as in the proof of
Proposition~\ref{propo4}, that $Q \fork_{\wtI} \wtP_1$.
\end{proof}

\medskip\noindent
{\bf 8.3.3~} We now switch back to the setting of
Proposition~\ref{propo23}.

Let $(P,Q,n)\in \cR (\wtI)$ with $P\subset P_s$. We say that $P$ is
{\it eventually $\wtI$-separated} from $\Qs$ if there exists an
$\wtI$-decomposition $(P_\al, Q_\al, n_\al)$ of $P$ such that $\Qs$
and $P_\al$ are $\wtI$-separated for every $\al$. We say that $P$ is
{\it eventually $\wtI$-transverse} to $\Qs$ if there exists an
$\wtI$-decomposition $(P_\al, Q_\al, n_\al)$ of $P$ such that $\Qs$
and $P_\al$ are $\wtI$-transverse for every $\al$. We say that $P$
is eventually {\it $\wtI$-$\Qs$-critical} if it is neither
eventually $\wtI$-separated from $\Qs$ nor eventually
$\wtI$-transverse to $\Qs$.

\begin{lema}\label{lema4}
If $P$ is eventually $\wtI$-transverse to $\Qs$ and $2
|P|^{1-\eta}\leqslant |\wtI|$ holds for some $t\in \wtI$, then $\Qs
\fork_{\wtI} P$ holds.
\end{lema}

\begin{proof}
This is an immediate consequence of Proposition~\ref{propo9}.
\end{proof}

\begin{lema}\label{lema5}
If $P$ is eventually $\wtI$-$\Qs$-critical, then $P$ is
$\wtI$-critical.
\end{lema}

\begin{proof}
If $P$ was $\wtI$-transverse, there would exist $(P_\al, Q_\al,
n_\al)\in \cR (\wtI)$ with $\Qs \cap Q_\al \ne \emptyset$ and $Q_\al
\fork_{\wtI} P$. But, then, from Proposition~\ref{propo2} if $Q_\al
\supset \Qs$ or Proposition~\ref{propo9} if $Q_\al \subset \Qs$, we
would deduce that $\Qs \fork_{\wtI} P$ also holds.
\end{proof}

\begin{lema}\label{lema6}
If $P$ is eventually $\wtI$-$\Qs$-critical and $|P| > |\wtI|^\be$
holds for some $t\in \wtI$, then some child of $P$ is also
eventually $\wtI$-$\Qs$-critical.
\end{lema}

\begin{proof}
By Lemma~\ref{lema5} and Corollary~\ref{coro4} in
Subsection~\ref{sub6.2}, $P$ is $\wtI$-decomposable. If all children
of $P$ were eventually $\wtI$-transverse to $\Qs$ (resp.~eventually
$\wtI$-separated from $\Qs$), we would put together the
corresponding $\wtI$-decompositions and obtain that $P$ is
eventually $\wtI$-transverse to $\Qs$ (resp.~eventually
$\wtI$-separated from $\Qs$). Therefore, we shall assume that some
child $P_0$ is eventually $\wtI$-transverse to $\Qs$, and some child
$P_1$ is eventually $\wtI$-separated from $\Qs$.

By contradiction, we assume that none of the children is eventually
$\wtI$-$\Qs$-critical and we will show that $\Qs$ is
$\wtI$-transverse. We construct an $\wtI$-decomposition
$(P_\al,Q_\al,n_\al)$ of $P_s$ such that every $P_\al$ is either
$\wtI$-separated from $\Qs$ or $\wtI$-transverse to $\Qs$ in the
following way.

Actually, it is sufficient to have an $\wtI$-decomposition such that
every $P_\al$ is either eventually $\wtI$-separated from $\Qs$ or
eventually $\wtI$-transverse to $\Qs$. Starting from the trivial
decomposition of $P_s$, we have at step $i$ an $\wtI$-decomposition
$(P^\pii_\al, Q^\pii_\al, n^\pii_\al)$. As long as there is one
$(P^\pii_\al, Q^\pii_\al, n^\pii_\al)$ with $P\subset P_\al^\pii$,
we observe that $P_\al^\pii$ is $\wtI$-critical and therefore
$\wtI$-decomposable and break it into its children to go to step
$i+1$.

After a finite number of steps, each $P_\al^\pii$ is either a child
of $P$ or disjoint from $P$. Comparing $P_\al^\pii$ with $P_0$ and
$P_1$, we conclude from Proposition~\ref{propo24} or from the remark
before this proposition that $P_\al^\pii$ is either eventually
$\wtI$-transverse to $\Qs$ or $\wtI$-separated from $\Qs$. Thus, we
have constructed the required $\wtI$-decomposition and the lemma is
proved.
\end{proof}

\begin{lema}\label{lema7}
If $P_0$, $P_1$ are eventually $\wtI$-$\Qs$-critical and disjoint,
then we have $|P_0|\leqslant C_0 |\wtI|$, $|P_1|\leqslant C_0
|\wtI|$ for all $t\in \wtI$.
\end{lema}

\begin{proof}
From Lemma~\ref{lema6}, we can find $(\whP_0, \whQ_0, \whn_0)$,
$(\whP_1, \whQ_1, \whn_1)$ in $\cR (\wtI)$ with $\whP_0 \subset
P_0$, $\whP_1 \subset P_1$, such that both $\whP_0$, $\whP_1$ are
eventually $\wtI$-$\Qs$-critical and we have
\begin{equation}
|\whP_0|\; \leqslant \; |\wtI|^\be,\qquad |\whP_1|\; \leqslant \;
|\wtI|^\be\;\;\; \text{ for all $\;t\in\wtI$}.\label{eq8.70}
\end{equation}
If we had, for all $t\in \wtI$ and $i=0$ or $1$,
\begin{equation}
\de (\Qs, \whP_i) \; \geqslant \; 2 |\wtI|,\label{eq8.71}
\end{equation}
then, from (\ref{eq8.62}) and (\ref{eq8.70}), one would have that
$\Qs \ovfork_{\wtI} \whP_i$ holds and $\whP_i$ would not be
eventually $\wtI$-$\Qs$-critical. We have, therefore,
\begin{equation}
\de_{LR} (\Qs, \whP_0) \; < \; 2 |\wtI|\;\;\; \text{ for some $\;t_0
\in \wtI\;$ and}\label{eq8.72}
\end{equation}
\begin{equation}
\de_{LR} (\Qs, \whP_1) \; < \; 2 |\wtI|\;\;\; \text{ for some $\;t_1
\in I$}.\label{eq8.73}
\end{equation}
In the same way, we must have
\begin{equation}
\de_{LR} (\Qs, \whP_0) \; \geqslant \; 0\;\;\; \text{ for some
$\;t'_0 \in \wtI$ and also}.\label{eq8.74}
\end{equation}
\begin{equation}
\de'_{LR} (\Qs, \whP_1) \; \geqslant \; 0\;\;\; \text{ for some
$\;t'_1 \in \wtI$}.\label{eq8.75}
\end{equation}

Let $(P',Q',n')$, $(P'_0,Q'_0,n'_0)$, $(P'_1,Q'_1,n'_1)$ be the
elements of $\cR (\wtI)$ such that $\Qs \subset Q'$, $\whP_0\subset
P'_0$, $\whP_1 \subset P'_1$, and
\begin{equation}
|Q'| \, < \, |\wtI|,\qquad |P'_0|\,< \, |\wtI|,\qquad |P'_1| \, < \,
|\wtI|,\quad \forall t\in \wtI\label{eq8.76}
\end{equation}
and which are maximal with these properties.

Then, $P'_0$ satisfies condition $(*s)$ of Subsection~\ref{sub7.6}:
this follows easily from (\ref{eq8.72}). Similarly, $P'_1$ satisfies
condition $(*s)$ and $Q'$ satisfies condition $(*u)$. We are
therefore allowed to apply the estimate of Corollary~\ref{coro11} in
Subsection~\ref{sub7.6} to conclude that
\begin{equation}
\max_{\wtI}\; \de_{LR} (Q',P'_i) - \min_{\wtI}\; \de_{LR} (Q',P'_i)
\; \leqslant \; 2 |\wtI|\label{eq8.77}
\end{equation}
for $i=0,1$, and similarly for $\de$, $\de_L$, $\de_R$. We either
have (\ref{eq8.63}) or (\ref{eq8.64}). Assume, for instance, that
(\ref{eq8.63}) holds. If $P_1\subset P'_1$, we have $|P_1|\leqslant
|\wtI|$ from (\ref{eq8.76}). If $P'_1\subset P_1$, we have
\begin{eqnarray}
\max_{\wtI}\; \de_L (\Qs,P_1) \; &\leqslant \; &\max_{\wtI} \; \de_L
(Q', P'_1)\label{eq8.78}\\
&\leqslant \; &\min_{\wtI} \; \de_L (Q', P'_1) + 2 |\wtI|\notag\\
&\leqslant \; &\min_{\wtI} \; (\de_L (\Qs, P_1) + c |Q'|) + 2 |\wtI|\notag\\
&\leqslant \; &c |\wtI|.\notag
\end{eqnarray}
Similarly, we have $|P_0|\leqslant |\wtI|$ from (\ref{eq8.76}) if $P_0 \subset P'_0$ and
otherwise
\begin{eqnarray}
\min_{\wtI}\; \de_{LR} (\Qs,P_0) \; &\geqslant \; &\min_{\wtI} \;
\de_{LR} (Q', P'_0) - c |\wtI|\label{eq8.79}\\
&\geqslant \; &\max_{\wtI} \; \de_{LR} (Q', P'_0) - (c+2)\,|\wtI|\notag\\
&\geqslant \; &-c |\wtI|.\notag
\end{eqnarray}
Observe that, for all $t\in \wtI$, the inequalities
\begin{equation}
c |\wtI| \; \geqslant \; \de_{L} (\Qs,P_1)\; > \; \de_{LR} (\Qs,P_0)
\; \geqslant \; -c |\wtI|\label{eq8.80}
\end{equation}
are also valid when $P_1\subset P'_1$ or $P_0 \subset P'_0$.
However, there is a constant $c>0$ such that, for any child $P^*_1$
of $P_1$, any child $P^*_0$ of $P_0$, we have, for all $t\in I$,
that
\begin{eqnarray}
\de_{L} (\Qs,P^*_1)\; &\geqslant \; &\de_{L} (\Qs,P_1) + c
|P_1|,\label{eq8.81}\\
\de_{LR} (\Qs,P^*_0)\; &\leqslant \; &\de_{LR} (\Qs,P_0) - c
|P_0|.\label{eq8.82}
\end{eqnarray}
As we also must have, for some children $P^*_1$, $P^*_0$, that
\begin{eqnarray}
\de_{L} (\Qs,P^*_1)\;  &\leqslant \; &c  |\wtI|,\label{eq8.83}\\
\de_{LR}(\Qs,P^*_0)\;  &\geqslant \; &-c |\wtI|,\label{eq8.84}
\end{eqnarray}
because these children are eventually $\wtI$-$\Qs$-critical, we
obtain the conclusion of the lemma.
\end{proof}

\medskip\noindent
{\bf 8.3.4~} Consider the set $\Pi$ of elements $(P,Q,n)\in \cR
(\wtI)$ which are eventually $\wtI$-$\Qs$-critical, satisfy
\begin{equation}
|P| \; \leqslant \; |\wtI|^{1+\tau}\label{eq8.85}
\end{equation}
for all $t\in \wtI$ and are maximal (in $P$) with respect to these
two properties.

\begin{lema}\label{lema8}
We have
\[
\# \; \Pi \; \leqslant \; \frac{1}{C_1}\; |\wtI|^{-\tau \dms},
\]
where $\dms = \dos + C\eta\,\tau^{-1}$ is as in the statement of
Proposition~\ref{propo24}.
\end{lema}

\begin{proof}
From Lemma~\ref{lema7}, there exists a unique element $(P_0, Q_0,
n_0)\in \cR (\wtI)$ with the following properties:
\begin{itemize}
\item[--] $P\subset P_0$ for all $(P,Q,n)\in \Pi$
\item[--] $|P_0| > C_0 |\wtI|$ for some $t\in\wtI$
\item[--] every child $P_1$ of $P_0$ which contains a rectangle $P$
with $(P,Q,n)\in \Pi$ satisfies $|P_1|\leqslant C_o |\wtI|$ for all
$t\in \wtI$.
\end{itemize}
As $P_0$ is eventually $\wtI$-$\Qs$-critical and $|P_0|> C_0 |\wtI|$
for some $t\in \wtI$, $P_0$ must be $\wtI$-critical. From
Proposition~\ref{propo20}, the number of children of $P_0$ is at
most $|\wtI|^{-c\eta}$.

For every $(P,Q,n)\in \Pi$, either $P$ is a child of $P_0$ or the
parent $\wtP$ of $P$ is contained in a child $P_1$ of $P_0$. Observe
that we have by definition of $\Pi$
\begin{equation}
|\wtP| \; > \; |\wtI|^{1+\tau}\;\;\; \text{ for some $\; t\in
\wtI$}.\label{eq8.86}
\end{equation}
As $\wtP$ satisfies condition $(*s)$ of Subsection~\ref{sub7.6}
since $P_0$ is $\wtI$-critical, we also have
\begin{equation}
|\wtP| \; \geqslant \; c^{-1}\, |\wtI|^{1+\tau}\;\;\; \text{ for all
$\; t\in \wtI$}.\label{eq8.87}
\end{equation}
As we also have
\begin{equation}
|P_1| \; \leqslant \; C_0\, |\wtI| \;\;\; \text{ for all $\; t\in
\wtI$},\label{eq8.88}
\end{equation}
the number of possible $\wtP$, given $P_1$, is from
Proposition~\ref{propo22} at most,
\begin{equation}
(C |\wtI|^\tau )^{-\dss}\label{eq8.89}
\end{equation}
with $\dss = \dos + \vep_0^{\frac{1}{5}\,\dos}$. The number of
children of $\wtP$ is at most $|\wtI|^{-c\eta}$. This gives a bound
for the cardinality of $\Pi$ by
\begin{equation}
|\wtI|^{-c\eta} + |\wtI|^{-2c\eta} (C |\wtI|^\tau)^{-\dss},
\label{eq8.90}
\end{equation}
in accordance with the statement of Lemma~\ref{lema8}.
\end{proof}

\medskip\noindent
{\bf 8.3.5~} {\it Proof of Proposition~\ref{propo23}}: By
Lemma~\ref{lema5} and Corollary~\ref{coro4} in
Subsection~\ref{sub6.2}, every $(P,Q,n)\in \cR (\wtI)$ such that $P$
is eventually $\wtI$-$\Qs$-critical and $|P|\geqslant |\wtI|^\be$
(for some $t\in \wtI$) is $\wtI$-decomposable.

Therefore, there exists an $\wtI$-decomposition $(P_\al, Q_\al,
n_\al)$ of $P_s$ such that every $(P_\al, Q_\al, n_\al)$ is either
eventually $\wtI$-separated from $\Qs$ or eventually
$\wtI$-transverse to $\Qs$ or an element of $\Pi$.

Let $I\subset \wtI$ be a candidate interval of the next level,
i.e.~$|I| = |\wtI|^{1+\tau}$, such that $\Qs$ is $I$-critical.

We claim that there exists $(P,Q,n)\in \Pi$ such that $P$ is
eventually $I$-$\Qs$-critical.

Indeed, every $(P_\al, Q_\al, n_\al)$ which is eventually
$\wtI$-transverse to $\Qs$ (resp.~eventually $\wtI$-separated from
$\Qs$) is a fortiori $I$-transverse to $\Qs$ (resp.~$I$-separated
from $\Qs$). If every $(P_\al, Q_\al, n_\al)\in \Pi$ was also either
eventually $I$-transverse to $\Qs$ or eventually $I$-separated from
$\Qs$, we would obtain a decomposition of $P_s$ which expresses that
$\Qs$ is $I$-transverse.

On the other hand, fix $(P,Q,n)\in \Pi$. We show that there are at
most $C_1$ candidates $I\subset \wtI$ such  that $(P,Q,n)$ is
$I$-critical. Together with Lemma~\ref{lema8}, this will imply the
statement of Proposition~\ref{propo23}.

Choose $(P',Q',n')\in \cR (\wtI)$, with $\Qs \subset Q'$, $|Q'| <
|\wtI|$ for all $t\in \wtI$, and maximal with this property as in
the proof of Lemma~\ref{lema7}. Then, $Q'$ satisfies condition
$(*u)$ of Subsection~\ref{sub7.6}; we already know that $P$
satisfies condition $(*s)$. Then, by Corollary~\ref{coro11} of
Subsection~\ref{sub7.6}, we have
\begin{equation}
\left\vert \frac{d}{dt}\; \de_{LR} (Q',P) - 1\right\vert \;
\leqslant \; C\,\vep_0^{\fud}.\label{eq8.91}
\end{equation}
But, if we have for all $t\in I$
\begin{equation}
\de_{LR} (Q',P) \; < \; 0,\label{eq8.92}
\end{equation}
then $Q'$ and $P$ are $I$-separated and a fortiori $\Qs$ and $P$ are
$I$-separated.

On the other hand, if we have
\begin{equation}
\de (Q',P) \; > \; 3\, |I|\;\;\; \text{ for all $\;t\in
I$},\label{eq8.93}
\end{equation}
it is easy to conclude from (\ref{eq8.62}) and
Proposition~\ref{propo9}, that $P$ is eventually $I$-transverse to
$\Qs$. But, by (\ref{eq5.11})--(\ref{eq5.14}), we have
\begin{equation}
\de_{LR} (Q',P) \; < \; \de\,(Q',P) + c\,|I|\label{eq8.94}
\end{equation}
and our claim then follows from (\ref{eq8.91}).
Proposition~\ref{propo23} is proved.\hfill $\square$

\newpage

\setcounter{section}{8}
\setcounter{equation}{0}

\section{Strong Regularity and Parameter Selection\label{sec9}}
\subsection{Partitions of the Critical Locus\label{sub9.1}}
As it was mentioned in Subsection~\ref{sub4.5}, regularity is a
rather qualitative property which is not appropriate for the
quantitative estimates needed for parameter selection. Therefore, we
will introduce later in this section, a stronger quantitative
property, that we call {\it strong regularity}.

This property is built up from several bounds: the number of
bicritical elements of a given size (including, of course, that
there are no ''fat'' bicritical elements), and also sizes of the
''critical locus'' (see below), which should be of approximate
dimension $\dos + \dou - 1$. These last bounds are more elementary
and will be taken care of first.

In all of Section~\ref{sec9}, we fix a parameter interval $\wtI$
which is strongly regular, i.e. satisfies all bounds that will be
stated shortly. We will prove at some point that the starting
interval satisfies this strong regularity. As mentioned before, and
as will be proved in Subsection~\ref{sub9.6}, strong regularity
implies $\be$-regularity for some $\be>1$. We denote by $I$ any
parameter interval contained in $\wtI$ of the next level, i.e. $|I|
= |\wtI|^{1+\tau}$; such a candidate interval is $\ovbe$-regular
(Corollary~\ref{coro9} in Subsection~\ref{sub6.7}), with $\ovbe =
\be (1+\tau)^{-1}$. The aim of this section is to estimate how many
of the candidates $I$ fail to be strongly regular.

We denote by $\cC_+ (\wtI)$ the set of $(P,Q,n)$ in $\cR (\wtI)$
such that $P$ is $\wtI$-critical, $|P|\leqslant |\wtI|^{1+\tau}$ for
all $t\in\wtI$, and $P$ is maximal with this property: the parent
$\wtP$ of $P$ satisfies $|\wtP| > |\wtI|^{1+\tau}$ for some $t\in
\wtI$.

Obviously, if $(P,Q,n)$, $(P',Q',n')$ are distinct elements in
$\cC_+ (\wtI)$, $P$ and $P'$ are disjoint. Moreover, if $(\whP,
\whQ, \whn)$ belongs to $\cR (\wtI)$, $\whP$ is $\wtI$-critical, and
$|\whP|\leqslant |\wtI|^{1+\tau}$ for all $t\in \wtI$, there is a
unique $(P,Q,n)\in \cC_+ (\wtI)$ such that $\whP \subset P$.

Exchanging $P$'s and $Q$'s, we define $\cC_- (\wtI)$ in a similar
way. The sets $\cC_+ (\wtI)$, $\cC_- (\wtI)$ correspond to the
$\wtI$-critical locus at the $|\wtI|^{1+\tau}$ scale. We will also
need to consider this at the $|\wtI|$ scale, as follows. We define
$\whcC_+ (\wtI)$ to be the set of $(P,Q,n)\in \cR (\wtI)$ such that
$P$ contains some $P'$ with $(P',Q',n')\in \cC_+ (\wtI)$,
$|P|\leqslant |\wtI|$ for all $t\in\wtI$, and $P$ is maximal with
this property. We define similarly $\whcC_- (\wtI)$.

We will need in the sequel to consider $\whI$-criticality (for some
parameter interval $\whI \supset \wtI$) for rectangles in $\cR
(\wtI)$ but not in $\cR (\whI)$. The following definition will be
useful.

\Def Let $I_\al\supset I$ be parameter intervals, and let
$(P,Q,n)\in \cR (I)$. We say that $P$ is {\it thin $I_\al$-critical}
if there exists $(P_\al,Q_\al,n_\al)\in\cR (I_\al)$ with $P\subset
P_\al$, $P_\al$ is $I_\al$-critical and
\[
|P_\al|^{1-\eta} \; \leqslant \; 2 |I_\al|
\]
for some $t\in I_\al$.

The main justification for this definition is the following lemma,
which is an immediate consequence of the structure theorem of
Subsection~\ref{sub6.7} and Proposition~\ref{propo11} in
Subsection~\ref{sub6.6}.

\begin{lema}\label{lema9}
Let $I$ be a parameter interval with parent $\wtI$, and let
$(P,Q,n)$ be an element of $\cR (I)$ which is not the restriction of
an element of $\cR (\wtI)$. Let $(P_i, Q_i, n_i)$, for $0\leqslant
i\leqslant k$, be the elements of $\cR (\wtI)$ given by the
structure theorem of Subsection~\ref{sub6.7}. Then $P_i$ is thin
$\wtI$-critical for $0 < i \leqslant k$, and $Q_i$ is thin
$\wtI$-critical for $0 \leqslant i < k$.
\end{lema}

\subsection{Size of the Critical Locus\label{sub9.2}}
We will state several inequalities related to the size of the sets
$\cC_+ (I)$, $\cC_- (I)$, $\whcC_+ (I)$, $\whcC_- (I)$. All these
inequalities are part of the definition of strong regularity: they
have to be satisfied by a strongly regular parameter interval. We
will then see that if these inequalities are satisfied by the parent
interval $\wtI$, most of the candidates $I\subset \wtI$ also satisfy
these inequalities.

Recall that $\dos$, $\dou$ are the transverse Hausdorff dimensions
of $W^s (K)$, $W^u (K)$ at $t=0$. In Proposition~\ref{propo23} of
Subsection~\ref{sub8.3}, we introduced $\dms = \dos + C\eta
\tau^{-1}$. Let also $\dum = \dou + C\eta \tau^{-1}$.

In Proposition~\ref{propo22} of Subsection~\ref{sub8.2}, we have
used $\dss := \dos + \vep_0^{\frac{1}{5}\,\dos} < \dms$. Define
similarly
\[
d^*_u \; = \; \dou + \vep_0^{\frac{1}{5}\,\dou}.
\]
We will control the cardinalities of $\cC_+ (I)$, $\cC_- (I)$,
$\whcC_+ (I)$, $\whcC_- (I)$, through:
\begin{alignat}{1}
\# \cC_+ (I) \; &\leqslant \; C \Bigl(\frac{|I|}{\vep_0}
\Bigr)^{1-\dms-\dum-\tau}\; \vep_0^{-\tau \dos},\tag*{(SR1)$_s$}\\
\# \whcC_+(I)\; &\leqslant \; C \Bigl(\frac{|I|}{\vep_0}
\Bigr)^{1-\dms-\dum-\tau},\tag*{(SR1)$_\whs$}\\
\# \cC_-(I)\; &\leqslant \; C \Bigl(\frac{|I|}{\vep_0}
\Bigr)^{1-\dms-\dum-\tau}\; \vep_0^{-\tau \dou},\tag*{(SR1)$_u$}\\
\# \whcC_-(I)\; &\leqslant \; C \Bigl(\frac{|I|}{\vep_0}
\Bigr)^{1-\dms-\dum-\tau}.\tag*{(SR1)$_\whu$}
\end{alignat}
We will also control a weighted cardinality through
\begin{alignat}{1}
\sum_{\cC_+ (I)} |Q|^{\dsu} \; &\leqslant \; C |Q_s|^{\dsu}
\Bigl(\frac{|I|}{\vep_0}\Bigr)^{1-\dum-\tau},\tag*{(SR2)$_s$}\\
\sum_{\whcC_+(I)} |Q|^{\dsu}\; &\leqslant \; C |Q_s|^{\dsu}
\Bigl(\frac{|I|}{\vep_0}\Bigr)^{1-\dum-\tau},\tag*{(SR2)$_\whs$}\\
\sum_{\cC_-(I)} |P|^{\dss} \; &\leqslant \; C |P_u|^{\dss}
\Bigl(\frac{|I|}{\vep_0}\Bigr)^{1-\dms-\tau},\tag*{(SR2)$_u$}\\
\sum_{\whcC_-(I)} |P|^{\dss} \; &\leqslant \; C |P_u|^{\dss}
\Bigl(\frac{|I|}{\vep_0} \Bigr)^{1-\dms-\tau}.\tag*{(SR2)$_\whu$}
\end{alignat}
The heuristics behind the second set of inequalities is the
following: in the mean, one expects that elements of $\cR (I)$ more
or less satisfy
\begin{equation}
|P|^{\dos} \; \sim \; |Q|^\dou\label{eq9.1}
\end{equation}
and, for $(P,Q,n)\in \whcC_+ (I)$ (for instance), one should have
\begin{equation}
|P| \; \sim \; |I|\label{eq9.2}
\end{equation}
which explains the relation between (SR1)$_\whs$ and (SR2)$_\whs$.

Let us check that the starting interval $\Io$ satisfies these eight
inequalities.

Then $\cC_+ (\Io)$ (resp.~$\cC_- (\Io)$) consists of the elements
$(P,Q,n)\in \cR (\Io)$ with $P\cap P_s\ne \emptyset$ (resp.~$Q\cap
Q_u \ne \emptyset$), $|P|\leqslant \vep_0^{1+\tau}$ for all $t\in
\Io$ (resp.~$|Q|\leqslant \vep_0^{1+\tau}$ for all $t\in \Io$) and
maximal with this property. Then, (SR1)$_s$ follows from
(\ref{eq8.38}) and (SR1)$_{u}$ is similar. Writing $(P,Q,n) = (P_s,
Q_s, n_s) \; * \; (P',Q',n')$, the inequality (SR2)$_s$ becomes
\begin{equation}
\sum_{\cC_+ (I)} |Q'|^{\dou + C\vep_0} \; \leqslant \;
C\label{eq9.3}
\end{equation}
which is a standard property of uniformly hyperbolic horseshoes. The
proof of (SR2)$_u$ is similar. The case of $\whcC_+ (\Io)$, $\whcC_-
(\Io)$ is even simpler. For the induction step, we have:

\begin{Propo}\label{propo25}
If the parent interval $\wtI$ is $\be$-regular and satisfies one of
the eight inequalities (SR) above, then all candidates $I\subset
\wtI$ satisfy the same inequality except perhaps for a proportion
not larger than $C |\wtI|^{\tau^2}$.
\end{Propo}

\medskip\noindent
{\bf Notation.~} Let $(P,Q,n)\in \cC_+ (\wtI)$. We denote by $Cr(P)$
the set of candidates $I\subset \wtI$ such that $P$ contains a thin
$I$-critical rectangle.

\begin{lema}\label{lema10}
For any $(P,Q,n)\in \cC_+ (\wtI)$, we have
\[
\# Cr (P) \; \leqslant \; C |\wtI|^{-\tau \dum}.
\]
\end{lema}

\begin{proof}
Let $(\Ps,\Qs,\ns)\in \cR (\wtI)$ be an element such that $\Ps
\subset P$ and
\begin{equation}
|\Ps| \; \leqslant \; \fud \,
|\wtI|^{(1+\tau)(1-\eta)^{-1}}\label{eq9.4}
\end{equation}
for all $t\in\wtI$. By Proposition~\ref{propo23} in
Subsection~\ref{sub8.3}, there are at most $|\wtI|^{-\tau\dum}$
candidates $I\subset \wtI$ such that $\Ps$ is $I$-critical. Let
$I\in Cr (P)$. By definition, there exists $(P_0,Q_0,n_0)\in \cR
(I)$ with $P_0\subset P$, $P_0$ $I$-critical and
\begin{equation}
|P_0|^{1-\eta} \; \leqslant \; 2 |I| \; \text{ for some } \; t_0\in
I.\label{eq9.5}
\end{equation}
But if $P_0$ is $I$-critical and (\ref{eq9.5}) holds, there must
exist $(P',Q',n') \in \cR (I)$ and $t_1\in I$ such that
\begin{equation}
0 \; \leqslant \; \de_{LR} (Q', P_0) \; \leqslant \; 2
|I|.\label{eq9.6}
\end{equation}
As $P_0$, $\Ps$ are contained in $P$ and $|P|\leqslant |I|$ for all
$t\in\wtI$, we also have
\begin{equation}
|\de_{LR} (Q', \Ps)| \; \leqslant \; C |I|.\label{eq9.7}
\end{equation}
Proceeding as in the Proof of Lemma~\ref{lema7},
Subsection~\ref{sub8.3}, we deduce from Corollary~\ref{coro11} in
Subsection~\ref{sub7.6}, that there exists $I'\subset \wtI$ at
distance $c |I|$ of $I$ such that $\Ps$ is $I'$-critical (except
perhaps when $I$ is very close to the boundary of $\wtI$). This
proves the lemma.
\end{proof}

\medskip\noindent
{\it Proof of Proposition~\ref{propo25}.~} We will deal with
(SR1)$_s$, (SR1)$_\whs$, (SR2)$_s$, (SR2)$_\whs$, the other four
being symmetric.

For each $(P_\al, Q_\al, n_\al)\in \cC_+ (\wtI)$, we consider the
elements $(P_{\al,i}, Q_{\al,i}, n_{\al,i})\in \cR (\wtI)$ which
satisfy $P_{\al,i} \subset P_\al$, $|P_{\al,i}|\leqslant
|\wtI|^{(1+\tau)^2}$ for all $t\in\wtI$ , and which are maximal with
respect to this property.

This gives an $\wtI$-decomposition of $P_\al$: indeed, it is easy to
see that for each $(P'_\al, Q'_\al, n'_\al)\in \cR (\wtI)$ with
$P'_\al \subset P_\al$ and $|P'_\al| > |\wtI|^{(1+\tau)^2}$ for some
$t\in \wtI$, $Q'_\al$ is $\wtI$-transverse and therefore $P'_\al$ is
$\wtI$-decomposable.

Using Propositions~\ref{propo20} and \ref{propo22}, we argue as in
the proof of Lemma~\ref{lema8}, Subsection~\ref{sub8.3}, to see that
for each $P_\al$, the number of $P_{\al,i}$ is not larger than
\begin{equation}
|\wtI|^{-c\eta}\, |\wtI|^{-\tau (1+\tau) \dss}.\label{eq9.8}
\end{equation}
with $\dss = \dos + \vep_0^{\frac{1}{5}\,\dos}$.

Obviously, if $(P,Q,n)\in \cC_+ (I)$, there exists $\al$, $i$ such
that $P\supset P_{\al,i}$ and $I$ must belong to $Cr (P_\al)$. We,
therefore, have
\begin{eqnarray}
\sum_{I\subset \wtI} \# \cC_+ (I) \; &\leqslant \; &C |\wtI|^{-c\eta
- \tau (1+\tau) \dss}\; \sum_{\cC_+ (\wtI)}\, \#Cr
(P_\al)\label{eq9.9}\\
&\leqslant \; &\# \cC_+ (\wtI)\; |\wtI|^{-c\eta-\tau(1+\tau) \dss}\;
\max \#Cr(P_\al).\notag
\end{eqnarray}
We have
\begin{equation}
\dms \; \geqslant \; \dss (1+\tau) + C\eta\tau^{-1}\label{eq9.10}
\end{equation}
and therefore, using also Lemma~\ref{lema10}:
\begin{equation}
\sum_{I\subset \wtI}\; \# \cC_+ (I) \; \leqslant \; C\,\# \cC_+
(\wtI) |\wtI|^{-\tau (\dms+\dum)}.\label{eq9.11}
\end{equation}
The induction step for (SR1)$_s$ follows immediately. In the same
way as in (\ref{eq9.8}), we obtain from Propositions~\ref{propo20}
and \ref{propo22} that
\begin{equation}
\#\cC_+ (\wtI) \; \leqslant \; |\wtI|^{-C\eta} \,
|\wtI|^{-\tau\dss}\; \#\whcC_+ (\wtI).\label{eq9.12}
\end{equation}
On the other hand, for every $(P,Q,n)\in \whcC_+ (I)$, there exists
a $(P_\al, Q_\al, n_\al) \in \cC_+ (\wtI)$ with $P_\al \subset P$
and $I\in Cr (P_\al)$. Therefore, we have
\begin{equation}
\sum_{I\subset \wtI}\; \# \whcC_+ (I) \; \leqslant \; \sum_{\cC_+
(\wtI)}\; \# Cr (P_\al).\label{eq9.13}
\end{equation}
Putting together (\ref{eq9.12}), (\ref{eq9.13}) and
Lemma~\ref{lema10}, we get
\begin{equation}
\sum \; \# \whcC_+ (I) \; \leqslant \; C |\wtI|^{-C\eta -
\tau\dss-\tau\dum}\, \# \whcC_+ (\wtI)\label{eq9.14}
\end{equation}
and we conclude as above.

Let us now consider (SR2)$_s$. We must sum $|Q|^\dsu$ over elements
$(P,Q,n)$ of $\cC_+ (I)$. Fix $(P_\al, Q_\al, n_\al)$ in $\cC_+
(\wtI)$ such that $I\in Cr (P_\al)$. When we consider the partial
sum over those $Q$ contained in $Q_\al$, we get, from
Proposition~\ref{propo22} in Subsection~\ref{sub8.2}.
\begin{equation}
\sum_{Q\subset Q_\al} \; |Q|^\dsu \; \leqslant \; C
|Q_\al|^\dsu.\label{eq9.15}
\end{equation}
As every $Q$ is contained in such a $Q_\al$, we will have:
\begin{equation}
\sum_{\cC_+ (I)} \; |Q|^\dsu \; \leqslant \; C\, \sum_{I\in Cr
(P_\al)} \,|Q_\al|^\dsu.\label{eq9.16}
\end{equation}
Summing then over candidates $I$ and using Lemma~\ref{lema10}, we
get
\begin{eqnarray}
\sum_{I} \, \sum_{\cC_+ (I)}\; |Q|^\dsu \; &\leqslant \; &C\,
\sum_{\cC_+ (\wtI)} \,|Q_\al|^\dsu\, \# Cr (P_\al)\label{eq9.17}\\
&\leqslant \; &C\,|\wtI|^{-\tau \dum}\; \sum_{\cC_+ (\wtI)}
\,|Q_\al|^\dsu,\notag
\end{eqnarray}
which allow us to conclude as above the induction step for
(SR2)$_s$. The argument for (SR2)$_\whs$ is similar and left to the
reader.\hfill $\square$

\subsection{Classes of Bicritical Rectangles\label{sub9.3}}
Now that the size of the critical locus is under control, we must
pay attention to the number of bicritical rectangles, which
represent the returns of the critical locus to itself under the
dynamics.

In order to have an appropriate induction scheme, we need to bound
the number of bicritical rectangles according to all width scales
and also according to the level of criticality (i.e., the distance
to critical locus) of both $P$ and $Q$. As we will see in the next
subsection, the number of bicritical elements experiments a ''phase
transition'' which is crucial for our argument but brings a lot of
complications.

Let $I$ be a candidate interval as above, and let $I_\al$, $I_\om$
be parameter intervals such that $I\subset I_\al \cap I_\om$. Let
also $x$ be a positive number.

\Def
We denote by $Bi_+ (I, I_\al, I_\om; x)$ the set of elements
$(P,Q,n)\in \cR (I)$ such that $P$ is thin $I_\al$-critical, $Q$ is
thin $I_\om$-critical and $|P|\geqslant x$ for some $t\in I$.

Similarly, $Bi_- (I, I_\al, I_\om; x)$ is the set of elements
$(P,Q,n)\in \cR (I)$ such that $P$ is thin $I_\al$-critical, $Q$ is
thin $I_\om$-critical and $|Q|\geqslant x$ for some $t\in I$.

We denote by $Bi_\pm^{new} (I, I_\al, I_\om; x)$ the set of elements
$(P,Q,n)\in Bi_\pm (I, I_\al, I_\om; x)$ that do not belong to $\cR
(\wtI)$, $\wtI$ being as above the parent interval of $I$.

We will estimate the cardinality of all sets $Bi_\pm (I, I_\al,
I_\om; x)$, by induction on the level of the parameter interval $I$.
The easy case is when $I$ is strictly smaller than $I_\al$ and
$I_\om$. In this case, no parameter selection is needed. An element
of $Bi_+ (I, I_\al, I_\om; x)$ is either in $Bi_+^{new} (I, I_\al,
I_\om; x)$ or in $Bi_+ (\wtI, I_\al, I_\om; x)$, where $\wtI$ is the
parent interval of $I$. We will see that the cardinality of
$Bi_+^{new} (I, I_\al, I_\om; x)$ can be estimated from the
induction hypothesis and the structure theorem of
Subsection~\ref{sub6.7}, and this cardinality is much smaller than
the cardinality of $Bi_+ (\wtI, I_\al, I_\om; x)$ which is
controlled by the induction hypothesis.

When $I$ is equal to $I_\al$, but strictly smaller than $I_\om$ (or
in the symmetric case $I = I_\om \subsetneqq I_\al$), a parameter
selection is needed in order to obtain satisfactory estimates. An
element of  $Bi_+ (I, I, I_\om; x)$ is either in $Bi_+^{new} (I, I,
I_\om; x)$, whose cardinality can be again estimated from the
structure theorem and the induction hypothesis, or it belongs to
$Bi_+ (\wtI, \wtI, I_\om; x)$. But in this last case, the vertical
rectangle $P$, which is known to be $\wtI$-critical, is only
$I$-critical for a small fraction of candidates $I$ which is
controlled by Proposition~\ref{propo23} in Subsection~\ref{sub8.3}.
Averaging, like in Subsection~\ref{sub9.2}, allows us to get the
required estimate.

By far the most difficult case occurs when $I=I_\al=I_\om$. When $x$
is large, we have a set of estimates, which is taken care of in the
same way as for $I=I_\al\ne I_\om$. But when $x$ is small, the
parameter selection process is much more subtle and will be
explained in Subsection~\ref{sub9.8}.

\subsection{Number of Bicritical Rectangles\label{sub9.4}}
In this subsection, we will state, and comment, the estimates for
the cardinalities of the sets $Bi_\pm$ introduced above. The rest of
Section~\ref{sec9} will then be devoted to the proof of these
estimates under appropriate parameter selection.

At this point, we have to break the symmetry between past and
future, $P$'s and $Q$'s, stable and unstable direction: the
estimates are indeed not symmetric, except when $\dos = \dou$, i.e.,
in the conservative case of area-preserving diffeomorphisms.

We will assume that $\dos\geqslant \dou$ (and $\dos + \dou > 1$).
The case $\dou \geqslant \dos$ is obviously symmetric.

For $I$, $I_\al$, $I_\om$, $x$ as above we want to have
\begin{alignat}{1}
\#\, Bi_+ (I, I_\al, I_\om; x)\; \leqslant\; CB,\tag*{(SR3)$_s$}
\end{alignat}
with
\begin{eqnarray}
B \; &= \; &\max (B_0, B_1),\label{eq9.18}\\
B_0 \; &= \; &\Bigl( \frac{x}{\vep_0 |P_u|}\Bigr)^{-\rho_0}\; \Bigl(
\frac{|I_\al|}{\vep_0}\Bigr)^{\sig_0+\sig_1}\; \Bigl(
\frac{|I_\om|}{\vep_0}\Bigr)^{\sig_0},\label{eq9.19}\\
B_1 \; &= \; &\Bigl( \frac{x}{\vep_0 |P_u|}\Bigr)^{-\rho_1}\; \Bigl(
\frac{|I_\al|}{\vep_0}\Bigr)^{\sig_1}\; \Bigl( \min\;
\Bigl(\frac{|I_\al|}{\vep_0}, \frac{|I_\om|}{\vep_0}\Bigr)
\Bigr)^{\sig_0}.\label{eq9.20}
\end{eqnarray}
Here $|P_u|$ denotes the supremum over $\Io$ of the width of $P_u$:
the exponents $\rho_0$, $\rho_1$, $\sig_0$, $\sig_1$ will be
specified more precisely later, but anyway they satisfy
\begin{eqnarray}
\rho_0 &= &\dos + o(1),\label{eq9.21}\\
\rho_1 &= &\frac{\dos}{\dos+\dou}\; (2\dos + \dou - 1) + o(1),\label{eq9.22}\\
\sig_0 &= &1 - \dos + o(1),\label{eq9.23}\\
\sig_1 &= &\dos - \dou + o(1).\label{eq9.24}
\end{eqnarray}
The meaning of the $o(1)$ terms in these formulas is that they
become arbitrarily small when $\tau \gg \eta \gg \vep_0$ are small
enough.

For the $Bi_-$ sets, we should have:
\begin{alignat}{1}
\#\, Bi_- (I, I_\al, I_\om; x)\; \leqslant\; CB',\tag*{(SR3)$_u$}
\end{alignat}
with
\begin{eqnarray}
B' \; &= \; &\max (B'_0, B'_1),\label{eq9.25}\\
B'_0 \; &= \; &\Bigl( \frac{x}{\vep_0 |Q_s|}\Bigr)^{-\rho'_0}\;
\Bigl( \frac{|I_\al|}{\vep_0}\Bigr)^{\sig_0+\sig_1}\; \Bigl(
\frac{|I_\om|}{\vep_0}\Bigr)^{\sig_0},\label{eq9.26}\\
B'_1 \; &= \; &\Bigl( \frac{x}{\vep_0 |Q_s|}\Bigr)^{-\rho'_1}\;
\Bigl( \frac{|I_\al|}{\vep_0}\Bigr)^{\sig_1}\; \Bigl( \min\;
\Bigl(\frac{|I_\al|}{\vep_0}, \frac{|I_\om|}{\vep_0}\Bigr)
\Bigr)^{\sig_0},\label{eq9.27}
\end{eqnarray}
\begin{eqnarray}
\rho'_0\; &= \; &\frac{\dou}{\dos} \; \rho_0 \; = \; \dou + o
(1),\label{eq9.28}\\
\rho'_1\; &= \; &\frac{\dou}{\dos} \; \rho_1 \; = \;
\frac{\dou}{\dos+\dou}\,(2\dos + \dou - 1) + o(1).\label{eq9.29}
\end{eqnarray}
Observe that the formulas (\ref{eq9.22}), (\ref{eq9.29}) for
$\rho_1$, $\rho'_1$ are {\it not} symmetric.

\Def A parameter interval is {\it strongly regular} if it satisfies
(SR3)$_s$, (SR3)$_u$ and the eight conditions (SR1), (SR2) of
Subsection~\ref{sub8.2}.

\begin{rem}
The definition will only be complete when we specify precisely the
exponents $\rho_0$, $\rho_1$, $\cdots$.
\end{rem}
We now comment on the inequalities above. First, observe that $B$
does not depend on $I$: this reflects the fact mentioned above that
$Bi_+^{new} (I, I_\al, I_\om; x)$ is small compared with $Bi_+
(\wtI, I_\al, I_\om; x)$.

From the formulas (\ref{eq9.21}), (\ref{eq9.22}), we have
\begin{equation}
\rho_1 \; < \; \rho_0.\label{eq9.30}
\end{equation}
Set
\begin{equation}
x_{cr} \; := \; \vep_0 |P_u|\, \Bigl( \max\;
\Bigl(\frac{|I_\al|}{\vep_0}, \, \frac{|I_\om|}{\vep_0}\Bigr)
\Bigr)^{\tfrac{\sig_0}{\rho_0-\rho_1}}. \label{eq9.31}
\end{equation}
Then, we have $B = B_0$ for $x\leqslant x_{cr}$ and $B=B_1$ for
$x\geqslant x_{cr}$: this is the ''phase transition'' mentioned
earlier. We have
\begin{eqnarray}
\rho_0 - \rho_1 \; &= \; &\frac{\dos (1-\dos)}{\dos+\dou} + o (1),\label{eq9.32}\\
\frac{\sig_0}{\rho_0-\rho_1} \; &= \; &\frac{\dos+\dou}{\dos} + o
(1)\; > \; 1.\label{eq9.33}
\end{eqnarray}
For $x = x_{cr}$, we have
\begin{equation}
B = B_{cr} := \Bigl( \frac{|I_\al|}{\vep_0}\Bigr)^{\sig_0+\sig_1}\,
\Bigl( \frac{|I_\om|}{\vep_0}\Bigr)^{\sig_0}\, \Bigl( \max\;
\Bigl(\frac{|I_\al|}{\vep_0},
\frac{|I_\om|}{\vep_0}\Bigr)\Bigr)^{-\tfrac{\rho_0\sig_0}{\rho_0-\rho_1}}.
\label{eq9.34}
\end{equation}
Assume $I_\al = I_\om$; we then have
\begin{equation}
B_{cr} = \Bigl( \frac{|I_\al|}{\vep_0}\Bigr)^{\sig_1+\sig_0\;
\tfrac{\rho_0-2\rho_1}{\rho_0-\rho_1}}\label{eq9.35}
\end{equation}
Here, the exponent satisfies
\begin{equation}
\sig_1 + \sig_0\; \frac{\rho_0-2\rho_1}{\rho_0-\rho_1} = 2 - 2\dos -
2\dou + o (1) < 0. \label{eq9.36}
\end{equation}
As $|I_\al| \leqslant \vep_0$, we have $B_{cr}\geqslant 1$. As $B$
is a decreasing function of $x$, we have $B < 1$ (in which case
(SR3)$_s$ means that the $Bi_+$ set is empty!) iff $B_1 < 1$ which
corresponds to
\begin{equation}
x > \ovx \; := \; \vep_0 |P_u| \, \Bigl(
\frac{|I_\al|}{\vep_0}\Bigr)^{\tfrac{\sig_0+\sig_1}{\rho_1}}.\label{eq9.37}
\end{equation}
The exponent here satisfies
\begin{equation}
\frac{\sig_0+\sig_1}{\rho_1} \; = \; \frac{1-\dou}{\dos}\;
\frac{\dos+\dou}{2\dos+\dou-1} + o (1)\label{eq9.38}
\end{equation}
We are finally able to justify the assumption (H4) of our Main
Theorem stated in Subsection~\ref{sub1.2}! Indeed, with
$\dos\geqslant \dou$, it means that
\begin{equation}
2 (\dos)^2 + (\dou)^2 + 2\dos\,\dou \; < \; 2\dos + \dou\tag{\bf H4}
\end{equation}
and this is exactly what is needed  to guarantee that
\begin{equation}
\frac{\sig_0+\sig_1}{\rho_1}\; > \; 1.\label{eq9.39}
\end{equation}
We will choose the constant $\be$ (related to the regularity
property) in order to have
\begin{equation}
1 \; < \; \be \; < \; \frac{\sig_0+\sig_1}{\rho_1}\label{eq9.40}
\end{equation}
and also
\begin{equation}
|P_u| \; < \; \vep_0^{\be-1}.\label{eq9.41}
\end{equation}
Then, from (\ref{eq9.37}), we will have  that
\begin{equation}
\ovx \; < \; |I_\al|^\be.\label{eq9.42}
\end{equation}
Summarizing: if (SR3)$_s$ holds, and if $(P,Q,n)\in \cR(I)$ is such
that $P$ and $Q$ are thin $I_\al$-critical, we must have
\begin{equation}
|P| \; < \; \ovx \; < \; |I_\al|^\be\label{eq9.43}
\end{equation}
for all $t\in I$. This is almost what we need for $\be$-regularity.
The full proof is given below.

The discussion for (SR3)$_u$ is similar; the critical threshold is
\begin{equation}
x'_{cr}\; := \; \vep_0 |Q_s|\; \Bigl( \max\;
\Bigl(\frac{|I_\al|}{\vep_0},
\frac{|I_\om|}{\vep_0}\Bigr)\Bigr)^{\tfrac{\sig_0}{\rho'_0-\rho'_1}},\label{eq9.44}
\end{equation}
with
\begin{eqnarray}
\rho'_0 - \rho'_1 &= &\frac{\dou (1-\dos)}{\dos+\dou} + o (1) =
\frac{\dou}{\dos}\; (\rho_0-\rho_1)\label{eq9.45}\\
\frac{\sig_0}{\rho'_0 - \rho'_1} &= &\frac{\sig_0}{\rho_0-\rho_1}\;
\frac{\dos}{\dou}\; = \; \frac{\dos+\dou}{\dou} + o (1)
> 1.\label{eq9.46}
\end{eqnarray}
When $I_\al = I_\om$, we have
\begin{equation}
B'_{cr}\; := \; \Bigl( \frac{|I_\al|}{\vep_0}\Bigr)^{\sig_1+\sig_0\;
\tfrac{\rho'_0-2\rho'_1}{\rho'_0-\rho'_1}} = B_{cr} \geqslant
1.\label{eq9.47}
\end{equation}
Thus, we have $B' < 1$ if
\begin{equation}
x \; > \; \ovx' \; := \; \vep_0 |Q_s| \; \Bigl(
\frac{|I_\al|}{\vep_0}\Bigr)^{\tfrac{\sig_0+\sig_1}{\rho'_1}}\label{eq9.48}
\end{equation}
We have here
\begin{equation}
\frac{\sig_0+\sig_1}{\rho'_1}\; = \; \frac{\sig_0+\sig_1}{\rho_1}\;
\frac{\dos}{\dou} \; > \; \be\label{eq9.49}
\end{equation}
and we choose $\be$ close enough to one to have
\begin{equation}
|Q_s| \; < \; \vep_0^{\be-1}.\label{eq9.50}
\end{equation}
Then, we are again able to conclude that, if (SR3)$_u$ holds, any
$(P,Q,n)\in \cR(I)$ such that both $P$ and $Q$ are thin
$I_\al$-critical must satisfy, for all $t\in I$
\begin{equation}
|Q| \; < \; \ovx' \; < \; |I_\al|^\be.\label{eq9.51}
\end{equation}

\begin{Propo}\label{propo26}
If a candidate interval satisfies (SR3)$_s$ and (SR3)$_u$, then it
is $\be$-regular. In particular, strong regularity implies
regularity.
\end{Propo}

\begin{proof}
We argue by induction on the level of the parameter interval. For
the starting interval $\Io$, we use (for the first time!) the
assumption (H1) of Subsection~\ref{sub1.2} that the periodic points
$p_s$, $p_u$ do not belong to the same periodic orbit. Then, if
$(P,Q,n)\in \cR (\Io)$ is such that $P\subset P_s$, $Q\subset Q_u$,
we must have, for all $t\in\Io$
\begin{equation}
|P| \; < \; \vep_0^\be,\qquad |Q| \; < \; \vep_0^\be\label{eq9.52}
\end{equation}
if we take $\be$ close enough to one. This proves that $\Io$ is
$\be$-regular (independently of (SR3)$_s$,  (SR3)$_u)$. Assume that
$I\ne \Io$ satisfies (SR3)$_s$, (SR3)$_u$ and that $(P,Q,n)\in
\cR(I)$ is $I$-bicritical. Assume also, for instance, that
\begin{equation}
\max_I \; |Q| \; \leqslant \; \max_I \; |P|\label{eq9.53}
\end{equation}
and, by contradiction that
\begin{equation}
\max_I \; |P| \; \geqslant \; |I|^\be.\label{eq9.54}
\end{equation}
From the proof of Corollary~\ref{coro9} in Subsection~\ref{sub6.7},
we know that $(P,Q,n)\in \cR (\wtI)$ ($\wtI$ being the parent of
$I$). As $(P,Q,n)$ is $\wtI$-bicritical, we must have
\begin{equation}
\max_I \; |P| \; < \; |\wtI|^\be.\label{eq9.55}
\end{equation}
Therefore, $P$ would be thin $I$-critical; similarly $Q$ would be
thin $I$-critical. But (SR3)$_s$ , (see (\ref{eq9.43}), says that such
a $(P,Q,n)$ satisfying (\ref{eq9.54}) does not exist.
\end{proof}

\begin{rem}
While there are only eight inequalities (SR1), (SR2) for each
parameter interval $I$, the inequalities (SR3)$_s$, (SR3)$_u$ form a
family parametrized not only by $I$, but also by the parameter
intervals $I_\al\supset I$ and $I_\om \supset I$ and the real number
$x>0$. Because each inequality, at least when $I=I_\al$ or
$I=I_\om$, is only obtained after parameter selection, we will
discretize the continuous variable $x$ by considering only the
values $x = 2^{-N}$, $N\geqslant 0$. There is still an infinite
number of inequalities, but we will see that they are trivially
satisfied if $N$ is large enough.
\end{rem}

\subsection{The Starting Interval\label{sub9.5}}
Our main purpose in this subsection is to prove the following fact.

\begin{Propo}\label{propo27}
The starting interval is strongly regular.
\end{Propo}

\begin{proof}
We have already checked in Subsection~\ref{sub8.2} that $\Io$
satisfies the eight inequalities (SR1), (SR2). We have therefore to
prove that (SR3)$_s$, (SR3)$_u$ are also satisfied. We clearly have
$I_\al = I_\om = \Io$. Then $Bi_+ (\Io, \Io, \Io; x)$ is the set of
$(P,Q,n)\in \cR (\Io)$ such that $P\subset P_s$, $Q\subset Q_u$ and
$|P|\geqslant x$ for some $t\in\Io$. Writing $(P,Q,n)$ as a simple
composition
\begin{equation}
(P,Q,n) \; = \; (P_s,Q_s,n_s)\; * \; (P',Q',n')\;*\; (P_u,Q_u,n_u))
\label{eq9.56}
\end{equation}
(here we use again assumption (H1)), we, then, use the standard
estimate (\ref{eq8.38}) to obtain:
\begin{equation}
\# Bi_+ (\Io, \Io, \Io; x) \; \leqslant \; C \Bigl( \frac{x}{\vep_0
|P_u|}\Bigr)^{-(\dos+C\vep_0)},\label{eq9.57}
\end{equation}
and similarly
\begin{equation}
\# Bi_- (\Io, \Io, \Io; x) \; \leqslant \; C \Bigl( \frac{x}{\vep_0
|Q_s|}\Bigr)^{-(\dou+C\vep_0)}.\label{eq9.58}
\end{equation}
Therefore, we obtain (SR3)$_s$ and (SR3)$_u$ if we have:
\begin{eqnarray}
\rho_0 \; &> \; &\dos + C \vep_0,\label{eq9.59}\\
\rho'_0 \; &> \; &\dou + C \vep_0,\label{eq9.60}
\end{eqnarray}
which is compatible with (\ref{eq9.21}), (\ref{eq9.28}).
\end{proof}
The other part of (SR3)$_s$, (SR3)$_u$ which can be taken care of
right now is the case where $x$ is extremely small.

In this case, we will just forget about the criticality conditions
for $P$ and $Q$ and bound the cardinality of $Bi_+ (I, I_\al, I_\om;
x)$ by the cardinality of the set of $(P,Q,n)\in\cR (I)$ for which
$|P|\geqslant x$ for at least some $t\in I$.

This cardinality was estimated in Proposition~\ref{propo22} of
Subsection~\ref{sub8.2}. Actually, the estimation was given for
fixed parameter but it is easy to check that the same proof gives
the same estimate as in Proposition~\ref{propo22} for the set we are
considering. One obtains
\begin{equation}
\# Bi_+ (I, \Io, \Io; x) \; \leqslant \; \Bigl(
\frac{x}{\vep_0}\Bigr)^{-\dss}\label{eq9.61}
\end{equation}
with $\dss = \dos + \vep_0^{\frac{1}{5}\;\dos}$ as above.

We want to have (for all $I_\al$, $I_\om \supset I$)
\begin{equation}
\Bigl( \frac{x}{\vep_0}\Bigr)^{-\dss}\; \leqslant \;
B_0,\label{eq9.62}
\end{equation}
which is equivalent to
\begin{equation}
\Bigl( \frac{x}{\vep_0}\Bigr)^{\rho_0-\dss}\; \leqslant \;
|P_u|^{\rho_0}\, \Bigl(
\frac{|I|}{\vep_0}\Bigr)^{2\sig_0+\sig_1}.\label{eq9.63}
\end{equation}
This will be satisfied if
\begin{equation}
x \; \leqslant \; x_{\min} \; := \;
|\wtI|^{C(\rho_0-\dss)^{-1}}.\label{eq9.64}
\end{equation}
So, we need to have
\begin{equation}
\rho_0 \; > \; \dss\label{eq9.65}
\end{equation}
which is compatible with (\ref{eq9.59}), (\ref{eq9.21}). The
exponent $C (\rho_0 - \dss)^{-1}$ will be very large. We have proved

\begin{Propo}\label{propo28}
The estimate (SR3)$_s$ is satisfied for all candidates $I$, all
$I_\al$, $I_\om \supset I$, as soon as $x\leqslant x_{\min}$. A
similar statement holds for (SR3)$_u$, with a threshold
\begin{equation}
x'_{\min} \; := \; |\wtI|^{C(\rho'_0-\dsu)^{-1}}\label{eq9.66}
\end{equation}
\end{Propo}

\subsection{New Bicritical Rectangles\label{sub9.6}}
We consider in this section the set $Bi_+^{new} (I, I_\al, I_\om;
x)$ of bicritical rectangles which were not defined over the parent
$\wtI$ of $I$: cf.~Subsection~\ref{sub9.3}. We assume that
$\wtI\subset I_\al \cap I_\om$. We apply to each element $(P,Q,n)$
in this set the structure theorem (Theorem~\ref{theo1}) of
Subsection~\ref{sub6.7}. We obtain an integer $k>0$, elements $(P_0,
Q_0, n_0), \cdots, (P_k, Q_k, n_k)$ of $\cR (\wtI)$ such that
\begin{equation}
(P,Q,n) \; \in \; (P_0, Q_0, n_0) \;\square \cdots \square\; (P_k,
Q_k, n_k).\label{eq9.67}
\end{equation}
Moreover, $P_i$ is $\wtI$-critical for $0 < i\leqslant k$ and $Q_i$
is $\wtI$-critical for $0\leqslant i < k$. On the other hand, $P_0$
is $I_\al$-critical because $P$ is $I_\al$-critical, and $Q_k$ is
$I_\om$-critical because $Q$ is $I_\om$-critical. Denote by $x_i =
2^{-n_i}$ the largest integral negative power of 2 such that
\begin{equation}
|P_i|\; \geqslant \; x_i\qquad \text{ for some
$\;t\in\wtI$}.\label{eq9.68}
\end{equation}

\begin{lema}\label{lema11}
We have
\begin{eqnarray*}
(P_0, Q_0, n_0) \; &\in \; &Bi_+ (\wtI, I_\al, \wtI; x_0),\\
(P_k, Q_k, n_k) \; &\in \; &Bi_+ (\wtI, \wtI, I_\om; x_k),\\
(P_i, Q_i, n_i) \; &\in \; &Bi_+ (\wtI, \wtI, \wtI; x_i)
\end{eqnarray*}
for $0 < i < k$.
\end{lema}

\begin{rem}
This lemma is the reason why we need to consider different levels of
criticality for $P$ and $Q$.
\end{rem}

\begin{proof}
For $0<i\leqslant k$, $P_i$ is $\wtI$-critical and it is even thin
$\wtI$-critical by Corollary~\ref{coro7} of Subsection~\ref{sub6.7}.
Similarly, for $Q_j$ when $0\leqslant j < k$. Finally, $P_0$ is thin
$P_\al$-critical because $P$ is thin $I_\al$-critical and $P_0$ is
the thinnest rectangle containing $P$ which is defined over $\wtI$
(which is contained in $I_\al$), so the thinnest rectangle
containing $P$ and defined over $I_\al$ also contains $P_0$.
Similarly, $Q_k$ is thin $I_\om$-critical.
\end{proof}

The widths $|P_i|$ are related to the width of $P$ by
Corollary~\ref{coro6} in Subsection~\ref{sub6.7} which gives
\begin{equation}
x\; \leqslant \; C^k \, |I|^{-\tfrac{k}{2}} \prod^k_0\,
x_i.\label{eq9.69}
\end{equation}
Let us write
\begin{eqnarray*}
\# (x_0) \; &:= \; &\# Bi_+ (\wtI, I_\al, \wtI; x_0),\\
\# (x_k) \; &:= \; &\# Bi_+ (\wtI, \wtI, I_\om; x_k),\\
\# (x_i) \; &:= \; &\# Bi_+ (\wtI, \wtI, \wtI; x_i)\; \text{for $0 <
i < k$.}
\end{eqnarray*}
Then, as each parabolic composition produces two elements, we have
\begin{equation}
\# Bi_+^{new} (I, I_\al, I_\om; x) \; \leqslant \; \sum_{k>0} 2^k
\sum_{x_0,\cdots,x_k}\; \prod^k_0 \# (x_i).\label{eq9.70}
\end{equation}
The term $\# (x_i)$ is estimated by the induction hypothesis
(SR3)$_s$ for $\wtI$. In view of the threshold (\ref{eq9.31}), we
divide $Bi_+^{new}$ into two parts. In the first, we put the
elements for which every $x_i$ is above the threshold $x_{i,cr}$
given by (\ref{eq9.31}). In the second, at least one of the $x_i$ is
below $x_{i,cr}$.

Let us consider the first part. Then all $\# (x_i)$ are estimated by
$B_1$ and we have
\begin{equation}
\prod^k_0 \#(x_i) \; \leqslant \; C^{k+1} \Bigl( \prod^k_0
x_i\Bigr)^{-\rho_1} (\vep_0 |P_u|)^{(k+1)\rho_1} \Bigl(
\frac{|I_\al|}{\vep_0} \Bigr)^{\sig_1} \Bigl( \frac{|\wtI|}{\vep_0}
\Bigr)^{(k+1) \sig_0+k \sig_1}.\label{eq9.71}
\end{equation}
In view of (\ref{eq9.69}), the right-hand side is bounded by
\begin{equation}
\Bigl( \frac{x}{\vep_0 |P_u|}\Bigr)^{-\rho_1} \Bigl(
\frac{|I_\al|}{\vep_0}\Bigr)^{\sig_1} \Bigl( \frac{|\wtI|}{\vep_0}
\Bigr)^{\sig_0} Z^k,\label{eq9.72}
\end{equation}
with
\begin{equation}
Z\; := \; \Bigl( C\vep_0 |P_u|\;|I|^{-\tfrac{1}{2}}\Bigr)^{\rho_1}
\Bigl( \frac{|\wtI|}{\vep_0}\Bigr)^{\sig_0+\sig_1}.\label{eq9.73}
\end{equation}
For all $0 \leqslant i\leqslant k$, from (\ref{eq9.31}) we have
\begin{equation}
x_{i,cr}\; \geqslant \; \vep_0 |P_u| \; \Bigl(
\frac{|\wtI|}{\vep_0}\Bigr)^{\tfrac{\sig_0}{\rho_0-\rho_1}}\label{eq9.74}
\end{equation}
(with equality when $i\ne 0,k$). Therefore, the number of
$(k+1)$-tuples $(n_0,\cdots,n_k)$ such that $2^{-n_i}\geqslant
x_{i,cr}$ for $0\leqslant i \leqslant k$ is at most $(C \log
|\wtI|^{-1})^{k+1} \leqslant (C\log |\wtI|^{-1})^{2k}$. We conclude
that the cardinality of the first part of $Bi_+^{new}$ is bounded by
\begin{equation}
\Bigl( \frac{x}{\vep_0 |P_u|} \Bigr)^{-\rho_1} \Bigl(
\frac{|I_\al|}{\vep_0}\Bigr)^{\sig_1} \Bigl(
\frac{|\wtI|}{\vep_0}\Bigr)^{\sig_0} \sum_{k>0} Z_1^k,\label{eq9.75}
\end{equation}
with
\begin{equation}
Z_1 \; := \; 2 (C \log |\wtI|^{-1})^2 Z.\label{eq9.76}
\end{equation}
In view of (\ref{eq9.73}), we have
\begin{equation}
Z_1 \; < \; \fud \, \Bigl( \frac{|\wtI|}{\vep_0}\Bigr)^{\sig_0 +
\sig_1 - \rho_1} |\wtI|^{\tfrac{1}{3}\,\rho_1},\label{eq9.77}
\end{equation}
with $\sig_0 + \sig_1 - \rho_1 > 0$ (cf.~(\ref{eq9.39})).

As $Z_1 < \fud$, the bound (\ref{eq9.75}) is smaller than
\begin{equation}
\Bigl( \frac{x}{\vep_0 |P_u|}\Bigr)^{-\rho_1} \Bigl(
\frac{|\wtI|}{\vep_0}\Bigr)^{2\sig_0 + \sig_1-\rho_1} \Bigl(
\frac{|I_\al|}{\vep_0}\Bigr)^{\sig_1} |\wtI|^{\tfrac{1}{3}\,\rho_1}.
\label{eq9.78}
\end{equation}
Thus, (\ref{eq9.78}) is a bound for the cardinality of the first
part of $Bi_+^{new} (I, I_\al, I_\om; x)$. Let us turn to the second
part. Let $J$ be the non-empty subset of indices $i\in \{
0,\cdots,k\}$ for which $x_i < x_{i,cr}$, and write $j = \# J$.

We first estimate the product $\prod_i \# (x_i)$. As $\rho_0 >
\rho_1$, we have from (\ref{eq9.69})
\begin{equation}
\prod_J x_i^{-\rho_0} \prod_{J^c} x_i^{-\rho_1} \; \leqslant \;
\Bigl (C^{-k} |I|^{\tfrac{k}{2}}\, x\Bigr)^{-\rho_0}.\label{eq9.79}
\end{equation}
As we also have $\wtI\subset I_\al$, $\wtI\subset I_\om$, we obtain
\begin{equation}
\prod_i \# (x_i) \; \leqslant \; C^{k+1} \Bigl( \frac{x}{\vep_0
|P_u|}\Bigr)^{-\rho_0} \Bigl(
\frac{|I_\al|}{\vep_0}\Bigr)^{\sig_0+\sig_1} \Bigl(
\frac{|I_\om|}{\vep_0}\Bigr)^{\sig_0} Y_0^{j-1} Y_1^k,\label{eq9.80}
\end{equation}
with
\begin{eqnarray}
Y_0 \; &= \; &(\vep_0 |P_u|)^{\rho_0-\rho_1} \Bigl(
\frac{|\wtI|}{\vep_0}\Bigr)^{\sig_0},\label{eq9.81}\\
Y_1 \; &= \; &\Bigl(C^{-1} |I|^{\tfrac{1}{2}}\Bigr)^{-\rho_0}
(\vep_0 |P_u|)^{\rho_1} \Bigl(
\frac{|\wtI|}{\vep_0}\Bigr)^{\sig_0+\sig_1}.\label{eq9.82}
\end{eqnarray}
We have
\begin{equation}
Y_0 \; < \; 1,\label{eq9.83}
\end{equation}
and we can rewrite $Y_1$ as
\begin{equation}
Y_1 \; = \;C^{\rho_0} \Bigl(
\frac{\wtI}{\vep_0}\Bigr)^{\sig_0+\sig_1-\fudt\,\rho_0 (1+\tau)} \;
\vep_0^{\rho_1 - \fudt\, \rho_0 (1+\tau)} \,
|P_u|^{\rho_1}.\label{eq9.84}
\end{equation}
The exponents satisfy
\begin{eqnarray}
\sig_0+\sig_1 - \fud\, \rho_0 (1+\tau) = 1 - \dou - \fud\, \dos + o
(1),\label{eq9.85}\\
\rho_1 - \fud\, \rho_0 (1+\tau) = \frac{\dos (3\dos + \dou - 2)}{2
(\dos + \dou)} + o (1).\label{eq9.86}
\end{eqnarray}
From (H4), we have $1-\dou - \fud\, \dos > 0$; from $\dos + \dou
\geqslant 1$, $\dos \geqslant \dou$, we have $3\dos + \dou - 2
\geqslant 0$. We can, therefore, find $\sig_2 > 0$ such that
\begin{equation}
Y_1 \; \leqslant \; |\wtI|^{2\sig_2}.\label{eq9.87}
\end{equation}
Consider now the number of $(k+1)$-tuples $(n_0,\cdots,n_k)$ for
this second part. From (\ref{eq9.69}), we have
\begin{equation}
\prod^k_0\, x_i \; \geqslant \; C^{-k}
x\,|I|^{\tfrac{k}{2}}\label{eq9.88}
\end{equation}
and we have only to consider the case
\begin{equation}
x\; \geqslant \; x_{\min} := |\wtI|^{C
(\rho_0-\dss)^{-1}}.\label{eq9.89}
\end{equation}
From Corollary~\ref{coro7} in Subsection~\ref{sub6.7}, we also have,
for $0<i\leqslant k$,
\begin{equation}
x_i\; < \;|\wtI|.\label{eq9.90}
\end{equation}
We conclude that the number of $(n_0,\cdots,n_k)$ is smaller than
\begin{equation}
\Bigl(C (\rho_0 - \dss)^{-1} \log
|\wtI|^{-1}\Bigr)^{k+1}.\label{eq9.91}
\end{equation}
Using $k+1\leqslant 2k$ for $k>0$, we carry this to (\ref{eq9.70})
to obtain a bound for the second part of $Bi_+^{new}$, which is
equal to
\begin{equation}
\Bigl( \frac{x}{\vep_0 |P_u|}\Bigr)^{-\rho_0} \Bigl(
\frac{|I_\al|}{\vep_0}\Bigr)^{\sig_0+\sig_1} \Bigl(
\frac{|I_\om|}{\vep_0} \Bigr)^{\sig_0} \sum_{k>0}\,
Y^k\label{eq9.92}
\end{equation}
with
\begin{equation}
Y \; = \; 2C (\rho_0 - \dss)^{-2} (\log |\wtI|^{-1})^2\,
Y_1.\label{eq9.93}
\end{equation}
Here, $(\rho_0 - \dss)^{-2}$ is large but independent of $\vep_0$,
which is always assumed to be as small as necessary. In view of
(\ref{eq9.87}), we have
\begin{equation}
\sum_{k>0} \, Y^k \; < \; |\wtI|^{\sig_2}.\label{eq9.94}
\end{equation}
We have now estimated the cardinalities of the two parts of
$Bi_+^{new}$. Taking $\sig_2$ smaller if necessary, we have
\begin{equation}
0 \; < \; \sig_2 \; \leqslant \; \frac{1}{3}\, \rho_1.\label{eq9.95}
\end{equation}
In (\ref{eq9.78}), as $\sig_0+\sig_1 - \rho_1 > 0$, we have,
\begin{equation}
\Bigl( \frac{|\wtI|}{\vep_0}\Bigr)^{2\sig_0+\sig_1-\rho_1}\;
\leqslant \; \Bigl( \min\; \frac{|I_\al|}{\vep_0},
\frac{|I_\om|}{\vep_0} \Bigr)^{\sig_0}.\label{eq9.96}
\end{equation}
Thus, we have proved

\begin{Propo}\label{propo29}
Assume that (SR3)$_s$ is satisfied by the parent interval $\wtI$.
Then, for all $x\geqslant x_{\min}$ and all candidates $I\subset
\wtI$, all $I_\al$, $I_\om$ containing $\wtI$, we have
\[
\# Bi_+^{new} (I, I_\al, I_\om; x) \; \leqslant \; 2B
|\wtI|^{\sig_2},
\]
with $B$ given by (\ref{eq9.18}).
\end{Propo}
The presence of the term $|\wtI|^{\sig_2}$ puts $Bi_+^{new}$ well
under control.

\begin{Coro}\label{coro13}
Assume that $\wtI\subset I_\al \cap I_\om$ and that (SR3)$_s$ holds
for $\whI = I_\al \cap I_\om$. Then, (SR3)$_s$ holds for all
candidates $I\subset \wtI$.
\end{Coro}

\begin{proof}
An element $(P,Q,n)\in Bi_+ (I, I_\al, I_\om; x)$ either belongs to
$Bi_+ (\whI, I_\al, I_\om; x)$ or to $Bi_+^{new} (\Is, I_\al, I_\om;
x)$ for some interval $\Is$ with $\whI \subsetneqq \Is \subset I$.
As the series $\sum_{\Is} |\Is|^{\sig_2}$ is bounded (actually very
small), the Corollary follows.
\end{proof}

The conclusions of Proposition~\ref{propo29} and
Corollary~\ref{coro13} also hold for the $Bi_-$ sets; let us review
briefly the proof of Proposition~\ref{propo29}.

In the formulas for $B'_0$, $B'_1$ the exponents $\sig_0$, $\sig_1$
are the same as for $B_0$, $B_1$ but we have to replace the
exponents $\rho_0$, $\rho_1$ by $\rho'_0 = \frac{\dou}{\dos}
\rho_0$, $\rho'_1 = \frac{\dou}{\dos} \rho_1$ (with $\dos\geqslant
\dou$).

Therefore, we have $\sig_0 + \sig_1 - \rho'_1 > 0$
(cf.~(\ref{eq9.77}), (\ref{eq9.78})), $\sig_0 + \sig_1 - \fud\,
\rho'_0 > 0$ (cf.~(\ref{eq9.85})), $\rho'_1 - \fud\, \rho'_0
\geqslant 0$ (cf.~(\ref{eq9.86}). Therefore,
Proposition~\ref{propo29} holds for the set $Bi_-^{new} (I, \Isal,
\Isom; x)$ with an appropriate choice of $\sig_2$ and
Corollary~\ref{coro13} holds for (SR3)$_u$.

The induction step in the case where $I$ is distinct from $I_\al$
and $I_\om$ is complete.

\subsection{The Case $I = I_\al \ne I_\om$\label{sub9.7}}
We consider here the set $Bi_+ (I, I, I_\om; x)$ with $\Isom \supset
\wtI$. We will estimate the size of this set for most candidates
$I\subset \wtI$.

Let $(P,Q,n)$ be an element of $Bi_+ (I,I,\Isom; x)$. Either it
belongs to $Bi_+^{new} (I,I,\Isom; x)$ or to $Bi_+ (\wtI, \wtI,
\Isom; x)$. In the second case, as $P$ in thin $I$-critical, there
exists $(\Ps, \Qs, \ns)\in \cC_+ (\wtI)$ such that $P\subset \Ps$
and, moreover, we have $I\in Cr (\Ps)$ by definition of this set
(Subsection~\ref{sub9.1}). We will denote by $Bi_+ (\Ps)$ the set of
$(P,Q,n)\in Bi_+ (\wtI,\wtI,\Isom; x)$ such that $P\subset \Ps$. We,
thus, have
\begin{equation}
\# Bi_+ (I,I,\Isom; x) \; \leqslant \; \# Bi_+^{new} (I,I,\Isom; x) +
\sum_{\substack{(\Ps,\Qs,\ns)\in\cC_+(\wtI)\\ I\in Cr (\Ps)}} \#
Bi_+ (\Ps)\label{eq9.97}
\end{equation}
As the sets $Bi_+ (\Ps)$ are disjoint, we have
\begin{equation}
\sum_{\cC_+ (\wtI)}\# Bi_+ (\Ps) \; \leqslant \; \# Bi_+
(\wtI,\wtI,\Isom; x).\label{eq9.98}
\end{equation}
We also have
\begin{equation}
Bi_+^{new} (I,I,\Isom; x) \; \subset \; Bi_+^{new} (I,\wtI,\Isom;
x).\label{eq9.99}
\end{equation}
We assume that (SR3)$_s$ holds for parameter intervals containing
$\wtI$. Write $\wtB = max (\wtB_0, \wtB_1)$ for the values of
(\ref{eq9.18})--(\ref{eq9.20}) with $I_\al = \wtI$. Thus, we have
\begin{equation}
\# Bi_+ (\wtI,\wtI,\Isom; x) \; \leqslant \; C\wtB,\label{eq9.100}
\end{equation}
and, from Proposition~\ref{propo29},
\begin{equation}
\# Bi_+^{new} (I,\wtI,\Isom; x) \; \leqslant \; 2\wtB
|\wtI|^{\sig_2}.\label{eq9.101}
\end{equation}
We now sum over candidates $I\subset \wtI$ the estimate
(\ref{eq9.97}). We obtain
\begin{equation}
\sum_{I\subset \wtI} \# Bi_+ (I,I,\Isom; x) \; \leqslant \; 2\wtB
|\wtI|^{\sig_2} |\wtI|^{-\tau} + \sum_{\cC_+ (\wtI)} \# Bi_+ (\Ps)
\# Cr (\Ps),\label{eq9.102}
\end{equation}
and then, using Lemma~\ref{lema10} and (\ref{eq9.100}), the same sum
is bounded by
\begin{equation}
2\wtB |\wtI|^{\sig_2-\tau} + C'\wtB |\wtI|^{-\tau\dum} \; \leqslant
\; C'' \wtB |\wtI|^{-\tau\dum}.\label{eq9.103}
\end{equation}
Write $B = \max (B_0,B_1)$ for the values of
(\ref{eq9.18})--(\ref{eq9.20}) with $I_\al = I$. The number of
candidates $I$ such that $\# Bi_+ (I, I, \Isom; x) > B$ is at most
\begin{equation}
C'' \wtB B^{-1} |\wtI|^{-\tau\dum}.\label{eq9.104}
\end{equation}
From formulas (\ref{eq9.18})--(\ref{eq9.20}), we see that
\begin{equation}
\wtB B^{-1}\; = \; |\wtI|^{-\tau (\sig_0+\sig_1)}.\label{eq9.105}
\end{equation}
The exponents $\sig_0$, $\sig_1$ will be chosen in order to have
\begin{equation}
\sig_0 + \sig_1 + \dum \; < \; 1-2\tau,\label{eq9.106}
\end{equation}
which is compatible with (\ref{eq9.23}), (\ref{eq9.24}). Then, we
obtain that (SR3)$_s$ holds except for a proportion of candidates no
greater than $C |\wtI|^{2\tau^2}$. This assertion has been proved
for a fixed value of $x$ and a fixed parameter interval $\Isom$. But
(in view of the constant $C$ in (SR3)$_s$), it is sufficient to
prove (SR3)$_s$ when $x = 2^{-n}$, $n$ a nonnegative integer, and
$x\geqslant x_{\min}$; and the number of intervals $\Isom \supset
\wtI$ is at most
\begin{equation}
C \tau^{-1} \, \log \, \Bigl( \frac{\log\,|\wtI|}{\log\,\vep_0}
\Bigr).\label{eq9.107}
\end{equation}
Thus, the total number of cases that we have to consider is much
smaller than $|\wtI|^{-\tau^2}$. We obtain

\begin{Propo}\label{propo30}
Assume that (SR3)$_s$ holds for parameter intervals containing
$\wtI$. Then, (SR3)$_s$ holds for $Bi_+ (I, I, \Isom; x)$, for all
$\Isom \supset \wtI$ and all $x$, except perhaps for a proposition
of candidates no greater than $|\wtI|^{-\tau^2}$.
\end{Propo}
Proposition~\ref{propo30} is also valid exchanging $I_\al$ and
$\Isom$. There are also two similar statements involving (SR3)$_u$
and $Bi_-$. We review briefly the slight difference in the proof.

For $Bi_- (I, I, \Isom; x)$ (with $\Isom\supset \wtI$), we proceed
exactly as for $Bi_+ (I, I, \Isom; x)$. For $Bi_+ (I, \Isal, I; x)$,
we use $\cC_- (\wtI)$ instead of $\cC_+ (\wtI)$ to subdivide the set
of bicritical elements.

We now have, by the dual version of Lemma~\ref{lema10}
\begin{equation}
\# Cr (\Qs) \; \leqslant \; C |\wtI|^{-\tau \dms}.\label{eq9.108}
\end{equation}
On the other hand, from formulas (\ref{eq9.18})--(\ref{eq9.20}), we
must replace (\ref{eq9.105}) by
\begin{equation}
\wtB B^{-1} \; = \; |\wtI|^{-\tau\sig_0}.\label{eq9.109}
\end{equation}
Thus, we will choose $\sig_0$ in order to have
\begin{equation}
\sig_0 + \dms \; < \; 1 - 2\tau.\label{eq9.110}
\end{equation}
Obviously, this is compatible with (\ref{eq9.23}), (\ref{eq9.24}),
(\ref{eq9.106}). The end of the argument is the same as before.

This completes the proof of the induction step except for the case
$I=\Isal=\Isom$.

In this last case, we can actually apply the same argument as above
to complete the induction step when $x$ is large. When $\Isal =
\Isom$ is equal to $I$ or $\wtI$, indeed we have
\begin{eqnarray}
B_0 (x) \; &= \; &\Bigl( \frac{x}{\vep_0|P_u|} \Bigr)^{-\rho_0}
\Bigl(\frac{|I|}{\vep_0}\Bigr)^{2\sig_0+\sig_1},\label{eq9.111}\\
B_1 (x) \; &= \; &\Bigl( \frac{x}{\vep_0|P_u|} \Bigr)^{-\rho_1}
\Bigl(\frac{|I|}{\vep_0}\Bigr)^{\sig_0+\sig_1},\label{eq9.112}\\
\wtB_0 (x) \; &= \; &\Bigl( \frac{x}{\vep_0|P_u|} \Bigr)^{-\rho_0}
\Bigl(\frac{|\wtI|}{\vep_0}\Bigr)^{2\sig_0+\sig_1},\label{eq9.113}\\
\wtB_1 (x) \; &= \; &\Bigl( \frac{x}{\vep_0|P_u|} \Bigr)^{-\rho_1}
\Bigl(\frac{|\wtI|}{\vep_0}\Bigr)^{\sig_0+\sig_1}.\label{eq9.114}
\end{eqnarray}
For $x\geqslant \wtx_{cr}$, we, therefore, have
\begin{equation}
\wtB (x) B (x)^{-1} \; = \;
|\wtI|^{-\tau(\sig_0+\sig_1)}\label{eq9.115}
\end{equation}
and the same argument as before applies.

For $x\leqslant x_{cr}$, we have
\begin{equation}
\wtB (x) B (x)^{-1} \; = \;
|\wtI|^{-\tau(2\sig_0+\sig_1)},\label{eq9.116}
\end{equation}
with $2\sig_0 + \sig_1 = 2 - \dos - \dou + o(1)$: this case will be
the object of the rest of Section~\ref{sec9}.

In the intermediate range $x_{cr}\leqslant x\leqslant \wtx_{cr}$, we
write $x = \wtx_{cr} y$, where we have, in view of (\ref{eq9.31}),
\begin{equation}
1 \; \geqslant\; y \; \geqslant \;
|\wtI|^{\tfrac{\tau\sig_0}{\rho_0-\rho_1}}.\label{eq9.117}
\end{equation}
We now have in this range:
\begin{eqnarray}
\wtB (x) B (x)^{-1} \; &= \; &\wtB_0 (x) B_1
(x)^{-1}\label{eq9.118}\\
&= \; &y^{\rho_1-\rho_0} |\wtI|^{-\tau (\sig_0+\sig_1)}.\notag
\end{eqnarray}
Therefore, the proportion of bad candidates will not be greater than
$C |\wtI|^{2\tau^2}$, as long as we have
\begin{equation}
y^{\rho_0-\rho_1} \; \geqslant\;
|\wtI|^{\tau(1-\sig_0-\sig_1-\dum-2\tau)}.\label{eq9.119}
\end{equation}
We recall that $\rho_0-\rho_1>0$ and that
$1-\sig_0-\sig_1-\dum-2\tau$ is positive according to
(\ref{eq9.106}), but it has to be small according to (\ref{eq9.23}),
(\ref{eq9.24}). We state

\begin{Propo}\label{propo31}
Assume that (SR3)$_s$ holds for all parameter intervals containing
$\wtI$. Then, it holds for $Bi_+ (I,I,I; x)$ in the range
\begin{equation}
x \; \geqslant\; \wtx_{cr}\,
|\wtI|^{\tfrac{\tau}{\rho_0-\rho_1}\,(1-\sig_0-\sig_1-\dum-2\tau)}\label{eq9.120}
\end{equation}
except for a proportion of candidates no greater than
$|\wtI|^{\tau^2}$.

For $Bi_- (I,I,I;x)$, the corresponding range is
\begin{equation}
x \; \geqslant\; \wtx'_{cr}\,
|\wtI|^{\tfrac{\tau}{\rho'_0-\rho'_1}\,(1-\sig_0-\sig_1-\dum-2\tau)}.\label{eq9.121}
\end{equation}
\end{Propo}

\subsection{The Case $I = I_\al = I_\om$: General Overview\label{sub9.8}}
We explain in this subsection the strategy which will be pursued in
the case $I = \Isal = \Isom$, $x$ small.

The exponent $2\sig_0+\sig_1$ in (\ref{eq9.116}) means that we have
take into account the criticality of both $P$ and $Q$ in the
selection process.

We start as before. An element $(P,Q,n)$ in $Bi_+ (I,I,I; x)$ is
either in $Bi_+^{new} (I,I,I; x)$ or in $Bi_+ (\wtI,\wtI, \wtI; x)$.
The cardinality of the first set is bounded as before by
\begin{equation}
\# Bi_+^{new} (I,I,I; x) \; \leqslant \; 2\wtB
|\wtI|^{\sig_2}.\label{eq9.122}
\end{equation}
On the other hand, let $(\Psal, \Qsal, \nsal)\in \cC_+ (\wtI)$ and
$(\Psom, \Qsom, \nsom)\in \cC_- (\wtI)$. We denote by $Bi_+ (\Psal,
Q_\om)$ the set of $(P,Q,n)\in Bi_+ (\wtI, \wtI,\wtI; x)$ such that
$P\subset P_\al$ and $Q\subset Q_\om$. These subsets are disjoint
and their union contain the elements $(P,Q,n)$ of $Bi_+ (I,I,I; x)$
which belong to $\cR (\wtI)$ because $P$ and $Q$ are thin
$I$-critical. Moreover, only those $P_\al$, $Q_\om$ such that $I\in
Cr (P_\al) \cap Cr (Q_\om)$ will actually appear.

Thus, in analogy with (\ref{eq9.97}), we have
\begin{equation}
\# Bi_+ (I,I,I; x) \; \leqslant \; 2\wtB |\wtI|^{\sig_2} + \sum \#
Bi_+ (\Psal, Q_\om),\label{eq9.123}
\end{equation}
where the sum is over pairs $(\Psal, \Qsal, \nsal)\in \cC_+ (\wtI)$,
$(\Psom, \Qsom, \nsom) \in \cC_- (\wtI)$ such that $I\subset Cr
(P_\al)$ and $I\subset Cr (Q_\om)$.

Unfortunately, we are not able to estimate the size of the
intersection $Cr (\Psal) \cap Cr (\Qsom)$ directly in a satisfactory
way: while these two sets seem ''independent'' in the naive sense,
it is another matter to translate this intuition in a quantitative
estimate in the style of Lemma~\ref{lema10}.

Instead of this approach, we will use some degree of independence,
when $x$ is small, of the variables $\Psal$ and $\Qsom$ in $Bi_+
(P_\al, Q_\om)$. To explain the technique, consider first the
unrealistic model case where we would have
\begin{equation}
\# Bi_+ (\Psal, Q_\om)\; = \; b_+ (\Psal) b_-
(\Qsom),\label{eq9.124}
\end{equation}
for some functions $b_+$, $b_-$ on $\cC_+ (\wtI)$, $\cC_- (\wtI)$,
respectively.

The formula (\ref{eq9.123}) now gives
\begin{equation}
\# Bi_+ (I,I,I; x) \; \leqslant \; 2\wtB |\wtI|^{\sig_2} + \phi_+
(I) \phi_- (I),\label{eq9.125}
\end{equation}
with
\begin{eqnarray}
\phi_+ (I) \; &= \; &\sum_{I\in Cr (\Psal)}\; b_+
(\Psal),\label{eq9.126}\\
\phi_- (I) \; &= \; &\sum_{I\in Cr (\Qsom)}\; b_-
(\Qsom).\label{eq9.127}
\end{eqnarray}
We now average {\it separately} $\phi_+$ and $\phi_-$. We obtain
\begin{eqnarray}
\sum_I \phi_+ (I) \; &\leqslant \; &\Bigl( \max_{\cC_+(\wtI)}\; \#
Cr
(\Psal)\Bigr) \sum_{\cC_+ (\wtI)}\; b_+ (\Psal),\label{eq9.128}\\
\sum_I \phi_- (I) \; &\leqslant \; &\Bigl( \max_{\cC_-(\wtI)}\; \#
Cr (\Qsom)\Bigr) \sum_{\cC_- (\wtI)}\; b_- (\Qsom),\label{eq9.129}
\end{eqnarray}
where $Cr (\Psal)$ and $Cr (\Qsom)$ are estimated by
Lemma~\ref{lema10} and we have
\begin{equation}
\sum b_+ (\Psal) \sum b_- (\Qsom) \; = \; \sum \# Bi_+ (\Psal,\Qsom)
\;\leqslant \; \wtB.\label{eq9.130}
\end{equation}
It is then sufficient to eliminate candidates for which either
$\phi_+$ or $\phi_-$ is much above its average value to be able to
conclude the proof.

As (\ref{eq9.124}) does not hold, we will find an appropriate
substitute as follows.

We will subdivide each class $Bi_+ (\Psal,\Qsom)$ into subclasses
$Bi_+ (\Psal,\Qsom,\ell)$; the index $\ell$ runs through a finite
large set $L$ dependent on $I$ and $x$ but independent on $\Psal$
and $\Qsom$. Moreover, we will have functions $b_+ (\Psal, \ell)$,
$b_- (\Qsom,\ell)$ on $\cC_+ (\wtI) \times L$, $\cC_- (\wtI)\times L$,
respectively, such that,
\begin{equation}
\# Bi_+ (\Psal, \Qsom, \ell) \; \leqslant \; b_+ (\Psal, \ell), b_-
(\Qsom, \ell).\label{eq9.131}
\end{equation}
We then set, for each $\ell\in L$:
\begin{eqnarray}
\phi_{+,\ell} (I)\; &= \; &\sum_{I\in Cr (\Psal)}\; b_+
(\Psal,\ell),\label{eq9.132}\\
\phi_{-,\ell} (I)\; &= \; &\sum_{I\in Cr (\Qsom)}\; b_-
(\Qsom,\ell).\label{eq9.133}
\end{eqnarray}
We average each of these functions to get, in view of
Lemma~\ref{lema10},
\begin{eqnarray}
\sum_I \; \phi_{+,\ell} (I)\; &\leqslant \; &C |\wtI|^{-\tau \dum}\;
b_+ (\ell),\label{eq9.134}\\
\sum_I \; \phi_{-,\ell} (I)\; &\leqslant \; &C |\wtI|^{-\tau \dms}\;
b_- (\ell),\label{eq9.135}
\end{eqnarray}
with
\begin{eqnarray}
b_+ (\ell) \; &= \; &\sum_{\cC_+ (\wtI)}\; b_+
(\Psal,\ell),\label{eq9.136}\\
b_- (\ell) \; &= \; &\sum_{\cC_- (\wtI)}\; b_-
(\Qsom,\ell).\label{eq9.137}
\end{eqnarray}
For each $\ell$, we will have
\begin{eqnarray}
\phi_{+,\ell} (I)\; &\leqslant \; &|\wtI|^{\tau(1-\dum-3\tau)}\;
b_+ (\ell)\label{eq9.138}\\
\phi_{-,\ell} (I)\; &\leqslant \; &|\wtI|^{\tau(1-\dms-3\tau)}\; b_-
(\ell)\label{eq9.139}
\end{eqnarray}
except for a proportion of candidates not greater than $C
|\wtI|^{3\tau^2}$. Set
\begin{equation}
\whB \; = \; \sum_L \; b_+ (\ell) b_- (\ell).\label{eq9.140}
\end{equation}
Because we need to eliminate candidates for each $\ell$, $L$ should
not be too large. We will have, see Proposition~\ref{propo33} in
Subsection~\ref{sub9.10}, that
\begin{equation}
\# L \; \leqslant \; C |\wtI|^{-\tau^2}.\label{eq9.141}
\end{equation}
Taking into account that we must eliminate candidates for each $x =
2^{-n} \geqslant x_{\min}$, the total proportion of the failed
candidates is at most $|\wtI|^{\tau^2}$. On the other hand, for the
surviving candidates, the discussion above gives
\begin{eqnarray}
\sum_{I\in Cr (\Psal) \cap Cr (\Qsom)} \# Bi_+ (\Psal,\Qsom)
\leqslant \; \sum_L \; \sum_{I\in Cr (\Psal)} \, \sum_{I\in Cr
(\Qsom)}\, b_+ (\Psal,\ell)\, b_- (\Qsom,\ell)\label{eq9.142}\\
= \sum_L \; \phi_{+,\ell} (I) \; \phi_{-,\ell} (I) \leqslant \;
|\wtI|^{\tau(2-\dms-\dum-6\tau)}\, \whB\notag
\end{eqnarray}
\begin{equation}
\# Bi_+ (I,I,I; x) \; \leqslant \; 2\wtB |\wtI|^{\sig_2} + \whB
|\wtI|^{\tau(2-\dms-\dum-6\tau)}.\label{eq9.143}
\end{equation}
We definitely have $\wtB |\wtI|^{\sig_2} \ll B$. In order to obtain
(SR3)$_s$, we need
\begin{equation}
\whB |\wtI|^{\tau (2-\dms-\dum-6\tau)} \; \leqslant \;
CB.\label{eq9.144}
\end{equation}
As $2-\dms-\dum-6\tau = 2\sig_0+\sig_1+o(1)$, we see from
(\ref{eq9.116}) that $\whB$ cannot be much larger than $\wtB$.

In the next three subsections, we will

--~~define precisely $L$ and the subclasses $Bi_+ (\Psal,
\Qsom,\ell)$;

--~~bound the cardinality of $L$ (to obtain (\ref{eq9.141}));

--~~obtain an appropriate estimate for $\whB$ (cf.~(\ref{eq9.144})).

\subsection{Subclasses of Bicritical Elements\label{sub9.9}}

{\bf 9.9.1~~Bound elements.~}

\Def Let $(\Psal,\Qsal,\nsal)\in\cC_+ (\wtI)$, $(\Psom, \Qsom,
\nsom)\in \cC_- (\wtI)$. An element $(P,Q,n)\in Bi_+ (\Psal,\Qsom)$
is {\it bound} if $n\leqslant \nsal + n_\om$. Otherwise, we say that
$(P,Q,n)$ is {\it free}. We will denote by $Bi_+ (\Psal,\Qsal,\loz)$
the subset of bound elements of $Bi_+ (\Psal,\Qsom)$.

Thus, $\loz$ is an element of $L$. On the other hand, free elements
will correspond to many elements of $L$. Recall that we have
$x\leqslant \wtx_{cr}$. When $x \ll \wtx_{cr}$, most elements are
free. When $x\gg \wtx_{cr}$, on the opposite, most elements are
bound.

\begin{Propo}\label{propo32}
For any $(\Psal,\Qsal,\nsal)$, $(\Psom,\Qsom,\nsom)\in \cR (\wtI)$,
and any $n\leqslant \nsal + \nsom$, there is at most one element
$(P,Q,n)\in \cR (\wtI)$ of length $n$ such that $P\subset \Psal$,
$Q\subset \Qsom$.
\end{Propo}

\begin{proof}
We argue by induction on the level of the parameter interval.

When $\wtI$ is the starting interval $\Io$, no parabolic composition
is involved and the result follows from usual symbolic dynamics: as
$n\leqslant \nsal + \nsom$, the word associated to a bound element
is determined by its initial and final parts.

Assume that the result holds for parameter intervals strictly larger
than $\wtI$. Denote by $\wtI_1$ the parent interval of $\wtI$.

Assume first that both $(\Psal,\Qsal,\nsal)$ and
$(\Psom,\Qsom,\nsom)$ belong to $\cR (\wtI_1)$. We claim that any
bound element also belongs to $\cR (\wtI_1)$, which allow us to
conclude the proof by the induction hypothesis. Indeed, if $(P,Q,n)$
satisfies $P\subset \Psal$, $Q\subset \Qsom$ and does not belong to
$\cR (\wtI_1)$, we apply the structure theorem in
Subsection~\ref{sub6.7}: it gives elements $(P_0,Q_0,n_0)$,
$(\Psk,\Qsk,\nsk)\in \cR (\wtI_1)$ such that $P_0$ is the thinnest
rectangle containing $P$ defined over $\wtI_1$, $\Qsk$ is the
thinnest rectangle containing $Q$ defined over $\wtI_1$, and
$n\geqslant n_0 + n_k + N_0$. Therefore, $n_0\geqslant \nsal$,
$n_k\geqslant \nsom$ and $n> \nsal + \nsom$.

We now consider the case when, for instance, $(\Psal,\Qsal,\nsal)$
does not belong to $\cR (\wtI_1)$. We now apply the structure
theorem of Subsection~\ref{sub6.7} to $(\Psal,\Qsal,\nsal)$ and also
to an element $(P,Q,n)\in\cR (\wtI)$ with $P\subset \Psal$, $Q
\subset \Qsom$, $n\leqslant \nsal + \nsom$. We obtain integers $0 <
j \leqslant k$, elements $(P_i, Q_i, n_i)\in \cR (\wtI_1)$ for
$0\leqslant i\leqslant k$ such that
\begin{equation}
(P,Q,n) \; \in \; (P_0, Q_0, n_0) \; \square \cdots \square\; (\Psk,
\Qsk, \nsk)\label{eq9.145}
\end{equation}
and also $(\wtP_j, \wtQ_j, \wtn_j)\in \cR (\wtI_1)$ such that
\begin{equation}
(\Psal, \Qsal, \nsal) \; \in \; (P_0, Q_0, n_0) \; \square \cdots
\square\; (\jsM P, \jsM Q, \jsM n)\; \square \; (\wtP_j, \wtQ_j,
\wtn_j).\label{eq9.146}
\end{equation}
There  also exists $m$ with $0\leqslant m\leqslant j$ and $(\wtP'_m,
\wtQ'_m, \wtn'_m)$ in $\cR (\wtI_1)$ such that

--~~either $m=j=k$, $(\wtP'_m, \wtQ'_m, \wtn'_m) = (\Psom,
\Qsom,\nsom)$;

--~~or $m< k$ and we have,
\begin{equation}
(\Psom, \Qsom, \nsom) \; \in \; (\wtP'_m, \wtQ'_m, \wtn'_m) \;
\square \cdots \square\; (\Psk, \Qsk, \nsk).\label{eq9.147}
\end{equation}
If $m<j$, the sequence $(P_i, Q_i, n_i)$ and the choice of the
result (out of two possibilities) in each parabolic composition are
completely determined by $\Psal$ and $\Qsom$: the assertion of the
proposition follows.

When $m=j$, the sequence $(P_i,Q_i,n_i)$ for $i\ne m = j$ is
determined by $\Psal$, $\Qsom$; but we also have $P_j \subset
\wtP_j$, $Q_j = Q_m \subset \wtQ_m$ and $n_j\leqslant \wtn_j +
\wtn'_m$, so by the induction hypothesis $(P_j, Q_j, n_j)$ is also
determined by $\Psal$, $\Qsom$. Again, the choices of the results in
the parabolic compositions are also determined by $\Psal$, $\Qsom$.
The proof of the proposition is complete.
\end{proof}

Recall that, by Proposition~\ref{propo13} in
Subsection~\ref{sub7.1}, we have
\begin{equation}
|P| \; \leqslant \; C\, \exp (-n^\ga)\label{eq9.148}
\end{equation}
for any $(P,Q,n)\in \cR (\wtI)$, any $t\in\wtI$, with $\ga = \log
\frac{3}{2}/\log 2$. Any bound element in $Bi_+ (\Psal, \Qsom,
\loz)$ must satisfy $|P|\geqslant x$ for some $t\in \wtI$, and we
therefore have
\begin{equation}
n \; \leqslant \; \Bigl( \log\,
\frac{C}{x}\Bigr)^{\tfrac{1}{\ga}}\label{eq9.149}
\end{equation}
For $x\geqslant x_{\min}$, we have
\begin{equation}
\Bigl( \log\, \frac{C}{x}\Bigr)^{\tfrac{1}{\ga}}\; \leqslant \;
(\log\, |\wtI|)^2.\label{eq9.150}
\end{equation}
We thus shall define
\begin{equation}
b_+ (\Psal, \loz) \; = \; b_- (\Qsom,\loz) \; = \; \log \,
|\wtI|^{-1},\label{eq9.151}
\end{equation}
and we will indeed have, from Proposition~\ref{propo32},
\begin{equation}
\# Bi_+ (\Psal, \Qsom, \loz) \; \leqslant \; b_+ (\Psal,\loz)\, b_-
(\Qsom, \loz).\label{eq9.152}
\end{equation}

{\bf 9.9.2~~Decomposition of a free element.~}

Let $(\Psal,\Qsal,\nsal)\in \cC_+ (\wtI)$, $(\Psom,\Qsom,\nsom)\in
\cC_- (\wtI)$ and let $(P,Q,n)\in Bi_+ (\Psal,\Qsom)$ be a {\it
free} element. We will analyze with respect to the structure theorem
of Subsection~\ref{sub6.7} the way in which $(P,Q,n)$,
$(\Psal,\Qsal,\nsal)$, $(\Psom,\Qsom,\nsom)$ have been created. This
will allow us in the sequel to define various subclasses of free
elements.

Denote by $\whI_0$ the largest parameter interval such that
$(P,Q,n)\in \cR (\whI_0)$. Elements $(P,Q,n)$ for which $\whI_0$ is
the starting interval $\Io$ are said to have depth $0$. They form a
first subclass of $Bi_+ (\Psal,\Qsom)$ denoted by $Bi_+
(\Psal,\Qsom,0)$.

We now assume that $\whI_0 \ne \Io$ and denote by $\wtI_0$ the
parent interval of $\whI_0$. We apply the structure theorem of
Subsection~\ref{sub6.7}. We obtain an integer $k>0$, elements
$(P_0,Q_0,n_0),\cdots,(\Psk,\Qsk,\nsk)$ in $\cR (\wtI_0)$ such that
\begin{equation}
(P,Q,n) \; \in \; (P_0, Q_0, n_0) \; \square \cdots \square\; (\Psk,
\Qsk, \nsk).\label{eq9.153}
\end{equation}
As in the proof of Proposition~\ref{propo32}, we find $0\leqslant
j\leqslant k$ and $(\wtP_j,\wtQ_j,\wtn_j)\in \cR (\wtI_0)$ such that
either $j=0$, $(\wtP_j, \wtQ_j, \wtn_j) = (\Psal,\Qsal,\nsal)$ (if
$(\Psal,\Qsal,\nsal) \in \cR (\wtI_0)$) or $j>0$ and
\begin{equation}
(\Psal,\Qsal,\nsal) \; \in \; (P_0, Q_0, n_0) \; \square \cdots
\square\; (\wtP_j, \wtQ_j, \wtn_j).\label{eq9.154}
\end{equation}
Similarly, we find $0\leqslant m\leqslant k$ and $(\wtP'_m, \wtQ'_m,
\wtn'_m)\in \cR (\wtI_0)$ such that either $m=k$, $(\wtP'_m,
\wtQ'_m, \wtn'_m) = (\Psom, \Qsom,\nsom)$ or $m< k$ and
\begin{equation}
(\Psom,\Qsom,\nsom) \; \in \; (\wtP'_m, \wtQ'_m, \wtn'_m) \; \square
\cdots \square\; (\Psk,\Qsk,\nsk).\label{eq9.155}
\end{equation}
We also must have $P_j\subset \wtP_j$, $Q_m\subset \wtQ'_m$.
Moreover, as $(P,Q,n)$ is free, we must have $j\leqslant m$ and,
when $j=m$, we must also have $n_j=n_m > \wtn_j + \wtn'_m$.

We say that $(P,Q,n)$ is {\it fully decomposed} if one has here
$j<m$ or $j=m$ and $(P_j,Q_j,n_j)\in \cR (\Io)$. Such elements are
said to have depth one.

Assume that $(P,Q,n)$ is not fully decomposed. Then, we have $j=m$,
$P_j\subset \wtP_j$, $Q_j\subset \wtQ'_j$ and the largest parameter
interval $\whI_1$ for which $(P_j,Q_j,n_j)\in \cR (\whI_1)$ is not
the starting interval $\Io$. We denote by $\wtI_1$ the parent
interval. We rewrite
\begin{eqnarray}
(P^1, Q^1, n^1) \; &:= \; &(P_j, Q_j, n_j),\label{eq9.156}\\
(P^1_\al, Q^1_\al, n^1_\al) \; &:= \; &(\wtP_j, \wtQ_j,
\wtn_j),\notag\\
(P^1_\om, Q^1_\om, n^1_\om) \; &:= \; &(\wtP'_j, \wtQ'_j,
\wtn'_j),\notag
\end{eqnarray}
and proceed with these elements as we did with $(P,Q,n)$,
$(\Psal,\Qsal,\nsal)$, $(\Psom,\Qsom,\nsom)$: we will find integers
$0\leqslant j_1\leqslant m_1 \leqslant k_1$ (with $k_1 > 0$),
elements $(P^1_i, Q^1_i, n^1_i)$ for $0\leqslant i \leqslant k_1$
and also $(P^2_\al, Q^2_\al, n^2_\al)$, $(P^2_\om, Q^2_\om,
n^2_\om)$, all in $\cR (\wtI_1)$, such that
\begin{eqnarray}
(P^1, Q^1, n^1) &\in &(P^1_0, Q^1_0, n^1_0)\; \square \cdots
\square\; (P^1_{k_1}, Q^1_{k_1}, n^1_{k_1}),\label{eq9.157}\\
(P^1_\al, Q^1_\al, n^1_\al) &\in &(P^1_0, Q^1_0, n^1_0)\; \square
\cdots \square\; (P^1_{j_{1}-1}, Q^1_{j_{1}-1}, n^1_{j_{1}-1}) \;
\square \; (P^2_\al, Q^2_\al, n^2_\al),\notag\\
(P^1_\om, Q^1_\om, n^1_\om) &\in &(P^2_\om, Q^2_\om, n^2_\om)\;
\square\; (P^1_{m_{1}+1}, Q^1_{m_{1}+1}, n^1_{m_{1}+1})\; \square
\cdots \square \; (P^1_{k_1}, Q^1_{k_1}, n^1_{k_1}).\notag
\end{eqnarray}
Again, we say that $(P^1, Q^1, n^1)$ is fully decomposed if either
$j_1< m_1$ or $j_1 = m_1$ and $(P^1_{j_1}, Q^1_{j_1}, n^1_{j_1})$ is
defined over the starting interval $\Io$; otherwise we set
\begin{equation}
(P^2, Q^2, n^2) \; := \; (P^1_{j_1}, Q^1_{j_1},
n^1_{j_1}),\label{eq9.158}
\end{equation}
and we go on. The sequence of parameter intervals $\whI_0 \subset
\whI_1 \subset \cdots$ is strictly increasing and therefore the
process will stop. We define inductively the {\it depth} of
$(P,Q,n)$ to be the depth of $(P^1, Q^1, n^1)$ plus one.

{\bf 9.9.3~~Size of the subclass of depth $0$}

We will define in this subsection $b_+ (\Psal,0)$, $b_- (\Qsom, 0)$
in order to have
\begin{equation}
\# Bi_+ (\Psal, \Qsom, 0)\; \leqslant \; b_+ (\Psal,0) \; b_-
(\Qsom, 0).\label{eq9.159}
\end{equation}
Let $(\Psal,\Qsal,\nsal)\in \cC_+ (\wtI)$, $(\Psom,\Qsom,\nsom)\in
\cC_- (\wtI)$ and let $(P,Q,n)\in Bi_+ (\Psal, \Qsom, 0)$. Then
$(P,Q,n)$ is obtained by a simple composition
\begin{equation}
(P,Q,n) \; = \; (\Psal,\Qsal,\nsal) \; * \; (P',Q',n') \; * \;
(\Psom, \Qsom, \nsom).\label{eq9.160}
\end{equation}
We have here, for all $t\in \wtI$
\begin{eqnarray}
|P| &\leqslant \; &C |P_\al|\; |P'| \; |P_\om|,\label{eq9.161}\\
|\Psal| &\leqslant \; &|\wtI|^{1+\tau}\;\; \text{(cf.
Subsection~\ref{sub9.1})},\label{eq9.162}
\end{eqnarray}
and also, for some $t_0\in \wtI$
\begin{equation}
|P| \; \geqslant \; x.\label{eq9.163}
\end{equation}
This gives, for this value $t_0$:
\begin{equation}
|P'| \; \geqslant \; C^{-1}\, |\wtI|^{-(1+\tau)}\, \Bigl(\max_{\wtI}
|\Psom|\Bigr)^{-1}\, x.\label{eq9.164}
\end{equation}
We observe that, as $(\Psom, \Qsom, \nsom)$ belongs to $\cR (\Io)$
and $|\Qsom|$ is of the order of $|\wtI|^{1+\tau}$, we have
\begin{equation}
\max_{\wtI} |\Psom| \; \leqslant \; C\; \min_{\wtI}\,
|\Psom|.\label{eq9.165}
\end{equation}
From the estimate (\ref{eq8.38}) in Subsection~\ref{sub8.2}, we can
thus define, as $\dos + C\vep_0 < \dss$,
\begin{equation}
b_+ (\Psal,0) \; =
\begin{cases}
(C |\wtI|^{1+\tau} x^{-1})^{\dss} \qquad &\text{ if $\;
(\Psal,\Qsal,\nsal)\in\cR (\Io)$},\\
0 \qquad &\text{ otherwise},
\end{cases}
\label{eq9.166}
\end{equation}
\begin{equation}
b_- (\Qsom,0) \; =
\begin{cases}
(\displaystyle{\min_{\wtI}} |\Psom|)^{\dss} \qquad &\text{ if $\;
(\Psom,\Qsom,\nsom)\in\cR (\Io)$},\\
0 \qquad &\text{}
\end{cases}
\label{eq9.167}
\end{equation}
Then, (\ref{eq9.159}) is satisfied.

{\bf 9.9.4~~Subclasses of higher depth}

Let $(\Psal,\Qsal,\nsal)\in \cC_+ (\wtI)$, $(\Psom, \Qsom, \nsom)\in
\cC_- (\wtI)$ and let $(P,Q,n)\in Bi_+ (\Psal,\Qsom)$ be an element
of depth $s>0$.

Let us first restate and extend somewhat the notations and the
setting of 9.9.2. We set
\begin{eqnarray}
(P^0, Q^0, n^0) \; &:= \; &(P,Q,n),\label{eq9.168}\\
(P^0_\al, Q^0_\al, n^0_\al) \; &:= \; &(P_\al,Q_\al,n_\al),\notag\\
(P^0_\om, Q^0_\om, n^0_\om) \; &:= \; &(P_\om,Q_\om,n_\om).\notag
\end{eqnarray}
We have

--~~a strictly increasing sequence of parameter intervals
\begin{equation}
\whI_0 \subset \whI_1 \subset \cdots \subset \whI_{s-1} \subset
\whI_s = I_0\label{eq9.169}
\end{equation}
with $\wtI \subset \whI_0$; we denote by $\wtI_r$ the parent
interval of $\whI_r$ for $0\leqslant r < s$;

--~~a sequence $(P^r, Q^r, n^r)$, $0\leqslant r\leqslant s$ such
that $(P^r, Q^r, n^r)$ belongs to $\cR (\whI_r)$ but not to $\cR
(\wtI_r)$ for $r<s$; also $(P^s, Q^s, n^s)\in \cR (\whI_{s-1})$;

--~~two sequences $(P^r_\al, Q^r_\al, n^r_\al)$, $(P^r_\om, Q^r_\om,
n^r_\om)$, $0\leqslant r\leqslant s$; for each $r<s$, resp.~$r=s$,
the two elements belonging to $\cR (\whI_r)$, resp.~$\cR
(\wtI_{s-1})$;

--~~two sequences $(P^r_+, Q^r_+, n^r_+)$, $(P^r_-, Q^r_-, n^r_-)$,
$0< r\leqslant s$; for each $r$, the two elements belonging to $\cR
(\whI_{r-1})$.

These data are related by the following properties: for each $0 < r
\leqslant s$, we have
\begin{eqnarray}
(P^{r-1}, Q^{r-1}, n^{r-1}) &\in &(P^r_-, Q^r_-, n^r_-)\; \square \;
(P^r, Q^r, n^r) \; \square \; (P^r_+, Q^r_+, n^r_+),\label{eq9.170}\\
(P^{r-1}_\al, Q^{r-1}_\al, n^{r-1}_\al) &\in &(P^r_-, Q^r_-,
n^r_-)\; \square \; (P^r_\al, Q^r_\al, n^r_\al),\label{eq9.171}\\
(P^{r-1}_\om, Q^{r-1}_\om, n^{r-1}_\om) &\in &(P^r_\om, Q^r_\om,
n^r_\om)\; \square \; (P^r_+, Q^r_+, n^r_+).\label{eq9.172}
\end{eqnarray}
The process stops at step $s$ because of one of the two following cases
occur

a)~~$(P^s, Q^s, n^s)$ does not belong to $\cR (\wtI_{s-1})$; then,
by the structure theorem of Subsection~\ref{sub6.7}, there exists an
integer $h>0$, elements $(P^s_0, Q^s_0, n^s_0) \cdots (P^s_h, Q^s_h,
n^s_h)$ in $\cR (\wtI_{s-1})$ with
\begin{equation}
(P^s, Q^s, n^s) \in (P^s_0, Q^s_0, n^s_0)\; \square \cdots \square
\; (P^s_h, Q^s_h, n^s_h)\label{eq9.173}
\end{equation}
and also
\begin{equation}
P^s_0 \subset P^s_\al,\qquad Q^s_h \subset Q^s_\om.\label{eq9.174}
\end{equation}

b)~~$(P^s, Q^s, n^s)$ belongs to $\cR (I_0)$; in this case we set
$h=0$.

We also observe that the parabolic compositions in (\ref{eq9.170})
through (\ref{eq9.172}) take place in $\cR (\whI_{r-1})$ but not in
$\cR (\wtI_{r-1})$; in (\ref{eq9.173}), they take place in $\cR
(\whI_{s-1})$ but not in $\cR (\wtI_{s-1})$.

A subclass $Bi_+ (\Psal,\Qsom,\ell)$, i.e. an element of $L$,
distinct from the two $(\loz, 0)$ that we already know is
determined by the following data

--~~the depth $s (>0)$;

--~~the sequence $\whI_0 \subset \cdots \subset \whI_s = \Io$;

--~~the integer $h\geqslant 0$;

--~~for each $0 < i < h$, the smallest integer $u_i$ such that
$|P^s_i|\geqslant 2^{-u_i} =: x_i$ for some $t\in\wtI$ when $h>1$;

--~~when $h>0$, the smallest integers $u_0$, $u_h$ such that
$|P_-|\geqslant 2^{-u_0} =: x_0$ for some $t_-\in \wtI$,
$|P_+|\geqslant 2^{-u_h} =: x_h$ for some $t_+\in\wtI$; here, the
elements $(P_-, Q_-, n_-)$, $(P_+, Q_+, n_+)$ are determined by
$P\subset P_-$, $Q\subset Q_+$ and
\begin{eqnarray}
(P_-, Q_-, n_-)\; \in \; (P^1_-, Q^1_-, n^1_-)\; \square \cdots
\square \; (P^s_-, Q^s_-, n^s_-)\; \square\; (P^s_0, Q^s_0,
n^s_0),\label{eq9.175}\\
(P_+, Q_+, n_+)\; \in \; (P^s_h, Q^s_h, n^s_h)\; \square \; (P^s_+,
Q^s_+, n^s_+)\; \square \cdots \square \; (P^1_+, Q^1_+,
n^1_+).\label{eq9.176}
\end{eqnarray}
Thus, we group together in a subclass $Bi_+ (\Psal,\Qsom,\ell)$ the
elements of $Bi_+ (\Psal,\Qsom)$ who share the same data; the
elements of $L$, distinct from $\loz$, $0$, are the sets of data
for which at least one subclass $Bi_+ (\Psal,\Qsom,\ell)$ is
non-empty, for some $(\Psal,\Qsal,\nsal)$ in $\cC_+ (\wtI)$,
$(\Psom, \Qsom, \nsom)$ in $\cC_- (\wtI)$.

The definition of the set $L$ is now complete.

{\bf 9.9.5~~Sizes of subclasses of higher depth}

The context and notations are the same as above. We want to define
$b_+ (\Psal,\ell)$ and $b_- (\Qsom,\ell)$ in order to satisfy
(\ref{eq9.131}) in Subsection~\ref{sub9.8}.

We first observe that $(\Psal,\Qsal,\nsal)$ determines $(P^1_-,
Q^1_-, n^1_-)$, $\cdots$, $(P^s_-, Q^s_-, n^s_-)$, $(P^s_\al,
Q^s_\al, n^s_\al)$ and the result of parabolic compositions between
these elements. Similarly, $(\Psom, \Qsom, \nsom)$ determines
$(P^1_+, Q^1_+, n^1_+)$, $\cdots$, $(P^s_+, Q^s_+, n^s_+)$,
$(P^s_\om, Q^s_\om, n^s_\om)$ and the result of parabolic
compositions between these elements. Therefore, the only ''freedom''
for the element $(P,Q,n)$ in the subclass $Bi_+ (\Psal,\Qsom,\ell)$
is through $(P^s, Q^s, n^s)$, and this freedom is constrained by the
relations $P^s \subset P^s_\al$, $Q^s \subset Q^s_\om$.

Consider first a subclass with $h=0$, i.e., $(P^s, Q^s, n^s) \in \cR
(\Io)$. The widths of the strips are related as follows: for every
$t\in \wtI$, we have
\begin{equation}
C^{-1}\; \frac{|P^s|}{|P^s_\al|\;|P^s_\om|} \; \leqslant \;
\frac{|P|}{|P_\al|\;|P_\om|} \; \leqslant \; C\;
\frac{|P^s|}{|P^s_\al|\;|P^s_\om|}.\label{eq9.177}
\end{equation}
This allows us to take, as in the case of depth $0$,
\begin{equation}
b_+ (\Psal,\ell) \; =
\begin{cases}
(C |\wtI|^{1+\tau} x^{-1})^{\dss} \qquad &\text{ if $\; (P^s_\al,
Q^s_\al, n^s_\al)\in\cR (\Io)$},\\
0 \qquad &\text{ otherwise},
\end{cases}
\label{eq9.178}
\end{equation}
\begin{equation}
b_- (\Qsom,\ell) \; =
\begin{cases}
|\Psom|^{\dss} \qquad &\text{ if $\; (P^s_\om,
Q^s_\om, n^s_\om)\in\cR (\Io)$},\\
0 \qquad &\text{}
\end{cases}
\label{eq9.179}
\end{equation}
Consider now a subclass with $h>0$, i.e., case a) in
Subsection~9.9.4 above.

By the structure theorem in Subsection~\ref{sub6.7}, see
Lemma~\ref{lema11} in Subsection~\ref{sub9.6}, for $0<i<h$, the
element $(P^s_i, Q^s_i, n^s_i)$ belongs to $Bi_+ (\wtI_{s-1},
\wtI_{s-1}, \wtI_{s-1}; x_i)$. From Corollary~\ref{coro6} in
Subsection~\ref{sub6.7}, we have
\begin{equation}
x\; \leqslant \; C \Bigl( C |\whI_{s-1}|^{-\fudt}\Bigr)^h\; x_0 x_1
\cdots x_h.\label{eq9.180}
\end{equation}
Thus, the data of every subclass must satisfy (\ref{eq9.180}).
Assuming that (\ref{eq9.180}) holds, we set $b_+ (\Psal,\ell) = 0$
if $(\Psal,\Qsal,\nsal)\not\in \cR (\whI_0)$. When
$(\Psal,\Qsal,\nsal)\in \cR (\whI_0)$, we set
\begin{equation}
b_+ (\Psal,\ell)\; = \; 2^h \Bigl( \prod_{0<i<h}\, \# Bi_+
(\wtI_{s-1}, \wtI_{s-1}, \wtI_{s-1}; x_i)\Bigr) \; \# Bi_+ (\Psal,
\wtI_{s-1}; x_0).\label{eq9.181}
\end{equation}
Here, $Bi_+ (\Psal, \wtI_{s-1}, x_0)$ is by definition the set of
elements $(P_-, Q_-, n_-)$ in $\cR (\wtI)$ such that $P_- \subset
\Psal$, $Q_-$ is thin $\wtI_{s-1}$-critical and $|P_-|\geqslant x_0$
for some $t\in\wtI$.

Similarly, when (\ref{eq9.180}) holds, we set $b_- (\Qsom, \ell) =
0$ if $(\Psom,\Qsom,\nsom)\not\in \cR (\whI_0)$. When $(\Psom,
\Qsom, \nsom) \in \cR (\whI_0)$, we set
\begin{equation}
b_- (\Qsom,\ell)\; = \; \# Bi_+ (\wtI_{s-1}, \Qsom;
x_h),\label{eq9.182}
\end{equation}
where now $Bi_+ (\wtI_{s-1}, \Qsom; x_h)$ is the set of elements
$(P_+, Q_+, n_+)$ in $\cR (\wtI)$ such that $Q_+ \subset \Qsom$,
$P_+$ is thin $\wtI_{s-1}$-critical and $|P_+|\geqslant x_h$ for
some $t\in \wtI$.

The factor $2^h$ in (\ref{eq9.181}) takes care of the possible
results of the ''free'' parabolic compositions, i.e., those
compositions which are not constrained by $(\Psal,\Qsal,\nsal)$ or
$(\Psom, \Qsom, \nsom)$.

The definition of $L$, $b_+$, $b_-$ is now complete, and relation
(9.131) is satisfied.

\subsection{The Size of the Index set $L$\label{sub9.10}}
It is not difficult from (\ref{eq9.180}) to see that the index set
$L$ is finite, but we need an explicit bound on its cardinality
(cf.~(\ref{eq9.141})).

\begin{Propo}\label{propo33}
The index set $L$ satisfies
\[
\# L \; \leqslant \; C |\wtI|^{-\tau^2}.
\]
\end{Propo}

\begin{proof}
In the first part of the proof, we fix the depth $s$ and the
sequence of intervals $\whI_0 \subset \cdots \subset \whI_{s-1}
\subset \whI_s = \Io$. There is one subclass with $h=0$ and we will
estimate the number of subclasses with $h>0$, i.e., the number of
$(h+1)$-tuples $(u_0, \cdots, u_h)$ such that (\ref{eq9.180}) is
satisfied; the integer $h$ itself is {\it not} fixed.

By Corollary~\ref{coro7} in Subsection~\ref{sub6.7}, we have
\begin{eqnarray}
x_i \; &< \; &|\wtI_{s-1}|^\be \;\; \text{ for $\; 0 < i <
h$},\label{eq9.183}\\
x_h \; &< \; &C |\wtI_{s-1}|^{(1-\eta)^{-1}}.\label{eq9.184}
\end{eqnarray}
As $P_-\subset \Psal$, we also have, for a non-empty subclass
\begin{equation}
x_0 \; < \; |\wtI|^{1+\tau}.\label{eq9.185}
\end{equation}
We rewrite (\ref{eq9.180}) as
\begin{equation}
\frac{x_0}{|\wtI|^{1+\tau}} \left( \prod_{0<i<h}\;
\frac{x_i}{|\wtI_{s-1}|^\be}\right) \; \frac{x_h}{C
|\wtI_{s-1}|^{\tfrac{1}{1-\eta}}}\; \geqslant\; \frac{x}{C
|\wtI|^{1+\tau}\,|\wtI_{s-1}|^{\be (h-1)+(1-\eta)^{-1}}}\, \Bigl(
C^{-1} \, |\whI_{s-1}|^{\fudt}\Bigr)^h.\label{eq9.186}
\end{equation}
Using $\be > 1$, and taking base-two logarithms, it is sufficient to
bound the number of non-negative integral solutions of
\begin{equation}
n_0 + \cdots + n_h \; \leqslant \; A_0 - A_1 h,\label{eq9.187}
\end{equation}
with
\begin{eqnarray}
A_0\; &= \; &\log_2 (|\wtI| x^{-1}),\label{eq9.188}\\
A_1\; &= \; &\frac{1}{3}\, \log_2  |\whI_{s-1}|^{-1}.\label{eq9.189}
\end{eqnarray}
As $x\leqslant \wtx_{cr}\ll |\wtI|$, both $A_0$ and $A_1$ are large;
by taking $A_0$ slightly larger and $A_1$ slightly smaller, we can
assume that both $A_0$, $A_1$ are integers. The number of
non-negative integral solutions of (\ref{eq9.187}) is then the
coefficient of $z^{A_0}$ in the power series for
\begin{equation}
\chi (z) \; := \; \sum_{h\geqslant 0} z^{A_1 h} (1-z)^{-h-2}\; = \;
(1-z)^{-1} (1-z-z^{A_1})^{-1}.\label{eq9.190}
\end{equation}
We estimate this coefficient by a Cauchy integral on the circle
$|z|=1 - 2 A_1^{-1} \log A_1$. On this circle, we have
\begin{eqnarray}
|z^{A_1}| \; &< \; &A_1^{-1},\label{eq9.191}\\
|\chi (z)| \; &< \; &\frac{1}{2}\; A_1^2 (\log\,
A_1)^{-2}.\label{eq9.192}
\end{eqnarray}
The number of solutions of (\ref{eq9.187}) is, therefore, not greater
than
\begin{equation}
A^2_1 (\log\, A_1)^{-2} (1- 2A_1^{-1} \log
\,A_1)^{-A_0}.\label{eq9.193}
\end{equation}
In view of (\ref{eq9.189}), this quantity is smaller than
\begin{equation}
(\log\, |\whI_{s-1}|)^2 \exp \left( C A_0 \; \frac{\log |\log
|\whI_{s-1}||}{|\log |\whI_{s-1}||}\right). \label{eq9.194}
\end{equation}
This is a bound for the number of subclasses with fixed depth $s$
and sequence $\whI_0 \subset \cdots \subset \whI_{s-1}$. We have now
to sum over these remaining data. Observe that (\ref{eq9.194})
depends on $|\whI_{s-1}|$, {\it not} on the depth $s$ and the
intervals $\whI_r$, $0\leqslant r < s-1$.

Fix an interval $\whI$ with $\wtI \subset \whI \subset \Io$, $\whI
\ne \Io$. Let $S (\whI)$ be the number of parameter intervals $\Is$
with $\wtI \subset \Is \subset \whI$, $\Is \ne \wtI$. Every $\Is$ in
this range may or may not be one of the $\whI_r$, for a sequence
$\whI_0 \subset \cdots \subset \whI_{s-1}$ terminating with
$\whI_{s-1} = \whI$; in other terms, there are exactly $2^{S(\whI)}$
such sequences (of various lengths). This means that the total
number of subclasses is bounded by
\begin{equation}
\sum_{\whI} \, 2^{S(\whI)} (\log\,|\whI|)^2 \exp \left( CA_0 \;
\frac{\log |\log\,|\whI||}{|\log\,|\whI||} \right).\label{eq9.195}
\end{equation}
We have here
\begin{equation}
\frac{\log |\log\,|\whI||}{|\log\,|\whI||} \; \leqslant \;
\frac{\log\,\log\,\vep_0^{-1}}{\log\,\vep_0^{-1}},\label{eq9.196}
\end{equation}
\begin{equation}
|\log\,|\wtI| \; = \; (1+\tau)^{S(\whI)} \;
\log\,|\whI|,\label{eq9.197}
\end{equation}
\begin{equation}
S (\whI) \; \leqslant \; 2\tau^{-1} \log \left(
\frac{\log\,|\wtI|^{-1}}{\log\,\vep_0^{-1}}\right) \; =: \;
S_{\max}.\label{eq9.198}
\end{equation}
The sum (\ref{eq9.195}) is thus bounded by
\begin{eqnarray}
C 2^{S_{\max}} (\log\,\vep_0^{-1})^2 \exp \left( CA_0 \; \frac{\log
\, \log\,\vep_0^{-1}}{\log\,\vep_0^{-1}}\right)\label{eq9.199}\\
\leqslant \; (\log\, |\wtI|^{-1})^{2\tau^{-1}}\, \exp \left( CA_0 \;
\frac{\log \, \log\,\vep_0^{-1}}{\log\,\vep_0^{-1}}\right).\notag
\end{eqnarray}
As $x\geqslant x_{\min} := |\wtI|^{C (\rho_0-\dss)^{-1}}$
(cf.~(\ref{eq9.64})), we have
\begin{equation}
A_0 \; \leqslant \; C (\rho_0 - \dss)^{-1} \, \log
|\wtI|^{-1}.\label{eq9.200}
\end{equation}
We choose the exponent $\rho_0$ in order to have
\begin{equation}
\rho_0 \; > \; \dss + \tau.\label{eq9.201}
\end{equation}
As $\vep_0$ can be chosen arbitrarily small, we have
\begin{equation}
C A_0 \; \frac{\log\,\log\,\vep_0^{-1}}{\log\,\vep_0^{-1}} \; < \;
\fud\, \tau^2 \,\log\, |\wtI|^{-1}.\label{eq9.202}
\end{equation}
We conclude, then, that with $\vep_0$ small enough, the term in
(\ref{eq9.199}) is indeed smaller than $|\wtI|^{-\tau^2}$.
\end{proof}

\subsection{The Size of $\whB$\label{sub9.11}}
According to the roadmap exposed in Subsection~\ref{sub9.8}, we have
now to estimate the quantity set in Subsection~\ref{sub9.8}
\begin{equation*}
\whB \; = \; \sum_L \; b_+ (\ell) b_- (\ell)\tag{9.140}
\end{equation*}
with
\begin{equation*}
b_+ (\ell) \; = \; \sum_{\cC_+ (\wtI)}\; b_+
(\Psal,\ell),\tag{9.136}
\end{equation*}
\begin{equation*}
b_- (\ell) \; = \; \sum_{\cC_- (\wtI)}\; b_-
(\Qsom,\ell).\tag{9.137}
\end{equation*}
Consider first the bound elements. In view of (\ref{eq9.150}), we
have:
\begin{eqnarray}
b_+ (\loz) \; = \; \# \cC_+ (\wtI) \, \log \, |\wtI|^{-1},\label{eq9.203}\\
b_- (\loz) \; = \; \# \cC_- (\wtI) \, \log \, |\wtI|^{-1}.\notag
\end{eqnarray}
Consider next the class of depth $0$, and also the classes of higher
depth with $h=0$: in view of (\ref{eq9.166})--(\ref{eq9.167}) and
(\ref{eq9.178})--(\ref{eq9.179}), we have in these cases
\begin{eqnarray}
b_+ (\ell) \; &\leqslant \; &(C |\wtI|^{1+\tau} x^{-1})^{\dss} \, \#
\cC_+ (\wtI),\label{eq9.204}\\
b_- (\ell) \; &\leqslant \; &\sum_{\cC_- (\wtI)}\, (\max_{\wtI} \,
|\Psom|)^{\dss}.\label{eq9.205}
\end{eqnarray}
Also, the number of such classes, according to the discussion in the
proof of Proposition~\ref{propo33} is not larger than
\begin{equation}
2^{S_{\max}} \; \leqslant \; \left( \frac{\log\,
|\wtI|^{-1}}{\log\,\vep_0^{-1}} \right)^{2\tau^{-1}}.\label{eq9.206}
\end{equation}
The remaining subclasses are more complicated! Formulas
(\ref{eq9.181}), (\ref{eq9.182}) suggest an induction. We thus
assume that (SR3)$_s$ is satisfied for all parameter intervals
containing $\wtI$. We have, for a class of depth $s>0$ with $h>0$:
\begin{equation}
b_+ (\ell) \; \leqslant \; 2^h \left( \prod_{0<i<h}\, (\# Bi_+
(\wtI_{s-1}, \wtI_{s-1}, \wtI_{s-1}; x_i))\right) \,\#
Bi_+(\wtI,\wtI, \wtI_{s-1}; x_0)\label{eq9.207}
\end{equation}
\begin{equation}
b_- (\ell) \; \leqslant \; \# Bi_+ (\wtI, \wtI_{s-1}, \wtI;
x_h).\label{eq9.208}
\end{equation}
Observe that, from (\ref{eq9.31}), the critical value $x_{cr}$ in
each of the $Bi_+$ sets above in the same and equal to
\begin{equation}
x_{cr} \; := \; \vep_0 |P_u| \Bigl( \frac{|\wtI_{s-1}|}{\vep_0}
\Bigr)^{\tfrac{\sig_0}{\rho_0-\rho_1}}.\label{eq9.209}
\end{equation}
As in Subsection~\ref{sub9.6}, we separate the subclasses into two
parts: those for which every $x_i$ is above the critical value
$x_{cr}$ and the others. In the first case, we have from (SR3)$_s$
\begin{equation}
\# Bi_+ (\wtI_{s-1}, \wtI_{s-1}, \wtI_{s-1}; x_i)\; \leqslant \; C
\Bigl( \frac{x_i}{\vep_0 |P_u|} \Bigr)^{-\rho_1} \, \Bigl(
\frac{|\wtI_{s-1}|}{\vep_0} \Bigr)^{\sig_0 + \sig_1},\text{ for $0 <
i < h$,} \label{eq9.210}
\end{equation}
\begin{equation}
\# Bi_+ (\wtI, \wtI, \wtI_{s-1}; x_0)\; \leqslant \; C \Bigl(
\frac{x_0}{\vep_0 |P_u|} \Bigr)^{-\rho_1} \, \Bigl(
\frac{|\wtI|}{\vep_0} \Bigr)^{\sig_0 + \sig_1},\label{eq9.211}
\end{equation}
\begin{equation}
\# Bi_+ (\wtI, \wtI_{s-1}, \wtI; x_h)\; \leqslant \; C \Bigl(
\frac{x_h}{\vep_0 |P_u|} \Bigr)^{-\rho_1} \, \Bigl(
\frac{|\wtI_{s-1}|}{\vep_0} \Bigr)^{\sig_1}\, \Bigl(
\frac{|\wtI|}{\vep_0} \Bigr)^{\sig_0}.\label{eq9.212}
\end{equation}
Multiplying these inequalities, we obtain, taking (\ref{eq9.180})
into account
\begin{equation}
b_+ (\ell)\, b_- (\ell) \; \leqslant \; A_2^h \, A_3\label{eq9.213}
\end{equation}
with
\begin{equation}
A_2 \; = \; 2 \Bigl( \frac{|\wtI_{s-1}|}{\vep_0}
\Bigr)^{\sig_0+\sig_1}\, \Bigl( C\vep_0 |P_u| \,
|\whI_{s-1}|^{-\fudt}\Bigr)^{\rho_1}\label{eq9.214}
\end{equation}
\begin{equation}
A_3 \; = \; C^{\rho_1} \Bigl( \frac{x}{\vep_0 |P_u|}
\Bigr)^{-\rho_1}\, \Bigl( \frac{|\wtI|}{\vep_0}
\Bigr)^{\sig_0+\sig_1}\label{eq9.215}
\end{equation}
In the second case, as $\rho_0 > \rho_1$, we have
\begin{equation}
\# Bi_+ (\wtI_{s-1}, \wtI_{s-1}, \wtI_{s-1}; x_i)\; \leqslant \; C
\Bigl( \frac{x_i}{\vep_0 |P_u|} \Bigr)^{-\rho_0} \, \Bigl(
\frac{|\wtI_{s-1}|}{\vep_0} \Bigr)^{\sig_0 + \sig_1},\label{eq9.216}
\end{equation}
\begin{equation}
\# Bi_+ (\wtI, \wtI, \wtI_{s-1}; x_0)\; \leqslant \; C \Bigl(
\frac{x_0}{\vep_0 |P_u|} \Bigr)^{-\rho_0} \, \Bigl(
\frac{|\wtI|}{\vep_0} \Bigr)^{\sig_0 + \sig_1},\label{eq9.217}
\end{equation}
\begin{equation}
\# Bi_+ (\wtI, \wtI_{s-1}, \wtI; x_h)\; \leqslant \; C \Bigl(
\frac{x_h}{\vep_0 |P_u|} \Bigr)^{-\rho_0} \, \Bigl(
\frac{|\wtI_{s-1}|}{\vep_0} \Bigr)^{\sig_1}\, \Bigl(
\frac{|\wtI|}{\vep_0} \Bigr)^{\sig_0}.\label{eq9.218}
\end{equation}
Moreover, comparing the $\sig$ exponents in (\ref{eq9.19}) and
(\ref{eq9.20}), we see that if $x_i \leqslant x_{cr}$, $0\leqslant i
\leqslant h$, the corresponding inequality is still true after
multiplying the right-hand side by
$\Bigl(\frac{|\wtI_{s-1}|}{\vep_0} \Bigr)^{\sig_0}$. As this happens
at least once, we get, by multiplying the three inequalities together:
\begin{equation}
b_+ (\ell)\, b_- (\ell) \; \leqslant \; \wtA_2^h \,
\wtA_3,\label{eq9.219}
\end{equation}
now with
\begin{equation}
\wtA_2 \; = \; 2 \Bigl( \frac{|\wtI_{s-1}|}{\vep_0}
\Bigr)^{\sig_0+\sig_1} \, \Bigl( C\vep_0 |P_u|\;
|\whI_{s-1}|^{-\fudt}\Bigr)^{\rho_0}\label{eq9.220}
\end{equation}
\begin{equation}
\wtA_3 \; = \; C^{\rho_0} \Bigl( \frac{x}{\vep_0 |P_u|}
\Bigr)^{-\rho_0} \, \Bigl(
\frac{|\wtI|}{\vep_0}\Bigr)^{2\sig_0+\sig_1}.\label{eq9.221}
\end{equation}
With $\wtB$ as in Subsection~\ref{sub9.8}, we have $\max (A_3,
\wtA_3) \leqslant C\wtB$.

We observe that in both (\ref{eq9.213}) and (\ref{eq9.219}), our
estimate for $b_+ (\ell)$, $b_- (\ell)$ depends on the class $\ell$
only through $\whI_{s-1}$ and $h$. We first sum over subclasses with
a fixed depth $s$ and sequence $\whI_0 \subset \cdots \subset
\whI_{s-1}$, using the same method of generating series as in the
proof of Proposition~\ref{propo33}. To deal with the two cases at
the same time, we first observe that
\begin{equation}
\sig_0 + \sig_1 - \fud\; \rho_0 (1+\tau) \; = \; 1 - \dou - \fud\;
\dos + o (1) \; > \; 0\label{eq9.222}
\end{equation}
under (H4), and a fortiori $\sig_0 + \sig_1 - \fud\, \rho_1 (1+\tau)
> 0$. Thus, $A_2$ and $\wtA_2$ are larger when $\wtI_{s-1}$ is
larger; the largest case is $\wtI_{s-1} = \Io$, which gives
\begin{equation}
\max (A_2, \wtA_2) \; \leqslant \; \whA_2 \; := \; 2 \Bigl(
C\vep_0^{\fud (1-\tau)} |P_u|\Bigr)^{\rho_1}.\label{eq9.223}
\end{equation}
We, then, set
\begin{eqnarray}
\chi_1 (z) \; &= \; &\sum_{h>0} \whA_2^h z^{A_1 h} \,
(1-z)^{-h-2}\label{eq9.224}\\
&= \; &\whA_2 z^{A_1} \, (1-z)^{-2} (1-z-\whA_2 z^{A_1})^{-1}.\notag
\end{eqnarray}
The (partial) sum of $b_+ (\ell) \, b_- (\ell)$ is, thus, not larger
than $C\wtB$ times the coefficient of $z^{A_0}$ in the power series
for $\chi_1 (z)$. Recall that $A_0$, $A_1$ were defined in
(\ref{eq9.188}), (\ref{eq9.189}).

We estimate this coefficient by Cauchy integration on the circle
$\{|z| = 1 - A_0^{-1} - \whA_2\}$, on which we have
\begin{eqnarray}
&&|1-z|^{-2} \; \leqslant \; (\whA_2 + A^{-1}_0)^{-2} \; \leqslant \; A^2_0,\label{eq9.225}\\
&&|1-z-\whA_2 z^{A_1}|^{-1} \; \leqslant \; A_0,\label{eq9.226}\\
&&|\chi_1 (z)| \; \leqslant \; \whA_2 A_0^3,\label{eq9.227}\\
&&|z^{-A_0}| \; \leqslant \; C (1 + \whA_2)^{A_0}.\label{eq9.228}
\end{eqnarray}
The (partial) sum of $b_+ (\ell) \, b_- (\ell)$ is therefore
dominated by
\begin{equation}
C (1 + \whA_2)^{A_0}\, \whA_2 \, A^3_0\, \wtB\label{eq9.229}
\end{equation}
We now have to sum over sequences $\whI_0 \subset \cdots \subset
\whI_{s-1}$ and depth $s$; but (\ref{eq9.229}) is independent of
these data and the same remarks as in the proof of
Proposition~\ref{propo33} apply. So, we finally obtain for the sum
of $b_+ (\ell)\, b_- (\ell)$ over subclasses with $s>0$ and $h>0$, a
bound by
\begin{equation}
C (1 + \whA_2)^{A_0}\, \whA_2 \, A^3_0\,
2^{S_{\max}}\,\wtB\label{eq9.230}
\end{equation}
with $S_{\max} := 2\tau^{-1} \log \Bigl(
\tfrac{\log\,|\wtI|^{-1}}{\log\,\vep_0^{-1}} \Bigr)$
(cf.~(\ref{eq9.198}).

We have here, from (\ref{eq9.200}), (\ref{eq9.201})
\begin{eqnarray}
A_0 \; &\leqslant \; &C\tau^{-1} \log \,|\wtI|^{-1},\label{eq9.231}\\
(1+\whA_2) A_0 \; &\leqslant \; &|\wtI|^{-C\tau^{-1}
\whA_2},\label{eq9.232}\\
\whA_2 \; &\leqslant \; &\vep_0^{\fudt\,\rho_1}.\label{eq9.233}
\end{eqnarray}
As $\vep_0$ can be made as small as we want with regard to $\tau$,
the term  in (\ref{eq9.230}) is bounded by
\begin{equation}
\vep_0^\sig \, |\wtI|^{-\vep_o^\sig}\, (\log\,
|\wtI|^{-1})^{2\tau^{-1}}\, \wtB\label{eq9.234}
\end{equation}
for some fixed $\sig>0$.

We summarize the calculations in this subsection in

\begin{Propo}\label{propo34}
The quantity $\whB = \sum_L b_+ (\ell) b_- (\ell)$ is bounded by
$\whB_1 + \whB_2 + \whB_3$, with
\begin{eqnarray*}
\whB_1 &= &(\# \cC_+ (\wtI))(\# \cC_- (\wtI))(\log |\wtI|^{-1})^2,\\
\whB_2 &= &\Bigl( \frac{\log |\wtI|^{-1}}{\log
\vep_0^{-1}}\Bigr)^{2\tau^{-1}} (C |\wtI|^{1+\tau} x^{-1})^{\dss}
(\# \cC_+ (\wtI)) \sum_{\cC_- (\wtI)} (\max_{\wtI}
|\Psom|)^{\dss},\\
\whB_3 &= &\vep_0^\sig |\wtI|^{-\vep_0^\sig} (\log
|\wtI|^{-1})^{2\tau^{-1}} \wtB.
\end{eqnarray*}
\end{Propo}

\subsection{End of the Induction Step for (SR3)$_s$\label{sub9.12}}
Of the two inequalities that were assumed in Subsection~\ref{sub9.8}
to make work the argument, the first has been the subject of
Proposition~\ref{propo33}. The second is
\begin{equation*}
\whB |\wtI|^{\tau (2-\dms-\dum-6\tau)} \; \leqslant \; CB\tag{9.144}
\end{equation*}
From Proposition~\ref{propo34}, it is sufficient to prove the same
inequality with $\whB$ replaced by $\whB_i$, $i=1,2,3$.

First consider $\whB_3$. For $x\leqslant \wtx_{cr}$, we have
\begin{equation}
\wtB = \wtB_0 = |\wtI|^{-\tau (2\sig_0+\sig_1})\, B_0 \; \leqslant \;
|\wtI|^{-\tau (2\sig_0+\sig_1)}\, B,\label{eq9.235}
\end{equation}
and, therefore,
\begin{equation}
\whB_3 |\wtI|^{\tau (2-\dms-\dum-6\tau)} \; \leqslant \; \vep_0^\sig
(\log\, |\wtI|^{-1})^{2\tau^{-1}}\,|\wtI|^\om\, B,\label{eq9.236}
\end{equation}
with $\om = \tau (2-\dms-\dum-6\tau-2\sig_0-\sig_1)-\vep_0^\sig$.

We choose the exponents $\sig_0$, $\sig_1$ in order to have
\begin{equation}
\om \; > \; \tau^2\label{eq9.237}
\end{equation}
which means
\begin{equation}
2\sig_0 + \sig_1 \; < \; 2-\dms-\dum-7\tau-\vep_0^\sig
\tau^{-1}\label{eq9.238}
\end{equation}
which is compatible with the previous conditions on $\sig_0$,
$\sig_1$.

As $\vep_0$ can be chosen arbitrarily small with respect to $\tau$,
(\ref{eq9.144}) holds for $\whB_3$.

Next consider $\whB_1$. We assume that the sizes of $\cC_+ (\wtI)$,
$\cC_- (\wtI)$ are controlled by (SR1)$_s$, (SR1)$_u$. Then, we have
\begin{equation}
\whB_1 |\wtI|^{\tau (2-\dms-\dum-6\tau)} \; \leqslant \; C
(\log\,|\wtI|^{-1})^2 \Bigl( \frac{|\wtI|}{\vep_0}
\Bigr)^{2-2\dms-2\dum-2\tau} \vep_0^{-\tau(\dos+\dou)} \,
|\wtI|^{\tau (2-\dms-\dum-6\tau)}.\label{eq9.239}
\end{equation}
The right-hand will be smaller than $CB_0$ as soon as
\begin{equation}
\Bigl( \frac{x}{\vep_0 |P_u|} \Bigr)^{\rho_0}\; \leqslant \; \Bigl(
\frac{|\wtI|}{\vep_0} \Bigr)^{2\sig_0+\sig_1-2+2\dms+2\dum+A\tau}\;
\vep_0^{A\tau}\; (\log\, |\wtI|^{-1})^{-2},\label{eq9.240}
\end{equation}
for some fixed constant $A>0$.

We would like (\ref{eq9.240}) to be a consequence of $x\leqslant
\wtx_{cr}$, but unfortunately this is only true if $\wtI$ is not too
large. Observe that
\begin{equation}
\rho^{-1}_0 (2\sig_0 + \sig_1 - 2 + 2\dms + 2\dum + A\tau) -
(\rho_0-\rho_1)^{-1} \sig_0 = o(1),\label{eq9.241}
\end{equation}
where $(\rho_0-\rho_1)^{-1} \sig_0$ is the exponent appearing in the
definition (\ref{eq9.31}) of $\wtx_{cr}$. Write
\begin{equation}
\sig\; := \; 2\sig_0 + \sig_1 - 2 + 2\dms + 2\dum + A\tau = \dms +
\dum + o(1).\label{eq9.242}
\end{equation}
We choose the exponent $\rho_1$ in order to have
\begin{equation}
\rho_1 > \rho_0 \Bigl( 1 - \frac{\sig_0}{\sig} \Bigr) +
\kappa,\label{eq9.243}
\end{equation}
with $\kappa>0$ small to guarantee that (\ref{eq9.22}) holds. Then, we
have
\begin{equation}
\rho_0^{-1} \sig \; < \; \frac{\sig_0}{\rho_0-\rho_1+\kappa}\; < \;
\frac{\sig_0}{(\rho_0-\rho_1)} - C^{-1} \kappa\label{eq9.244}
\end{equation}
and therefore (\ref{eq9.240}) holds as soon as
\begin{equation}
x\; \leqslant \; \wtx_{cr} \Bigl( \frac{|\wtI|}{\vep_0}
\Bigr)^{-C^{-1}\,\kappa} \vep_0^{A\tau} (\log\,
|\wtI|^{-1})^{-\tfrac{2}{\rho_0}}\label{eq9.245}
\end{equation}
Keeping $\kappa\gg \tau$, the right-hand side is larger than $\wtx_{cr}$
when $|\wtI| < \vep_0^{1+C\kappa^{-1}\tau}$, $C$ fixed large enough.
Thus, we are able to conclude that (\ref{eq9.144}) holds for
$\whB_1$ except in the range
\begin{eqnarray}
\vep_0 \; &\geqslant \; |\wtI| \; &\geqslant \;
\vep_0^{1+C \kappa^{-1}\tau},\label{eq9.246}\\
\wtx_{cr} \; &\geqslant \; x \; &\geqslant \; \wtx_{cr}\,
\vep_0^{2A\tau}.\notag
\end{eqnarray}
We shall deal directly with this case below. Before that, we
consider (\ref{eq9.144}) for $\whB_2$. In this case, we assume that
(SR1)$_s$ and (SR2)$_u$ hold. A small calculation shows that
(\ref{eq9.144}) holds as soon as
\begin{equation}
\Bigl( \frac{x}{\vep_0 |P_u|} \Bigr)^{\rho_0-\dss}\; \leqslant \;
\Bigl( \frac{|\wtI|}{\vep_0} \Bigr)^{\whsig} \;\vep_0^{\whA \tau}\,
\Bigl(\frac{\log\,|\wtI|^{-1}}{\log\,\vep_0^{-1}}\Bigr)^{-2\tau^{-1}}
\label{eq9.247}
\end{equation}
with
\begin{equation}
\whsig\; = \; -2 + 2\dms + \dum - \dss + 2\sig_0 + \sig_1 + \whA
\tau\label{eq9.248}
\end{equation}
and $\whA$ a fixed positive constant. Here, both $\rho_0 - \dss$ and
$\whsig$ are $o(1)$.

With $\kappa$ as above, i.e., $\kappa = o(1)$, $\tau=o(\kappa)$, we choose the
exponents in order to have
\begin{eqnarray}
\rho_0 \; &> \; \dss + \kappa^2,\label{eq9.249}\\
\whsig \; &< \; -2\kappa,\label{eq9.250}
\end{eqnarray}
which corresponds to
\begin{equation}
2\sig_0 + \sig_1 \; < \; 2-2\dms - \dum + \dss - \whA \tau - 2\kappa.
\label{eq9.251}
\end{equation}
If $(\wtI,x)$ is not in the range
\begin{eqnarray}
\vep_0 \; &\geqslant \; |\wtI| \; &\geqslant \; \vep_0^{1+C\kappa^{-1}\tau},\label{eq9.252}\\
\wtx_{cr} \; &\geqslant \; x \; &\geqslant \; \wtx_{cr}
\vep_0^{C\kappa^{-1}\tau},\notag
\end{eqnarray}
then, (\ref{eq9.247}) follows from $x\leqslant \wtx_{cr}$.

The final step in the inductive proof of (SR3)$_s$ is, therefore,
the proof of

\begin{Propo}\label{propo35}
Assume that
\[
|\wtI| \; \leqslant \; \vep_0^{1+C\tau \kappa^{-1}}
\]
and that (SR1)$_\whs$, (SR1)$_\whu$, (SR2)$_\whu$ hold for some
candidate $I\subset \wtI$. Then (SR3)$_s$ holds for $I$ in the range
\[
\wtx_{cr} \; \geqslant \; x \; \geqslant\; \wtx_{cr} \,
\vep_0^{C\tau \kappa^{-1}}.
\]
\end{Propo}

\begin{proof}
Any element $(P,Q,n)$ in $Bi_+ (I,I,I; x)$ satisfies $P\subset P_s$,
$Q\subset Q_u$. Using (H1), this implies, for all $t\in I$:
\begin{equation}
|P| \; \leqslant \; C\vep_0 |P_u|.\label{eq9.253}
\end{equation}
It, then, follows from Corollary~\ref{coro6} in
Subsection~\ref{sub6.7} that for $x\geqslant \wtx_{cr}\,
\vep_0^{C\tau \kappa^{-1}}$, one must have $(P,Q,n)\in \cR (\Io)$.

As $P$ is thin $I$-critical, there exists $(\Psal, \Qsal, \nsal)\in
\whcC_+ (I)$ with $P\subset \Psal$. Similarly, there exists $(\Psom,
\Qsom, \nsom) \in \whcC_- (I)$ with $Q\subset \Qsom$; see
Subsection~\ref{sub9.1} for the definitions of $\whcC_- (I)$,
$\whcC_+ (I)$.

Let us estimate the number of possible $(P,Q,n)$ for fixed $\Psal$,
$\Qsom$.

If $n\leqslant \nsal + \nsom$, there is at most one for each value
of $n$. If $n\geqslant \nsal + \nsom$, we can write (as $(P,Q,n)\in
\cR (\Io)$), $(P,Q,n)$ as a simple composition
\[
(P,Q,n) \; = \; (\Psal, \Qsal, \nsal)\; * \; (P', Q', n') \; * \;
(\Psom, \Qsom, \nsom)
\]
and conclude that there are no more than
\begin{equation}
C \, \left( \frac{x}{|I|\,|\Psom|} \right)^{-\dos}\label{eq9.254}
\end{equation}
possible $P'$ (at this scale, the dependence of dimension on the
parameter is not relevant).

The total number is, therefore, at most
\begin{equation}
C\,\log \Bigl( \frac{\vep_0 |P_u|}{x} \Bigr)\, (\# \whcC_- (I))\,
(\# \whcC_+ (I)) + C \, (\# \whcC_+ (I)) \Bigl( \frac{x}{|I|}
\Bigr)^{-\dss} \sum_{\whcC_-(I)}\, |\Psom|^{\dss}.\label{eq9.255}
\end{equation}
In view of (SR1)$_\whs$, (SR1)$_\whu$, (SR2)$_\whu$, this is not
greater than
\begin{equation}
C \Bigl( \log\,\frac{\vep_0 |P_u|}{x} \Bigr)\, \Bigl(
\frac{|I|}{\vep_0} \Bigr)^{2-2\dms-2\dum-2\tau} + C \, \Bigl(
\frac{|I|}{\vep_0} \Bigr)^{2-2\dms-\dum-2\tau+\dss}\, \Bigl(
\frac{x}{\vep_0 |P_u|}\Bigr)^{-d_s^*}\label{eq9.256}
\end{equation}
and this should be smaller than $CB_0$ with
\begin{equation}
B_0 \; = \; \Bigl( \frac{x}{\vep_0 |P_u|}\Bigr)^{-\rho_0} \, \Bigl(
\frac{|I|}{\vep_0} \Bigr)^{2\sig_0 + \sig_1}\label{eq9.257}
\end{equation}
and
\begin{equation}
\frac{x}{\vep_0 |P_u|} \; \leqslant \; \Bigl( \frac{|I|}{\vep_0}
\Bigr)^{\tfrac{\sig_0}{\rho_0-\rho_1}}.\label{eq9.258}
\end{equation}
The second part of (\ref{eq9.256}) is as required as $\rho_0 >
\dss$ (cf.~(\ref{eq9.65}) and we ask that
\begin{equation}
2\sig_0 + \sig_1 \; < \; 2 - 2\dms - \dum - 2\tau +
\dss.\label{eq9.259}
\end{equation}
In the first part of (\ref{eq9.256}), we bound the logarithmic term
by a small negative power of $\frac{x}{\vep_0 |P_u|}$ and we use
(\ref{eq9.258}) to conclude.

This ends the proof of (SR3)$_s$ in this case.
\end{proof}
The induction step for (SR3)$_s$  is complete. We have proved

\begin{Theo}\label{theo2}
Assume that all parameter intervals which contain $\wtI$ are strongly
regular. Then, all candidates but a proportion not larger than
$|\wtI|^{\tau^2}$ satisfy (SR3)$_s$.
\end{Theo}

\subsection{The Induction Step for (SR3)$_u$\label{sub9.13}}
We have already explained for (SR3)$_u$ the cases where $\Isal$ or
$\Isom$ contains $\wtI$, and the case where $I = \Isal = \Isom$ and
$x$ is large (cf.~Proposition~\ref{propo31}).

The case $I = \Isal = \Isom$, $x$ small which has been treated for
(SR3)$_s$ in the Subsections~\ref{sub9.8}--\ref{sub9.13} is
completely similar for (SR3)$_u$. It is only not completely
symmetric because we have assumed that $\dos \geqslant \dou$ and
thus the formulas for the exponents are not symmetric. So, one has
only to be careful with the inequalities involving the exponents.
For instance, in (\ref{eq9.222}) we had
\[
\sig_0 + \sig_1 - \fud\; \rho_0 (1+\tau)\; > \; 0.
\]
As $\rho'_0 = \frac{\dou}{\dos}\leqslant \rho_0$, we still have
\[
\sig_0 + \sig_1 - \fud\; \rho'_0 (1+\tau)\; > \; 0.
\]
Checking everything in this way is rather tedious, and we leave this
to the reader.

\begin{rem}
At several points, we have asked that the exponents $\rho_0$,
$\rho_1$, $\rho'_0$, $\rho'_1$, $\sig_0$, $\sig_1$ of (SR3)$_s$,
(SR3)$_u$ should satisfy some inequalities; one could worry whether
these inequalities are compatible between themselves (it is easy to
check that each is compatible with (\ref{eq9.21}) through
(\ref{eq9.24}) and (\ref{eq9.28})  through (\ref{eq9.29}). But we
are always bounding $\rho_0$, $\rho_1$, $\rho'_0$, $\rho'_1$ from
below and $\sig_0$, $\sig_1$ from above, hence the compatibility is
obvious.
\end{rem}

\newpage

\setcounter{section}{9}
\setcounter{equation}{0}

\section{The Well-Behaved Part of the Dynamics for Strongly Regular Parameters\label{sec10}}
\subsection{Prime Elements and Prime Decomposition\label{sub10.1}}
In the last two sections, we fix a strongly regular parameter, i.e.
the intersection of a decreasing sequence $(\Ism)_{m\geqslant 0}$ of
strongly regular parameter intervals.

The sequence $\cR (\Ism)$ is increasing and we set
\begin{equation}
\cR \; = \; \bigcup_{m\geqslant 0} \; \cR (\Ism).\label{eq10.1}
\end{equation}

\Def An element $(P,Q,n)\in \cR$ is {\it prime} if $n>0$ and it
cannot be written as a {\it simple} composition of two shorter
elements.

Obviously, for any  $(a,a')\in\cB$, the element $(P_{aa'}, Q_{aa'},
1)$ is prime. Such elements are called {\it trivial} primes. Non
trivial primes are those of length bigger than~1.

There are only finitely many trivial primes. On the other hand,
there are typically countably many non trivial ones.

\begin{Propo}\label{propo36}
Any element $(P,Q,n)\in \cR$ with $n>0$ can be uniquely written as a
simple composition of a finite sequence of prime elements.
\end{Propo}

\begin{proof}
The existence of such a decomposition is clear. We have to show it
is unique. Assume on the opposite that we can write
\begin{eqnarray}
(P,Q,n) &= &(\Psu,\Qsu,\nsu)\; * \cdots * \; (\Psr, \Qsr,
\nsr)\label{eq10.2}\\
&= &(\Psul,\Qsul,\nsul)\; * \cdots * \; (P'_s , Q'_s , n'_s)\notag
\end{eqnarray}
It is sufficient to show that $(\Psu,\Qsu,\nsu) =
(\Psul,\Qsul,\nsul)$. This is true if $\nsu = \nsul$. Assume for
instance that $\nsu < \nsul$. Then we have $P\subset \Psul \subset
\Psu$ with $\Psul \ne \Psu$. By Proposition~\ref{propo8} in
Subsection~\ref{sub6.5}, we can write
\begin{equation}
(\Psul,\Qsul,\nsul) \; = \; (\Psu, \Qsu, \nsu) \; * \; (\wtP, \wtQ,
\wtn)\label{eq10.3}
\end{equation}
which contradicts the fact that $(\Psul,\Qsul,\nsul)$ is prime.
\end{proof}

\begin{rem}
In the prime decomposition
\begin{equation}
(P,Q,n) \; = \; (\Psu, \Qsu, \nsu) \; * \; (\Psr, \Qsr,
\nsr),\label{eq10.4}
\end{equation}
$\Psu$ can be characterized as the thinnest prime rectangle
containing $P$.
\end{rem}
We will denote by $\cP$ the set of prime elements of $\cR$. We
denote by $\cRs$ the set of elements of $\cR$ of length $>0$.

Let $(P,Q,n)$ be an element of $\cRs$ and let
\begin{equation}
(P,Q,n) \; = \; (\Psu, \Qsu, \nsu) \; * \cdots * \; (\Psr, \Qsr,
\nsr),\label{eq10.5}
\end{equation}
be its prime decomposition. We define
\begin{eqnarray}
\Tim ((P,Q,n)) \; &= \; &(\Psd,\Qsd,\nsd) \; * \cdots * \; (\Psr,
\Qsr,\nsr),\label{eq10.6}\\
\TiM ((P,Q,n)) \; &= \; &(\Psu,\Qsu,\nsu) \; * \cdots * \; (\Psru,
\Qsru,\nsru),\notag
\end{eqnarray}
if $r>1$. When $(P,Q,n)$ is prime, with $P\subset R_a$ and $Q\subset
R_{a'}$, we set
\begin{eqnarray}
\Tim ((P,Q,n)) \; &= \; &(R_{a'},R_{a'},0)\label{eq10.7}\\
\TiM ((P,Q,n)) \; &= \; &(R_a,R_a,0).\notag
\end{eqnarray}
For $S = (P,Q,n)\in \cR$, we write $S * \cR$, resp.~$\cR * S$, for
the set of elements which can be written as $(P,Q,n) * (P',Q',n')$,
resp.~$(P',Q',n') * (P,Q,n)$, for some $(P',Q',n')\in \cR$. We have
partitions
\begin{equation}
\cRs \; = \; \bigsqcup_\cP \; S * \cR \; = \; \bigsqcup_\cP \; \cR *
S. \label{eq10.8}
\end{equation}
Moreover, for any $S\in \cP$, the restriction of $\Tim$,
resp.~$\TiM$, to $S*\cR$, resp.~$\cR * S$, is a bijection onto
$\cR$, inverse of $S' \mapsto S * S'$, resp.~$S' \mapsto S' * S$.

\subsection{Number of Factors in a Prime Decomposition\label{sub10.2}}
We write $r(S)$ for the number of factors in the prime decomposition
of an element $S$ of $\cR$ (setting $r(S) = 0$ if $S$ has length
$0$). Let $(P,Q,n)$, $(P',Q',n')$ be elements of $\cR$ such that
$P'$ is a child of $P$. When $P'$ is a simple child, it is obtained
by simple composition of $P$ with an element of length 1 and we have
\begin{equation}
r (P',Q',n') \; = \; r (P,Q,n) + 1.\label{eq10.9}
\end{equation}

\begin{Propo}\label{propo37}
If $P'$ is a non-simple child of $P$, we have
\[
r (P',Q',n') \; \leqslant \; r (P,Q,n)
\]
\end{Propo}

\begin{proof}
Let
\begin{equation}
(P,Q,n) \; = \; (\Psu,\Qsu,\nsu) * \cdots * (\Psr, \Qsr,
\nsr)\label{eq10.10}
\end{equation}
be the prime decomposition of $(P,Q,n)$.

Let $(\whP,\whQ,\whn)$ be the element of $R$ such that
\begin{equation}
(P',Q',n') \; \in \; (P,Q,n) \; \square \; (\whP, \whQ,
\whn),\label{eq10.11}
\end{equation}
(cf.~Proposition~\ref{propo5} in Subsection~\ref{sub6.4}). There
exists $m\geqslant 0$ such that $Q$ and $\whP$ are
$\Ism$-transverse. Define
\begin{eqnarray}
(P^i,Q^i,n^i) \; &= \; &(P_i,Q_i,n_i) * \cdots * (\Psr, \Qsr,
\nsr)\label{eq10.12}\\
&= \; &(\Tim)^{i-1} \, (P,Q,n), \; 1\leqslant i\leqslant r.\notag
\end{eqnarray}
We have an increasing sequence
\begin{equation}
Q = Q^1 \subset Q^2 \subset \cdots \subset Q^r \; =
Q_r.\label{eq10.13}
\end{equation}
Let $\rl$ be the largest integer in $\{1,\cdots,r\}$ such that
$Q^\rl$ and $\whP$ are $\Ism$-transverse for some $m\geqslant 0$
(and then for all large enough $m$). Define $(\wtP, \wtQ, \wtn) \in
\cR$ by the condition $Q' \subset \wtQ$ and
\begin{equation}
(\wtP, \wtQ, \wtn) \; \in \; (P^\rl, Q^\rl, n^\rl) \; \square \;
(\whP, \whQ, \whn).\label{eq10.14}
\end{equation}
We then have
\begin{equation}
(P', Q', n') \; = \; (\Psu, \Qsu, \nsu) \; * \cdots * \; (P_{r'-1},
Q_{r'-1}, n_{r'-1}) \; * \; (\wtP, \wtQ, \wtn).\label{eq10.15}
\end{equation}

The assertion of the proposition is, thus, a consequence of
\begin{lema}\label{lema12}
$(\wtP, \wtQ, \wtn)$ is prime.
\end{lema}

{\it{Proof.}} Assume by contradiction that we can write
\begin{equation}
(\wtP, \wtQ, \wtn) \; = \; (\wtP_1, \wtQ_1, \wtn_1) \; * \; (\wtP_2,
\wtQ_2, \wtn_2)\label{eq10.16}
\end{equation}
with $\wtn_1$, $\wtn_2 > 0$. Define
\begin{equation}
(\Psul, \Qsul, \nsul) \; = \; (\Psu,\Qsu,\nsu) \; * \cdots * \;
(P_{r'-1}, Q_{r'-1}, n_{r'-1}) * (\wtP_1, \wtQ_1,
\wtn_1).\label{eq10.17}
\end{equation}
We have $P'\subset \Psul$, $P'\ne \Psul$, hence $P\subset \Psul$. We
also have
\begin{equation}
(P', Q', n') \; = \; (\Psul,\Qsul,\nsul) \; * \; (\wtP_2, \wtQ_2,
\wtn_2).\label{eq10.18}
\end{equation}
By Proposition~\ref{propo8} in Subsection~\ref{sub6.5}, there exists
$(\Ps, \Qs, \ns)$ such that
\begin{equation}
(P, Q, n) \; = \; (\Psul,\Qsul,\nsul) \; * \; (\Ps, \Qs,
\ns).\label{eq10.19}
\end{equation}
We claim that there exists $j\in \{1,\cdots,r\}$ such that
\begin{equation}
(\Psul, \Qsul, \nsul) \; = \; (\Psu,\Qsu,\nsu) \; * \cdots * \;
(\Psj, \Qsj, \nsj).\label{eq10.20}
\end{equation}
Indeed, if $j$ is the smallest integer such that
\begin{equation}
n'_1 \; \leqslant \; n_1 + \cdots + n_j,\label{eq10.21}
\end{equation}
it follows from Proposition~\ref{propo8} that we can write
\begin{equation}
(\Psu, \Qsu, \nsu) \; * \cdots * \; (\Psj, \Qsj, \nsj) = (\Psul,
\Qsul, \nsul)\; * \; (\ovP, \ovQ, \ovn)\label{eq10.22}
\end{equation}
for some $(\ovP,\ovQ,\ovn)\in \cR$. By the same proposition, we can
also write
\begin{equation}
(\Psj, \Qsj, \nsj) \; = \; (\ovP', \ovQ', \ovn') \; * \; (\ovP,
\ovQ, \ovn)\label{eq10.23}
\end{equation}
for some $(\ovP', \ovQ', \ovn')\in \cR$. We have $\ovn < n_j$ by the
definition of $j$, hence $\ovn' > 0$. As $(\Psj, \Qsj, \nsj)$ is
prime, we have $\ovn = 0$ which implies our claim.

The integer $j$ satisfies $r'\leqslant j \leqslant r$. We have
\begin{equation}
(\wtP_2, \wtQ_2, \wtn_2) \; \in \; (\jim P, \jim Q, \jim n) \;
\square \; (\whP, \whQ, \whn),\label{eq10.24}
\end{equation}
which contradicts the definition of $r'$.

This concludes the proof of the lemma and also of the proposition.
\end{proof}

\subsection{A Weighted Estimate on the Number of Children\label{sub10.3}}
We present in this subsection a variation over the estimates in
Subsection~\ref{sub8.1}, which will be important in the definition
of a transfer operator.

We fix a constant $\ka\in (0,1)$ close to 1, but independent of
$\vep_0$. We set
\begin{equation}
\dMs = \dos - C\vep_0,\label{eq10.25}
\end{equation}
with a constant $C$ sufficiently large so $\dMs$ is smaller than the
transverse Hausdorff dimension of the stable foliation $W^s(K)$ for
the parameter that we are dealing with.

For $S = (P,Q,n)$, we set
\begin{equation}
||P|| \; = \; |P|^{\dMs}\; \ka^{r(S)},\label{eq10.26}
\end{equation}
(we will also write $r(P)$ instead of $r(S)$).

\begin{Propo}\label{propo38}
For any $m\geqslant 1$, any $(P,Q,n)\in \cR$, we have
\[
\sum_{P'} ||P'|| \; \leqslant \; C\, \ka^{\tfrac{m}{2}} \, ||P||
\]
where the sum in the left-hand side is over elements $(P',Q',n')$
such that $P'$ is a descendent of the $m^{th}$ generation of $P$.
\end{Propo}
We will first state a Lemma, then prove the proposition from the
Lemma, and finally prove the Lemma.

\begin{lema}\label{lema13}
Let $\vep_1 > 0$. If $\vep_0$ is small enough, we have
\[
\sum_{P'} ||P'|| \; \leqslant \; \vep_1 \, ||P||
\]
for all $(P,Q,n)\in \cR$, where the sum in the left-hand side is
over non-simple children of $P$.
\end{lema}

\medskip\noindent{\it Proof of the Proposition.}

Let $m_0\geqslant 1$ be an integer to be determined later. Consider
all chains
\begin{equation}
P = P^0 \supset P^1 \supset \cdots \supset P^{m_0} =
P'\label{eq10.27}
\end{equation}
where $P$ is given and $\iim P$ is a child of $P^i$. If $\iim P$ is,
for each $i$, a simple child of $P^i$, one has $r(P') = m_0 + r (P)$
and the corresponding part of the sum in Proposition~\ref{propo38}
satisfies
\begin{equation}
\sum \; ||P'|| \; \leqslant \; C\, \ka^{m_0}\, ||P||\label{eq10.28}
\end{equation}
(as long as $m_0 = o (\vep_0^{-1})$). We choose $m_0$ such that in
(\ref{eq10.28}) we have
From the lemma above, it follows that for every $\wtP$, we have
\begin{equation}
2C\, \ka^{m_0} \; \leqslant \; \ka^{\tfrac{m_0}{2}}.\label{eq10.29}
\end{equation}
\begin{equation}
\sum \; ||\wtP'|| \; \leqslant \; C\, ||\wtP||\label{eq10.30}
\end{equation}
where the sum is over all children of $\wtP$. Using the lemma again,
when we sum over chains such that $\ism P$ is a non-simple child of
$P_i$ for some $i$, we obtain
\begin{equation}
\sum \; ||\wtP'|| \; \leqslant \; m_0\, C^{m_0-1}\, \vep_1
||P||.\label{eq10.31}
\end{equation}
Taking $\vep_1$ small enough, we obtain
\begin{equation}
\sum \; ||P'|| \; \leqslant \; \ka^{\tfrac{m_0}{2}}\,
||P||\label{eq10.32}
\end{equation}
where the sum is now over all chains. The proposition follows
immediately from (\ref{eq10.32}) and (\ref{eq10.30}).

\medskip\noindent
{\it Proof of Lemma~\ref{lema13}.} Let $(P,Q,n)\in \cR$. Any
non-simple child $P'$ of $P$ is obtained as
\begin{equation}
(P',Q',n')\; \in \; (P,Q,n) \; \square \; (\Psu, \Qsu, \nsu)
\label{eq10.33}
\end{equation}
and we denote by $\wtP_1$, the parent of $P_1$.  One has
\begin{equation}
|P'| \; \leqslant \; C \, |P|\; |\Psu| \, \de (Q, \Psu)^{-\fudt}.
\label{eq10.34}
\end{equation}
Therefore, we will have
\begin{equation}
||P||^{-1} \sum \, ||P'|| \; \leqslant \; C \, \sum |P_1|^{\dMs}\;
\ka^{r(P')-r(P)} \de^{-\fudt\,\dMs}.\label{eq10.35}
\end{equation}
From the proof of Proposition~\ref{propo37}, there is an increasing
sequence
\begin{equation}
Q \; = \; Q^1 \subset Q^2 \subset \cdots \subset
Q^{r(P)}\label{eq10.36}
\end{equation}
such that $r (P')$ is the largest integer for which $Q^r$ and $P_1$
are $\Ism$-transverse for large enough $m$. If
\begin{equation}
|Q^r| \; \leqslant \; \de (Q, P_1)^2\label{eq10.37}
\end{equation}
then, by Proposition~\ref{propo10}, $Q^r$ and $P_1$ are
$\Ism$-transverse for large $m$ and thus $r (P')\geqslant r$. On the
other hand, there exists $\ka^* \in (0,1)$ such that
\begin{equation}
|Q^r| \; \leqslant \; \ka^* \, |Q^{r+1}|\label{eq10.38}
\end{equation}
for $r< r(P)$. We infer that
\begin{equation}
r(P) - r(P') \; \leqslant \; C\,\log (\de
(Q,P_1))^{-1}.\label{eq10.39}
\end{equation}
Therefore, if $\ka$ is close enough to 1, we have
\begin{equation}
\ka^{r(P')-r(P)} \; \leqslant \; (\de
(Q,P_1))^{-\tfrac{1}{6}\,\dMs}\label{eq10.40}
\end{equation}
and the right-hand side of (\ref{eq10.35}) is bounded by
\begin{equation}
C\; \sum \; |P_1|^{\dMs}\, (\de
(Q,P_1))^{-\tfrac{2}{3}\,\dMs}\label{eq10.41}
\end{equation}
Using (R7), this is smaller than
\begin{equation}
C\; \sum \; |P_1|^{\tfrac{1}{3}\,\dMs}.\label{eq10.42}
\end{equation}
To estimate this sum, we first fix the parent $\wtP_1$ and sum over
children $P_1$; it follows from Proposition~\ref{propo21} that the
corresponding sum is bounded by $C |\wtP_1|^{\frac{1}{3}\,\dMs}$ and
(\ref{eq10.42}) is not greater than
\begin{equation}
C\; \sum \; |\wtP_1|^{\tfrac{1}{3}\,\dMs}.\label{eq10.43}
\end{equation}
For each integer $m$, let us count now how many $\wtP_1$ may satisfy
\begin{equation}
2^{-m} \, |P_s|\; \geqslant \; |\wtP_1| \; \geqslant \; 2^{-m-1} \;
|P_s|.\label{eq10.44}
\end{equation}
As $Q$ is transverse to $P_1$ but not to $\wtP_1$, we must have, by
Proposition~\ref{propo10}:
\begin{equation}
\de (Q, \wtP_1) \; \leqslant \; C\, 2^{-m(1-\eta)}\label{eq10.45}
\end{equation}
which shows that there are no more than $C\,2^{m\eta}$ such
$\wtP_1$'s. This implies that the sum (\ref{eq10.43}) is at most of
order $C\,\vep_0^{\frac{1}{3}\,\dMs}$, which yields the statement of
the lemma.\hfill $\square$

\medskip\noindent
{\bf Remark.}

{\it 1. In the lemma, the value $\dMs = \dos - C\,\vep_0$ that has
been used to define $||P||$ is irrelevant. The assertion of the
lemma is still true if we replace $\dMs$ by any positive number
bounded away from $0$.

2. In the proposition, we can replace $\dMs$ by a slightly lower
value (assuming, as usual, that $\vep_0$ is arbitrarily small): if
we take
\begin{equation}
\whd_s^- \; = \; \dos - o (\log\, \ka^{-1}),\label{eq10.46}
\end{equation}
the same argument works and the result of the proposition is still
valid.}

\begin{Coro}\label{coro14}
Let $\vep_1 > 0$. If $\vep_0$ is small enough, we have
\[
\sum \, n |P|^{\dMs}\; < \; \vep_1
\]
where the sum in the left-hand side is over non-trivial primes
$(P,Q,n)$.
\end{Coro}

\begin{proof}
Let $(P,Q,n)$ be a non trivial prime. We have, by
Proposition~\ref{propo13} of Subsection~\ref{sub7.1}
\begin{equation}
n \; \leqslant \; \Bigl( \log\,(C |P|^{-1})
\Bigr)^{\tfrac{\log\,2}{\log 3/2}}.\label{eq10.47}
\end{equation}
Consequently, we have
\begin{equation}
n \, |P|^{\dMs} \; \leqslant \; |P|^{\whd^-_s},\label{eq10.48}
\end{equation}
with $\whd^-_s$ as in (\ref{eq10.46}) and $\dos - \whd_s^-$ being
independent of $\vep_0$ (as $P$ is a non trivial prime, one has
$|P|<\vep_0$). Observe also that the thinnest $(\wtP, \wtQ, \wtn)
\in \cR (\Io)$ with $P\subset \wtP$ satisfies $\wtQ \subset Q_u$
hence $|\wtP| \leqslant \vep_0^\al$ for some fixed positive $\al$.

We apply the proposition, taking the Remark~2 above into account and
using $\whd^-_s$ instead of $\dMs$; we obtain
\begin{eqnarray}
\sum \, n |P|^{\dMs} \; &\leqslant \; &\sum
|P|^{\whd^-_s}\label{eq10.49}\\
&= \; &\ka^{-1} \, \sum ||P||\notag\\
&\leqslant \; &C\,\ka^{-1} (1-\ka^{\fud})^{-1} \, \sum
||\wtP||,\notag
\end{eqnarray}
where, in the last sum, $(\wtP, \wtQ, \wtn)$ runs through the
elements of $\cR (\Io)$ with $|\wtP|$ of the order of $\vep_0^\al$.
We have
\begin{eqnarray}
\sum \, ||\wtP|| \; &\leqslant \; &C\, \sum
\; \vep_0^{\al \whd^-_s}\, \ka^{C^{-1}\,\log\, \vep_0^{-1}}\label{eq10.50}\\
&\leqslant \; &C\,\vep_0^{\al(\dMs-\dos) + C^{-1}\,\log\,
\ka^{-1}},\notag
\end{eqnarray}
and the exponent is positive from (\ref{eq10.46}). Putting this into
(\ref{eq10.49}) yields Corollary~\ref{coro14}.
\end{proof}

\begin{rem}
We do not know the range of $d$ for which
\[
\sum_\cP\, |P|^d
\]
is convergent. The Corollary shows at least that there are
relatively few primes: the sum
\[
\sum_{\cR (\Io)} \, |P|^{\dMs}\;\; \text{ is convergent.}
\]
\end{rem}

\subsection{Stable Curves\label{sub10.4}}
Let $(\Psk, \Qsk, \nsk)_{k\geqslant 0}$ be a sequence of elements of
$\cR$ such that $\ksm P$ is strictly contained in $\Psk$ for
$k\geqslant 0$. Let $R_a$ be the rectangle of the Markov partition
which contains $P_0$. The vertical part of the boundary of $\Psk$ is
the union of two graphs $\{ x_a = \vphi_k^\pm (y_a)\}$, and the
$\vphi_k^\pm$ are uniformly bounded in the $C^2$ topology. Moreover,
there exists $\kas\in (0,1)$ such that $|\ksm P|\leqslant \kas
|\Psk|$ for all $k\geqslant 0$. It follows that both sequences
$\vphi_k^\pm$ converge in the $C^1$ topology to the same limit
$\vphi_\infty$, which is of class $C^{1+\lip}$, where $\lip$ stands
for Lipschitz. We state this as

\begin{Propo}\label{propo39}
The intersection $\bigcap\limits_{k\geqslant 0} \Psk$ is the graph
$\{ x_a = \vphi_\infty (y_a)\}$ of a $C^{1+\lip}$ function.
Moreover, the $C^{1+\lip}$ norm of $\vphi_\infty$ is bounded
independently of the sequence $(\Psk)_{k\geqslant 0}$.
\end{Propo}

\Def

1.~~A {\it stable curve} is the intersection $\om =
\bigcap\limits_{k\geqslant 0} \Psk$ of a decreasing sequence of
vertical-like rectangles as above. An {\it unstable curve} is the
intersection $\om' = \bigcap\limits_{k\geqslant 0} \Qskl$ of a
decreasing sequence of horizontal-like strips.

2.~~The {\it set} of stable curves, resp.~unstable curves, is
denoted by $\cR^\infty_+$, resp.~$\cR^\infty_-$. The {\it union} of
stable curves, resp.~unstable curves, is denoted by
$\wtcR^\infty_+$, resp.~$\wtcR^\infty_-$.

3.~~Any stable curve $\om \subset R_a$ has a {\it canonical defining
sequence} characterized by the following conditions: $P_0 = R_a$
and, for each $k$, $\ksm P$ is a child of $\Psk$.

4.~~Two stable curves are equal or disjoint. Hence there is a {\it
canonical projection}
\[
\pi: \; \wtcR^\infty_+ \; \mapsto \; \cR^\infty_+.
\]
We will now define dynamics on a part of the sets $\cR^\infty_+$,
$\wtcR^\infty_+$.

Let $\cN_+$ be the set of stable curves $\om$ which are contained in
infinitely many prime elements and let $\cD_+$ be the complementary
subset in $\cR^\infty_+$. For $(P,Q,n)\in \cP$, denote by
$\cR^\infty_+ (P)$ the set of stable curves $\om \in \cD_+$ such
that $P$ is the thinnest prime containing $\om$.

We, thus, have partitions
\begin{eqnarray}
\cR^\infty_+ \; &= \; &\cN_+ \; \bigsqcup \; \cD_+,\label{eq10.51}\\
\cD_+        \; &= \; &\bigsqcup_\cP \; \cR_+^\infty
(P).\label{eq10.52}
\end{eqnarray}
We denote by $\wtcN_+$, $\wtcD_+$, $\wtcR^\infty_+ (P)$ the
respective pre-images by $\pi$.

Let $(P,Q,n)\in \cP$, $\om\in \cR^\infty_+ (P)$. For any $(\Psk,
\Qsk, \nsk)$ with $\om \subset \Psk \subset P$, we can write
(cf.~Remark after Proposition~\ref{propo36})
\begin{equation}
(\Psk, \Qsk, \nsk) = (P,Q,n) \; * \; (\Pskl, \Qskl, \nskl)
\label{eq10.53}
\end{equation}
for some $(\Pskl, \Qskl, \nskl)\in \cR$; we have
\begin{equation}
\Tim (\Psk, \Qsk, \nsk)\; = \; (\Pskl, \Qskl, \nskl)\label{eq10.54}
\end{equation}
and we define $\om' = \Tim (\om)$ to be the stable curve obtained by
the intersection of the $\Pskl$ when $\Psk$ decrease to $\om$.  We
have
\begin{equation}
g^n (\Psk) \subset \Pskl,\qquad g^n (\om) \subset
\om'\label{eq10.55}
\end{equation}
and we also define
\begin{equation}
\wtT^+ \left / \right. \wtcR^\infty_+ (P) = g^n \left / \right.
\wtcR^\infty_+ (P). \label{eq10.56}
\end{equation}
We, thus, have a commutative diagram
\begin{equation}
\begin{CD}
\wtcD_+      @>\wtT^+>>    \wtcR^\infty_+\\
@VV{\pi}V                  @VV{\pi}V\\
\cD_+        @>\Tim>>      \cR^\infty_+
\end{CD} \label{eq10.57}
\end{equation}
We observe that for $(P,Q,n)\in \cP$ with $Q\subset R_a$, the image
$\Tim (\om)$ of any $\om\in \cR^\infty_+ (P)$ is contained in $R_a$.

Conversely, let $(P,Q,n)\in \cP$ with $Q\subset R_a$ and let
$\om'\in \cR^\infty_+$, $\om'\subset R_a$. For any $(\Pskl, \Qskl,
\nskl)$ with $\om\subset \Pskl$, we define $(\Psk, \Qsk, \nsk)$ by
(\ref{eq10.53}); the intersection $\om$ of the $\Pskl$s, when
$\Pskl$ decrease to $\om'$, is the unique stable curve in
$\cR^\infty_+ (P)$ such that $\Tim (\om) = \om'$.

Thus, $\Tim$ induces a bijection from $\cR^\infty_+ (P)$ on the set
$\cR^\infty_+ (a)$ of stable curves contained in $R_a$. For $\om\in
\cR^\infty_+ (P)$, we have
\begin{equation}
\wtT^+ (\om) \; = \; \om' \cap Q.\label{eq10.58}
\end{equation}

\subsection{Topology and Geometry of $\cR^\infty_+$ and $\wtcR_+^\infty$\label{sub10.5}}
Each stable curve is a compact subset of $R = \cup R_a$. Therefore,
$\cR^\infty_+$ may be viewed as a subset of the set of non empty
compact subsets of $R$ endowed with the Hausdorff topology. The
topology induced on $\cR^\infty_+$ can also be viewed directly: for
any $\om = \cap \Psk$ in $\cR^\infty_+$, a basis of neighbourhoods
of $\om$ is obtained by considering for each $k$ the set $\Vsk$ of
stable curves contained in $\Psk$.

Equipped with this topology, $\cR^\infty_+$ is a Cantor set. Each
$\cR^\infty_+ (P)$, $P\in \cP$, is a closed subset, and also a
Cantor set. The restriction of $\Tim$ to each $\cR^\infty_+ (P)$ is
a homeomorphism onto $\cR^\infty_+ (a)$ (with $Q\subset R_a$).

However, the subset $\cN_+$ may be dense and the map $\Tim$ in
general is not continuous on the whole of $\cD_+$. We will see in
the sequel that $\cN_+$ is, in some appropriate sense, negligible.

For each $\om \in \cR^\infty_+ (a)$, we denote by $\vphi_\om$ the
$C^{1+\lip}$ map such that $\om = \{ x_a = \vphi_\om (y_a)\}$; for
each $a\in \AL$, each $y^0_a\in I^u_a$, the map
\begin{eqnarray}
\phi_{y^0_a} \; : \; &\cR^\infty_+ (a) \; \mapsto \;
I^s_a\label{eq10.59}\\
&\om \; \mapsto \; \vphi_\om (y^0_a)\notag
\end{eqnarray}
is a homeomorphism onto its image. Letting $y^0_a$ vary, we get an
homeomorphism from $\cR^\infty_+ (a) \times I^u_a$ onto
$\wtcR^\infty_+ (a)$.

The regularity of the partial foliation $\wtcR^\infty_+ (a)$ is
given by

\begin{Propo}\label{propo40}
There exists $C>0$ such that, for all $a\in \AL$, all distinct
$\om$, $\om'\in \cR^\infty_+ (a)$, all $y$, $y'\in I^u_a$,
(cf.~Subsection~\ref{sub2.1}) we have
\begin{equation}
\left\vert \log\; \frac{\vphi_\om(y)-\vphi_{\om'}(y)}{\vphi_\om (y')
- \vphi_{\om'} (y')} \right\vert \; \leqslant \; C\, |y-y'|.
\label{eq10.60}
\end{equation}
In particular, the homeomorphisms $\phi_{y'}\circ \phi^{-1}_y$ are
bi-Lipschitzian, uniformly in $y$, $y'$.
\end{Propo}

\begin{proof}
The calculations that support the proof will be found in Appendix~B.
Write $\om = \cap \Psk$, $\om' = \cap \Pskl$, where
$(\Psk)_{k\geqslant 0}$, $(\Pskl)_{k\geqslant 0}$ are the canonical
sequences associated with $\om$, $\om'$. We write $\vphi_k$,
$\vphi'_k$ for $\vphi^+_k$, $\vphi^{'+}_k$. Let $\ell$ be the
largest integer such that $P_\ell = P'_\ell$. There are two cases.

\medskip\noindent
{\bf Case~1.} $P_{\ell+1}$ or $P'_{\ell+1}$ is a simple child of
$P_\ell$. Let $(A,B)$ be the implicit representation for $(P_\ell,
Q_\ell, n_\ell)$. We have in this case, for all $y\in I^u_a$
\begin{equation}
C^{-1} |P_\ell| \; \leqslant \; |\vphi_\om (y) - \vphi_{\om'} (y)|\;
\leqslant \; C |P_\ell|.\label{eq10.61}
\end{equation}
To prove (\ref{eq10.60}), we observe (cf.~Appendix~B) that we can
write, for $k>\ell$,
\begin{eqnarray}
\vphi_k (y) \; &= \; A (y, \psi_k (y)),\label{eq10.62}\\
\vphi'_k (y) \; &= \; A (y, \psi'_k (y)),\notag
\end{eqnarray}
with $\psi_k$, $\psi'_k$ uniformly $C^1$-bounded and
\begin{equation}
|\psi_k (y) - \psi'_k (y)| \; \geqslant \; C > 0.\label{eq10.63}
\end{equation}
We then have,
\begin{equation}
\vphi'_k (y) - \vphi_k (y) = \int^1_0 A_x \Bigl(y, (1-t) \psi_k (y)
+ t\, \psi'_k (y)\Bigr) dt [\psi'_k (y) - \psi_k
(y)].\label{eq10.64}
\end{equation}
From (\ref{eq10.64}) we deduce (\ref{eq10.60}) with $\vphi_\om$,
$\vphi_{\om'}$ replaced by $\vphi_k$, $\vphi'_k$. We then let $k$
goes to $+\infty$ to obtain (\ref{eq10.60}) for $\vphi_\om$,
$\vphi_{\om'}$.

\medskip\noindent
{\bf Case~2.} $P_{\ell+1}$ and $P'_{\ell+1}$ are non-simple children
of $P_\ell$. In this case (\ref{eq10.61}) does not hold. The
computation in Appendix~B shows that
\begin{equation}
\left\vert \log\; \frac{\vphi_{\ell+1} (y) - \vphi'_{\ell+1}(y)}
{\vphi_{\ell+1} (y') - \vphi'_{\ell+1} (y')} \right\vert \;
\leqslant \; C\, |y-y'|. \label{eq10.65}
\end{equation}
By the same proof as in Case~1, we have
\begin{equation}
\left\vert \log\; \frac{\vphi_{\ell+1} (y) - \vphi_\om (y)}
{\vphi_{\ell+1} (y') - \vphi_\om (y')} \right\vert \; \leqslant \;
C\, |y-y'|,\label{eq10.66}
\end{equation}
\begin{equation}
\left\vert \log\; \frac{\vphi'_{\ell+1} (y) - \vphi_\om' (y)}
{\vphi'_{\ell+1} (y') - \vphi_\om' (y')} \right\vert \; \leqslant \;
C\, |y-y'|.\label{eq10.67}
\end{equation}
As we have
\begin{eqnarray}
|\vphi_{\ell+1}(y) - \vphi'_{\ell+1} (y)| \; &\leqslant \; &C
|\vphi_\om (y) - \vphi_{\om'} (y)|,\label{eq10.68}\\
|\vphi_{\ell+1}(y) - \vphi_\om (y)| \; &\leqslant \; &C
|\vphi_\om (y) - \vphi_{\om'} (y)|,\label{eq10.69}\\
|\vphi'_{\ell+1}(y) - \vphi_{\om'} (y)| \; &\leqslant \; &C
|\vphi_\om (y) - \vphi_{\om'} (y)|,\label{eq10.70}
\end{eqnarray}
the inequality in the proposition follows.
\end{proof}

The result of Proposition~\ref{propo40} implies that the transverse
Hausdorff dimension $d_s = d_s (g)$ of $\qqim\wtcR$ is well-defined,
being equal to the Hausdorff dimension of $\phi_y (\qqim\cR (a))$
for any $a\in \AL$, $y\in I^u_a$. We have just proved that it does
not depend on $y$. That it does not depend on $a$ is seen as
follows: for $(a,a')\in \cB$, $g$ sends $\qqim\wtR \cap P_{aa'}$
into $\qqim\wtR \cap R_{a'}$; the transverse Hausdorff dimension of
$\qqim\wtR \cap R_{a'}$ is therefore not smaller than that of
$\qqim\wtR \cap R_a$; as this is true for all $(a,a')\in \cB$, the
conclusion follows.

We will also identify below in this section the transverse Hausdorff
dimension $d_s$ through a transfer operator in the classical manner
of Bowen, Ruelle and Sinai.

Another partial control on the geometry of $\qqim\wtR$ is obtained
via the usual estimate on the graph transform, in the setting of
uniformly hyperbolic dynamics. Let $\om$, $\whom$ be two stable
curves and $j>0$. Assume that $\om$ and $\whom$ belong to the domain
of $(\Tim)^j$. We say that $\om$ and $\whom$ belong to the same
component of the domain of $(\Tim)^j$ if, for each $0\leqslant i <
j$, there exists a prime $P_i$ such that $(\Tim)^i (\om)$ and
$(\Tim)^i (\whom)$ belong to $\qqim\cR (P_i)$.

\begin{Propo}\label{propo41}
There exists $\te_0 \in (0,1)$ such that, if $\om$, $\whom$ belong
to the same component of the domain of $(\Tim)^j$, we have
\[
|D\vphi_\om (y) - D\vphi_\whom (y)| \; \leqslant \; C\, \te_0^j
\]
for all $y$.
\end{Propo}

\begin{proof}
We leave the proof of this very standard estimate to the reader.
\end{proof}

\subsection{Transverse Dilatation\label{sub10.6}}
This subsection is a preparation for the definition of a transfer
operator in the next subsection. The weight function in this
transfer operator is, up to a coboundary term, given by a transverse
dilatation.

Let $(P,Q,n)\in \cP$, $\om = \{ x = \vphi_\om (y)\}$ a stable curve
in $\cRim (P)$, $\om' = \Tim (\om)$ its image. Let $(A,B)$ be the
implicit representation of $(P,Q,n)$.

For $z = (\vphi_\om (y), y)\in \om$, let
\begin{equation}
v_\om (z) = \frac{\pa}{\pa y} + D\vphi_\om (y) \; \frac{\pa}{\pa
x}\label{eq10.71}
\end{equation}
be the normalized tangent vector to $\om$ at $z$.

The matrix of $D\wtT_+$ at $z$, computed in the bases
$(\frac{\pa}{\pa x}, v_\om (z))$ at $z$, $(\frac{\pa}{\pa x},
v_{\om'} (z'))$ at $z' = \wtT^+ (z)$, is lower triangular; the first
diagonal coefficient is
\begin{equation}
A^{-1}_x (y,x')\; \Bigl( 1 - B_x (y,x') \, D \vphi_{\om'}
(y')\Bigr)\label{eq10.72}
\end{equation}
We denote by $\wtb (z)$ the logarithm of the absolute value of this
coefficient. As $\vphi_\om$ is $C^{1+\lip}$ uniformly in $\om$ and
$g^n: P \mapsto Q$ has bounded distortion, we have, for all $z$,
$z^*\in \om$
\begin{equation}
|\wtb (z) - \wtb (\zs)| \; \leqslant \; C |z-\zs|.\label{eq10.73}
\end{equation}
Let $j>0$ and let $\om$, $\whom$ be stable  curves which belong to
the same component of the domain of $(\Tim)^j$. Let $z$, $\whz$ be
point of $\om$, $\whom$, respectively, with the same $y$ coordinate.
It follows from Proposition~\ref{propo41} that one has
\begin{equation}
|\wtb (z) - \wtb (\whz)| \; \leqslant \; C \,
\te^j_0.\label{eq10.74}
\end{equation}
We also have, from the definition of $\wtb$:
\begin{equation}
|\wtb (z) + \log \, |P|| \; \leqslant \; C.\label{eq10.75}
\end{equation}
We want to get rid of the dependence of $\wtb$ on the $y$ coordinate
along $\om$ by adding a coboundary term. We define
\begin{equation}
\qqim \cD\; = \; \bigcap_{j\geqslant 0}\; \dom\, (\Tim)^j \; = \;
\bigcap_{j\geqslant 0} \; (\Tim)^{-j} \,(\cD_+).\label{eq10.76}
\end{equation}
For each $a$, fix some $y^0_a\in I^u_a$. Then, for $\om\in
\qqim\cD$, $\om \subset R_a$, $z\in\om$, define
\begin{equation}
\De b (z) = \sum_{i\geqslant 0}\; \Bigl( \wtb ((\wtT^+)^i (z)) -
\wtb ((\wtT^+)^i (z^0))\Bigr),\label{eq10.77}
\end{equation}
where $z^0$ is the point on $\om$ with vertical coordinate equal to
$y^0_a$.

From the cone condition, we have, for $i\geqslant 0$:
\begin{equation}
||(\wttm)^i (z) - (\wtT^+)^i(z^0)||\; \leqslant \; C\la^{-i}.
\label{eq10.78}
\end{equation}
The series defining $\De b$ is uniformly convergent from
(\ref{eq10.73}), (\ref{eq10.78}), and $\De b$ is bounded on
$\qqim\wtcD := \pi^{-1} (\qqim\cD)$.

Write $z^1$ for the point on $T_+ (\om)\in R_{a'}$ with vertical
coordinate $y^0_{a'}$. We have
\begin{equation}
\De b (z) - \De b (\wttm (z))\; = \; \wtb (z) - b
(\om),\label{eq10.79}
\end{equation}
with
\begin{equation}
b (\om) = \wtb (z^0) + \sum_{i\geqslant 0} \Bigl[ \wtb
((\wttm)^{i+1} (z^0)) - \wtb ((\wttm)^i (z^1))\Bigr].\label{eq10.80}
\end{equation}
We call $b$ the (logarithmic) mean transverse dilatation.

\begin{Propo}\label{propo42}
The mean transverse dilatation $b$, which differs from $\wtb$ on
$\qqim\cD$ by the coboundary of the bounded function $\De b$,
satisfies
\begin{equation}
|b (\om) - b (\whom)|\; \leqslant \; C \te^j_1\label{eq10.81}
\end{equation}
if $\om$, $\whom$ belong to the same component of the domain of
$(\Tim)^j$. Here $\te_1$ is a fixed constant in $(0,1)$ larger than
$\te_0$.
\end{Propo}

\begin{proof}
We have only to prove (\ref{eq10.81}). Let $z^0\in \om$, $z^1\in
\Tim (\om)$ as above and let $\whz^0\in \whom$, $\whz^1\in \Tim
(\whom)$ be similarly defined. We have, for $i\geqslant 0$:
\begin{eqnarray}
|\wtb ((\wttm)^{i+1} (z^0)) - \wtb ((\wttm)^i (z^1))| \; &\leqslant
\; &C \la^{-i},\label{eq10.82}\\
|\wtb ((\wttm)^{i+1} (\whz^0)) - \wtb ((\wttm)^i (\whz^1))| \;
&\leqslant \; &C \la^{-i},\label{eq10.83}
\end{eqnarray}
From (\ref{eq10.74}), we also have
\begin{equation}
|\wtb (z^0) - \wtb (\whz^0)|\; \leqslant \; C
\te^j_0.\label{eq10.84}
\end{equation}
For $0\leqslant i < j-1$, we compare $\wtb ((\wttm)^{i+1} (z^0))$
and $\wtb ((\wttm)^{i+1} (\whz^0))$ as follows: let $\whz^{i+1}$ be
the point on $T_+^{i+1} (\whom)$ with the same $y$-coordinate as
$(\wttm)^{i+1} (z^0)$; one has
\begin{eqnarray}
\left\vert \wtb (\whz^{i+1}) - \wtb ((\wttm)^{i+1}
(z^0))\right\vert \; &\leqslant \; &C\,\te_0^{j-i-1},\label{eq10.85}\\
||\whz^{i+1} - (\wttm)^{i+1} (\whz^0)|| \; &\leqslant \; &C
\la^{i+1-j}.\label{eq10.86}
\end{eqnarray}
Dealing in the same way with the terms involving $z^1$, $\whz^1$, we
obtain (\ref{eq10.81}) with
\begin{equation}
\te_1 \; = \; \max \, (\te_0^{1/2}, \la^{-1/2})\label{eq10.87}
\end{equation}
proving the proposition.
\end{proof}

\subsection{Definition of a Transverse Operator\label{sub10.7}}
As $\De b$ is bounded, it follows from (\ref{eq10.75}) that
\begin{equation}
|b (\om) + \log \, |P|| \; \leqslant \; C\label{eq10.88}
\end{equation}
for all $(P,Q,n)\in \cP$, $\om\in \qqim\cR (P)$.

It is then a consequence of Corollary~\ref{coro14} in
Subsection~\ref{sub10.3} that the series
\begin{equation}
\sum_{\Tim \om'=\om}\; \exp\, (-db (\om'))\label{eq10.89}
\end{equation}
over pre-images $\om'$ of a given stable curve $\om$ is converging,
uniformly in $\om$, for $d\geqslant \dMs$. We will, therefore, define
a transfer operator $\Ld$ for $d\geqslant \dMs$ as follows: for a
bounded function $h$ defined on $\qqim\cD$, for $\om\in\qqim\cD$, we
set
\begin{equation}
\Ld\, h (\om) \; = \; \sum_{\Tim \om'=\om}\; \exp\, (-db (\om')) \,
h (\om').\label{eq10.90}
\end{equation}
We can also view this sum over pre-images as a sum over inverse
branches of $\Tim$, which are in one-to-one correspondence with the
primes $(P,Q,n)$ such that $Q$ and $\om$ belong to the same
rectangle $R_a$. Accordingly, we split the series in two parts: a
finite sum corresponding to the trivial primes,
(cf.~Subsection~\ref{eq10.1}), which we denote by $\Lod$ and which
is defined for all values of $d$, and a perturbative term which we
denote by $\De \Ld$. The formula (\ref{eq10.90}) defines a bounded
operator from the space of bounded functions on $\qqim \cD$ into itself, but to have nice spectral
properties we need, as usual, to restrict to spaces of slightly more
regular functions.

Let $\te$ be a constant with
\begin{equation}
\te_1\; < \;\te \;< \;1\label{eq10.91}
\end{equation}
where $\te_1$ comes from Proposition~\ref{propo42} and satisfies
$\te_1 > \la^{-1}$ (cf.~\ref{eq10.87}). Denote by $E$ the space of
bounded functions $h$ on $\qqim\cD$ which satisfy, for some constant
$C > 0$,
\begin{equation}
|h (\om) - h (\whom)| \; \leqslant \; C \te^j\label{eq10.92}
\end{equation}
whenever $\om$, $\whom$ belong to the same component of the domain
of $(\Tim)^j$. We denote by $||h||_\infty$ the usual norm on bounded
functions, by $|h|_E$ the best possible $C$ and set
\begin{equation}
||h ||_E \; = \; \max \, (|h|_E, ||h||_\infty).\label{eq10.93}
\end{equation}
It is clear that $E$ is a Banach space.

\begin{Propo}\label{propo43}
For $d\geqslant \dMs$, $\Ld$ restricts to a bounded operator on $E$.
Moreover, the norm of the perturbative part $\De\Ld$ is as small as
we want if $\vep_0$ is small enough.
\end{Propo}

\begin{proof}
Let $h\in E$, $\om$, $\whom\in\qqim\cD$, $j>0$.  Assume that $\om$,
$\whom$ belong to the same component of the domain of $(\Tim)^j$.
Let $(P,Q,n)$ be a prime such that $Q$, $\om$, $\whom$ belong to the
same rectangle $\Ra$, and let $\om_1$, $\whom_1$ be the inverse
images of $\om$, $\whom$ by $\Tim$ corresponding to this inverse
branch. By the definition of $|\,\;|_E$, we have
\begin{equation}
|h (\om_1) - h (\whom_1)| \; \leqslant \; |h|_E \; \te^{j+1}.
\label{eq10.94}
\end{equation}
From Proposition~\ref{propo42}, we have
\begin{equation}
|b (\om_1) - b (\whom_1)| \; \leqslant \; C\te_1^{j+1}.
\label{eq10.95}
\end{equation}
It follows from (\ref{eq10.88}) that
\begin{equation}
|\exp (-db (\om_1)) - \exp (-db (\whom_1))| \; \leqslant \; C(d) \,
|P|^d\, \te_1^{j+1}. \label{eq10.96}
\end{equation}
Putting together (\ref{eq10.94}) and (\ref{eq10.96}), we have
\begin{equation}
|h (\om_1) \exp (-db (\om_1)) - h (\whom_1) \exp (-db (\whom_1))| \;
\leqslant \; C(d) \, |P|^d\, (\te^{j+1} |h|_E + \te^{j+1}_1
||h||_\infty).\label{eq10.97}
\end{equation}
Summing over (non trivial) primes yields for $d\geqslant \dMs$:
\begin{eqnarray}
|\De \Ld h|_E \, &< \, &\vep_1 ||h||_E,\label{eq10.98}\\
|\Ld h|_E \, &< \, &C ||h||_E,\label{eq10.99}
\end{eqnarray}
where $\vep_1$ can be made arbitrarily small if $\vep_0$ is small
enough, according to Corollary~\ref{coro14}. The same estimates (for
$d\geqslant \dMs$) for $||\De\Ld h||_\infty$ and $||\Ld h||_\infty$
are easier and can be seen directly. The proposition follows.
\end{proof}

\subsection{Spectral Properties of the Transfer Operator\label{sub10.8}}
Let us denote by $\qqim\cR (K)$ the set of stable curves $\om$ which
are intersections of a sequence of rectangles belonging to $\cR
(\Io)$; these stable curves are precisely those which meet the
initial horseshoe $K$.

Observe that $\qqim\cR (K) \subset \qqim\cD$. Denote by $E_K$ the
space of bounded functions $h$ on $\qqim\cR (K)$ which satisfy
\begin{equation}
|h (\om) - h (\whom)| \; \leqslant \; C\te^j,\label{eq10.100}
\end{equation}
whenever $\om$, $\whom$ belong to the same component of the domain
of $(\Tim)^j$. Define $|h|_{E_K}$, $||h||_{E_K}$ as above, which
makes $E_K$ a Banach space.

Let $h\in E$; the restriction of $h$ to $\qqim\cR (K)$ belongs to
$E_K$ and we have
\begin{equation}
||h / \qqim\cR (K)||_{E_K} \; \leqslant \; ||h||_E.\label{eq10.101}
\end{equation}
The formula for $\Lod$ defines a bounded operator, still denoted by
$\Lod$, on $E_K$ and we have a commutative diagram
\begin{equation}
\begin{CD}
E           @>>{\Lod}>      E\\
@V{r}VV                     @VV{r}V\\
E_K         @>>{\Lod}>      E_K
\end{CD} \label{eq10.102}
\end{equation}
where $r: E \mapsto E_K$ is the restriction operator. The bounded
operator $\Lod: E_K \mapsto E_K$ is the subject of the classical
theory by Bowen, Ruelle, Sinai for uniformly hyperbolic systems.

Let us recall some standard results of this theory.

\medskip\noindent
a)~~There is a direct sum invariant decomposition
\begin{equation}
E_K = \Bbr\, \hdl \; \oplus \; \Hdl\label{eq10.103}
\end{equation}
depending analytically on the parameter $d$, such that $\hdl$ is a
positive eigenfunction, with associated eigenvalue $\ladl > 0$, and
such that
\begin{equation}
sp \Bigl( \Lod / \Hdl \Bigr) \subset \{ |z|\; < \; \ladl
\}.\label{eq10.104}
\end{equation}

\medskip\noindent
b)~~There exists a (unique) probability measure $\mudl$ on $\qqim\cR
(K)$ such that
\begin{equation}
\Hdl \; = \; \{ h\in E_K,\; \int hd\mudl = 0\}.\label{eq10.105}
\end{equation}
One normalizes $\hdl$ to have $\int \hdl d\mudl = 1$. Then, the
probability measure $\nudl = \hdl \mudl$ is invariant under the
restriction of $\Tim$ to $\qqim\cR (K)$ (observe that $\wttm$ on
$\qqim\wtR (K)$ is just the restriction of~$g$).

Let $\Eio$ be the kernel of the restriction operator $r: E \mapsto
E_K$. It is invariant under $\Lod$.

\begin{lema}\label{lema14}
One has, for all $d\in\Bbr$.
\[
sp \Bigl( \Lod / \Eio\Bigr) \subset \{ |z|\; \leqslant \; \te
\ladl\}.
\]
\end{lema}

\begin{proof}
Let $h\in \Eio$, $j\geqslant 0$. We have
\begin{equation}
(\Lod)^j h (\om) = \sum\nolimits^0_{(\Tim)^j (\om')=\om} \; h (\om')
\exp (-db^\pjj (\om')),\label{eq10.106}
\end{equation}
where the symbol $\sum^0$ indicates that we only consider inverse
branches of $\Tim$ associated with trivial primes. The notation
$b^\pjj$ denotes the Birkhoff sum
\begin{equation}
b^\pjj (\om') \; = \; \sum_{0\leqslant i < j}\; b ((\Tim)^i
(\om')).\label{eq10.107}
\end{equation}
We observe that in the sum in (\ref{eq10.106}), each $\om'$ belongs
to the same component of the domain of $(\Tim)^j$ as a stable curve
in $\qqim\cR (K)$. As $h$ belongs to $\Eio$, this implies that for
such a $\om'$ we have
\begin{equation}
|h (\om')| \; \leqslant \; |h|_E \; \te^j.\label{eq10.108}
\end{equation}

On the other hand, we have
\begin{equation}
\sum\nolimits^0 \exp (-d b^\pjj (\om')) \; \leqslant \; C
\la_d^{'j}, \label{eq10.109}
\end{equation}
and it follows that
\begin{equation}
||(\Lod)^j \, h||_\infty \; \leqslant \; C\la_d^{'j} \te^j ||h||_E.
\label{eq10.110}
\end{equation}

Let $\whom\in \qqim\cR$ belong to the same component of the domain
of $(\Tim)^\ell$ as $\om$. Denote by $\whom'$ the inverse image of
$\whom$ associated to the same sequence of trivial primes as $\om'$.
We have
\begin{equation}
|h (\om') - h (\whom')| \; \leqslant \; |h|_E \; \te^{j+\ell},
\label{eq10.111}
\end{equation}
and, from Proposition~\ref{propo42}
\begin{equation}
|b^\pjj (\om') - b^\pjj (\whom')| \; \leqslant \; C\,
\te_1^\ell.\label{eq10.112}
\end{equation}
Using also (\ref{eq10.108}) and (\ref{eq10.109}), we obtain
\begin{equation}
\left\vert (\Lod)^j h (\om) - (\Lod)^j h (\whom)\right\vert \;
\leqslant \; C \te^{j+\ell} \la_d^{'j} ||h||_E,\label{eq10.113}
\end{equation}
which implies the statement of the Lemma.
\end{proof}

We deduce from Lemma~\ref{lema14} that there is a unique function in
$E$, still denoted by $\hdl$, which restricts to $\hdl$ on $\qqim\cR
(K)$ and satisfies
\begin{equation}
\Lod\, (\hdl) \; = \; \ladl \, \hdl.\label{eq10.114}
\end{equation}
Moreover, defining a supplementary hyperplane by
\begin{equation}
\Hdll \; = \; r^{-1} (\Hdl) \oplus \Eio,\label{eq10.115}
\end{equation}
we have that $\Hdll$ is invariant under $\Lod$ and
\begin{equation}
sp \Bigl( \Lod / \Hdll\Bigr) \; \subset \; \{ |z| \; \leqslant \;
\ladll\},\label{eq10.116}
\end{equation}
where $\ladll < \ladl$ is independent of $\vep_0$.

Using Proposition~\ref{propo43}, we now consider $\Ld$ itself,
assuming that $\vep_0$ is small enough and $d\geqslant \dMs$.

As the norm of the perturbation part $\De \Ld$ is arbitrarily small,
we conclude that $\Ld$ has a positive eigenfunction $\hd$, with
associated eigenvalue $\lad$ arbitrarily close to $\ladl$, and an
invariant supplementary hyperplane $\Hd$ satisfying
\begin{equation}
sp \Bigl( \Ld / \Hd \Bigr) \; \subset \; \{ |z| \; < \;
\lad\}.\label{eq10.117}
\end{equation}
Moreover, $\hd$, $\lad$ and $\Hd$ depend analytically on $d$ for
$d>\dMs$ because $\Ld$ does.  We check that
\begin{equation}
\hd \; \geqslant \; c^{-1} > 0.\label{eq10.118}
\end{equation}

Indeed, the sequence $h^\pnn = \la_d^{-n} L_d^n (1)$ converge to a
positive multiple of $\hd$. We have
\begin{equation}
h^\pnn (\om) \; = \; \sum_{(\Tim)^n (\om')=\om} \exp (-d b^\pnn
(\om')).\label{eq10.119}
\end{equation}
Let $\om$, $\whom$ be elements of $\qqim\cD$ in the same rectangle
$\Ra$; let $\om'$, $\whom'$ be pre-images of $\om$, $\whom$ by
$(\Tim)^n$ associated with the same sequence of primes. We have
(cf.~\ref{eq10.112})
\begin{equation}
|b^\pnn (\om') - b^\pnn (\whom')| \; \leqslant \; C,\label{eq10.120}
\end{equation}
and it follows that
\begin{equation}
C^{-1} \; \leqslant \; (h^\pnn (\om'))^{-1} h^\pnn (\om) \;
\leqslant \; C.\label{eq10.121}
\end{equation}
This implies (\ref{eq10.118}). One normalizes $\hd$ in order to have
\begin{equation}
\hd \; = \; \lim_{n\to +\infty}\; \la^{-n}_d \, L^n_d\,
(1).\label{eq10.122}
\end{equation}
Denote then by $\mud$ the linear form on $E$ with kernel $\Hd$
normalized by $\mud (\hd) = 1$. We have, for all $h\in E$
\begin{equation}
\lim_{n\to\infty} \; \la_d^{-n} L^n_d\, h\; = \; \mud (h)
\hd.\label{eq10.123}
\end{equation}
As $\Ld$ is a positive operator, $\mud$ is positive. Observe also
that for all $(P,Q,n)\in\cR$, the characteristic function $\chi_P$
(equal to 1 if $\om\subset P$, $0$ otherwise) belongs to  $E$ and
satisfies $L^n \chi_P > 0$ everywhere for some $n>0$. Therefore,
there exists a unique probability measure on $\qqim\cR$, still
denoted by $\mud$, which coincides with $\mud$ on the intersection
of $E$ with $C (\qqim\cR)$.

\subsection{The Gibbs Measure\label{sub10.9}}
From the defining property (\ref{eq10.123}) of $\mud$, we have, for
all $h\in E$
\begin{equation}
\mud (\Ld\, h) \; = \; \lad \, \mud (h).\label{eq10.124}
\end{equation}
We will now check the classical Jacobian property for $\mud$.

Let $(P,Q,n)\in\cR$. Let $P_1\supset P$ be the thinnest prime
containing $P$. The application $\Tim$ is injective on the set
$\qqim\cR (P)$ formed by the $\om\in \qqim\cR (P_1)$ which are
contained in $P$. The image of this set by $\Tim$ is exactly the set
of stable curves contained in $P'$, with $\Tim (P,Q,n) =
(P',Q',n')$.

Let $h$ be a function in $E$ which vanishes outside $\qqim\cR (P)$.
Then, $\Ld\,h$ vanishes on any curve not contained in $P'$, and
satisfies
\begin{equation}
\Ld h (\Tim \om) \; = \; h (\om)\, \exp (-d b(\om))\label{eq10.125}
\end{equation}
for $\om\in \qqim\cR (P)$. The relation (\ref{eq10.124}) for $h$ is
the Jacobian property.

Consider in particular the case where $h$ is the characteristic
function of $\qqim\cR (P)$. We then obtain
\begin{equation}
\lad\,\mud (\qqim\cR (P)) \; = \; \int_{P'} \exp (-d b ((\Tim)^{-1}
\om)) d\mud (\om).\label{eq10.126}
\end{equation}

We now will specify the value of $d$ by asking that
\begin{equation}
\lad \; = \; 1.\label{eq10.127}
\end{equation}
Indeed, we have, for $d\geqslant \dMs$
\begin{equation}
\frac{\pa}{\pa d}\; \lad \; < \; 0.\label{eq10.128}
\end{equation}
This follows from the formula for $\Ld$ and the fact that $b^\pnn$
increases at least linearly with $n$ (cf.~(\ref{eq10.107})).

We also see easily that
\begin{equation}
\lim_{d\to +\infty}\; \lad \; = \; 0.\label{eq10.129}
\end{equation}
Finally, we have
\begin{equation}
\la_{\dMs} \; > \; 1.\label{eq10.130}
\end{equation}
Indeed, $\dMs$ was chosen in order to be smaller than the transverse
Hausdorff dimension of $W^s (K)$. This means that the eigenvalue
$\la'_{\dMs}$ for $\L^0_{\dMs}$ on $E_K$ satisfies $\la'_{\dMs} >
1$. As $\De L_{\dMs}$ is also a nonnegative operator, we have
$\la_{\dMs} \geqslant \la'_{\dMs}$. Therefore, (\ref{eq10.130})
holds, and it follows from (\ref{eq10.128})--(\ref{eq10.130}) that
(\ref{eq10.127}) holds for a unique value of $d$. We will denote
this value by $d_s$. We shall indeed see that $d_s$ is the
transverse Hausdorff dimension of $\qqim\wtcR$ which we were able to
define in Subsection~\ref{sub10.5}.

We just write $\mu$ for the measure $\mu_{d_s}$ and $h^*$ for the
eigenfunction $h_{d_s}$.

\begin{Propo}\label{propo44}
For any $(P,Q,n)\in \cR$, we have
\[
C^{-1} |P|^{d_s} \; \leqslant \; \mu (\{\om \subset P\})\; \leqslant
\; C |P|^{d_s}.
\]
\end{Propo}

\begin{proof}
Let
\begin{equation}
(P,Q,n)\; = \; (P_1, Q_1, n_1)\; * \cdots * \; (P_r,Q_r,n_r)
\label{eq10.131}
\end{equation}
be the prime decomposition of $(P,Q,n)$. If $\om\in\qqim\cD$
satisfies $(\Tim)^i (\om) \in \qqim\cR (P_{i+1})$ for $0\leqslant i
< r$, we have, from (\ref{eq10.88}), (see the definition of
$b^{(r)}$ in (\ref{eq10.107})):
\begin{equation}
C^{-1} |P|^\ds \; \leqslant \; \exp (-\ds b^\prr (\om)) \; \leqslant
\; C |P|^\ds.\label{eq10.132}
\end{equation}
It, then, follows from the Jacobian property that
\begin{eqnarray}
\mu ( \{\om\subset P\})\; &\geqslant \; &C^{-1} |P|^\ds \, \mu
(\{\om\in \qqim\cR (a)\}),\label{eq10.133}\\
&\geqslant \; &C^{-1} |P|^\ds,\notag
\end{eqnarray}
where $\Ra$ is the rectangle containing $Q$.

For the opposite inequality, we have also to take into account the
other inverse branches of $T_+^r$ when we estimate $L^r_{\ds}
(\chi_P)$, where $\chi_P$  is the characteristic function of
$\{\om\subset P\}$. For $0\leqslant i\leqslant r$, let
\begin{equation}
(P^i,Q^i,n^i)\; = \; (P_{i+1}, Q_{i+1}, n_{i+1})\; * \cdots * \;
(P_r,Q_r,n_r)\label{eq10.134}
\end{equation}
(with $(P^r, Q^r, n^r)\; = \; (\Ra, \Ra, 0)$). We have
\begin{equation}
L_{\ds} \chi_P \; = \; \chi^1_P + \De \chi^1_P\label{eq10.135}
\end{equation}
where
\begin{equation}
\chi^1_P (\om^1) =
\begin{cases}
0 &\text{ if $\;\om^1 \not\subset P^1$}\\
\exp (-\ds b (\om^0) &\text{ if $\;\om^1 = \Tim (\om^0)$ for some
$\om^0 \in \qqim\cR (P)$}
\end{cases}
\label{eq10.136}
\end{equation}
and
\begin{equation}
\De \chi^1_P \; \leqslant \; C\sum \; |P^*_1|^{\ds},\label{eq10.137}
\end{equation}
where the sum runs over prime elements $(P^*_1, Q^*_1, n^*_1)$ with
$P^*_1$ contained in $P$ and distinct from $P$.  By
Proposition~\ref{propo38} in Subsection~\ref{sub10.3}, we obtain
\begin{equation}
\mu (\De \chi^1_P) \; \leqslant \; C |P|^{\ds}\,
\kappa^{\tfrac{r-1}{2}}.\label{eq10.138}
\end{equation}
If $r>1$, we write similarly
\begin{equation}
\Ld \chi^1_P \; = \; \chi^2_P + \De \chi^2_P,\label{eq10.139}
\end{equation}
where $\chi^2_P$ is associated with the inverse branch defined by
the prime $P_2$ and vanishes outside $P^2$. The perturbative term
satisfies
\begin{equation}
\De \chi^2_P \; \leqslant \; C |P_1|^{\ds} \sum \;
|P^*_2|^{\ds},\label{eq10.140}
\end{equation}
where the sum now is over primes $P^*_2$ contained in $P^1$ and
distinct from $P^1$. Proposition~\ref{propo38} now gives
\begin{eqnarray}
\mu (\De \chi^2_P) \; &\leqslant \; &C |P_1|^\ds \, |P^1|^\ds
\kappa^{\tfrac{r-2}{2}}\label{eq10.141}\\
&\leqslant \; &C |P|^\ds \, \kappa^{\tfrac{r-2}{2}}.\notag
\end{eqnarray}
We iterate this process. At the last step, we have from
(\ref{eq10.132})
\begin{equation}
\mu (\chi^r_P) \; \leqslant \; C |P|^{\ds}.\label{eq10.142}
\end{equation}
The contribution of the perturbative terms is bounded by
\begin{equation}
\mu \Bigl( \sum^r_{1} \,\chi^i_P\Bigr) \; \leqslant \; C |P|^{\ds}\,
\sum^r_1 \, \kappa^{\tfrac{r-i}{2}}\; \leqslant \; C
|P|^\ds.\label{eq10.143}
\end{equation}
\end{proof}

\begin{Coro}\label{coro15}
The transverse Hausdorff dimension of $\qqim \wtcR$ is $\leqslant
\ds$. More precisely, for any curve $\ga$ which is transverse to
$\qqim\wtcR$, the Hausdorff measure in dimension $\ds$ of the
intersection of $\ga$ with $\qqim \wtcR$ is finite.
\end{Coro}
We will see below that the transverse Hausdorff dimension is equal
to $\ds$.

\begin{proof}
Let $\de > 0$, choose a finite collection of disjoint rectangles
$P_i$ with $|P_i|\leqslant \de$ for each $i$ and $\qqim\wtcR \subset
\cup P_i$. We have
\begin{eqnarray}
1 \; = \; \sum \mu (P_i) &\geqslant \; &C^{-1} \sum\; |P_i|^\ds\label{eq10.144}\\
&\geqslant \; &C^{-1} \sum\; [\diam (\ga \cap P_i)]^\ds\notag
\end{eqnarray}
and the statement of the Corollary follows.
\end{proof}
The following statement shows that the dynamics $\Tim$ is only
undefined on a small set.

\begin{Propo}\label{propo45}
The transverse Hausdorff dimension  of the set $\qqim\wtcR -
\qqim\wtcD$ is $\leqslant \dMs < \ds$.
\end{Propo}

\begin{proof}
We have
\begin{equation}
\qqim\wtcR - \qqim\wtcD = \bigcup_{n\geqslant 0} \; (\Tim)^{-n}
(\cN_+)\label{eq10.145}
\end{equation}
As each $(\Tim)^n$ has countably many inverse branches which are
Lipschitzian, it is sufficient to prove that the transverse
Hausdorff dimension of $\cN_+$ is $\leqslant \dMs$. But by the
definition of $\cN_+$, for any $\de>0$, the union of prime
rectangles $P$ with $|P|<\de$ contains $\cD$.
Proposition~\ref{propo45} then follows from Corollary~\ref{coro14}.
\end{proof}

\subsection{Transverse Hausdorff Dimension of $\wtcR^\infty_+$\label{sub10.10}}
\begin{Propo}\label{propo46}
The transverse Hausdorff dimension of $\qqim\wtcR$ is the number
$\ds$ characterized by $\la_{\ds} = 1$.
\end{Propo}

\begin{rem}
We have already seen that the Hausdorff measure in dimension $\ds$
of the intersection of $\qqim\wtcR$ with a transverse curve is
always finite. We do not know whether it is positive or always zero.
\end{rem}

\begin{proof}
Let $\ga$ be a smooth horizontal-like curve in some $\Ra$. We denote
by $[\ga]$ the set of stable curves which meet $\ga$. We will show
that
\begin{equation}
\mu ([\ga]) \; \leqslant \; C (\diam\,\ga)^\ds\; (\log
(\diam\,\ga)^{-1})^{C_0}.\label{eq10.146}
\end{equation}
This, being true for all such $\ga$, clearly implies that the
transverse Hausdorff dimension of $\qqim\wtcR$ is $\geqslant \ds$,
which is sufficient to prove the proposition in view of
Corollary~\ref{coro15}.

Clearly, we may assume that $\mu ([\ga])>0$. Define
$(P_0,Q_0,n_0)\in\cR$ to be the element such that $P_0$ is the
thinnest rectangle containing any stable curve in $[\ga]$. There are
at least two children of $P_0$ which contain a stable curve in
$[\ga]$. If one of these children is simple, we must have
\begin{equation}
\diam\,\ga \; \geqslant \; C^{-1} |P_0|\label{eq10.147}
\end{equation}
and also
\begin{equation}
\mu ([\ga]) \; \leqslant \; \mu (\{ \om\subset P_0\}) \; \leqslant
\; C |P_0|^\ds\label{eq10.148}
\end{equation}
by Proposition~\ref{propo44}, which gives the required estimate (and
even better). This case is said to have complexity $0$. In the
remaining case, denote by $P_{0,i}$ the (non-simple) children of
$P_0$ which contain a stable curve in $[\ga]$. Each $P_{0,i}$ is
obtained by parabolic composition:
\begin{equation}
(P_{0,i}, Q_{0,i}, n_{0,i}) \; \in \; (P_0, Q_0, n_0) \; \square \;
(P^*_{0,i}, Q^*_{0,i}, n^*_{0,i})\label{eq10.149}
\end{equation}
and the widths are related through
\begin{equation}
C^{-1} \; \leqslant \; |P_{0,i}|\;|P_0|^{-1} |P^*_{0,i}|^{-1} \de
(Q_0, P^*_{0,i})^{\fud} \; \leqslant \; C.\label{eq10.150}
\end{equation}
Let $\de_0 = \sup\limits_i \, \de (Q_0, P^*_{0,i})$; let $\ga_1$ be
any horizontal-like curve with the following property: a stable
curve meets $\ga_1$ if and only if it is contained in some
$P^*_{0,i}$. For $\ell \geq 0$ , denote by $\gamma_{1, \ell}$ a piece of $\gamma_1$ with the property that a stable curve meets $\gamma_{1, \ell}$ iff it is contained in some $P^*_{0,i}$ with
\begin{equation}
\de_0 \, 2^{-\ell-1}\; \leqslant \; \de (Q_0, P^*_{0,i}) \;
\leqslant \; \de_0 2^{-\ell}\label{eq10.151}
\end{equation}
(if there is no such $P^*_{0,i}$, take $\ga_{1,\ell} = \emptyset$).
We can now write
\begin{eqnarray}
\mu ([\ga]) \; &\leqslant \; &C \sum_i \, |P_{0,i}|^\ds \; \text{
 (from Proposition~44)}\label{eq10.152}\\
&\leqslant \; &C\,|P_0|^\ds\; \sum_i\,|P^*_{0,i}|^\ds\,\de (Q_0, P^*_{0,i})^{-\fudt\,\ds}\notag\\
&\leqslant \; &C\,|P_0|^\ds\;\de_0^{-\fudt\,\ds} \,
\sum_{\ell\geqslant 0}\,2^{\tfrac{\ell ds}{2}}\, \sum\nolimits^{(\ell)} \; |P^*_{0,i}|^\ds\notag\\
&\leqslant \; &C\,|P_0|^\ds\;\de_0^{-\fudt\,\ds} \,
\sum_{\ell\geqslant 0}\,2^{\tfrac{\ell ds}{2}}\, \mu
([\ga_{1,\ell}])\notag
\end{eqnarray}
again by Proposition~\ref{propo44}. We have written $\sum^{(\ell)}$
for the partial sum over those $P^*_{0,i}$ satisfying
(\ref{eq10.151}). Assume for some constant $A>0$, that we have, for
each $\ell\geqslant 0$:
\begin{equation}
\mu ([\ga_{1,\ell}])\; \leqslant \; A (\diam\,
\ga_{1,\ell})^\ds.\label{eq10.153}
\end{equation}
Observe that we have
\begin{eqnarray}
\diam \, \ga_{1,\ell} \; &\leqslant \; &C 2^{-\ell} \, \diam\,
\ga_1,\label{eq10.154}\\
\diam \, \ga  \; &\geqslant \; &C^{-1}\,\de_0^{-\fudt} \,|P_0|\,
\diam\, \ga_1.\label{eq10.155}
\end{eqnarray}
Making use of (\ref{eq10.153})--(\ref{eq10.155}) in (\ref{eq10.152})
yields
\begin{equation}
\mu ([\ga]) \; \leqslant \; A C (\diam\,\ga)^\ds,\label{eq10.156}
\end{equation}
which is of the same form as (\ref{eq10.153}), but with a worse
constant $A C$ instead of $A$.

To obtain (\ref{eq10.146}), it is thus sufficient to define a
complexity index $c (\ga)\in \Bbn$ which satisfies
\begin{eqnarray}
c (\ga) \; &\leqslant \; &c\, \log\, \log\, |P_0|^{-1},\label{eq10.157}\\
c (\ga) \; &= \; &1 + \sup_\ell\, c\,
(\ga_{1,\ell}),\label{eq10.158}
\end{eqnarray}
the case of complexity $0$ having already been defined and dealt
with.

We want to use (\ref{eq10.158}) to give an inductive definition of
$c(\ga)$. This will work if the $\ga_{1,\ell}$ are in some sense
"simpler" than $\ga$. If all $\ga_{1,\ell}$ have complexity $0$, we
just set $c(\ga)=1$. Assume therefore that some $\ga_{1,\ell}$ has
complexity $>0$. This means that there exists an element
$(P_{1,\ell}, Q_{1,\ell}, n_{1,\ell})$ with the following
properties:

-- each $P^*_{0,i}$ related to $\ga_{1,\ell}$ through
(\ref{eq10.151}) is contained in some non-simple child of
$P_{1,\ell}$;

-- at least two non-simple children of $P_{1,\ell}$ contain some
$P^*_{0,i}$.

For any parameter interval $I$ containing the given parameter value,
$P_{1,\ell}$ is $I$-critical (cf.~Proposition~\ref{propo5} in
Subsection~\ref{sub6.4}). We now distinguish two cases.

\indent {\bf Case~1.~} $P_0$ is also $I$-critical (for any $I$ as
above).

Let $\Is$ be the largest parameter interval for which we have
\begin{equation}
|P_0| \; > \; |\Is|^\be.\label{eq10.159}
\end{equation}
Then, $(P_0,Q_0,n_0)$ cannot be $\Is$-bicritical as $\Is$ is
$\be$-regular, hence $Q_0$ cannot be $\Is$-critical. This implies
\begin{equation}
\de_0\, 2^{-\ell} \; \geqslant \; C^{-1} |\Is| \; \geqslant \;
C^{-1} |P_0|^{\be^{-1}(1+\tau)}.\label{eq10.160}
\end{equation}
As $P_{1,\ell}$ is not transverse to $Q_0$ (because $P_{0,i}$ was a
child of $P_0$), we must have from Proposition~\ref{propo10} in
Subsection~\ref{sub6.6} that
\begin{equation}
|P_{1,\ell}|\; \geqslant \; C^{-1} (\de_0
2^{-\ell})^{(1-\eta)^{-1}}.\label{eq10.161}
\end{equation}
Comparing with (\ref{eq10.160}), this guarantees that
\begin{equation}
|P_0|\; \ll \; |P_{1,\ell}|^{\fudt\,(1+\be)}.\label{eq10.162}
\end{equation}
This means indeed that every $P_{1,\ell}$ (such that the complexity
of $\ga_{1,\ell}$ is $>0$) is indeed simpler than $P_0$ and allows
us to use (\ref{eq10.158}) to define inductively $c(\ga)$. Observe
that the hypothesis of case~1 is always satisfied by the
$P_{1,\ell}$. The inequality (\ref{eq10.157}) follows from
(\ref{eq10.162}).

\indent {\bf Case~2.~} $P_0$ is $I$-transverse for $I$ small enough.

From Case~1, we have already defined the complexity indices
$c(\ga_{1,\ell})$ using (\ref{eq10.158}) and again we define
$c(\ga)$ by (\ref{eq10.158}). We have to check (\ref{eq10.157}) in
this case. This will hold if we have, for each $\ell$,
\begin{equation}
\log\, \log |P_{1,\ell}|^{-1} \; \leqslant \; c\, \log\, \log\,
|P_0|^{-1}.\label{eq10.163}
\end{equation}
But (\ref{eq10.161}) still holds. We also have
\begin{eqnarray}
|Q_0| \; &\ll \; &\de_0 2^{-\ell}\; \text{ (from
(R7))},\label{eq10.164}\\
\log\, |Q_0|^{-1} \; &\leqslant \; &C n_0,\label{eq10.165}\\
\log\, \log \, |P_0|^{-1} \; &\geqslant \;
&\frac{\log\,\frac{3}{2}}{\log\,2} \, \log\, n_0 -
C,\label{eq10.166}
\end{eqnarray}
from Proposition~\ref{propo13} of Subsection~\ref{sub7.1}. Putting
this together, we obtain (\ref{eq10.163}). The proof of the
proposition is now complete.
\end{proof}

\subsection{Invariant Measures\label{sub10.11}}
From the Gibbs measure $\mu$, which is not invariant but has the
Jacobian property, we define a $\Tim$-invariant measure $\nu$ on
$\qqim\cD$ by
\begin{equation}
d\nu\; = \; h^* d\mu.\label{eq10.167}
\end{equation}
The measure $\nu$ is actually a probability measure on $\qqim\cR$ by
Proposition~\ref{propo44} and Corollary~\ref{coro14} (see proof of
Proposition~\ref{propo45}), which implies that
\begin{equation}
\mu (\qqim\cR - \qqim\cD) \; = \; 0.\label{eq10.168}
\end{equation}
As $\hs$ is bounded and bounded away from $0$
(cf.~(\ref{eq10.118})), the statement of Proposition~\ref{propo44}
is also valid for $\nu$ instead of $\mu$. To check that $\nu$ is
indeed $\Tsm$-invariant, we first observe that, if $h_0$, $h_1\in
E$, the product $h_0 h_1$ also belongs to $E$; indeed we have
\begin{equation}
|h_0 h_1|_E \; \leqslant \; ||h_0||_\infty |h_1|_E + |h_0|_E
||h_1||_\infty.\label{eq10.169}
\end{equation}
In particular, for any $h\in E$, $h\hs$ also belongs to $E$. Let
$h\in E$. We write
\begin{eqnarray}
\int h (\Tim\om) d\nu (\om) \; &= \; &\int h (\Tim\om) \hs (\om)
d\mu (\om)\label{eq10.170}\\
&= \; &\sum_\cP \int h (\Tim\om) \hs(\om) \chi^*_P (\om) d\mu
(\om),\notag
\end{eqnarray}
where $\chi^*_P$ is the characteristic function of $\qqim\cR (P)$.
The Jacobian property gives
\begin{equation}
\int h (\Tim\om) \hs (\om) \chi^*_P (\om) d\mu (\om) = \int h (\om)
\hs (\om') \exp (-\ds b (\om')) \, d\mu (\om)\label{eq10.171}
\end{equation}
where $\om'$ is the image of $\om$ under the inverse branch of
$\Tsm$ associated with $P$. Summing over $P$ and using that $\hs$ is
$L_{\ds}$-invariant gives
\begin{equation}
\int h (\Tim\om) d\nu (\om)\; = \; \int h (\om) d\nu
(\om).\label{eq10.172}
\end{equation}
But $E\cap C (\qqim\cR)$ is dense in the space of continuous
functions $C (\qqim\cR)$; the invariance of $\nu$ follows.

Let us now check that the invariant measure $\nu$ is ergodic. Let
$A\subset \qqim\cR$ be a $\Tim$-invariant Borel subset with $\nu
(A)>0$ and $A^c$ its complement. Let $\vep>0$. We will prove that
there exists $a\in\AL$ such that
\begin{equation}
\nu (A \cap \qqim\cR (a))\; \geqslant \; (1-\vep)\; \nu (\qqim\cR
(a)).\label{eq10.173}
\end{equation}
As $\vep>0$ is arbitrary, this easily implies $\nu(A)=1$.

As $\nu (A)>0$, we can find $(P,Q,n)$ such that
\begin{equation}
\nu (\{ \om\subset P\} \cap A^c)\; \leqslant \vep'\,\nu (\{\om \subset
P\}),\label{eq10.174}
\end{equation}
where $\vep'\vep^{-1}$ is small. Let $r$ be the number of factors in
the prime decomposition of $(P,Q,n)$. Up to a set of measure $0$, we
have
\begin{equation}
\{ \om\subset P\} \; = \; \bigcup_{0\leqslant j\leqslant r}\;
\bigcup_{P_j}\, (\Tim)^{-j} (\qqim\cR (P_j))\; \text{ mod}\;
0\label{eq10.175}
\end{equation}
where $P_j$ runs through prime elements satisfying $P_j\subset
(\Tim)^j (P)$ and $(\Tim)^{-j}$ is the inverse branch of $(\Tim)^j$
whose image contains $P$. From (\ref{eq10.174}), there exists
$0\leqslant j\leqslant r$ and $P_j$ such that
\begin{equation}
\nu (A^c \cap (\Tim)^{-j} (\qqim\cR (P_j))\; \leqslant \; \vep'\,\nu
((\Tim)^{-j} (\qqim\cR (P_j))).\label{eq10.176}
\end{equation}
We apply the Jacobian property, taking (\ref{eq10.120}) into account
to get (\ref{eq10.173}) with $\vep = C\vep'$. We have proved that
$\nu$ is ergodic. We summarize:

\begin{Propo}\label{propo47}
The measure $d\nu = \hs\,d\mu$ is $\Tim$-invariant, ergodic. It
satisfies, for all $(P,Q,n)\in \cR$:
\[
C^{-1} |P|^\ds \; \leqslant \; \nu (\{\om\subset P\})\; \leqslant \;
C |P|^\ds
\]
and $\nu (\qqim\cD) = 1$.
\end{Propo}
We will now lift $\nu$ to obtain a $\wttm$-invariant probability
measure on $\qqim\wtcR$.

\begin{Propo}\label{propo48}
There exists a unique probability measure $\wtnu$ on $\qqim\wtcR$
which is $\wttm$-invariant and projects onto $\nu$ under $\pi$. It
is ergodic.
\end{Propo}

\begin{proof}
The arguments are standard.

\noindent {\bf Existence.~} Denote by $\cM (\nu)$ the set of
probability measures on $\qqim\cR$ which project onto $\nu$. This is
a compact set for the weak topology, invariant under $\wttm$ because
$\nu$ is $\Tim$-invariant. One obtains a $\wttm$-invariant measure
in $\cM (\nu)$ by taking any $\wtnu_0\in \cM (\nu)$ and choosing a
weak limit of a subsequence of
\begin{equation}
\frac{1}{n}\; \sum_0^{n-1} \, [(\wttm)^j]^*
(\wtnu_0).\label{eq10.177}
\end{equation}

\noindent {\bf Uniqueness.~} The set of fixed points for the action
of $\wttm$ on $\cM (\nu)$ is thus non-empty. It is also compact and
convex. If it has more than one point, it has at least two distinct
extremal points $\wtnu_0$, $\wtnu_1$. As $\nu$ is ergodic, $\wtnu_0$
and $\wtnu_1$ are also ergodic. Still by the ergodicity of $\nu$,
some stable curve $\om$ must meet the basins of both $\wtnu_0$ and
$\wtnu_1$. But stable curves are contracted exponentially fast under
positive iteration by $\Tim$; we should thus have $\wtnu_0 =
\wtnu_1$, a contradiction.

We have already said that $\wtnu$ is ergodic.
\end{proof}

Finally, we want to ''spread'' the $\wttm$-invariant measure $\wtnu$
in order to obtain a $g$-invariant measure $\sig$. Let $\La = \La_g$
as in the Introduction (cf.~Subsection~\ref{eq1.2}).

We first observe that the support of $\wtnu$ is contained into
$\La\cap\qqim\wtcR$: if $N\subset \qqim\wtcR$ is compact and
disjoint from $\La$, then $N$ is disjoint from the image of
$(\wttm)^j$ if $j$ is large enough, hence $\wtnu (N) = 0$.

Let now $h$ be a continuous, and thus bounded, function on $\La$.
For $x\in \La\cap \qqim\cD$, we write
\begin{equation}
\wttm (x) \; = \; g^{N(x)}(x),\label{eq10.178}
\end{equation}
where $N(x)=n$ if $x\in\qqim\wtcR (P)$ with $(P,Q,n)\in\cP$. We
define:
\begin{equation}
Sh (x) \; = \; \sum_{0\leqslant j < N(x)}\; h
(g^j(x)).\label{eq10.179}
\end{equation}
The function $Sh$ is defined $\wtnu$-almost everywhere. It
satisfies:
\begin{equation}
|Sh (x)| \; \leqslant \; ||h||_\infty N(x).\label{eq10.180}
\end{equation}
By Proposition~\ref{propo47} and Corollary~\ref{coro14} in
Subsection~\ref{sub10.3}, the function $N$ is $\wtnu$-integrable. We
have therefore defined an operator
\begin{equation}
S:\;\; C (\La) \; \mapsto \; L^1 (\wtnu).\label{eq10.181}
\end{equation}
where $C(\La)$ stands for the space of continuous functions on $\La
= \La_g$.

We define a finite measure $\sig$ on $\La$ by
\begin{equation}
\int h d\sig \; = \; \int Sh\,d\wtnu,\label{eq10.182}
\end{equation}
for $h\in C(\La)$. From the definition of $Sh$, we have
\begin{equation}
S (h \circ g) \; = \; Sh + h \circ \wttm - h.\label{eq10.183}
\end{equation}
Thus, the $\wttm$-invariance of $\wtnu$ implies that $\sig$ is
$g$-invariant. It is ergodic. The Lyapunov exponents of $\wttm$ for
$\wtnu$ are non-zero because $\wttm$ is uniformly hyperbolic. To get
the Lyapunov exponents of $g$ for $\sig$ we have only to change
time, which is possible since $N$ is $\wtnu$-integrable.

In the next and last section, we will see that in some appropriate
geometric sense, the measure $\sig$ captures ''most'' of the
dynamics on $\La$, and therefore can be considered as a naturally
defined geometric invariant measure on $\La$.

We end this section by observing that everything that has been done
for $\Tim$ and positive iteration in Section~\ref{sec10}, can also
be done for $\TiM$ and negative iteration, leading to another
naturally defined geometric invariant measure $\sig^-$ on $\La$.

\newpage

\setcounter{section}{10}
\setcounter{equation}{0}

\section{Some Further Geometric Properties of the Invariant Set\label{sec11}}
In this final section we pursue the geometric study of the invariant
set $\La = \La_g$ in two directions. First, we will describe in a
rather precise way, both from a dynamical and a geometrical point of
view, the intersection of an unstable curve in $\cR^\infty_-$, as
defined in Subsection~\ref{sub10.4}, with the invariant set $\La$.
In the second part of the section, we prove that $\La$ is a
saddle-like invariant set in the measure-theoretical sense: both its
stable and unstable sets have Lebesgue measure $0$; thus, no
attractors are present in $\La$.

\subsection{One-Dimensional Analysis of the Invariant Set\label{sub11.1}}
Let $\oms\in \cR^\infty_-$ be an unstable curve as defined in
Subsection~\ref{sub10.4}. Let $(\Pssk, \Qssk, \nssk)_{k\geqslant 0}$
be the canonical sequence associated to $\oms$ (cf.~definition also
in Subsection~\ref{sub10.4}). We have
\begin{equation}
\oms \; = \; \bigcap_{k\geqslant 0}\; \Qssk,\label{eq11.1}
\end{equation}
where $\Qsso$ is a rectangle $\Ra$ and $Q^*_{k+1}$ is a child of
$\Qssk$ for each $k\geqslant 0$. We want to analyze the intersection
$\oms \cap \La$. In Section~\ref{sec10}, we have analyzed the set
$\wtcR^\infty_+$ and we know, in particular, that $\oms \cap \La$
contains the subset $\oms \cap \wtcR^\infty_+$; this last subset has
Hausdorff dimension $d_s$ characterized in terms of the transfer
operator studied in Section~\ref{sec10}; in particular, this
dimension is independent of $\oms$.

Let us summarize the results of our analysis in this section.

\begin{Theo}\label{theo3}
The intersection $\oms \cap \La$ is the disjoint union of

--~~a, at most countable, family of Cantor sets $\La_i (\oms)$,

--~~a, at most countable, set $Cr (\oms)$,

--~~an exceptional set $\cE (\oms)$,

with the following properties

(i)~~For each $i$, there exists a piece $\oms (i)$ of $\oms$
containing $\La_i (\oms)$, an unstable curve $\oms_i$ and an integer
$n_i$ such that
\begin{eqnarray}
g^{n_i} (\oms (i)) \; &= \; &\om^*_i,\label{eq11.2}\\
g^{n_i} (\La_i (\oms)) \; &= \; &\om^*_i \cap
\wtcR^\infty_+.\label{eq11.3}
\end{eqnarray}
In particular, there is a special index $i=0$ for which $n_0 = 0$,
$\oms (0) = \om_0^* = \oms$, $\La_0 (\oms) = \oms \cap
\wtcR^\infty_+$.

(ii)~~For every point $c\in Cr (\oms)$, there exists a stable curve
$\om_+ (c) \in \cR^\infty_+$, an unstable curve $\om_- (c) \in
\cR^\infty_-$, a positive integer $n(c)$ such that $g^{n(c)} (c)$ is
a quadratic tangency point between $\om^+ (c)$ and $g^{N_0}
(\om^-(c) \cap L_u)$.

(iii)~~The Hausdorff dimension of $\cE (\oms)$ is not greater than
\begin{equation}
(\dos + \dou - 1) \; \frac{2\dos}{2\dou+\dos} + o (1)\label{eq11.4}
\end{equation}
where the $o(1)$ term is small provided $\tau$ is small enough.
Consequently, the Hausdorff dimension of $\oms \cap \La$ is equal to
$d_s$.

(iv)~~Every point $x\in \cE (\oms)$ is the intersection of a
decreasing sequence of pieces $(\oms (i_n (x)))_{n\geqslant 0}$.
\end{Theo}

\bigskip\noindent {\bf Remark.}

{\it 1.~~The structure will be made more precise in the next
subsections. We have tried here to extract the most significant
features of our analysis.

2.~~Even with $\dos+\dou > 1$, it may happen that $\La$ is a
uniformly hyperbolic horseshoe; then, the family $(\La_i (\oms))_i$
is finite, $Cr (\oms)$ and $\cE (\oms)$ are empty. When $\La$ is not
uniformly hyperbolic, the family $(\La_i (\oms))_i$ is countable and
$\cE (\oms)$ is a Cantor set; it is not clear in this case if $Cr
(\oms)$ can be empty.}

\subsection{Parabolic Cores\label{sub11.2}}
Let $(P,Q,n)\in \cR$, $\cR$ as in Subsection~\ref{sub10.1}.

\Def The {\it parabolic core} of $P$, denoted by $c (P)$, is the set
of points of $W^s (\La, \whR)$ which belong to $P$ but not to any
child of $P$. The parabolic core of $Q$, denoted by $c (Q)$, is the
set of points of $W^u (\La,\whR)$ which belong to $Q$ but not to any
child of $Q$.

We have partitions
\begin{eqnarray}
R \cap W^s (\La, \whR) \; &= \; &\bigsqcup_\cR \; c (P) \sqcup \; \wtcR^\infty_+,\label{eq11.5}\\
R \cap W^u (\La, \whR) \; &= \; &\bigsqcup_\cR \; c (Q) \sqcup \;
\wtcR^\infty_-.\label{eq11.6}
\end{eqnarray}
If $\Ra$ is the rectangle which contains $\oms$, we also have
\begin{equation}
\oms \cap \La \; = \; \bigsqcup_{P\subset \Ra}\, (\oms \cap c(P)) \;
\sqcup \; (\oms \cap \wtcR^\infty_+).\label{eq11.7}
\end{equation}
The parabolic core is empty if and only if $P$ is $I$-decomposable
for a small enough parameter interval containing the given strongly
regular parameter value. In particular, $c(P)$ is empty if $Q$ is
$I$-transverse. Thus, the union in (\ref{eq11.5}), (\ref{eq11.7})
can be restricted to those $(P,Q,n)\in \cR$ such that $Q$ is
$I$-critical for all $I$.

We will denote by $C (\oms)$ the set of elements $(P,Q,n)\in \cR$
such that $c (P) \cap \oms$ is not empty. For any $(P,Q,n)\in C
(\oms)$, $Q$ is $I$-critical for all $I$.

\subsection{Decomposition of $c(P)\cap \oms$\label{sub11.3}}
Let $(P,Q,n)\in C (\oms)$. For $k\geqslant 0$, set
\begin{eqnarray}
(\Psk, \Qsk, \nsk) \; &= \; &(\Pssk, \Qssk, \nssk) \; * \; (P,Q,n),\label{eq11.8}\\
\om^*_P \; &= \; &\bigcap_{k\geqslant 0}\; \Qsk.\label{eq11.9}
\end{eqnarray}
The unstable curve $\om^*_P$ is contained in $Q$ and we have
\begin{equation}
g^n (\oms \cap c (P)) \; \subset\; \om^*_P \cap L_u.\label{eq11.10}
\end{equation}
We define a tree $\cA (\oms, P)$ as follows. The vertices are the
rectangles $P' \subset P_s$ with the following property: for any
parameter interval $I$ (containing the given parameter value, say
$t$), for any $\Qsk \supset \om^*_P$, $\Qsk$ and $P'$ are not
$I$-separated, and $\Qsk$ and the parent of $P'$ are $I$-critically
related.

We connect two vertices by an (oriented) edge if one is the parent
of the other. We say that a vertex $P'$ is {\it critical} if, for
all $I$ and $\Qsk \supset \om^*_P$, $\Qsk$ and $P'$ are
$I$-critically related. Otherwise, we say that $P'$ is transverse.
The parent of a vertex is always a critical vertex, except if this
vertex is $P_s$, the {\it root} of the tree. When $P'$ is a
transverse vertex, the smallest integer $k$ such that $\Qsk$, $P$
are $I$-transverse for $I$ small enough is called the {\it level} of
$P'$.

Let $P'$ be a critical vertex; then, for every parameter interval
$I\ni t$, $P'$ is $I$-critical and, therefore, decomposable.

Let $P'$ be a transverse vertex of level $0$. We have $Q_0 = Q$.
Therefore, the parabolic composition $(P,Q,n)\; \square \;
(P',Q',n')$ is well defined and produces two children of $P$.

Let $P'$ be a transverse vertex of level $k>0$. For all $m\geqslant
k$, the parabolic composition $(\Psm, \Qsm, \nsm) \; \square \;
(P',Q',n')$ is well-defined
and produces two elements $(P^\pm_m, Q^\pm_m, n^\pm_m)$. The
formulas
\begin{eqnarray}
\om^*_{P,P',+} \; &:= \; &\cap\; Q^+_m,\label{eq11.11}\\
\om^*_{P,P',-} \; &:= \; &\cap\; Q^-_m,\notag
\end{eqnarray}
define unstable curves $\om^*_{P,P',\pm}$ contained in $Q'$. We also
define pieces $\oms (P, P', \pm)$ through
\begin{eqnarray}
g^{n_{P,P'}} (\oms (P, P', \pm)) \; &= \; &\om^*_{P,P',\pm},\label{eq11.12}\\
n_{P,P'} \; &:= \; &n + n' + N_0.\label{eq11.13}
\end{eqnarray}

\begin{lema}\label{lema15}
Let $x$ be a point in $\oms \cap c (P)$, $y = g^{n+N_0} (x)$. Either
$y$ belong to a transverse vertex of level $>0$ or it belongs to an
infinite decreasing sequence of critical vertices.
\end{lema}

\begin{proof}
We have $g^n (x)\in L_u$ (cf.~(\ref{eq11.10})), $y\in L_s \subset
P_s$, and $P_s$ is the root and a critical vertex of the tree $\cA
(\oms, P)$. We assume that the first possibility in the statement of
the lemma does not hold and construct, starting with $P_s$, a
sequence of critical vertices containing $y$.

Assume that $y$ belongs to a critical vertex $P'$. As $P'$ is
indecomposable and $y\in W^s (\La)$, $y$ belongs to some child
$P'_1$ of $P'$. This rectangle is a vertex of the tree: otherwise,
$\Qsk$ and $P'_1$ would be $I$-separated if $I$ and $\Qsk$  are thin
enough, and then $g^{N_0} (g^n (\oms) \cap L_u) \cap P'_1$ (which
contains $y$) would be empty. The vertex $P'_1$ cannot be transverse
of level $0$ because, as remarked above, the parabolic composition
of $(P,Q,n)$ and $(P'_1, Q'_1, n'_1)$ would produce a child of $P$
containing $x$, contradicting the hypothesis that $x\in c (P)$.
Finally $P'_1$ cannot be transverse of level $>0$ by hypothesis. It
must be a critical vertex, and the induction step is complete.
\end{proof}

\begin{Propo}\label{propo49}
There is at most one point $x\in \oms \cap c (P)$ such that $y =
g^{n+N_0} (x)$ belongs to a decreasing sequence of critical
vertices. When such a point exists, the intersection of this
decreasing sequence of vertices is a stable curve which intersects
$g^{N_0} (L_u \cap \om^*_P)$ at $y$ as a quadratic tangency point.
\end{Propo}

\begin{proof}
Let $x$ be a point in $\oms \cap c(P)$ such that $y = g^{n+N_0} (x)$
belongs to a decreasing sequence $(P'_\ell)_{\ell\geqslant 0}$ of
critical vertices. Denote by $\om_+$ the stable curve which is the
intersection of these critical vertices. For all parameter intervals
$I$, all $k\geqslant 0$, $\ell\geqslant 0$, $\Qsk$ and $P'_\ell$ are
$I$-critically related. This implies that
\begin{equation}
\lim_{\substack{k\to +\infty\\ \ell\to +\infty}} \; \de (\Qsk,
P'_\ell) \; = \; 0.\label{eq11.14}
\end{equation}
For large $k$ and $\ell$, let $\ga_k$ (resp.~($\ga'_\ell$) be the
image in $\Qsk$ (resp.~the inverse image in $P'_\ell$) of the
intersection of $\Psk$ with an horizontal curve (resp.~the
intersection of $Q'_\ell$ with a vertical curve). By
(\ref{eq11.14}), the distance between the vertical-like curve
$\ga'_\ell$ and the tip of the parabolic-like curve $g^{N_0}
(\ga_k)$ goes to zero as $k$, $\ell$ go to $+\infty$. Passing to the
limit, we see that $\om_+$ has a tangency with $g^{N_0} (\om^*_P
\cap L_u)$. This tangency is quadratic in the following sense
(cf.~also the remark after the end of the proof): First, $g^{N_0}
(\om^*_P \cap L_u)$ is contained, with the exception of the tangency
point, in one of the components of $P_s - \om_+$; moreover, the
angle between the tangent lines to $\om_+ (x)$, $g^{N_0} (L_u \cap
\om^*_P)$ at points on these curves at the same distance and on the
same side of the tangency point is of the same order as this
distance to the tangency point. This is a consequence of the uniform
estimates (\ref{eq3.21}),  (\ref{eq3.22}) in
Subsection~\ref{sub3.5}.

As $\om_+$ and $g^{N_0} (L_u \cap \om^*_P)$ meet at only one point,
this point must be $y$. If $x'$ is a point with the same property as
$x$, and we construct $\om'_+$ in the same way as $\om_+$, we must
have $\om_+ = \om'_+$ because otherwise $g^{N_0} (L_u \cap \om^*_P)
\cap \om_+$ or $g^{N_0} (L_u \cap \om^*_P) \cap \om'_+$ is empty.
But, then, we have $y' := g^{n+N_0} (x') = y$ and $x' = x$.
\end{proof}

\begin{rem}
Calculations involving partial derivatives of higher order for the
maps $(A,B)$, which implicitly represent elements of $\cR$, show
that stable curves and unstable curves are actually of class
$C^\infty$, with uniform estimates in the $C^k$ topology for all
$k$. Then, quadratic tangency can be taken in the usual sense.
However, the calculations involved, especially when considering
parabolic composition, are quite long and not very interesting; we
decided to stick to the $C^{1+\lip}$ regularity class, where the
notion of ''quadratic'' tangency, as explained in the proof of
Proposition~\ref{propo49}, still makes sense.
\end{rem}

It is easy to see exactly  when a point $x\in\oms \cap c (P)$ with
the property specified in Proposition~\ref{propo49} does exist: a
necessary and sufficient condition is that the tree $\cA (\oms, P)$
is infinite. In this case, the point $x$ will be a point of the set
$Cr (\oms)$ in the statement of Theorem~\ref{theo3} and the point $y
= g^{n+N_0} (x)$ is said to be {\it critical}.

Summarizing what we have established so far, two cases may happen:

1)~~The tree $\cA (\oms,P)$ is finite. Then, the intersection $\oms
\cap c (P)$ is the finite disjoint union of the sets
\begin{equation}
\oms (P, P',  \pm) \; \cap \; \La\label{eq11.15}
\end{equation}
where $P'$ runs through the vertices of the tree which are
transverse of level $>0$. The image under $g^{n_{P,P'}}$ of the set
(\ref{eq11.15}) is the intersection $\om^*_{P,P',\pm} \cap \La$.

2)~~The tree $\cA (\oms, P)$ is infinite. Then, the intersection
$\oms\cap c (P)$ is the countable disjoint union of the sets $\oms
(P, P', \pm) \cap \La$ as above and a single point $x\in Cr (\oms)$.
The point $x = x_P$ is the limit of the pieces $\oms (P, P', \pm)$
(whose diameters goes to $0$ as $|P'|$ goes to $0$).

\subsection{The Structure of $\om^* \cap \La$\label{sub11.4}}
We are now ready to prove all the statements in Theorem~\ref{theo3},
stated above in Subsection~\ref{sub11.1}, with the exception of (iii)
(the estimate on the Hausdorff dimension of $\cE (\oms)$).

The structure of $\oms \cap \La$ that we are looking for, which is
roughly described in Theorem~\ref{theo3}, is obtained by iterating
the partition (\ref{eq11.7}) and the decomposition of $\oms \cap
c(P)$ described in Subsection~\ref{sub11.3}.

At the first step, we have partitioned $\oms \cap \La$ into the
following subsets:

--~~the intersection $\oms \cap \wtcR^\infty_+$; points in this set
ar said of type~I;

--~~for each $(P,Q,n)\in C (\oms)$ such that $\cA (P, \oms)$ is
infinite, a point $x_P$ such that $y_P = g^{n+N_0} (x_P)$ is
critical; such points $x_P$ are said of type~II;

--~~for each $(P,Q,n)\in C (\oms)$, each vertex $(P',Q',n')$ of $\cA
(\oms, P)$ which is transverse of level bigger than $0$, each
$\vep\in \{ +, -\}$, the intersection $\oms (P, P', \vep) \cap \La$;
the image of this set under $g^{n_{P,P'}}$ is the intersection
$\om^*_{P,P',\vep} \cap \La$ of another unstable curve with $\La$.

The intersection $\om^*_{P,P',\vep} \cap \La$ will be analyzed in
the same way that $\oms \cap \La$.

Consider a point $z_0\in \oms \cap \La$. If it is of type I, it
belongs to the set $\La_0 (\oms) := \oms \cap \wtcR^\infty_+$ of the
statement of Theorem~\ref{theo3}. If it is of type~II, it belongs to
$Cr (\oms)$. Assume now that it is of type~III. Then, it belongs to
some $\oms (P,P',\vep) \cap \La$ as above. Define
\begin{equation}
z_1 \; = \; g^{n_{P,P'}} \, (z_0),\label{eq11.16}
\end{equation}
which belongs to $\om^*_{P,P',\vep} \cap \La =: \om^*_1$. This point
may in turn be of type I, II, III with respect to $\om^*_1$. The
process stops if $z_1$ is of type I or II; if $z_1$ is of type III,
it belongs to some piece $\om^*_1 (P_1, P'_1, \vep_1)$; we define
\begin{equation}
z_2 \; = \; g^{n_{P_1,P'_1}} \, (z_1),\label{eq11.17}
\end{equation}
which belongs to $\om^*_2 \cap \La$, with
\begin{equation}
\om^*_2 \; := \; g^{n_{P_1,P'_1}} \, (\om^*_1 (P_1, P'_1,
\vep)).\label{eq11.18}
\end{equation}
Iterating this process lead to one of three possible outcomes:

1)~~the $z_k$'s are defined and of type III for all $k\geqslant 0$;
the corresponding initial points $z_0$ form the set $\cE (\oms)$.

2)~~the $z_k$'s are defined for $0\leqslant k\leqslant \ell$ and
$z_\ell$ is of type I, i.e. it belongs to $\wtcR^\infty_+$; let
$(\Psk, \Pskl, \vepk)$ for $0\leqslant k < \ell$ be the data
involved in the definitions of the $\zk$'s. We collect together the
initial points $\zo$'s with the same set of data; such a set form
one of the Cantor sets $\La_i (\oms)$ in Theorem~\ref{theo3}.

3)~~the $\zk$'s are defined for $0\leqslant k\leqslant \ell$ and
$\zell$ is of type II. Then $\zo$ belongs to the set $Cr (\oms)$.

We have now completely defined the partition of $\oms \cap \La$
described in Theorem~\ref{theo3}. The properties (i), (ii), (iv)
follow immediately from the definitions.

\subsection{Hausdorff Dimension of the Exceptional Set $\cE (\oms)$\label{sub11.5}}
The self-similar structure apparent in the definition of $\cE
(\oms)$ is the key to obtain an estimate of the dimension of this
set. More specifically, we have
\begin{equation}
\cE (\oms) \; = \; \bigsqcup_{(P,P',\vep)} \, g^{-n_{P,P'}}\, (\cE
(\oms_{P,P',\vep})),\label{eq11.19}
\end{equation}
where $\vep\in \{+,-\}$, $P$ runs through $C (\oms)$ and $P'$
through vertices of $\cA (\oms, P)$ which are transverse of level
$>0$.

\begin{lema}\label{lema16}
The maps
\[
g^{n_{P,P'}}:\; \oms (P,P',\vep) \; \to \; \oms_{P,P',\vep}
\]
have uniformly bounded distortion.
\end{lema}

\begin{proof}
Let $k$ be an integer larger than the level of the transverse vertex
$P'$. Then, the parabolic composition of $(\Psk, \Qsk, \nsk)$
(cf.~(\ref{eq11.8})) and $(P',Q',n')$ is defined and produces an
element $(\Pskl, \Qskl, \nskl)$ such that $\Qskl$ contains
$\om^*_{P,P',\vep}$. Let $\ga^*_k$ be an horizontal segment in
$\Pssk$, $\ga_k$ its image under $g^{\nssk}$, $\gaskl$ the image of
$\ga^*_k \cap \Pskl$ under $g^{\nskl}$.

The affine-like maps
\begin{equation}
g^{\nssk}\; : \; \Pssk \to \Qssk,\qquad\qquad g^{\nskl}\; : \; \Pskl
\to \Qskl,\label{eq11.20}
\end{equation}
have bounded distortion, hence the one-dimensional map
\begin{equation}
g^{\nssk}\; \circ \; (g^{\nskl})^{-1}\; : \; \gaskl \to
\gask\label{eq11.21}
\end{equation}
have also uniformly bounded distortion. Letting $k$ go to $+\infty$,
$\gaskl$ converge to $\om^*_{P,P',\vep}$ and $\gask$ to $\oms$ in
the $C^{2-\vep}$-topology for all $\vep > 0$. The statement of the
lemma follows.
\end{proof}

\begin{lema}\label{lema17}
Let
\[
\de (\om^*_P, P') \; = \; \lim_{k\to +\infty} \; \de (\Qsk, P')
\]
We have
\[
C^{-1} \; \leqslant \; \frac{\diam\,\oms (P,P',\vep)}{|P|\,|P'|(\de
(\om^*_P,P'))^{-\fudt}} \; \leqslant \; C
\]
\end{lema}

\begin{proof}
As in the proof of Lemma~\ref{lema16}, we write
\begin{equation}
g^{n_{P,P'}}= g^{\nskl} \circ (g^{n^*_k})^{-1}.\label{eq11.22}
\end{equation}
From the estimate (\ref{eq3.27}) for parabolic composition in
Subsection~\ref{sub3.5}, we have
\begin{equation}
C^{-1} \; \leqslant \; \frac{\diam\, \gassk}{|\Psk|,|P'|(\de
(\Qsk,P'))^{-\fudt}} \; \leqslant \; C\label{eq11.23}
\end{equation}
We also have, from the estimates on simple composition
\begin{eqnarray}
C^{-1} \; &\leqslant \; &\frac{\diam\,\oms
(P,P',\vep)\,|\Pssk|}{\diam\,\gassk} \; \leqslant \; C,\label{eq11.24}\\
C^{-1} \; &\leqslant \; &\frac{|\Psk|}{|P|\,|\Pssk|} \; \leqslant \;
C.\label{eq11.25}
\end{eqnarray}
Multiplying these three inequalities yields the Lemma.
\end{proof}

Let us introduce
\begin{equation}
\chi (d) \; = \; \sum_{(P,P',\vep)} \, [\diam\, \oms (P,P',\vep)]^d.
\label{eq11.26}
\end{equation}
If we are able, for some value of $d$, to show that the series
defining $\chi$ is convergent and $\chi (d)$ is small, then by
(\ref{eq11.19}) and Lemma~\ref{lema16}, we will deduce that the
Hausdorff dimension of $\cE (\oms)$ is $\leqslant d$.

In order to study $\chi$, we will first fix $P$ in $C (\oms)$ and
sum over $(P',\vep)$. As $\vep$ takes only two values, and in view
of Lemma~\ref{lema17}, we define, for $P\in C (\oms)$:
\begin{equation}
\chi_P (d) \; = \; \sum_{P'} \, |P'|^d \, \de (\om^*_P,
P')^{-\fudt\,d}.\label{eq11.27}
\end{equation}
We will then have
\begin{equation}
\chi (d) \; \leqslant \; C\, \sum_{P} \, |P|^d \, \chi_P
(d).\label{eq11.28}
\end{equation}
In the sum (\ref{eq11.27}), $P'$ is a transverse vertex of level
$>0$, and we therefore must have
\begin{equation}
\de (\om^*_P, P') \; \leqslant \; \de_{\max} \; := \; \min (\vep_0,
C |Q|^{1-\eta}).\label{eq11.29}
\end{equation}
In the series (\ref{eq11.27}), we first sum over those $P'$ such
that
\begin{equation}
2^{-\ell}\; \de_{\max} \; \geqslant \; \de (\om^*_P, P') \;
\geqslant \; 2^{-\ell-1}\, \de_{\max}\label{eq11.30}
\end{equation}
for some fixed $\ell\geqslant 0$. This allow us to write
\begin{equation}
\chi_P (d) \; \leqslant \; C \, \de_{\max}^{-\fudt\,d} \,
\sum_{\ell\geqslant 0} \, 2^{\tfrac{\ell d}{2}}
\Bigl(\sum\nolimits^{(\ell)} |P'|^d\Bigr),\label{eq11.31}
\end{equation}
where $\sum^{(\ell)}$ means that $P'$ is constrained by
(\ref{eq11.30}). We divide $\sum^{(\ell)}$ into two parts.

In the first part, denoted by $\sum^{(\ell)}_1$, we consider only
those $P'$ such that its parent $\wtP'$ satisfies
\begin{equation}
|\wtP'| \; \leqslant \; 2^{-\ell}\, \de_{\max}.\label{eq11.32}
\end{equation}
To estimate $\sum^{(\ell)}_1 |P'|^d$, first observe that, with $d$
bounded away from $0$, it follows from Proposition~\ref{propo21} in
Subsection~\ref{sub8.1} that the sum of $|P'|^d$ over children of a
fixed parent $\wtP'$ is bounded by $C |\wtP'|^d$. We must therefore
bound $\sum^{(\ell)}_1 |\wtP'|^d$.

Also, as $\wtP'$ is a critical vertex, $\wtP'$ cannot be very thin:
from Proposition~\ref{propo10} in Subsection~\ref{sub6.6}, we have
\begin{equation}
|\wtP'| \; \geqslant \; C^{-1} [\de (\om^*_P, P')]^{(1-\eta)^{-1}}
\; \geqslant \; C^{-1} (\de_{\max}
2^{-\ell})^{(1-\eta)^{-1}}.\label{eq11.33}
\end{equation}
Finally, the number of $\wtP'$ with $|\wtP'|$ of order
$2^{-m-\ell}\, \de_{\max}$ is at most $C2^m$ and the integer $m$
here is restricted by (\ref{eq11.33}) to the range
\begin{equation}
1\; \leqslant \; 2^m \; \leqslant \; C (\de_{\max} 2^{-\ell})^{-\eta
(1-\eta)^{-1}}.\label{eq11.34}
\end{equation}
We, therefore, obtain for $d$ bounded away from 0 and 1,
\begin{eqnarray}
\sum\nolimits^{(\ell)}_1 \, |P'|^d \; &\leqslant \; &C \sum\nolimits^{(\ell)}_1 \,
|\wtP'|^d\label{eq11.35}\\
&\leqslant \; &C \, \de^d_{\max} \, 2^{-\ell d}\, \sum_m \,
2^{m(1-d)}\notag\\
&\leqslant \; &C \, (\de_{\max} \, 2^{-\ell})^{d-\eta}.\notag
\end{eqnarray}
In the second part of $\sum^{(\ell)}$, denoted by $\sum^{(\ell)}_2$,
we have on the opposite
\begin{equation}
|\wtP'|\; > \; 2^{-\ell} \, \de_{\max}.\label{eq11.36}
\end{equation}
As $|\wtP'| > \de (\om^*_P, P')$, the number of possibilities for
$\wtP'$ is now bounded. As each $P'$ is a transverse vertex, we must
have (by (R7))
\begin{equation}
|P'|\; \leqslant \; C (2^{-\ell} \,
\de_{\max})^{(1-\eta)^{-1}}.\label{eq11.37}
\end{equation}
In particular, from (\ref{eq11.36}),  (\ref{eq11.37}), $P'$ is a
non-simple child of $\wtP'$. From Proposition~\ref{propo21} in
Subsection~\ref{sub8.1}, the number of $P'$ with $|P'|$ of order
$2^{-m}\vep_0$ is at most $2^{Cm\eta}$.

We have
\begin{equation}
\sum\nolimits^{(\ell)}_2 \, |P'|^d \; \leqslant \; \vep_0^d \,
\sum_m 2^{-m(d-C\eta)} \leqslant \;C \vep_0^{-C\eta} (2^{-\ell}
\de_{\max})^{d-C\eta}.\label{eq11.38}
\end{equation}
Putting (\ref{eq11.35}) and (\ref{eq11.38}) together yields
\begin{equation}
\sum^\ell |P'|^d \; \leqslant \; C (\de_{\max}\,
2^{-\ell})^{d-C'\eta}\label{eq11.39}
\end{equation}
and introducing this in (\ref{eq11.31}) allow us to estimate
$\chi_P$:
\begin{equation}
\chi_P (d) \; \leqslant \; C \de_{\max}^{\fudt\, d
-C'\eta}.\label{eq11.40}
\end{equation}
Finally, we obtain
\begin{equation}
\chi (d) \; \leqslant \; C \sum_{C (\oms)} |P|^d [\min (\vep_0,
|Q|)]^{\fudt\, d - C\eta}.\label{eq11.41}
\end{equation}
We do not know exactly the set $C (\oms)$, but we know that if
$(P,Q,n)\in C (\oms)$, the parabolic core $c (P)$ is non-empty and
$Q$ must be $I$-critical for all parameter intervals $I$ containing
the given parameter value.

We use H\"older's inequality to separate the $P$ and $Q$ in
(\ref{eq11.41}): for any $p$, $q > 1$ such that
\begin{equation}
\frac{1}{p} + \frac{1}{q} \; = \; 1,\label{eq11.42}
\end{equation}
we have
\begin{equation}
\chi (d) \; \leqslant \; C \chi_+ (d)^{\tfrac{1}{p}} \chi_-
(d)^{\tfrac{1}{q}}\label{eq11.43}
\end{equation}
where
\begin{eqnarray}
\chi_+ (d) \; = \; \sum_{Q \text{ critical}}\,
|P|^{dp},\label{eq11.44}\\
\chi_- (d) \; = \; \sum_{Q \text{ critical}}\, \min (\vep_0 \, , \,
|Q|)^{\left(\fudt\,d - C\eta\right)\,q}.\label{eq11.45}
\end{eqnarray}
We will choose $d$, $p$, $q$ (satisfying (\ref{eq11.42}) in order to
have $\chi_+ (d)$ bounded and $\chi_- (d)$ small (when $\vep_0$ is
small). For such a choice, we can conclude that the Hausdorff
dimension of $\cE (\oms)$ is $\leqslant d$.

We now use that the parameter value is strongly regular, more
precisely that the eight estimates (SR1), (SR2) of
Subsection~\ref{sub9.2} on the size of the critical locus are
satisfied.

It is not difficult to deduce from (SR1)$_\whu$ that, if
\begin{equation}
\fud\, dq \; > \; \dos + \dou - 1\label{eq11.46}
\end{equation}
then, $\chi_- (d)$ will be small.

From (SR2)$_\whu$, one can also deduce that if
\begin{equation}
dp\; \frac{\dou}{\dos} \; > \; \dos + \dou - 1\label{eq11.47}
\end{equation}
then $\chi_+$ will be bounded. The relations (\ref{eq11.42}),
(\ref{eq11.46}), (\ref{eq11.47}) are compatible exactly when
\begin{equation}
d \; > \; (\dos + \dou - 1) \;
\frac{2\dos}{2\dou+\dos}.\label{eq11.48}
\end{equation}
We observe that the right hand size is always $<\dos$. This ends the
proof of Theorem~\ref{theo3}.

\begin{rem}
The inequalities (\ref{eq11.46}), (\ref{eq11.47}) should be
understood in the following sense: the difference between the left
and right-hand sides is much larger than $\tau$ (which is itself
much larger than $\eta$).
\end{rem}

\subsection{The Stable and Unstable Sets of $\La$\label{sub11.6}}
Our goal at the end of this final section is to prove that the
invariant set $\La$ is a saddle-like object in the following
measure-theoretical sense:

\begin{Theo}\label{theo4}
For a strongly regular parameter, both the stable set $W^s (\La)$
and the unstable set $W^u (\La)$ have Lebesgue measure $0$.
\end{Theo}
The situation is symmetrical and we will deal with the stable set.

We have:
\begin{equation}
W^s (\La) \; = \; \bigsqcup_{n\geqslant 0}\, g^{-n} (W^s (\La,\whR) \cap R).
\label{eq11.49}
\end{equation}
Therefore, it is sufficient to show that $W^s (\La,\whR) \cap R$ has Lebesgue
measure 0. We write
\begin{equation}
R \cap W^s (\La, \whR) \; = \; \bigsqcup_{n\geqslant 0}\, \Bigl( W^s (\La, \whR) \cap R
\cap g^{-n} (\wtcR^\infty_+)\Bigr) \sqcup \cE^+,\label{eq11.50}
\end{equation}
with
\begin{equation}
\cE^+ \; = \; \{ z\in W^s (\La, \whR) \cap R \, , g^n (z) \not\in \wtcR^\infty_+ \;
\text{ for all } \; n\geqslant 0\}.\label{eq11.51}
\end{equation}
We have seen in Section~\ref{sec10} that $\wtcR^\infty_+$ is
Lipschitzian with transverse Hausdorff dimension $\ds$. Therefore,
the Hausdorff dimension of $\wtcR^\infty_+$ is $1+\ds$ and its
Lebesgue measure is $0$. The same is true of $g^{-n}
(\wtcR^\infty_+)$. We have to prove that the Lebesgue measure of
$\cE^+$ is equal to $0$.

\subsection{Decomposition of $\cE^+$\label{sub11.7}}
By the definition of $\cE^+$ and of the parabolic cores, we can
write
\begin{equation}
\cE^+ \; = \; \bigsqcup_{P_0} \, \cE^+ (P_0),\label{eq11.52}
\end{equation}
where
\begin{equation}
\cE^+  (P_0) \; = \; \cE^+ \cap  c(P_0)\label{eq11.53}
\end{equation}
and $(P_0, Q_0, n_0)$ runs through the set $\cC_-$ of elements of
$R$ with $c (P_0) \ne\emptyset$. In particular, $Q_0$ is
$I$-critical for all $I$ containing the given parameter value.

For any such $P_0$, we have
\begin{eqnarray}
g^{n_0} (\cE^+ (P_0)) &\subset &Q_0 \cap L_u \cap
\cE^+\label{eq11.54}\\
g^{n_0+N_0} (\cE^+ (P_0)) &\subset &L_s \cap \cE^+.\label{eq11.55}
\end{eqnarray}
For $P_1\in \cC_-$, define
\begin{equation}
\cE^+  (P_0, P_1) \; = \; \{ z\in \cE^+ (P_0), g^{n_0+N_0} (z) \in c
(P_1)\}.\label{eq11.56}
\end{equation}
We have a partition
\begin{equation}
\cE^+  (P_0) \; = \; \bigsqcup_{P_1} \, \cE^+ (P_0,
P_1).\label{eq11.57}
\end{equation}
At step $k$, we have a partition
\begin{equation}
\cE^+  \; = \; \bigsqcup_{P_0,\cdots,P_k} \, \cE^+ (P_0,\cdots,
P_k)\label{eq11.58}
\end{equation}
where the $(P_i, Q_i, n_i)$ runs through $\cC_-$. We write
\begin{eqnarray}
m_0 \; &= \; &n_0,\label{eq11.59}\\
m_1 \; &= \; &n_0 + N_0 + n_1\notag\\
&\vdots &\notag\\
m_j \; &= \; &n_0 + N_0 + n_1 + N_0 + \cdots + n_{j-1} + N_0 +
n_j\notag\\
&= \; &m_{j-1} + N_0 + n_j.\notag
\end{eqnarray}
For $0\leqslant j\leqslant k$, we have
\begin{eqnarray}
g^{m_j} (\cE^+ (P_0,\cdots,P_k)) &\subset &Q_j \cap L_u \cap \cE^+,\label{eq11.60}\\
g^{m_j+N_0} (\cE^+ (P_0,\cdots,P_k)) &\subset &L_s \cap
\cE^+.\label{eq11.61}
\end{eqnarray}
We define, for $P_{k+1}\in \cC_-$
\begin{equation}
\cE^+  (P_0,\cdots,\Psk, P_{k+1}) \; = \; \{ z\in \cE^+ (P_0,\cdots,
P_k), g^{m_k+N_0}(z)\in c (P_{k+1})\}\label{eq11.62}
\end{equation}
and we have
\begin{equation}
\cE^+  (P_0,\cdots,\Psk) \; = \; \bigsqcup_{P_{k+1}} \, \cE^+
(P_0,\cdots, P_k, P_{k+1}).\label{eq11.63}
\end{equation}
However, in order to have $\cE^+ (P_0,\cdots,\Psk) \ne \emptyset$
strong restrictions on the $P_i$ must take place. We have already
mentioned that $(P_i, Q_i, n_i)\in \cC_-$. This is the only
restriction on $(P_0, Q_0, n_0)$. But, from (\ref{eq11.55}), $P_1$
must meet $P_s$ and we also know that $Q_1$ is critical. As the
parameter is regular, we must have
\begin{equation}
\max (|P_1|, |Q_1|) \; \leqslant \; \vep_0^\be.\label{eq11.64}
\end{equation}
Assume that $\cE^+ (P_0,\cdots,P_{k+1})$ is not-empty. We already
know that $Q_{k+1}$ is critical. It is also true that, for any
parameter interval $I$ containing the given parameter value, $\Qsk$
and $P_{k+1}$ cannot be $I$-transverse: if they were, parabolic
composition would produce children of $\Psk$ whose union contains
$\cE_+ (P_0,\cdots,P_{k+1})$, in contradiction with the definition
of the parabolic core of $\Psk$.

Let $I$ be the smallest parameter interval satisfying
\begin{equation}
|I|\; \geqslant \; 2 (\max\, |\Qsk|, \,
|P_{k+1}|)^{1-\eta}.\label{eq11.65}
\end{equation}
As $\Qsk$, $P_{k+1}$ are not $I$-transverse, it follows from
Proposition~\ref{propo9} in Subsection~\ref{sub6.6} that $P_{k+1}$
is $I$-critical. Therefore, we must have
\begin{equation}
\max\, (|P_{k+1}|, |Q_{k+1}|) \; < \; |I|^\be\label{eq11.66}
\end{equation}
which implies
\begin{equation}
\max\, (|P_{k+1}|, |Q_{k+1}|) \; \leqslant \; C
|\Qsk|^{\wtbe},\label{eq11.67}
\end{equation}
with $\wtbe = \be (1-\eta)(1+\tau)^{-1}$. Taking $\whbe < \wtbe$ but
close to $\be$ and $\vep_0$ sufficiently small, (\ref{eq11.67}) and
(\ref{eq11.64}) give
\begin{equation}
\max\, (|P_j|, |Q_j|) \; \leqslant \;
\vep_0^{\whbe^j}.\label{eq11.68}
\end{equation}

\subsection{Size and Area of Parabolic Cores\label{sub11.8}}
\begin{Propo}\label{propo50}
Let $(P,Q,n)\in \cC_-$. With Leb standing for Lebesgue measure, we
have
\begin{eqnarray}
\diam (g^n (c(P))) \; &\leqslant \; &C
|Q|^{\fudt\,(1-\eta)}\label{eq11.69}\\
\leb (g^n (c(P))) \; &\leqslant \; &C
|Q|^{\tfrac{3}{2}-\fudt\,\eta}\label{eq11.70}\\
\leb (c(P)) \; &\leqslant \; &C
|P|\,|Q|^{\fudt\,(1-\eta)}.\label{eq11.71}
\end{eqnarray}
\end{Propo}

\begin{rem}
A posteriori, $c(P)$, which is contained in $W^s (\La, \whR)$, will have
zero Lebesgue measure. However, we estimate here the diameter and
Lebesgue measure of a larger set, as will be apparent in the proof.
\end{rem}

\begin{proof}
We start with a general observation on an affine-map with implicit
representation $(A,B)$. The Jacobian of the map is the product
$A_x^{-1} B_y$. The distortion of Lebesgue measure under the map,
which is produced by the oscillation of the logarithm of the
Jacobian, is, therefore, controlled by the distortion of the
affine-like map in the sense of Subsection~\ref{sub3.2}. In
particular, the distortion of Lebesgue measure by the restriction of
iterates corresponding to the elements of $\cR$ is uniformly
bounded.

Thus, the third inequality (\ref{eq11.71}) in the proposition is a
consequence of the second. On the other hand, as $g^n (c(P)) \subset
Q$, the second inequality (\ref{eq11.70}) is an obvious consequence
of the first. We have, therefore, only to prove (\ref{eq11.69}). Set
\begin{equation}
Z \; = \; g^{n+N_0} (c(P)).\label{eq11.72}
\end{equation}
This set is contained in $g^{N_0} (Q\cap L_u) \cap W^s (\La, \whR)$,
and, a fortiori, in $P_s$. Let $R(Z)$ be the set of $(\Ps, \Qs,
\ns)\in \cR$ with $\Ps \cap Z \ne \emptyset$, $\Ps\subset P_s$.

By the definition of the parabolic core, $Q$ and $\Ps$ cannot be
$I$-transverse (for any parameter interval $I$). Assume that $|\Ps|
> 2 |Q|$; it follows from Proposition~\ref{propo10} in
Subsection~\ref{sub6.6} that $\Ps$ is $I$-critical (for any $I$);
this implies that $\Ps$ is $I$-decomposable if $I$ is small enough.

We conclude that
\begin{equation}
Z \subset \bigcup_{\substack{(\Ps,\Qs,\ns)\in R(Z)\\
|\Ps|\leqslant 2|Q|}} \; [\Ps \cap g^{N_0} (Q\cap
L_u)].\label{eq11.73}
\end{equation}
We replace $\Ps \cap g^{N_0} (Q\cap L_u)$ by the larger Jordan
domain $V (\Ps, Q)$ defined as follows (see figure~9): the boundary
of $V (\Ps,Q)$ is made of one arc in the boundary of $\Ps$ and one
arc in the boundary of $g^{N_0} (Q\cap L_u)$ and the interior of $V
(\Ps,Q)$ meets both $\Ps$ and $g^{N_0} (Q\cap L_u)$.

\begin{center}
\includegraphics{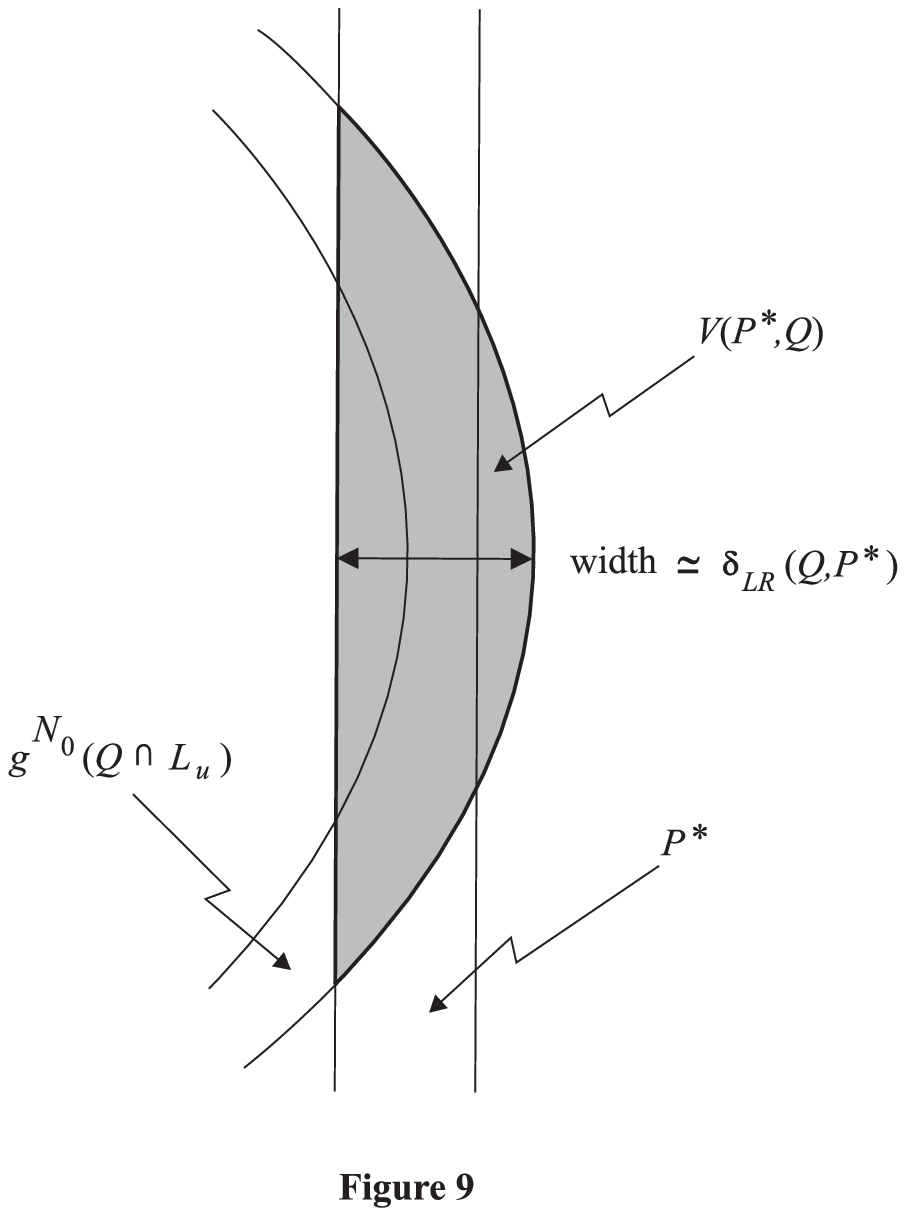}
\end{center}

Observe that, if $(P^*_1, Q^*_1, n^*_1)$,  $(P^*_2, Q^*_2, n^*_2)$
belong to $R(Z)$, we have either $V (P^*_1, Q) \subset V (P^*_2, Q)$
or the opposite inclusion. Moreover, the parabolic geometry of the
picture (cf.~(\ref{eq3.22})) gives the estimate
\begin{equation}
\diam \, V (\Ps, Q) \; \leqslant \; C (\de_{LR} (Q,
\Ps))^{\fudt}.\label{eq11.74}
\end{equation}
But, for $|\Ps| \leqslant 2 |Q|$, the fact that $Q$ and $\Ps$ are
not $I$-transverse (for any $I$) implies that
\begin{equation}
\de_{LR} (Q, \Ps)) \; \leqslant \; C |Q|^{1-\eta}.\label{eq11.75}
\end{equation}
We conclude that
\begin{eqnarray}
\diam\, Z \; &\leqslant \; &\sup_{\substack{(\Ps,\Qs,\ns)\in R(Z)\\
|\Ps|\leqslant 2|Q|}} \; \diam \, V (\Ps, Q)\label{eq11.76}\\
&\leqslant \; &C |Q|^{\fudt\, (1-\eta)}\notag
\end{eqnarray}
Taking the image by the fixed map $g^{-N_0}$ yields (\ref{eq11.69}).
\end{proof}

\subsection{Proof of Theorem~4\label{sub11.9}}
We will estimate first the Lebesgue measure of each domain $\cE^+
(P_0, \cdots, \Psk)$. We have
\begin{equation}
g^{m_{k-1} + N_0} \, (\cE^+ (P_0, \cdots, \Psk)) \subset c
(\Psk).\label{eq11.77}
\end{equation}
We now use that both the fixed map $g^{N_0}$ and the affine-like
iterates $g^{n_j}: P_j \to Q_j$ (for $0\leqslant j < k$) have
uniformly bounded distortion with respect to Lebesgue measure. We
are, therefore, able to deduce from (\ref{eq11.70}) in
Proposition~\ref{propo50} that
\begin{equation}
\leb (\cE^+ (P_0,\cdots, \Psk)) \; \leqslant \; C^{k+1}
|\Qsk|^{\tfrac{3}{2}-\tfrac{1}{2}\,\eta} \; \prod^k_0 \;
\frac{|P_j|}{|Q_j|}.\label{eq11.78}
\end{equation}
By (\ref{eq11.67}), we have $|P_{j+1}| \ll |Q_j|$ for $j>0$ \, and it is
easy to check that this still holds for $j=0$ (using (\ref{eq11.64})
if $|Q_0|\geqslant \vep_0$, (\ref{eq11.67}) otherwise). It then
follows from (\ref{eq11.78}) that we have (for $k>0$)
\begin{equation}
\leb (\cE^+ (P_0,\cdots, \Psk)) \; \ll \; |P_0|\, |\Qsk|^{\fudt\,
(1-\eta)}.\label{eq11.79}
\end{equation}
To obtain the estimate for $\cE^+$, we have to sum over sequences
$(P_0,\cdots,\Psk)$. We first estimate, when $(P_0, Q_0, n_0)$ and
$(\Psk, \Qsk, \nsk)$ are fixed, how many admissible sequences have
these two extremities.

The element $(P_{k-1}, Q_{k-1}, n_{k-1})$ must satisfy
\begin{equation}
|Q_{k-1}|\; \geqslant \; C^{-1} \max
(|\Psk|\,|\Qsk|)^{1/\widetilde\beta}\label{eq11.80}
\end{equation}
and also that $Q_{k-1}$, $\Psk$ are not $I$-transverse (for any
$I$). This implies that there is at most
\begin{equation}
C |Q_k|^{-\eta/\widetilde\beta}\label{eq11.81}
\end{equation}
possibilities for $(P_{k-1}, Q_{k-1}, n_{k-1})$. Iterating leads to
a total number of sequences (with $P_0$, $Q_k$ fixed) which is at
most
\begin{equation}
C |Q_k|^{-C\eta}.\label{eq11.82}
\end{equation}
Therefore, we obtain
\begin{equation}
\leb\, (\cE^+) \; \leqslant \; \sum_{P_0,Q_k}\,
|P_0|\,|Q_k|^{\fudt-C\eta}.\label{eq11.83}
\end{equation}
The sum $\sum |P_0|$ is obviously bounded. The sum $\sum_{Q_k}
|Q_k|^{\fud-C\eta}$ is arbitrarily small (when $k$ is large).
Indeed, we have from (\ref{eq11.68})
\begin{equation}
|Q_k| \; \leqslant \; \vep_0^{\whbe^k}.\label{eq11.84}
\end{equation}
On the other hand, recall that $\Qsk$ is $I$-critical (for any $I$).
As the parameter is strongly regular, we know, for instance, from
(SR1)$_\whu$, that the series
\begin{equation}
\sum_{Q \text{ critical}} \, |Q|^d\label{eq11.85}
\end{equation}
is convergent if $d > \dos + \dou - 1 + o(1)$. As the maximal value
taken by $\dos + \dou - 1$ under hypothesis (H4) is $\frac{1}{5}$,
we indeed have
\begin{equation}
\lim_{k\to +\infty} \, \sum_{Q_k} |Q_k|^{\fudt-C\eta} \; = \;
0,\label{eq11.86}
\end{equation}
and this concludes the proof of Theorem~\ref{theo4} and, thus, of
the Main Theorem in the paper.

We can sum up the results in Sections~\ref{sec10} and \ref{sec11} by
rephrasing our main result as follows:

\begin{Theo}\label{theo5}
Assume (H1)--(H4). Then, for most $g\in \cU_+$, $\La_g \subset W^s
(\La_g)$ and $\La_g \subset W^u (\La_g)$ carry geometric invariant
measures, à la Sinai-Ruelle-Bowen \cite{Si,Ru,BR}, with non-zero
Lyapunov exponents. Both $W^s (\La_g)$ and $W^u (\La_g)$ have
Lebesgue measure zero and thus $\La_g$ carries no attractors nor
repellors.
\end{Theo}

\newpage

\setcounter{equation}{0}
\appendix

\section*{Appendix~A\\
Composition Formulas for Affine-Like Maps}\label{app-A}

We mostly recall in this appendix the formulas for the simple and
parabolic compositions of implicitly defined affine-like maps.

We follow closely \cite{PY2}. The main difference with \cite{PY2} is
that we consider maps depending on a parameter $t$, and we are
interested also in some partial derivatives with respect to the
parameter.

\medskip\noindent
{\bf A.1~Simple Composition.~} Here we consider a map $F_t: (x_0,
y_0) \mapsto (x_1, y_1)$ implicitly defined by
\begin{equation}
\begin{cases}
x_0 = A (y_0, x_1, t)\\
y_1 = B (y_0, x_1, t).
\end{cases}\tag{A.1}
\end{equation}
and a map $F'_t: (x_1, y_1) \mapsto (x_2, y_2)$ implicitly defined
by
\begin{equation}
\begin{cases}
x_1 = A' (y_1, x_2, t)\\
y_2 = B' (y_1, x_2, t).
\end{cases}\tag{A.2}
\end{equation}
The composition $F''_t = F'_t \circ F_t$ is implicitly defined by
\begin{equation}
\begin{cases}
x_0 = A'' (y_0, x_2, t)\\
y_2 = B'' (y_0, x_2, t)
\end{cases}\tag{A.3}
\end{equation}
and we want to relate the partial derivatives of $A''$, $B''$ to
those of $A$, $B$, $A'$, $B'$. Set
\begin{equation}
\De\, := \,1 - A'_y (y_1, x_2, t) B_x (y_0, x_1, t).\tag{A.4}
\end{equation}
When we solve the system (A.1), (A.2) for $x_1$, $y_1$, we obtain
\begin{equation}
\begin{cases}
x_1 = X (y_0, x_2, t)\\
y_1 = Y (y_0, x_2, t)
\end{cases}\tag{A.5}
\end{equation}
where the partial derivatives of $X$, $Y$ are given by
\begin{equation}
\begin{cases}
X_x = A'_x\,\dmu \\
X_y = A'_y B_y\,\dmu \\
X_t = (A'_t + A'_y B_t)\,\dmu \\
Y_x = A'_x B_x\,\dmu \\
Y_y = B_y\,\dmu \\
Y_t = (B_t + A'_t B_x)\,\dmu.
\end{cases}\tag{A.6}
\end{equation}
We have
\begin{equation}
\begin{cases}
A'' (y_0, x_2, t) = A (y_0, X, t)\\
B'' (y_0, x_2, t) = B' (Y, x_2, t),
\end{cases}\tag{A.7}
\end{equation}
which gives
\begin{equation}
\begin{cases}
A''_x = A_x A'_x\,\dmu\\
B''_y = B'_y B_y\,\dmu,
\end{cases}\tag{A.8}
\end{equation}
\begin{equation}
\begin{cases}
A''_y = A_y + A_x X_y\\
B''_x = B'_x + B'_y Y_x,
\end{cases}\tag{A.9}
\end{equation}
\begin{equation}
\begin{cases}
A''_t = A_t + A_x X_t\\
B''_t = B'_t + B'_y Y_t.
\end{cases}\tag{A.10}
\end{equation}
Next, from (A.4), we have
\begin{equation}
\begin{cases}
-\De_x = B_{xx} X_x A'_y + B_x A'_{xy} + B_x A'_{yy} Y_x\\
-\De_y = A'_{yy} Y_y B_x + A'_y B_{xy} + A'_y B_{xx} X_y\\
-\De_t = B_{xt} A'_y + B_{xx} X_t A'_y + B_x A'_{yt} + B_x A'_{yy}
Y_t.
\end{cases}\tag{A.11}
\end{equation}
Taking logarithmic derivatives in (A.8) gives
\begin{equation}
\pa_x \log |A''_x| = \pa_x \log |A'_x| + Y_x \pa_y \log |A'_x| + X_x
\pa_x \log |A_x| - \De_x \dmu,\tag{A.12}
\end{equation}
\begin{equation}
\pa_y \log |A''_x| = \pa_y \log |A_x| + X_y \pa_x \log |A_x| + Y_y
\pa_y \log |A'_x| - \De_y \dmu,\tag{A.13}
\end{equation}
\begin{equation}
\pa_t \log |A''_x| = \pa_t \log |A_x| + \pa_t \log |A'_x| + X_t
\pa_x \log |A_x| + Y_t \pa_y \log |A'_x| - \De_t \dmu,\tag{A.14}
\end{equation}
\begin{equation}
\pa_y \log |B''_y| = \pa_y \log |B_y| + X_y \pa_x \log |B_y| + Y_y
\pa_y \log |B'_y| - \De_y \dmu,\tag{A.15}
\end{equation}
\begin{equation}
\pa_x \log |B''_y| = \pa_x \log |B'_y| + Y_x \pa_y \log |B'_y| + X_x
\pa_x \log |B_y| - \De_x \dmu,\tag{A.16}
\end{equation}
\begin{equation}
\pa_t \log |B''_y| = \pa_t \log |B'_y| + \pa_t \log |B_y| + Y_t
\pa_y \log |B'_y| + X_t \pa_x \log |B_y| - \De_t \dmu.\tag{A.17}
\end{equation}
Taking derivatives in (A.9) gives
\begin{equation}
\begin{cases}
A''_{yy} = A_{yy} + 2A_{xy} X_y + A_{xx} X^2_y + A_x X_{yy}\\
B''_{xx} = B'_{xx} + 2B'_{xy} Y_x + B'_{yy} Y^2_x + B'_y Y_{xx},
\end{cases}\tag{A.18}
\end{equation}
\begin{equation}
\begin{cases}
A''_{yt} = A_{yt} + X_t A_{xy} + X_y A_{xt} + X_t X_y A_{xx} + A_x X_{yt}\\
B''_{xt} = B'_{xt} + Y_t B'_{xy} + Y_x B'_{yt} + Y_t Y_x B'_{yy} +
B'_y Y_{xt},
\end{cases}\tag{A.19}
\end{equation}
where the partial derivatives of $X$, $Y$ are obtained from (A.6):
\begin{equation}
X_{yy} = B_y \dmu (A'_{yy} Y_y + A'_y \pa_y \log |B_y| + A'_y X_y
\pa_x \log |B_y| - A'_y \De_y \dmu),\tag{A.20}
\end{equation}
\begin{equation}
Y_{xx} = A'_x \dmu (B_{xx} X_x + B_x \pa_x \log |A'_x| + B_x Y_x
\pa_y \log |A'_x| - B_x \De_x \dmu),\tag{A.21}
\end{equation}
\begin{equation}
X_{yt} = B_y \dmu (A'_{yy} Y_t + A'_{yt} + A'_y \pa_t \log |B_y| +
A'_y X_t \pa_x \log |B_y| - A'_y \De_t \dmu),\tag{A.22}
\end{equation}
\begin{equation}
Y_{xt} = A'_x \dmu (B_{xx} X_t + B_{xt} + B_x \pa_t \log |A'_x| +
B_x Y_t \pa_y \log |A'_x| - B_x \De_t \dmu).\tag{A.23}
\end{equation}

\medskip\noindent
{\bf A.2~Parabolic Composition.~} We have now a fold map $G_t = G_+
\circ G_0 \circ G_-$:
\[
(x_u, y_u) \overset{G_-}{\longrightarrow} (w, y_u)
\overset{G_0}{\longrightarrow} (x_s, w)
\overset{G_+}{\longrightarrow} (x_s, y_s),
\]
with
\begin{equation}
\begin{cases}
y_s = Y_s (w, x_s, t)\\
x_u = X_u (w, y_u, t),
\end{cases}\tag{A.24}
\end{equation}
\begin{equation}
w^2 = \te (y_u, x_s, t).\tag{A.25}
\end{equation}
We also have an affine like map $F_0: (x_0, y_0) \mapsto (x_u, y_u)$
implicitly defined by
\begin{equation}
\begin{cases}
x_0 = A_0 (y_0, x_u, t)\\
y_u = B_0 (y_0, x_u, t),
\end{cases}\tag{A.26}
\end{equation}
and another affine-like map $F_1: (x_1, y_1) \mapsto (x_1, y_1)$
implicitly defined by
\begin{equation}
\begin{cases}
x_s = A_1 (y_s, x_1, t)\\
y_1 = B_1 (y_s, x_1, t).
\end{cases}\tag{A.27}
\end{equation}
We assume that (PC1), (PC2) in Subsection~\ref{sub3.5} are
satisfied. As we have seen in \cite{PY2} and
Subsection~\ref{sub3.5}, the first step is to write
\begin{equation}
\begin{cases}
x_u = X (w, y_0, t)\\
y_s = Y (w, x_1, t),
\end{cases}\tag{A.28}
\end{equation}
where the partial derivatives of $X$, $Y$ are given by
\begin{equation}
\begin{cases}
X_w = X_{u,w}\,\De^{-1}_0\\
X_y = X_{u,y} B_{0,y}\,\De^{-1}_0\\
X_t = (X_{u,t} + X_{u,y} B_{0,t})\,\De^{-1}_0\\
Y_w = Y_{s,w}\,\duu\\
Y_x = Y_{s,x} A_{1,x}\,\duu\\
Y_t = (Y_{s,t} + Y_{s,x} A_{1,t})\,\duu,
\end{cases}\tag{A.29}
\end{equation}
\begin{equation}
\begin{cases}
\De_0 = 1 - X_{u,y} B_{0,x}\\
\De_1 = 1 - Y_{s,x} A_{1,y}.
\end{cases}\tag{A.30}
\end{equation}
We set
\begin{equation}
\begin{cases}
\ovY (w, y_0, t)\, := \,B_0 (y_0, X, t)\\
\ovX (w, x_1, t)\, := \,A_1 (Y, x_1, t),
\end{cases}\tag{A.31}
\end{equation}
\begin{equation}
C (w, y_0, x_1, t)\, := \,w^2 - \te (\ovX, \ovY, t).\tag{A.32}
\end{equation}
The partial derivatives are given by
\begin{equation}
\begin{cases}
\ovY_w = B_{0,x} X_w\\
\ovY_y = B_{0,y} + B_{0,x} X_y = B_{0,y}\,\De^{-1}_0\\
\ovY_t = B_{0,t} + B_{0,x} X_t = (B_{0,t} + B_{0,x} X_{u,t})\, \duz,
\end{cases}\tag{A.33}
\end{equation}
\begin{equation}
\begin{cases}
\ovX_w = A_{1,y} Y_w\\
\ovX_x = A_{1,x} + A_{1,y} Y_x = A_{1,x}\,\duu\\
\ovX_t = A_{1,t} + A_{1,y} Y_t = (A_{1,t} + A_{1,y} Y_{s,t})\, \duu,
\end{cases}\tag{A.34}
\end{equation}
\begin{equation}
\begin{cases}
-C_w = -2w + \te_x \ovX_w + \te_y \ovY_w\\
-C_x = \te_x \ovX_x\\
-C_y = \te_y \ovY_y\\
-C_t = \te_x \ovX_t + \te_y \ovY_t + \te_t.
\end{cases}\tag{A.35}
\end{equation}
We solve
\begin{equation}
C(w, y_0, x_1, t) = 0\tag{A.36}
\end{equation}
to define
\begin{equation}
w = W (y_0, x_1, t)\tag{A.37}
\end{equation}
(there are two solutions $W^+$ and $W^-$).

The corresponding branch of the parabolic composition is implicitly
defined by
\begin{equation}
\begin{cases}
x_0 = A_0 (y_0, X (W, y_0, t), t) =: A (y_0, x_1, t)\\
y_1 = B_1 (Y (W, x_1, t), x_1, t) =: B (y_0, x_1, t).
\end{cases}\tag{A.38}
\end{equation}
The partial derivatives of $A$, $B$, $W$ are given by
\begin{equation}
\begin{cases}
A_x = A_{0,x} X_w W_x\\
A_y = A_{0,y} + A_{0,x} (X_y + X_w W_y)\\
A_t = A_{0,t} + A_{0,x} (X_t + X_w W_t)
\end{cases}\tag{A.39}
\end{equation}
\begin{equation}
\begin{cases}
B_y = B_{1,y} Y_w W_y\\
B_x = B_{1,x} + B_{1,y} (Y_x + Y_w W_x)\\
B_t = B_{1,t} + B_{1,y} (Y_t + Y_w W_t)
\end{cases}\tag{A.40}
\end{equation}
\begin{equation}
\begin{cases}
W_x = -C_x C_w^{-1}\\
W_y = -C_y C_w^{-1}\\
W_t = -C_t C_w^{-1}.
\end{cases}\tag{A.41}
\end{equation}
Substituting (A.29), (A.41), (A.35), (A.34) in the formulas
(A.39)--(A.40) leads to
\begin{equation}
\begin{cases}
A_x = A_{0,x} A_{1,x} C^{-1}_w \te_x X_{u,w} \duz \duu\\
B_y = B_{1,y} B_{0,y} C^{-1}_w \te_y Y_{s,w} \duz \duu
\end{cases}\tag{A.42}
\end{equation}
\begin{equation}
\begin{cases}
A_y = A_{0,y} + A_{0,x} B_{0,y} \duz (X_{u,y} + X_{u,w} \te_y \duz C^{-1}_w)\\
B_x = B_{1,x} + B_{1,y} A_{1,x} \duu (Y_{s,x} + Y_{s,w} \te_x \duu
C^{-1}_w)
\end{cases}\tag{A.43}
\end{equation}
\begin{equation}
\begin{cases}
A_t = A_{0,t} + A_{0,x} \duz [ X_{u,t} + X_{u,y} B_{0,t} + X_{u,w}
C^{-1}_w (\te_t + \te_x \ovX_t + \te_y \ovY_t)]\\
B_t = B_{1,t} + B_{1,y} \duu [ Y_{s,t} + Y_{s,x} A_{1,t} + Y_{s,w}
C^{-1}_w (\te_t + \te_x \ovX_t + \te_y \ovY_t)]
\end{cases}\tag{A.44}
\end{equation}
Taking the logarithmic derivatives in the first formula of (A.39),
we obtain
\begin{equation}
\pa_x \log |A_x| = W_x X_w \pa_x \log |A_{0,x}| + W_x \pa_w \log
|X_w| + \pa_x \log |W_x|,\tag{A.45}
\end{equation}
\begin{alignat}{2}
\pa_y \log |A_x| &= &\pa_y \log |A_{0,x}| + \pa_x \log |A_{0,x}|
(X_y + X_w W_y)\tag{A.46}\\
&+ &\pa_y \log |X_w| + W_y \pa_w \log |X_w| + \pa_y \log
|W_x|,\notag
\end{alignat}
\begin{alignat}{2}
\pa_t \log |A_x| &= &\pa_t \log |A_{0,x}| + \pa_x \log |A_{0,x}|
(X_t + X_w W_t)\tag{A.47}\\
&+ &\pa_t \log |X_w| + W_t \pa_w \log |X_w| + \pa_t \log
|W_x|.\notag
\end{alignat}
From the second formula in (A.39), one gets
\begin{alignat}{2}
&A_{yy} &= &A_{0,yy} + 2 A_{0,xy} (X_y + X_w W_y) + A_{0,xx} (X_y +
X_w W_y)^2\tag{A.48}\\
& &+ &A_{0,x} (X_{yy} + 2 X_{wy} W_y + X_{ww} W^2_y + X_w
W_{yy}),\notag
\end{alignat}
\begin{alignat}{2}
&A_{yt} &= &A_{0,yt} + A_{0,xy} (X_y + X_w W_t) + A_{0,xt} (X_y +
X_w W_y)\tag{A.49}\\
& &+ &A_{0,xx} (X_t + X_w W_t) (X_y + X_w W_y)\notag\\
& &+ &A_{0,x} (X_{yt} + X_{wy} W_t + X_{wt} W_y + X_{ww} W_y W_t +
X_w W_{yt}).\notag
\end{alignat}
The symmetric formulas for $B$ are
\begin{equation}
\pa_y \log |B_y| = W_y Y_w \pa_y \log |B_{1,y}| + \pa_y \log |W_y| +
W_y \pa_w \log |Y_w|,\tag{A.50}
\end{equation}
\begin{alignat}{2}
&\pa_x \log |B_y| &= &\pa_x \log |B_{1,y}| + \pa_x \log
|Y_w|\tag{A.51}\\
& &+ &\pa_x \log |W_y| + \pa_y \log |B_{1,y}| (Y_x + Y_w W_x) + W_x
\pa_w \log |Y_w|,\notag
\end{alignat}
\begin{alignat}{2}
\pa_t \log |B_y| &= &\pa_t \log |B_{1,y}| + \pa_y \log |B_{1,y}|
(Y_t + Y_w W_t)\tag{A.52}\\
&+ &\pa_t \log |Y_w| + W_t \pa_w \log |Y_w| + \pa_t \log
|W_y|,\notag
\end{alignat}
\begin{alignat}{2}
&B_{xx} &= &B_{1,xx} + 2 B_{1,xy} (Y_x + Y_w W_x) + B_{1,yy} (Y_x +
Y_w W_x)^2\tag{A.53}\\
& &+ &B_{1,y} (Y_{xx} + 2Y_{wx} W_x + Y_{ww} W_x^2 + Y_w
W_{xx}),\notag
\end{alignat}
\begin{alignat}{2}
&B_{xt} &= &B_{1,xt} + B_{1,xy} (Y_t + Y_w W_t) + B_{1,yt} (Y_x +
Y_w W_x)\tag{A.54}\\
& &+ &B_{1,yy} (Y_x + Y_w W_x)(Y_t + Y_w W_t)\notag\\
& &+ &B_{1,y} (Y_{xt} + Y_{wx} W_t + Y_{wt} W_x + Y_{ww} W_x W_t +
Y_w W_{xt}).\notag
\end{alignat}
In formulas (A.45)--(A.54), the partial derivatives of order 2 of
$W$ are obtained from (A.41):
\begin{equation}
\begin{cases}
W_{xx} = -C^{-1}_w (C_{ww} W^2_x + 2C_{wx} W_x + C_{xx})\\
W_{xy} = -C^{-1}_w (C_{ww} W_x W_y + C_{wx} W_y + C_{wy} W_x + C_{xy})\\
W_{yy} = -C^{-1}_w (C_{ww} W^2_y + 2C_{wy} W_y + C_{yy})\\
W_{xt} = -C^{-1}_w (C_{ww} W_x W_t + C_{wx} W_t + C_{wt} W_x + C_{xt})\\
W_{yt} = -C^{-1}_w (C_{ww} W_y W_t + C_{wy} W_t + C_{wt} W_y +
C_{yt})
\end{cases}\tag{A.55}
\end{equation}
The partial derivatives of order 2 of $C$ are obtained from (A.35):
\begin{equation}
-C_{ww} = -2 + \te_x \ovX_{ww} + \te_y \ovY_{ww} + \te_{xx} \ovX^2_w
+ 2\te_{xy} \ovX_w \ovY_w + \te_{yy} \ovY^2_w,\tag{A.56}
\end{equation}
\begin{equation}
\begin{cases}
-C_{wx} = \te_{xx} \ovX_w \ovX_x + \te_{xy} \ovY_w \ovX_x + \te_x \ovX_{wx}\\
-C_{wy} = \te_{yy} \ovY_w \ovY_y + \te_{xy} \ovX_w \ovY_y + \te_y
\ovY_{wy}
\end{cases}\tag{A.57}
\end{equation}
\begin{equation}
-C_{wt} = \te_{xt} \ovX_w + \te_{xx} \ovX_w \ovX_t + \te_{xy}
(\ovX_w \ovY_t + \ovY_w \ovX_t) + \te_{yy} \ovY_w \ovY_t + \te_{yt}
\ovY_w + \te_x \ovX_{wt} + \te_y \ovY_{wt},\tag{A.58}
\end{equation}
\begin{equation}
\begin{cases}
-C_{xx} = \te_{xx} \ovX^2_x + \te_x \ovX_{xx}\\
-C_{xy} = \te_{xy} \ovX_x \ovY_y\\
-C_{yy} = \te_{yy} \ovY^2_y + \te_y \ovY_{yy},
\end{cases}\tag{A.59}
\end{equation}
\begin{equation}
\begin{cases}
-C_{xt} = \te_{xt} \ovX_x + \te_{xx} \ovX_x \ovX_t + \te_{xy} \ovX_x \ovY_t + \te_x \ovX_{xt}\\
-C_{yt} = \te_{yt} \ovY_y + \te_{yy} \ovY_y \ovY_t + \te_{xy} \ovY_y
\ovX_t + \te_y \ovY_{yt}.
\end{cases}\tag{A.60}
\end{equation}
The partial derivatives of order 2 of $\ovX$, $\ovY$ are obtained
from (A.34):
\begin{equation}
\begin{cases}
\ovX_{ww} = A_{1,yy} Y^2_w + A_{1,y} Y_{ww}\\
\ovX_{wx} = A_{1,xy} Y_w + A_{1,yy} Y_w Y_x + A_{1,y} Y_{wx}\\
\ovX_{wt} = A_{1,yt} Y_w + A_{1,yy} Y_w Y_t + A_{1,y} Y_{wt}\\
\ovX_{xx} = A_{1,xx} + 2A_{1,xy} Y_x + A_{1,yy} Y^2_x + A_{1,y} Y_{xx}\\
\ovX_{xt} = A_{1,xt} + A_{1,xy} Y_t + A_{1,yt} Y_x + A_{1,yy} Y_x
Y_t + A_{1,y} Y_{xt}.
\end{cases}\tag{A.61}
\end{equation}
\begin{equation}
\begin{cases}
\ovY_{ww} = B_{0,xx} X^2_w + B_{0,x} X_{ww}\\
\ovY_{wy} = B_{0,xy} X_w + B_{0,xx} X_w X_y + B_{0,x} X_{wy}\\
\ovY_{wt} = B_{0,xt} X_w + B_{0,xx} X_w X_t + B_{0,x} X_{wt}\\
\ovY_{yy} = B_{0,yy} + 2B_{0,xy} X_y + B_{0,xx} X^2_y + B_{0,x} X_{yy}\\
\ovY_{yt} = B_{0,yt} + B_{0,xy} X_t + B_{0,xt} X_y + B_{0,xx} X_y
X_t + B_{0,x} X_{yt}.
\end{cases}\tag{A.62}
\end{equation}
Finally, from (A.29), we obtain
\begin{equation}
\begin{cases}
X_{ww} = \duz (X_{u,ww} + 2X_{u,wy} \ovY_w + X_{u,yy} \ovY^2_w +
X_{u,y} B_{0,xx} X^2_w),\\
X_{wy} = \duz (X_{u,wy} \ovY_y + X_{u,yy} \ovY_w \ovY_y + X_{u,y}
X_w (B_{0,xy} + B_{0,xx} X_y)),\\
X_{yy} = \duz (X_{u,yy} \ovY^2_y + X_{u,y} (B_{0,yy} + 2B_{0,xy} X_y
+ B_{0,xx} X^2_y))\\
X_{wt} = \duz [X_{u,wt} + X_{u,wy} \ovY_t + X_{u,yy} \ovY_w \ovY_t +
X_{u,yt} \ovY_w\\
\qquad + X_{u,y} X_w (B_{0,xt} + B_{0,xx} X_t)],\\
X_{yt} = \duz [X_{u,yt} \ovY_y + X_{u,yy} \ovY_y \ovY_t\\
\qquad + X_{u,y} (B_{0,yt} + B_{0,xy} X_t + B_{0,xt} X_y + B_{0,xx}
X_y X_t)],\\
\end{cases}\tag{A.63}
\end{equation}

\begin{equation}
\begin{cases}
Y_{ww} = \duu (Y_{s,ww} + 2Y_{s,wx} \ovX_w + Y_{s,xx} \ovX^2_w +
Y_{s,x} A_{1,yy} Y^2_w),\\
Y_{wx} = \duu (Y_{s,wx} \ovX_x + Y_{s,xx} \ovX_w \ovX_x + Y_{s,x}
Y_w (A_{1,xy} + A_{1,yy} Y_x)),\\
Y_{xx} = \duu (Y_{s,xx} \ovX^2_x + Y_{s,x} (A_{1,xx} + 2A_{1,xy} Y_x
+ A_{1,yy} Y^2_x))\\
Y_{wt} = \duu [Y_{s,wt} + Y_{s,wx} \ovX_t + Y_{s,xx} \ovX_w \ovX_t +
Y_{s,xt} \ovX_w\\
\qquad + Y_{s,x} Y_w (A_{1,yt} + A_{1,yy} Y_t)],\\
Y_{xt} = \duu [Y_{s,xt} \ovX_x + Y_{s,xx} \ovX_x \ovX_t\\
\qquad + Y_{s,x} (A_{1,xt} + A_{1,xy} Y_t + A_{1,yt} Y_x + A_{1,yy}
Y_x Y_t)].\\
\end{cases}\tag{A.64}
\end{equation}

\newpage

\setcounter{equation}{0}

\section*{Appendix~B\\
On the Lipschitz Regularity of $\wtR^\infty_+$}

{\bf B.1~} In this appendix, we will perform some calculations and
obtain some estimates that uphold the proof of
Proposition~\ref{propo40} in Subsection~\ref{sub10.5}, showing that
the holonomy maps of the partial foliation $\wtcR^\infty_+$ are
uniformly Lipschitz.

As always, we have to consider two cases, dealing separately with
affine-like maps and parabolic composition, the later case being
more involved than the first.

\medskip\noindent
{\bf B.2~The Case of Affine-Like Maps.~} We consider here, with the
notations of Subsection~\ref{sub3.1}, an affine-like map $F$ with
implicit representation $(A,B)$: the domain of $F$ is a vertical
strip $P$ in a rectangle $I^s_0 \times I^u_0$, its image is a
horizontal strip in a rectangle $I^s_1 \times I^u_1$, the respective
coordinates being $x_0$, $y_0$, $x_1$, $y_1$.

In $I^s_1 \times I^s_1$, we are given a one-parameter family of
vertical-like curves which are graphs
\begin{equation}
\om (s) = \{ x_1 = \vphi (y_1, s)\}, \text{ with }\; \Bigl\vert
\frac{\pa\vphi}{\pa y_1}\Bigr\vert \; \text{ not too large}.
\tag{B.1}
\end{equation}
We assume that, for all $y_1$, $s$
\begin{equation}
\frac{\pa \vphi}{\pa s}\; (y_1, s) > 0\tag{B.2}
\end{equation}
and that we have a uniform bound
\begin{equation}
\Bigl\vert \frac{\pa}{\pa y_1}\;\log\, \frac{\pa\vphi}{\pa s}\,
(y_1, s)\Bigr\vert\; \leqslant \; T.\tag{B.3}
\end{equation}
We define the vertical-like curves in $P$ by
\begin{equation}
\Om (s) = F^{-1} (\om (s) \cap Q).\tag{B.4}
\end{equation}
We assume that $F$ satisfies a cone condition and we have a good
control on the distortion of $F$. Under these hypotheses, $\Om (s)$
will be a graph
\begin{equation}
\Om (s) = \{ x_0 = \Phi (y_0, s)\}\tag{B.5}
\end{equation}
and we want to obtain the analogue of (B.3):
\begin{equation}
\biggl\vert \frac{\pa}{\pa y_0}\;\log\, \Bigl\vert\frac{\pa\Phi}{\pa
s}\, (y_0, s)\Bigr\vert\biggr\vert\; \leqslant\; T'.\tag{B.6}
\end{equation}
The relation between $\vphi$ and $\Phi$ is as follows in the
equation
\begin{equation}
y_1 = B (y_0, \vphi (y_1, s)),\tag{B.7}
\end{equation}
we solve for $y_1$ to define
\begin{equation}
y_1 = \psi (y_0, s).\tag{B.8}
\end{equation}
We then have
\begin{equation}
\Phi (y_0, s) = A \Bigl(y_0, \vphi (\psi (y_0, s),
s)\Bigr).\tag{B.9}
\end{equation}
From (B.7), (B.8), we get
\begin{equation}
\frac{\pa\psi}{\pa s} = B_x \; \frac{\pa\vphi}{\pa s} \Bigl( 1 -
B_x\; \frac{\pa\vphi}{\pa y_1}\Bigr)^{-1}\tag{B.10}
\end{equation}
(where, as always, the functions must be taken at the appropriate
arguments); then, from (B.9) we get
\begin{equation}
\frac{\pa\Phi}{\pa s} = A_x \; \frac{\pa\vphi}{\pa s} \Bigl( 1 -
B_x\; \frac{\pa\vphi}{\pa y_1}\Bigr)^{-1}.\tag{B.11}
\end{equation}
We must now take the logarithmic derivative of the right-hand term
with respect to $y_0$. It is the sum of three terms:
\begin{equation}
Z_1 = \pa_y \log |A_x| + \pa_x \log |A_x| \; \frac{\pa\vphi}{\pa
y_1} \; \frac{\pa\psi}{\pa y_0},\tag{B.12}
\end{equation}
\begin{equation}
Z_2 = \frac{\pa}{\pa y_1} \, \Bigl( \log \; \frac{\pa\vphi}{\pa
s}\Bigr)\; \frac{\pa\psi}{\pa y_0},\tag{B.13}
\end{equation}
\begin{equation}
Z_3 = \Bigl( 1 - B_x \; \frac{\pa\vphi}{\pa y_1}\Bigr)^{-1} \Bigl\{
\frac{\pa^2\vphi}{\pa y^2_1}\; \frac{\pa\psi}{\pa y_0}\,B_x +
\frac{\pa\vphi}{\pa y_1} \Bigl( B_y \pa_x \log |B_y| + B_{xx} \;
\frac{\pa\vphi}{\pa y_1}\; \frac{\pa\psi}{\pa y_0}\Bigr)\Bigr\}.
\tag{B.14}
\end{equation}
In these formulas, we have from (B.7), (B.8)
\begin{equation}
\frac{\pa\psi}{\pa y_0}  = B_y \, \Bigl( 1 - B_x \;
\frac{\pa\vphi}{\pa y_1}\Bigr)^{-1}.\tag{B.15}
\end{equation}
We see that we obtain (B.6) with
\begin{equation}
T' \; \leqslant \; a\,T + C.\tag{B.16}
\end{equation}
Here, the constant $C$ is bounded in terms of
$\Vert(1-B_x\,\tfrac{\pa\vphi}{\pa y_1})^{-1}\Vert_\infty$, the
distortion of $F$, the cone condition, $\Vert \tfrac{\pa\vphi}{\pa
y_1}\Vert_\infty$ and $\Vert \tfrac{\pa^2\vphi}{\pa
y^2_1}\Vert_\infty$. The constant $a$ is given by:
\begin{equation}
a \; = \; \Vert \frac{\pa\psi}{\pa y_0}\Vert_\infty \; \leqslant \;
C \Vert B_y \Vert_\infty\tag{B.17}
\end{equation}
and, therefore, $a$ will be $< \tfrac{1}{2}$ if $|Q| =
||B_y||_\infty$ is sufficiently small.

\begin{rem}
The control of $\tfrac{\pa \Phi}{\pa y_0}$ is guaranteed by the cone
condition and the control of $\tfrac{\pa^2\Phi}{\pa y^2_0}$ by the
distortion of $F$.
\end{rem}
Thus, the mixed second derivative does not ``explode'' by taking
pre-images by affine-like maps.\hfill $\square$

\medskip\noindent
{\bf B.3~The Parabolic Case.} We now use the setting and notations
of Subsection~\ref{sub3.5}. We have intervals $I^s_0$, $I^u_0$,
$I^s_u$, $I^u_u$, $I^s_s$, $I^u_s$ with respective coordinates
$x_0$, $y_0$, $x_u$, $y_u$, $x_s$, $y_s$.

We have an affine-like map $F$ with domain $P\subset I^s_0 \times
I^u_0$, image $Q \subset I^s_u\times I^u_u$, and implicit
representation $(A,B)$. We also have a folding map $G$ with domain
$L_u \subset I^s_u \times I^u_u$, image $L_s \subset I^s_s \times
I^u_s$. The map $G$ is implicitly defined by the system
\begin{equation}
\begin{cases}
y_s = Y_s (w, x_s)\\
x_u = X_u (w, y_u)\\
w^2 = \te (y_u, x_s)
\end{cases}\tag{B.18}
\end{equation}
with $Y_s$, $X_u$, $\te$ as in Subsections~\ref{sub2.3} and
\ref{sub3.5}. The affine-like map $F$ should satisfy (cf.~(PC1) in
Subsection~\ref{sub3.5}, (R4) in Subsection~\ref{sub5.3}):
\begin{equation}
|B_x| \; < \; b,\qquad\qquad |B_{xx}| \; < \; b.\tag{B.19}
\end{equation}
with $b\ll 1$. As above we are given a one-parameter family of
vertical-like curves which are graphs
\begin{equation}
\om (s) = \{ x_s = \vphi (y_s, s)\},\tag{B.20}
\end{equation}
but these curves are now very close to being exactly vertical
\begin{equation}
\Bigl\vert \frac{\pa\vphi}{\pa y_s}\Bigr\vert \; < \; b,\qquad
\Bigl\vert \frac{\pa^2\vphi}{\pa y^2_s}\Bigr\vert \; < \;
b.\tag{B.21}
\end{equation}
We also assume as above that
\begin{equation}
\frac{\pa\vphi}{\pa s}\; (y_s, s)\; > \; 0\tag{B.22}
\end{equation}
and that we have some uniform bound
\begin{equation}
\Bigl\vert \frac{\pa}{\pa s}\, \log \, \frac{\pa\vphi}{\pa
s}\Bigr\vert \; \leqslant \; T.\tag{B.23}
\end{equation}
We then define the curves $\Om^\pm (s)$ as the connected components
of $F^{-1} (Q \cap G^{-1} (\om (s) \cap L_s))$.

They should be graphs
\begin{equation}
\Om^\pm (s) = \{ x_0 = \Phi^\pm (y_0, s)\}\tag{B.24}
\end{equation}
and we want some uniform control
\begin{equation}
\biggl\vert \frac{\pa}{\pa y_0}\, \log \,
\Bigl\vert\frac{\pa\Phi}{\pa s}\Bigr\vert\biggr\vert \; \leqslant \;
T'.\tag{B.25}
\end{equation}
As in Subsection~\ref{sub3.5} (formulas (\ref{eq3.14}),
(\ref{eq3.15}) we eliminate $y_u$ to write
\begin{equation}
x_u \; = \; X (w, y_0).\tag{B.26}
\end{equation}
Similarly, in the equation
\begin{equation}
y_s \; = \; Y_s (w, \vphi (y_s, s)),\tag{B.27}
\end{equation}
we solve for $y_s$ to write
\begin{equation}
y_s \; = \; Y (w, s).\tag{B.28}
\end{equation}
We then set
\begin{alignat}{3}
&\ovX (w,s) &:= &\vphi (Y (w, s), s),\tag{B.29}\\
&\ovY (w, y_0) &:= &B (y_0, X (w, y_0)),\tag{B.30}\\
&C (w, y_0, s) &:= &w^2 - \te (\ovY (w, y_0), \ovX (w,s)).\tag{B.31}
\end{alignat}
Then, one should solve $C = 0$ to get
\begin{align}
&w = W^\pm (y_0, s),\tag{B.32}\\
&\Phi (y_0, s) = A (y_0, X (W (y_0,s), y_0)),\tag{B.33}
\end{align}
with $W = W^\pm$, $\Phi = \Phi^\pm$.

The relations (B.26) through (B.33) are the equations which allow us
to calculate $\tfrac{\pa}{\pa y_0}\, \log\, \vert\tfrac{\pa\Phi}{\pa
s}\vert$.

From (B.33), we have
\begin{equation}
\frac{\pa\Phi}{\pa s} = A_x X_w \, \frac{\pa W}{\pa s}.\tag{B.34}
\end{equation}
From the implicit definition of $X$
\begin{equation}
X (w, y_0) = X_u \Bigl(w, B_0 (y_0, X (w, y_0))\Bigr),\tag{B.35}
\end{equation}
we get
\begin{equation}
X_w = (1 - X_{u,y} B_x)^{-1} X_{u,w}.\tag{B.36}
\end{equation}
One has also
\begin{alignat}{3}
&\frac{\pa W}{\pa s} &= &-\frac{\pa C}{\pa s}\; C^{-1}_w,\tag{B.37}\\
&\frac{\pa C}{\pa s} &= &-\te_x \Bigl( \frac{\pa\vphi}{\pa s} +
&\frac{\pa\vphi}{\pa y}\; \frac{\pa Y}{\pa s}\Bigr).\tag{B.38}
\end{alignat}
From (B.27), (B.28), one obtains
\begin{equation}
\frac{\pa Y}{\pa s} = Y_{s,x}\; \frac{\pa\vphi}{\pa s} \Bigl( 1 -
Y_{s,x} \; \frac{\pa\vphi}{\pa y}\Bigr)^{-1},\tag{B.39}
\end{equation}
and putting this in (B.38) leads to
\begin{equation}
\frac{\pa C}{\pa s} = -\te_x\; \frac{\pa\vphi}{\pa s} \Bigl( 1 -
Y_{s,x} \; \frac{\pa\vphi}{\pa y}\Bigr)^{-1}.\tag{B.40}
\end{equation}
Therefore, introducing (B.36), (B.37), (B.40) into (B.34), we have
\begin{equation}
\frac{\pa\Phi}{\pa s} = A_x X_{u,w} (1 - X_{u,y} B_x)^{-1} C^{-1}_w
\te_x \, \frac{\pa\vphi}{\pa s}\, \Bigl( 1 - Y_{s,x}\,
\frac{\pa\vphi}{\pa y}\Bigr)^{-1}.\tag{B.41}
\end{equation}
We now take the logarithmic derivative of this product of seven
terms with relation to $y_0$. We obtain a sum of seven terms:
\begin{alignat}{3}
&Z_1\; &= \;&\pa_y \log |A_x| + \pa_x \log |A_x| \; (X_y + X_w
W_y),\tag{B.42}\\
&Z_2\; &= \;&(\pa_w \log |X_{u,w}|) W_y + \pa_y \log |X_{u,w}|
(\ovY_y + \ovY_w W_y),\tag{B.43}\\
&Z_3\; &= \;&(1 - X_{u,y} B_x)^{-1} \Bigl\{ X_{uy} (B_{xy} + B_{xx}
(X_y + X_w W_y))\tag{B.44}\\
& &+\; &B_x (X_{u,yw} W_y + X_{u,yy} (\ovY_y + \ovY_w W_y))
\Bigr\},\notag\\
&Z_4\; &=\; &C^{-1}_w \Bigl[ \te_x \ovX_w W_y + \te_y (\ovY_y +
\ovY_w W_y)\Bigr],\tag{B.45}\\
&Z_5\; &=\; &\pa_y \log |\te_x| (\ovY_y + \ovY_w W_y) + \pa_x \log
|\te_x| \ovX_w W_y,\tag{B.46}\\
&Z_6\; &= \; &\Bigl( \frac{\pa}{\pa y_s}\,\log\, \frac{\pa\vphi}{\pa
s}\Bigr) \, Y_w W_y,\tag{B.47}\\
&Z_7\; &= \;&\Bigl( 1 - Y_{s,x} \, \frac{\pa\vphi}{\pa
y_s}\Bigr)^{-1} \Bigl\{ Y_{s,x}\,\frac{\pa^2 \vphi}{\pa y^2_s}\, Y_w
W_y + \frac{\pa\vphi}{\pa y_s} (Y_{s,xw} W_y + Y_{s,xx} \ovX_w
W_y)\Bigr\}.\tag{B.48}
\end{alignat}
The terms in these expressions that have not yet been introduced
are:
\begin{alignat}{3}
&X_y\; &= \;&X_{u,y} B_y (1 - X_{u,y} B_x)^{-1},\tag{B.49}\\
&W_y\; &=\; &-C_y C_w^{-1} = \te_y \ovY_y C^{-1}_w,\tag{B.50}\\
&\ovY_y\; &=\; &B_y + B_x X_y,\tag{B.51}\\
&\ovY_w\; &=\; &B_x X_w,\tag{B.52}\\
&\ovX\;_w &=\; &\frac{\pa\vphi}{\pa y_s}\, Y_w,\tag{B.53}\\
&Y_w\; &=\; &Y_{s,w} \Bigl( 1 - Y_{s,x}\, \frac{\pa\vphi}{\pa
y_s}\Bigr)^{-1}.\tag{B.54}
\end{alignat}
We obtain the estimates
\begin{alignat}{3}
&|X_w|\; &\leqslant \; &C,\; |Y_w|\; \leqslant \; C,\tag{B.55}\\
&|\ovX_w| \; &\leqslant \; &Cb,\; |\ovY_w|\; \leqslant \; Cb,\tag{B.56}\\
&|X_y|\; &\leqslant \; &C |Q|,\; |\ovY_y|\; \leqslant \; C |Q|,\tag{B.57}\\
&|W_y|\; &\leqslant \; &C |Q|\,|C^{-1}_w|,\tag{B.58}\\
&|Z_i|\; &\leqslant \; &C + C |Q|\,|C^{-1}_w|\; \text{ for
$i=1,2,3,5,7$;}\tag{B.59}\\
&|Z_4|\; &\leqslant \; &C |Q|\,|C^{-1}_w| + Cb |Q|\,|C^{-1}_w|^2;\tag{B.60}\\
&|Z_6|\; &\leqslant \; &C |Q|\,|C^{-1}_w|\, T.\tag{B.61}
\end{alignat}
Therefore, we obtain (B.25) with
\begin{equation}
T' = C |Q|\,||C^{-1}_w||_\infty\; T + C + C|Q|\,||C^{-1}_w||_\infty
+ Cb |Q|\,||C^{-1}_w||^2_\infty.\tag{B.62}
\end{equation}
Regarding $C^{-1}_w$, the same considerations as in
Subsection~\ref{sub3.5} apply. One checks that $C_{ww}$ is close to
two, hence the value of $|C^{-1}_w|$ at the argument where $C$
vanishes is of the order of $|\ovC|^{-\fud}$ where
\begin{equation}
\ovC (y_0, s) = \min_w \, C (w, y_0, s).\tag{B.63}
\end{equation}
In the context of Section~\ref{sec10}, we have
\begin{equation}
|Q| \; \ll \; |\ovC|\; \ll \; 1,\tag{B.64}
\end{equation}
and, therefore,
\begin{equation}
T' \; \leqslant \; C |Q|^{\fudt} T + C.\tag{B.65}
\end{equation}
This shows that $T'$ does not ''blow up''. \hfill $\square$

\medskip\noindent
{\bf B.4~} In the context of Proposition~\ref{propo40} in
Subsection~\ref{sub10.5}, Case~2 in the proof, we had an element
$(P_\ell, Q_\ell, n_\ell)\in \cR$ and two non-simple children
$P_{\ell+1}$, $P'_{\ell+1}$ of $P_\ell$. Let $\om_{\ell+1}$,
$\om'_{\ell+1}$ be components of the vertical part of the boundary
of $P_{\ell+1}$, $P'_{\ell+1}$ respectively. We wanted to have
(cf.~(\ref{eq10.65}))
\begin{equation}
\Bigl\vert \log\, \frac{\vphi_{\ell+1}(y) -
\vphi'_{\ell+1}(y)}{\vphi_{\ell+1} (y') - \vphi'_{\ell+1} (y')}
\Bigr\vert \; \leqslant \; C |y-y'|\tag{B.66}
\end{equation}
where $\om_{\ell+1} = \{ x = \vphi_{\ell+1} (y)\}$, $\om'_{\ell+1} = \{
x = \vphi'_{\ell+1} (y)\}$.

To prove (\ref{eq10.66}), we will imbed both $\om_{\ell+1}$ and
$\om'_{\ell+1}$ in a one-parameter family of curves
\begin{equation}
\Om (s) = \{ x = \Phi (y,s)\}\tag{B.67}
\end{equation}
with
\begin{equation}
\Om (s_0) = \om_{\ell+1},\qquad \Om (s_1) = \om'_{\ell+1}\tag{B.68}
\end{equation}
\begin{equation}
\frac{\pa}{\pa s}\; \Phi \; > \; 0\tag{B.69}
\end{equation}
\begin{equation}
\Bigl\vert \frac{\pa}{\pa y}\, \log\, \frac{\pa\Phi}{\pa
s}\Bigr\vert \; \leqslant \; C. \tag{B.70}
\end{equation}
Indeed, in this case, we can write
\begin{equation}
\frac{\vphi_{\ell+1}(y)-\vphi'_{\ell+1}(y)}{\vphi_{\ell+1}(y')-\vphi'_{\ell+1}(y')}
= \frac{\int^{s_1}_{s_0}\,\tfrac{\pa\Phi}{\pa s}\,(y,s)
ds}{\int^{s_1}_{s_0}\,\tfrac{\pa\Phi}{\pa s}\,(y',s) ds}, \tag{B.71}
\end{equation}
with, from (B.70), for any $s$
\begin{equation}
e^{-C|y-y'|}\, \leqslant \,\Bigl( \frac{\pa \Phi}{\pa s}
(y,s)\Bigr)^{-1}\; \frac{\pa \Phi}{\pa s}\, (y',s)\, \leqslant \,
e^{C|y-y'|},\tag{B.72}
\end{equation}
which yields (B.66).\hfill $\square$

\newpage

\setcounter{equation}{0}

\section*{Appendix~C\\
A Toy Model for the Transversality Relation}

\medskip\noindent
{\bf C.1~} Our goal in this appendix is to explain why the
complicated definition of the transversality relation in
Subsection~\ref{sub5.4}, is in some way ''natural'', if we require
some useful properties for the proof of our Main Theorem. The toy
model that we are considering is an abstract one. It is much simpler
than the real situation of Section~\ref{sec5} because the sets in
which the relation takes place are well defined to begin with: in
Section~\ref{sec5}, we need to know the transversality relation in
order to construct the classes $\cR (I)$.

\medskip\noindent
{\bf C.2~} A partially ordered set $X$ is a {\it forest} if, for any
$x_0\in X$, the set $\{ x\geqslant x_0\}$ is finite and totally
ordered. A {\it tree} is a forest with a single maximal element. Let
$X_1,\cdots, X_n$ be forests and let $A$ be a subset of $X = X_1
\times \cdots \times X_n$ (one should think of $A$ as the graph of
an {\it n-ary} relation). We say that $A$ is {\it hereditary} if
whenever $x = (x_1, \cdots, x_n)$, $y = (y_1, \cdots, y_n)$ are such
that and $y_i \leqslant x_i$ for all $1\leqslant i\leqslant n$
(abbreviated as $y\leqslant x$), then $y\in A$ if $x\in A$.

Two points $x = (x_1, \cdots, x_n)$, $y = (y_1,\cdots,y_n)$ of $X$
are {\it coordinate-wise comparable} ($c$-{\it comparable} for
short) if for each $i\in \{1,\cdots,n\}$, we have $x_i \geqslant
y_i$ or $x_i\leqslant y_i$. In this case, we set
\begin{equation}
x \vee y = (\max (x_i, y_i))_{1\leqslant i\leqslant n}.\tag{C.1}
\end{equation}
The set $A$ is concave if, whenever $x$, $y\in A$ are
$c$-comparable, then the point $x\vee y$ also belongs to $A$.

The intersection of hereditary, resp.~concave, subsets of $X$ is
hereditary, resp.~concave. It follows that any subset $A\subset X$
is contained in a smallest concave hereditary subset, called the
$c.h$-envelope and denoted by $\whA$.

\Exem When the number of factors $n=1$, any subset is concave: the
$c.h$-envelope of $A\subset X_1$ is the set of $x\in X_1$ such that
$x\leqslant y$ for some $y\in A$.

\medskip\noindent
{\bf C.3~} We construct the $c.h$-envelope when $n=2$.

\medskip\noindent
{\bf Proposition.~}{\it Let $X_1$, $X_2$ be forests and $A$ be a
subset of $X = X_1\times X_2$. Let $A_1$ be the set of $x\in X$ such
that $x = y \vee z$ for some $c$-comparable $y$, $z\in A$. The
$c.h$-envelope of $A$ is equal to the set $A_2$ of $t=(t_1,t_2)$
such that $t_1\leqslant x_1$, $t_2\leqslant x_2$ for some $x = (x_1,
x_2)\in A_1$.}

\medskip\noindent
{\it Proof.~} It is clear that the set $A_2$ defined in the
proposition is hereditary and it is contained in the $c.h$-envelope
$\whA$ of $A$. We have to prove that $A_2$ is concave. We first
prove the

\begin{lema}
If $y$, $z\in A_1$ are $c$-comparable, $y\vee z$ also belongs to
$A_1$.
\end{lema}

\begin{proof}
By the definition of $A_1$, we can write $y = y' \vee y''$, $z = z'
\vee z''$ with $y'$, $y''$, $z'$, $z''$ in $A$, $y'$ and $y''$
$c$-comparable, $z'$ and $z''$ $c$-comparable.  We may assume that
$y_1\leqslant z_1$ and $y_2\geqslant z_2$, and also $z_1 = z'_1$,
$y_2 = y'_2$; then, we have $z'_1 = z_1\geqslant y_1 \geqslant y'_1$
and $y'_2 = y_2 \geqslant z_2 \geqslant z'_2$, hence $y'$, $z'$ are
$c$-comparable with $y'\vee z' = y \vee z$.
\end{proof}

\medskip\noindent
{\bf End of the Proof of the Proposition.~} Let $t'$, $t''\in A_2$
be $c$-comparable and let $x'$, $x''\in A_1$ be such that
$x'_i\geqslant t'_i$, $x''_i \geqslant t''_i$ for $i=1,2$. As $X_1$
and $X_2$ are forests, $x'$ and $x''$ are $c$-comparable. From the
lemma, $x'\vee x''$ belongs to $A_1$; then $t'\vee t''$ belongs to
$A_2$.\hfill $\square$

\medskip\noindent
{\bf C.4~} For $n\geqslant 3$, the situation is more complicated, as
the two examples below indicate.

\medskip\noindent
{\bf Example~1.~} Let $X_1$, $X_2$, $X_3$ be forests and let $x$,
$y$, $z$ be three points of $X = X_1 \times X_2 \times X_3$ such
that
\begin{alignat}{3}
&x_1 &\geqslant &y_1 \geqslant z_1,\tag{C.2}\\
&y_2 &\geqslant &x_2,\;\; y_2\; \geqslant z_2,\tag{C.3}\\
&z_3 &\geqslant &x_3,\;\; z_3\; \geqslant y_3.\tag{C.4}
\end{alignat}
Let $A = \{ x,y,z\}\subset X$. If we define, as in the proposition
above,
\begin{equation}
A_1 = \{ u\vee v,\;\; u, v\in A,\;\; u,v\;\;
c\text{-comparable}\}\tag{C.5}
\end{equation}
and if we assume that $x_2$, $z_2$ are {\it not} comparable and
$x_3$, $y_3$ are {\it not} comparable, then we have
\begin{equation}
A_1 = \{ x, y, z, \;\; y\vee z = (y_1, y_2, z_3)\}.\tag{C.6}
\end{equation}
On the other hand, the point $w = (x_1, y_2, z_3) = x \vee (y \vee
z)$ certainly belongs to the $c.h$-envelope of $A$, but does not
satisfy $w_i\leqslant u_i\; (i=1,2,3)$ for any $u\in A_1$. This
example shows that the analogue of the proposition above is false
for $n=3$.

\medskip\noindent
{\bf Example~2.~} Let $X_1$, $X_2$, $X_3$ forests and let $x$, $y$,
$z\in X = X_1\times X_2 \times X_3$ such that
\begin{alignat}{3}
x_1 &\geqslant &y_1,\qquad x_1 \geqslant z_1\tag{C.7}\\
y_2 &\geqslant &x_2,\qquad y_2 \geqslant z_2\tag{C.8}\\
z_3 &\geqslant &x_3,\qquad z_3 \geqslant y_3\tag{C.9}
\end{alignat}
but none of the pairs $(y_1, z_1)$, $(x_2, z_2)$, $(x_3, y_3)$ is
made of comparable elements. Let $A = \{x,y,z\}$. The sets
$\{u\leqslant x\}$, $\{ v\leqslant y\}$, $\{ w\leqslant z\}$ are
disjoint and their union is the $c.h$-envelope of $A$: any
$u\leqslant x$, $v\leqslant y$ cannot be $c$-comparable; otherwise,
as $X_3$ is a forest and $x$, $y$ are larger than $u\wedge v = (\min
(u_i, v_i))$, $x_3$ and $y_3$ would be comparable. On the other
hand, if $u'\leqslant x$, $u''\leqslant x$ and $u'$, $u''$
$c$-comparable, then $u'\vee u''\leqslant x$.

\medskip\noindent
{\bf C.5~} We have the following partial result:

\medskip\noindent
{\bf Proposition.~}{\it Let $X_1,\cdots,X_n$ be forests and let $A$
be a subset of $X = X_1\times \cdots \times X_n$. Let $A_1$ be the
set of elements $x\in X$ for which there exists $x^1$, $x^2, \cdots,
x^n$ in $A$ with $x_j = x^j_j \geqslant x^i_j$ for all $1\leqslant
i,j\leqslant n$. Let $A_2$ be the set of elements $y\in X$ such that
$y\leqslant x$ for some $x$ in $A_1$. Then $A_2$ contains the
$c.h$-envelope of $A$.}

\begin{rem}
Example~2 above shows that $A_2$ can be strictly larger than the
$c.h$-envelope. Example~1 shows that the straightforward
generalization of the case $n=2$ does not work.
\end{rem}

\medskip\noindent
{\it Proof of the Proposition.~} It is very similar to the proof of
the proposition in C.3 above and left to the reader. \hfill
$\square$

\medskip\noindent
{\bf C.6~} We will now see how the definition of the transversality
relation in Subsection~\ref{sub5.4} is a natural consequence of the
proposition above. As observed earlier, an essential difference with
the toy model is that the transversality relation is used to
construct the classes $\cR (I)$. So, let us just try to define the
relation for the starting class $\cR (\Io)$ associated to the
initial horseshoe $K$. We would have:
\begin{alignat}{3}
&X_1 &=  &\{ (P,Q,n) \in \cR (\Io),\;\; Q\subset Q_u\},\tag{C.10}\\
&X_2 &=  &\{ (P',Q',n') \in \cR (\Io),\;\; P\subset P_s\},\tag{C.11}
\end{alignat}
and $X_3$ is the set of parameter intervals. All sets are partially
ordered by inclusion (of the $Q$'s for $X_1$, of the $P$'s for
$X_2$), and are obviously {\it trees}, with respective roots $(P_u,
Q_u, n_u)$, $(P_s, Q_s, n_s)$, $\Io$.

We start from an intuitive definition of transversality: for
$(P,Q,n)\in X_1$, $(P', Q', n')\in X_2$, $I\in X_3$, we write
\[
Q\; \whforki\;  P'
\]
if for all $t\in I$ we have
\begin{equation}
\de (Q, P') \; \geqslant \; 2 \max (I, |Q|^{1-\eta},
|P'|^{1-\eta}).\tag{C.12}
\end{equation}
(The number $\eta$ in the exponent is necessary in order to keep the
distortions under control.)

The corresponding subset of $X_1 \times X_2 \times X_3$ is
\begin{equation}
A = \{ (Q, P', I),\;\;\; Q\,\whforki\,P'\}.\tag{C.13}
\end{equation}
This set is hereditary but it is not, a priori, concave. The
concavity property (Propositions~\ref{propo4} and \ref{propo7} in
Section~\ref{sec6}) is very useful in many places. So, we wish to
replace $A$ by a larger set which is hereditary and concave. If we
apply the recipe of the proposition in C.5, we are led first to
define a set $A_1$ and then a set $A_2$.  According to the
proposition, $A_1$ should be the set of $(Q, P', I)$ for which there
exist $Q_2$, $Q_3\subset Q$, $P'_1$, $P'_3\subset P'$, $I_1$,
$I_2\subset I$ satisfying
\begin{equation}
Q \;\whfork_{I_1}\;P'_1,\qquad Q_2 \;\whfork_{I_2}\;P',\qquad Q_3
\;\whforki\;P'_3.\tag{C.14}
\end{equation}
As $I_1$, $I_2$ can be chosen arbitrarily small, we can replace them
with single values $t_1$, $t_2\in I$ of the parameter; the three
conditions in (C.14) become:

--~~there exists $P'_1\subset P'$, $t_1\in I$ such that
\begin{equation}
\de (Q, P'_1) \; \geqslant \; 2 \max (|Q|^{1-\eta},
|P'_1|^{1-\eta})\tag{C.15}
\end{equation}

--~~there exists $Q_2\subset Q$, $t_2\in I$ such that
\begin{equation}
\de (Q_2, P') \; \geqslant \; 2 \max (|Q_2|^{1-\eta},
|P'|^{1-\eta})\tag{C.16}
\end{equation}

--~~there exists $Q_3\subset Q$, $P'_3\subset P'$ such that
\begin{equation}
\de (Q_3, P'_3) \; \geqslant \; 2 \max (|I|, |Q_3|^{1-\eta},
|P'_3|^{1-\eta})\tag{C.17}
\end{equation}
for all $t\in I$.

As $P'_1$, $Q_2$, $Q_3$, $P'_3$ may be chosen arbitrarily thin, it
is natural to replace (C.15), (C.16), (C.17) by
\begin{alignat}{3}
&\de (Q, P'_1) \;&\geqslant \; &2 |Q|^{1-\eta} \; \text{ for
$t=t_1$};\tag*{(C.15)$'$}\\
&\de (Q_2, P') \;&\geqslant \; &2 |P'|^{1-\eta} \; \text{ for
$t=t_2$};\tag*{(C.16)$'$}\\
&\de (Q_3, P'_3) \;&\geqslant &2 |I| \; \text{ for all $t\in
I$}.\tag*{(C.17)$'$}
\end{alignat}
Finally, the largest $t$ value of $\de (Q, P'_1)$ that one can hope
for (by choosing $P'_1\subset P'$ appropriately) is $\de_R (Q, P')$;
similarly, the largest value of $\de (Q_2, P')$ that one can hope
for is $\de_L (Q, P')$ and the largest value of $\de (Q_3, P'_3)$
that one can hope for is $\de_{LR} (Q, P')$. Notice that we need
anyway to eliminate $P'_1$, $Q_2$, $P'_3$, $Q_3$ from the definition
because in $\cR (I)$ (instead of $\cR (\Io)$), elements are
constructed inductively and thinner rectangles are constructed at
the end. Replacing $\de (Q, P'_1)$ by $\de_R (Q, P')$, $\de (Q_2,
P')$ by $\de_L (Q, P')$ and $\de (Q_3, P'_3)$ by $\de_{LR} (Q, P')$,
we obtain the three conditions, (T1), (T2),(T3) in
Subsection~\ref{sub5.4}. This defines $\ovfork$.

The last step is to go from $A_1$ to $A_2$, taking the hereditary
envelope of $A_1$, which corresponds exactly to the transition from
$\ovfork$ to $\fork$ in Subsection~\ref{sub5.4}.

\newpage

\bibliographystyle{plain}
\bibliography{main}
\end{document}